\documentclass{book}
\usepackage{amsmath,amssymb,amscd,stmaryrd,bm,rotating,trfsigns,setspace,enumitem}
\usepackage{amsthm} 
\usepackage{mathtools,extarrows}
\usepackage[french]{babel}
\usepackage[all]{xy}
\usepackage{enumitem} 
\usepackage{xcolor}
\usepackage{hyperref}

\makeatletter
\def\revddots{\mathinner{\mkern1mu\raise\p@
		\vbox{\kern7\p@\hbox{.}}\mkern2mu
		\raise4\p@\hbox{.}\mkern2mu\raise7\p@\hbox{.}\mkern1mu}}
\makeatother

\makeatletter
\newcommand*{\relrelbarsep}{.386ex}
\newcommand*{\relrelbar}{%
	\mathrel{%
		\mathpalette\@relrelbar\relrelbarsep
	}%
}
\newcommand*{\@relrelbar}[2]{%
	\raise#2\hbox to 0pt{$\m@th#1\relbar$\hss}%
	\lower#2\hbox{$\m@th#1\relbar$}%
}
\providecommand*{\rightrightarrowsfill@}{%
	\arrowfill@\relrelbar\relrelbar\rightrightarrows
}
\providecommand*{\leftleftarrowsfill@}{%
	\arrowfill@\leftleftarrows\relrelbar\relrelbar
}
\providecommand*{\xrightrightarrows}[2][]{%
	\ext@arrow 0359\rightrightarrowsfill@{#1}{#2}%
}
\providecommand*{\xleftleftarrows}[2][]{%
	\ext@arrow 3095\leftleftarrowsfill@{#1}{#2}%
}
\makeatother

\usepackage{scalerel,stackengine}
\stackMath
\newcommand\reallywidehat[1]{%
	\savestack{\tmpbox}{\stretchto{%
			\scaleto{%
				\scalerel*[\widthof{\ensuremath{#1}}]{\kern-.6pt\bigwedge\kern-.6pt}%
				{\rule[-\textheight/2]{1ex}{\textheight}}%WIDTH-LIMITED BIG WEDGE
			}{\textheight}% 
		}{0.5ex}}%
	\stackon[1pt]{#1}{\tmpbox}%
}

\def\dar[#1]{\ar@<2pt>[#1]\ar@<-2pt>[#1]}

\newtheoremstyle{propo}% name
{1.5em}%         Space above, empty = `usual value'
{1.5em}%         Space below
{\itshape
	\interlinepenalty=10%    Discourage breaking between lines.
	\predisplaypenalty=100% Never break right before a display equation.
	\postdisplaypenalty=100%  Discourage breaking after
}% Body font
{}%         Indent amount (empty = no indent, \parindent = para indent)
{\bfseries}% Thm head font
{.\,--}%        Punctuation after thm head
{\newline}% Space after thm head: \newline = linebreak
{}%         Thm head spec
\theoremstyle{propo}

\newcounter{thm}[section]
\renewcommand{\thethm}{\Roman{chapter}.\arabic{section}.\arabic{thm}}
\newenvironment{thm}[1][]{\refstepcounter{thm}\par\medskip
	\noindent \textbf{Th\'{e}or\`{e}me~\thethm. -- #1} 
	{\color{white}Grothendieck}
	
	\begin{em}}{\end{em}\medskip}
\newenvironment{lem}[1][]{\refstepcounter{thm}\par\medskip
	\noindent \textbf{Lemme~\thethm. -- #1} 
	{\color{white}Grothendieck}
	
	\begin{em}}{\end{em}\medskip}
\newenvironment{defn}[1][]{\refstepcounter{thm}\par\medskip
	\noindent \textbf{D\'{e}finition~\thethm. -- #1} 
	{\color{white}Grothendieck}
	
	\begin{em}}{\end{em}\medskip}
\newenvironment{prop}[1][]{\refstepcounter{thm}\par\medskip
	\noindent \textbf{Proposition~\thethm. -- #1} 
	{\color{white}Grothendieck}
	
	\begin{em}}{\end{em}\medskip}

\newenvironment{cor}[1][]{\refstepcounter{thm}\par\medskip
	\noindent \textbf{Corollaire~\thethm. -- #1} 
	{\color{white}Grothendieck}
	
	\begin{em}}{\end{em}\medskip}

\newenvironment{lemqed}[1][]{\refstepcounter{thm}\par\medskip
	\noindent \textbf{Lemme~\thethm. -- #1} 
	{\color{white}Grothendieck}
	
	\begin{em}}{\end{em}{\color{white}Grothendieck}\null\hfill\qedsymbol\medskip}
\newenvironment{propqed}[1][]{\refstepcounter{thm}\par\medskip
	\noindent \textbf{Proposition~\thethm. -- #1} 
	{\color{white}Grothendieck}
	
	\begin{em}}{\end{em}{\color{white}Grothendieck}\null\hfill\qedsymbol\medskip}
\newenvironment{corqed}[1][]{\refstepcounter{thm}\par\medskip
	\noindent \textbf{Corollaire~\thethm. -- #1} 
	{\color{white}Grothendieck}
	
	\begin{em}}{\end{em}{\color{white}Grothendieck}\null\hfill\qedsymbol\medskip}
\newenvironment{remark}[1][]{\par\medskip
	\noindent \textbf{Remarque~ : #1} 
	{\color{white}Grothendieck}
	
	\begin{rm}}{\end{rm}\medskip}

\newenvironment{remarks}[1][]{\par\medskip
	\noindent \textbf{Remarques~ : #1} 
	{\color{white}Grothendieck}
	
	\begin{rm}}{\end{rm}\medskip}

\newenvironment{remarkqed}[1][]{\par\medskip
	\noindent \textbf{Remarque~ : #1} 
	{\color{white}Grothendieck}
	
	\begin{rm}}{\end{rm}{\color{white}Grothendieck}\null\hfill\qedsymbol\medskip}

\newenvironment{remarksqed}[1][]{\par\medskip
	\noindent \textbf{Remarques~ : #1} 
	{\color{white}Grothendieck}%\vspace{-0.05cm}
	
	\begin{rm}}{\end{rm}{\color{white}Grothendieck}\nopagebreak\null\hfill\qedsymbol\medskip}

\newenvironment{ex}[1][]{\par\medskip
	\noindent \textbf{Exemple~ : #1} 
	{\color{white}Grothendieck}
	
	\begin{rm}}{\end{rm}\medskip}

\newenvironment{exs}[1][]{\par\medskip
	\noindent \textbf{Exemples~ : #1} 
	{\color{white}Grothendieck}
	
	\begin{rm}}{\end{rm}\medskip}

\newenvironment{exqed}[1][]{\par\medskip
	\noindent \textbf{Exemple~ : #1} 
	{\color{white}Grothendieck}%\vspace{-0.05cm}
	
	\begin{rm}}{\end{rm}{\color{white}Grothendieck}\nopagebreak\null\hfill\qedsymbol\medskip}

\newenvironment{exsqed}[1][]{\par\medskip
	\noindent \textbf{Exemples~ : #1} 
	{\color{white}Grothendieck}%\vspace{-0.05cm}
	
	\begin{rm}}{\end{rm}{\color{white}Grothendieck}\nopagebreak\null\hfill\qedsymbol\medskip}

\newenvironment{refqed}[1][]{\par\medskip
	\noindent \textbf{R\'{e}f\'{e}rence~ : #1} 
	{\color{white}Grothendieck}%\vspace{-0.05cm}
	
	\begin{rm}}{\end{rm}{\color{white}Grothendieck}\nopagebreak\null\hfill\qedsymbol\medskip}

\newenvironment{demo}[1][]{\par\medskip
	\noindent \textbf{D\'{e}monstration~ : #1} 
	{\color{white}Grothendieck}%\vspace{-0.05cm}
	
	\begin{rm}}{\end{rm}{\color{white}Grothendieck}\nopagebreak\null\hfill\qedsymbol\medskip}

\newenvironment{demolem}[1][]{\par\medskip
	\noindent \textbf{D\'{e}monstration du lemme~ : #1} 
	{\color{white}Grothendieck}%\vspace{-0.05cm}
	
	\begin{rm}}{\end{rm}{\color{white}Grothendieck}\nopagebreak\null\hfill\qedsymbol\medskip}

\newenvironment{demoprop}[1][]{\par\medskip
	\noindent \textbf{Démonstration de la proposition~ : #1} 
	{\color{white}Grothendieck}%\vspace{-0.05cm}
	
	\begin{rm}}{\end{rm}{\color{white}Grothendieck}\nopagebreak\null\hfill\qedsymbol\medskip}

\newenvironment{suitedemoprop}[1][]{\par\medskip
	\noindent \textbf{Suite de la démonstration de la proposition~ : #1} 
	{\color{white}Grothendieck}%\vspace{-0.05cm}
	
	\begin{rm}}{\end{rm}{\color{white}Grothendieck}\nopagebreak\null\hfill\qedsymbol\medskip}

\newenvironment{demosansqed}[1][]{\par\medskip
	\noindent \textbf{Démonstration~ : #1} 
	{\color{white}Grothendieck}
	
	\begin{rm}}{\end{rm}\medskip}

\newenvironment{esquissedemo}[1][]{\par\medskip
	\noindent \textbf{Esquisse de démonstration~ : #1} 
	{\color{white}Grothendieck}%\vspace{-0.05cm}
	
	\begin{rm}}{\end{rm}{\color{white}Grothendieck}\nopagebreak\null\hfill\qedsymbol\medskip}
\newenvironment{esquissedemolem}[1][]{\par\medskip
	\noindent \textbf{Esquisse de démonstration du lemme~ : #1} 
	{\color{white}Grothendieck}%\vspace{-0.05cm}
	
	\begin{rm}}{\end{rm}{\color{white}Grothendieck}\nopagebreak\null\hfill\qedsymbol\medskip}

\newenvironment{esquissedemopropsansqed}[1][]{\par\medskip
	\noindent \textbf{Esquisse de démonstration de la proposition~ : #1} 
	{\color{white}Grothendieck}
	
	\begin{rm}}{\end{rm}\medskip}

\newenvironment{esquissedemosansqed}[1][]{\par\medskip
	\noindent \textbf{Esquisse de démonstration~ : #1} 
	{\color{white}Grothendieck}
	
	\begin{rm}}{\end{rm}\medskip}

\newenvironment{principedemosansqed}[1][]{\par\medskip
	\noindent \textbf{Principe de démonstration~ : #1} 
	{\color{white}Grothendieck}
	
	\begin{rm}}{\end{rm}\medskip}

\newenvironment{pourdemo}[1][]{\par\medskip
	\noindent \textbf{Pour la démonstration~ : #1} 
	{\color{white}Grothendieck}%\vspace{-0.05cm}
	
	\begin{rm}}{\end{rm}{\color{white}Grothendieck}\nopagebreak\null\hfill\qedsymbol\medskip}

\newenvironment{indicationdemo}[1][]{\par\medskip
	\noindent \textbf{Indication de démonstration~ : #1} 
	{\color{white}Grothendieck}%\vspace{-0.05cm}
	
	\begin{rm}}{\end{rm}{\color{white}Grothendieck}\nopagebreak\null\hfill\qedsymbol\medskip}

\newenvironment{principedemo}[1][]{\par\medskip
	\noindent \textbf{Principe de démonstration~ : #1} 
	{\color{white}Grothendieck}%\vspace{-0.05cm}
	
	\begin{rm}}{\end{rm}{\color{white}Grothendieck}\nopagebreak\null\hfill\qedsymbol\medskip}

\newenvironment{justifdefi}[1][]{\par\medskip
	\noindent \textbf{Justification de la définition~ : #1} 
	{\color{white}Grothendieck}%\vspace{-0.05cm}
	
	\begin{rm}}{\end{rm}{\color{white}Grothendieck}\nopagebreak\null\hfill\qedsymbol\medskip}

\newenvironment{justifdefisansqed}[1][]{\par\medskip
	\noindent \textbf{Justification de la définition~ : #1} 
	{\color{white}Grothendieck}
	
	\begin{rm}}{\end{rm}\medskip}

\newenvironment{listeisansmarge}[1][]{\begin{enumerate}[label=(\roman*),leftmargin=0pt]}{\end{enumerate}}

\newenvironment{listeimarge}[1][]{\begin{enumerate}[label=(\roman*)]}{\end{enumerate}}

\allowdisplaybreaks
\begin{document}
	
	\title{Engendrement de topologies de Grothendieck, 
		
		démontrabilité et opérations sur les sous-topos}
	
	\author{\Large Olivia Caramello$^{(*)}$ et Laurent Lafforgue$^{(**)}$}
	
	\date{
		(*) {\normalsize Università di Insubria, Como, Italie} \\
		{\normalsize Istituto Grothendieck, Corso Statuto 24, Mondovi CN, Italie} \\
		{\normalsize \&} \\
		(**) {\normalsize Centre Lagrange, 103 rue de Grenelle, 75007 Paris, France} \\
		{\normalsize Centre de recherche de Huawei, 18 quai du Point du Jour,} \\
		{\normalsize  92100 Boulogne-Billancourt, France}
	}
	
	\maketitle
	
	\section*{Préface}
	
	Cette monographie étudie systématiquement une notion centrale de la théorie des topos, celle de sous-topos, des points de vue géométrique, algébrique et logique. Cet invariant fondamental peut être calculé efficacement de différentes manières, correspondant aux différentes perspectives que l'on peut avoir sur le topos concerné : ainsi, lorsque le topos est représenté comme la catégorie des faisceaux sur un site, ses sous-topos correspondent aux topologies de Grothendieck sur la catégorie sous-jacente au site qui raffinent la topologie donnée, tandis que si le topos est présenté comme classifiant une théorie, ils correspondent aux ``quotients'' de cette théorie, c'est-à-dire aux manières d'enrichir cette théorie en lui ajoutant des axiomes supplementaires formulés dans le même langage.
	
	Cela rend la notion de sous-topos particulièrement bien adaptée à une application de la technique des ``ponts'' de \cite{Memoire} et \cite{TST}. Certains de ces ``ponts'' ont déjà été étudiés dans \cite{TST}, et ont notamment conduit à des équivalences entre les systèmes de déduction classiques pour les théories géométriques et des nouveaux ``systèmes de démonstration'' pour les cribles, dont les règles d'inférence correspondent aux axiomes des topologies de Grothendieck et dans lesquels des calculs qui seraient difficilement réalisables dans les systèmes classiques deviennent possibles :
	$$
	\resizebox{\textwidth}{!}
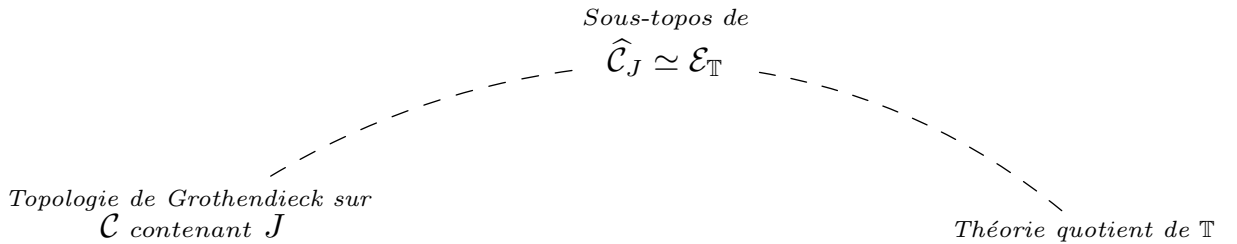
{
		\xymatrix {
			& & \stackrel{\stackrel{\scriptstyle Sous\textup{-}topos\textup{ }de}{\textstyle \textup{ }}}{  \widehat{\mathcal C}_J \simeq {\cal E}_{\mathbb T}} \ar@/^12pt/@{--}[drr] & & \\
			\stackrel{\scriptstyle \textit{Topologie de Grothendieck sur }}{ \displaystyle{ \mathcal{C} } \hspace{0.1cm}{\scriptstyle contenant } \textrm{ } \displaystyle{J} } \ar@/^12pt/@{--}[urr] & & & & {\scriptstyle \textit{Théorie quotient de } {\mathbb T}}
	}}
	$$
	
	Ces équivalences permettent notamment de traduire les problèmes de démontrabilité dans les théories géométriques en termes de problèmes d'engendrement de topologies de Grothendieck. Plus précisément, tous les problèmes de démontrabilité se traduisent constructivement en la question de savoir si certains cribles appartiennent à la topologie de Grothendieck engendrée par une certaine famille de cribles, une tâche qui peut s'avérer plus facile à traiter, comme déjà illustré dans \cite{TST} à travers la discussion de plusieurs exemples. Il existe même une formule fermée pour calculer la topologie de Grothendieck engendrée par une famille donnée de cribles : une prémière version de cette formule avait déjà été donnée dans \cite{TST}, et une version simplifiée est présentée dans la section 1 du chapitre III. De plus, une expression alternative en termes d'arbres est fournie dans la section 2 du même chapitre.
	
	Dans ce travail, nous approfondissons ces équivalences et les décrivons avec plus de détails dans des cas notables qui sont pertinents pour des applications et des implémentations informatiques. De plus, nous présentons plusieurs nouveaux résultats sur le calcul de certaines opérations sur les familles de cribles, les topologies de Grothendieck et les sous-topos.
	
	Une autre de nos motivations pour entreprendre cette étude systématique des sous-topos est leur pertinence en relation avec le développement d'un apprentissage profond topossique, basé sur la représentation d'éléments de réalité (tels que des images, des textes ou tout autre type de données) en termes de \emph{sous-topos} plutôt que de \emph{vecteurs} (comme c'est le cas dans l'apprentissage profond usuel), et sur les opérations géométriques et algébriques sur les sous-topos qui peuvent remplacer celles entre vecteurs dans ce cadre. Alors que l'apprentissage profond classique est basé sur la métrique euclidienne, utilisée pour mesurer la distance entre les vecteurs, les correspondances géométriques entre topos fournissent une notion abstraite de distance par rapport à laquelle on peut developper une théorie de l'approximation. D'autre part, les correspondances induisent des opérations d'image directe et d'image réciproque sur les sous-topos. Afin de jeter les bases d'un tel développement, nous étudions en détail ces opérations, ainsi que d'autres constructions impliquant des sous-topos, en fournissant dans chaque cas des formules explicites pour les calculer.
	
	Le contenu du livre peut être résumé comme suit.
	
	Le chapitre I commence par des rappels, dans la section 1, des préliminaires pertinents de théorie des topos, notamment la relation entre cette théorie et celle des espaces topologiques ordinaires, la sémantique catégorique des théories (géométriques) du premier ordre et la construction des topos classifiants : les résultats rappelés sont bien connus et le lecteur pourra trouver des dévéloppements complets à leur sujet dans les références classiques \cite{SGA1}, \cite{Elephant} et \cite{CategoricalLogic}. Ensuite, dans la section 2, la notion de sous-topos est étudiée à la fois du point de vue géométrique des sites et du point de vue logique des théories, ce qui conduit aux équivalences découlant des ``ponts'' mentionnées ci-dessus : nous formulons un théorème général permettant de transformer les problèmes de démontrabilité en problèmes d'engendrement de topologies de Grothendieck et discutons de son application dans des cas particuliers notables. Dans la section 3, nous nous intéressons aux opérations sur les sous-topos : l'ensemble des sous-topos d'un topos de Grothendieck donné possède un ordre naturel défini par la relation d'inclusion, qui en fait une algèbre de co-Heyting. Des descriptions explicites des opérations de Heyting corresponantes sur les topologies de Grothendieck et les théories quotient avaient déjà été données dans \cite{TST}, mais nous en fournissons ici des descriptions alternatives. De plus, nous étudions les opérations d'image directe et d'image réciproque des sous-topos le long des morphismes de topos, en donnant des descriptions géométriques et logiques de celles-ci. Nous montrons notamment que, si le morphisme est localement connexe (au sens que sa composante d'image inverse a un adjoint à gauche qui est compatible avec les changements de base), alors l'opération induite d'image réciproque sur les sous-topos a un adjoint à gauche, que nous appelons l'opération d'\emph{image directe extraordinaire}. Nous considérons également les correspondances entre topos et les opérations sur les sous-topos qu'elles induisent. 
	
	Le chapitre II est une étude abstraite des topologies de Grothendieck, des opérations de fermeture sur les monomorphismes dans un topos et des sous-topos, visant à faire apparaître naturellement ces notions comme des points fixes de certaines connexions de Galois. Nous commençons par rappeler une méthode générale, déjà largement utilisée par Birkhoff (cf. \cite{Galois}), pour engendrer des connexions de Galois à partir de relations entre ensembles ou classes, ainsi que la caractérisation générale des points fixes des connexions de Galois. Nous introduisons ensuite trois relations pertinentes : une entre les cribles et les préfaisceaux, une entre les monomorphismes dans un topos et les objets du topos (représentant une version plus intrinsèque et une extension de la première), et une entre les cribles et les plongements de préfaisceaux. Les équivalences qui en résultent entre les points fixes des connexions de Galois induites donnent respectivement l'équivalence bien connue entre les topologies de Grothendieck et les sous-topos d'un topos de préfaisceaux, l'équivalence (plus générale et intrinsèque) entre les opérations de fermeture sur les monomorphismes d'un topos donné et les sous-topos de ce topos, et une équivalence entre les topologies de Grothendieck et certaines ``propriétés de fermeture'' sur les plongements de préfaisceaux. Une conséquence technique particulièrement intéressante de ces connexions de Galois est qu'elles permettent de réaliser concrètement des opérations de ``fermeture'' ou d'``engendrement'' comme les composés des deux foncteurs formant la connexion. Par exemple, la topologie de Grothendieck engendrée par une famille donnée $J$ de cribles est précisément $G(F(J))$, où $G$ et $F$ sont les deux foncteurs de la connexion de Galois, tandis que le plus petit sous-topos contenant une famille donnée $I$ de préfaisceaux n'est autre que $F(G(I))$. Ces réinterprétations ont à leur tour des conséquences remarquables sur certaines propriétés fondamentales de ces constructions, qui seraient peu visibles ou plus difficiles à établir par d'autres moyens. Par exemple, on peut déduire, comme simples conséquences formelles des propriétés de ces connexions de Galois, que les intersections de sous-topos peuvent être calculées directement en termes des cribles qui appartiennent à l'une ou l'autre des topologies de Grothendieck associées (sans qu'il soit nécessaire de calculer la topologie de Grothendieck engendrée par leur réunion). On déduit également qu'un plongement de préfaisceaux est fermé par rapport à une famille de cribles qui est stable par image réciproque si et seulement s'il est fermé par rapport à la topologie de Grothendieck engendrée par cette famille. Ce dernier fait jouera un rôle clé au chapitre III pour établir une formule fermée caractérisant les cribles couvrant de la topologie de Grothendieck engendrée par une famille donnée de cribles.
	
	Le chapitre III est consacré à l'étude du phénomène d'engendrement des topologies de Grothendieck et à son lien avec les problèmes de démontrabilité en logique géométrique du premier ordre. La section 1 présente la formule déjà mentionnée pour la topologie de Grothendieck engendrée  par une famille donnée de crible qui est stable par image réciproque. Cette formule, qui représente un raffinement de la formule déjà obtenue dans \cite{TST} (cf. proposition 4.1.1 dans cet ouvrage), est remarquablement simple : les cribles couvrants pour la topologie engendrée par une telle famille $J$ sont précisément les cribles dont la fermeture est le crible maximal, c'est-à-dire ceux tels que tout crible $J$-fermé les contenant est le crible maximal. Cette formule est ensuite appliquée pour obtenir une description explicite de la topologie de Grothendieck engendrée par une famille de topologies, ainsi que pour le problème du calcul des produits finis de topos. En effet, nous montrons que les topologies de Grothendieck qui présentent le produit (fini) de topos de faisceaux sur des sites coincident avec les topologies engendrées par les images réciproques, le long des projections canoniques du produit, des topologies qui présentent les facteurs; elles peuvent donc être calculées à l'aide de notre formule. Dans la dernière partie de la section, nous établissons un corollaire montrant que le produit des topos de faisceaux sur des espaces topologiques qui sont tous localement compacts - à l'exception éventuelle d'un seul - est équivalent au topos des faisceaux sur l'espace produit. Dans la section 2, nous obtenons une formule alternative pour calculer la topologie de Grothendieck engendrée par une famille donnée de (pré)cribles, formulée en termes d'arbres et de ``multi-recouvrements''. Pour obtenir cette formule, qui \emph{a priori} nécessiterait une induction transfinie (comme le croyait Johnstone – cf. les remarques suivant la preuve de la proposition 2.1.9 dans \cite{Elephant}), nous exploitons notamment une propriété de stabilité de l'opération de multi-composition de précribles par rapport aux images réciproques, que nous formulons explicitement dans le langage des arbres. Il est intéressant de comparer cette formule avec celle de la section 1, car elles sont de natures assez différentes : la dernière est plus utile dans les situations où la topologie est engendrée à partir de \emph{précribles}, alors que la première est plus utile lorsque la topologie est engendré à partir d'une famille de \emph{cribles}. Cela reflète une ``dichotomie'' intrigante dans les topologies de Grothendieck qui apparaissent naturellement dans la pratique mathématique : certaines d'entre elles, apparaissant principalement dans des contextes logiques (comme, par exemple, la topologie dense ou atomique, ou la topologie de De Morgan) ont pour recouvrements des cribles qui peuvent être caractérisés par une ``condition fermée'', tandis que d'autres, qui apparaissent principalement dans des contextes géométriques (comme, par exemple, la topologie de Zariski ou la topologie étale), sont présentées en termes de \emph{prétopologies}, c'est-à-dire en spécifiant des familles génératrices de \emph{précribles}. La section 3 est consacrée à l'étude des ``ponts'' entre les topologies de Grothendieck et les théories ``quotient'' induits par la notion invariante de sous-topos lorsque l'on dispose d'une équivalence entre le topos des faisceaux sur un certain site et le topos classifiant d'une théorie géométrique. Comme expliqué ci-dessus, ces ``ponts'' permettent d'interpréter les problèmes de démontrabilité en logique géométrique en termes de problèmes d'engendrement pour les topologies de Grothendieck, et réciproquement. Nous présentons un théorème formalisant une méthode très générale pour étudier les problèmes de démontrabilité en termes de topologies de Grothendieck, en nous concentrant sur le contexte dans lequel les calculs peuvent être effectués le plus efficacement, à savoir celui des sites cartésiens de présentation pour le topos classifiant donné. Nous donnons une description explicite des catégories syntactiques cartésiennes des théories sans axiomes et observons qu'une telle description se simplifie considérablement si la signature donnée ne contient pas de symboles de fonction. En effet, la construction dans ce cas est purement syntactique, c'est-à-dire ne nécessite pas de considérer des relations d'équivalence non triviales sur les termes. Cela conduit à introduire une construction permettant de transformer toute théorie géométrique en une théorie sur une signature relationnelle qui possède le même topos classifiant : en effet, nous montrons que les catégories syntactiques géométriques (resp. cohérentes, resp. régulières) des deux théories sont équivalentes, ce qui les rend \emph{bi-interprétables}.
	
	Le chapitre IV est une étude de l'operation d'image réciproque des sous-topos le long de morphismes localement connexes. Dans la section 1, nous établissons un théorème fournissant une description compacte de la topologie de Grothendieck présentant un sous-topos image réciproque en termes de la topologie présentant le sous-topos d'origine. Ce résultat peut être considéré comme une généralisation topossique du théorème de Giraud qui caractérise les topologies image réciproque le long d'une fibration. Rappelons en effet que les fibrations, munies de topologies de Giraud, induisent des morphismes de topos localement connexes - un résultat déjà prouvé dans \cite{Denseness}. Toutes les connaissances nécessaires sur les fibrations et les morphismes localement connexes sont rappelées dans la première partie de la section, à l'intention du lecteur qui ne serait pas familier avec ces notions ; en particulier, un résultat crucial de \cite{Localfibrations}, caractérisant les morphismes localement connexes comme les morphismes essentiels dont l'image essentielle est une fibration, est rappelé et expliqué en détail. Ensuite, dans la section 2, une description explicite de l'opération d'image directe extraordinaire pour les sous-topos induite par un morphisme localement connexe est fournie, ainsi que des caractérisations explicites des topologies de Grothendieck correspondantes dans le cas d'un morphisme localement connexe induit par une fibration. L'existence de cette opération, qui est adjointe à l'opération d'image réciproque, garantit que cette dernière préserve les réunions arbitraires. En exploitant ce fait, ainsi que la possibilité de factoriser tout morphisme de topos en une inclusion suivie d'un morphisme localement connexe, nous déduisons que l'opération d'image réciproque des sous-topos \emph{le long de n'importe quel morphisme de topos} préserve les réunions finies de sous-topos, un nouveau résultat qui sera certainement utile dans des applications.
	
	Des travaux supplémentaires, qui seront inclus dans une version future de cette monographie, porteront sur l'application des méthodes introduites ici au calcul des produits orientés et fibrés de topos.
	
	\tableofcontents
	
	\chapter{Introduction : Pourquoi les sous-topos?}\label{chap1}

\section{Pourquoi les topos?}\label{sec1.1}

\subsection{Les topos comme vaste généralisation des espaces topologiques}\label{ssec1.1.1}

Grothendieck a remarqué qu'associer à tout espace topologique $X$ la catégorie ${\mathcal E}_X$ des faisceaux d'ensembles sur $X$ permet de plonger le monde des espaces topologiques dans un monde beaucoup plus vaste, celui des topos, où toutes les notions de base de la topologie restent définies et les intuitions fondamentales de nature géométrique continuent à s'appliquer.

Un topos est par définition une catégorie qui peut être présentée à équivalence près comme la catégorie des faisceaux d'ensembles sur un ``site'', c'est-à-dire sur une catégorie petite ou essentiellement petite munie d'une ``topologie''. La notion de site est le contexte le plus vaste dans lequel la notion de faisceau garde un sens.

La notion de topologie sur une catégorie ${\mathcal C}$ localement petite peut s'exprimer en termes de ``cribles'' ou en termes de ``familles couvrantes'' de morphismes. Sous sa seconde forme, elle a un sens sur n'importe quelle catégorie:

\begin{defn} \label{defn1.1.1} 

Une topologie $J$ sur une catégorie ${\mathcal C}$ consiste en une notion de recouvrement des objets $X$ de ${\mathcal C}$ par des familles de morphismes de ${\mathcal C}$
$$
\left( X_i \xrightarrow{ \ x_i \ } X \right)_{i \in I} \qquad \mbox{dites alors $J$-couvrantes,}
$$
qui satisfait les conditions suivantes:

%\smallskip

$\left\lmoustache\begin{matrix}
(0) &\mbox{Pour qu'une famille de morphismes $\left( X_i \xrightarrow{ \ x_i \ } X \right)_{i \in I}$ soit $J$-couvrante,} \hfill \\
&\mbox{il suffit qu'existe une famille $J$-couvrante} \hfill \\
&\left( Y_{i'} \xrightarrow{ \ y_{i'} \ } X \right)_{i' \in I'} \\
&\mbox{dont chaque élément $Y_{i'} \xrightarrow{ \ y_{i'} \ } X$ se factorise à travers l'un au moins des} \hfill \\
&\mbox{morphismes $x_i$ en $Y_{i'} \longrightarrow X_i \xrightarrow{ \ x_i \ } X$.} \hfill \\
(1) &\mbox{Pour tout objet $X$ de ${\mathcal C}$, le morphisme $X \xrightarrow{ \ {\rm id}_X \ } X$ est $J$-couvrant.} \hfill \\
(2) &\mbox{Pour tout morphisme $X' \xrightarrow{ \ x \ } X$ de ${\mathcal C}$ et toute famille $J$-couvrante} \hfill \\
&\left( X_i \xrightarrow{ \ x_i \ } X \right)_{i \in I} \, , \\
&\mbox{il existe une famille $J$-couvrante $\left( X'_{i'} \xrightarrow{ \ x'_{i'} \ } X' \right)_{i' \in I'}$ dont} \hfill 
\end{matrix}\right.
$

$\left\rmoustache\begin{matrix}
&\mbox{chaque élément $x'_{i'} : X'_{i'} \to X'$ s'inscrit dans un carré commutatif} \hfill \\
&\xymatrix{
X'_{i'} \ar[d]_{x'_{i'}} \ar[r] &X_i \ar[d]^{x_i} \\
X' \ar[r]^x &X
} \\
&\mbox{pour au moins un des morphismes $x_i : X_i \to X$.} \hfill \\
(3) &\mbox{Pour tout $J$-recouvrement d'un objet $X$ par des objets $X_i$, $i \in I$,} \hfill \\
&\left( X_i \xrightarrow{ \ x_i \ } X \right)_{i \in I_{ \ }}  \\
&\mbox{et pour tout recouvrement de chaque objet $X_i$, $i \in I$,} \hfill \\
&\left( X_{i,j} \xrightarrow{ \ x_{i,j} \ } X_i \right)_{j \in I_i} \, ,  \\
&\mbox{la famille des composés} \hfill \\
&\left( X_{i,j} \xrightarrow{ \ x_{i,j} \ } X_i \xrightarrow{ \ x_i \ } X \right)_{i \in I \atop j \in I_i} \, ,  \\
&\mbox{forme un $J$-recouvrement de $X$.} \hfill
\end{matrix}\right.$
\end{defn}

\begin{remark}
Si ${\mathcal C}$ est une catégorie petite ou essentiellement petite, la propriété (0) implique que la notion de famille $J$-couvrante
$$
\left( X_i \xrightarrow{ \ x_i \ } X \right)_{i \in I}
$$
dépend seulement du ``crible'' engendré par cette famille, constitué des morphismes
$$
X' \longrightarrow X
$$
qui se factorisent à travers l'un au moins des $x_i$ en $X' \longrightarrow X_i \xrightarrow{ \ x_i \ } X$.
\end{remark}

\begin{exsqed}
	\begin{listeisansmarge}
\item Si $X$ est un espace topologique, on dispose de la petite catégorie ${\mathcal C}_X$ constitu\'ee des ouverts $U$ de $X$ et des morphismes d'inclusion entre ouverts $U' \hookrightarrow U$. 

La topologie ordinaire $J_X$ de $X$ est celle pour laquelle une famille d'inclusions
$$
\left( U_i \xymatrix{\ar@{^{(}->}[r] & } U \right)_{i \in I}
$$
est couvrante si et seulement si
$$
U = \bigcup_{i \in I} U_i \, .
$$

\item La topologie canonique d'un topos ${\mathcal E}$ est définie en demandant qu'une famille de morphismes de ${\mathcal E}$
$$
\left( E_i \xrightarrow{ \ e_i \ } E \right)_{i \in I}
$$
soit couvrante si et seulement si elle est globalement épimorphique.
	\end{listeisansmarge}
\end{exsqed}

\begin{defn}\label{defn1.1.2}
	\begin{listeimarge}
\item  Un préfaisceau d'ensembles sur une catégorie ${\mathcal C}$ est un foncteur contravariant de ${\mathcal C}$ vers la catégorie $\mbox{\rm Ens}$ des ensembles
$$
P : {\mathcal C}^{\rm op} \longrightarrow {\rm Ens}.
$$

\item Un faisceau d'ensembles sur une catégorie ${\mathcal C}$ munie d'une topologie $J$ est un préfaisceau
$$
F : {\mathcal C}^{\rm op} \longrightarrow {\rm Ens}
$$
tel que, pour toute famille couvrante
$$
\left( X_i \xrightarrow{ \ x_i \ } X \right)_{i \in I} \, ,
$$
toute famille d'éléments
$$
m_i \in F(X_i) \, , \quad i \in I \, ,
$$
qui est ``compatible'' au sens que, pour tout carré commutatif de ${\mathcal C}$
$$
\xymatrix{
&X_{i_1} \ar[rd]^{x_{i_1}} \\
X' \ar[ru]^{a} \ar[rd]_b &&X \\
&X_{i_2} \ar[ru]_{x_{i_2}}
}
$$
on a
$$
F(a)(m_{i_1}) = F(b) (m_{i_2}) \quad \mbox{dans} \quad F(X') \, ,
$$
provient par les applications $F(x_i) : F(X) \to F(X_i)$, $i \in I$, d'un unique élément $m \in F(X)$.
	\end{listeimarge}
\end{defn}

\begin{remarksqed}
		\begin{listeisansmarge}
\item Si ${\mathcal C}$ est une catégorie localement petite, elle est munie du foncteur de Yoneda
$$
\begin{matrix}
y : {\mathcal C} &\longrightarrow &\widehat{\mathcal C} = [{\mathcal C}^{\rm op} , {\rm Ens}] \\
\hfill X &\longmapsto &{\rm Hom} (\bullet,X) \hfill
\end{matrix}
$$
vers la catégorie $\widehat{\mathcal C}$ des préfaisceaux sur ${\mathcal C}$, c'est-à-dire des foncteurs ${\mathcal C}^{\rm op} \to {\rm Ens}$.

D'après le lemme de Yoneda, le foncteur
$$
y : {\mathcal C} \longrightarrow \widehat{\mathcal C}
$$
est pleinement fidèle.

Les objets de $\widehat{\mathcal C}$, c'est-à-dire les préfaisceaux, qui sont isomorphes à l'image par $y$ d'un objet de ${\mathcal C}$ (nécessairement unique à unique isomorphisme près) sont dits ``représentables''.

\item La remarque précédente s'applique en particulier aux topos, qui sont des catégories localement petites.

On montre de plus qu'un préfaisceau sur un topos ${\mathcal E}$ est représentable si et seulement si il est un faisceau pour la topologie canonique de ${\mathcal E}$. 
\end{listeisansmarge}
\end{remarksqed}

\begin{defn}\label{defn1.1.3}
		\begin{listeimarge}
\item On appelle ``site'' la donnée $({\mathcal C},J)$ d'une catégorie ${\mathcal C}$ petite ou essentiellement petite et d'une topologie $J$ de ${\mathcal C}$.

On note alors
$$
\widehat{\mathcal C}_J
$$
la sous-catégorie pleine de $\widehat{\mathcal C}$ constituée des $J$-faisceaux.

\item On appelle ``topos'' les catégories ${\mathcal E}$ qui admettent des équivalences
$$
{\mathcal E} \xrightarrow{ \ \sim \ } \widehat{\mathcal C}_J
$$
avec des catégories de faisceaux sur des sites $({\mathcal C},J)$.

\item On appelle morphisme de topos
$$
f : {\mathcal E}' \longrightarrow {\mathcal E}
$$
les paires de foncteurs adjoints
$$
f = \left( {\mathcal E} \xrightarrow{ \ f^* \ } {\mathcal E}' \, , \ {\mathcal E}' \xrightarrow{ \ f_* \ } {\mathcal E} \right)
$$
dont la composante adjointe à gauche
$$
f^* : {\mathcal E} \longrightarrow {\mathcal E}'
$$
respecte non seulement les colimites arbitraires mais aussi les limites finies.

\item On appelle transformations de morphismes de topos
$$
\left( {\mathcal E}' \xrightarrow{ \ f \ } {\mathcal E} \right) \xrightarrow{ \ \sigma \ } \left( {\mathcal E}' \xrightarrow{ \ g\ } {\mathcal E} \right)
$$
les transformations naturelles de foncteurs
$$
\sigma : f^* \longrightarrow g^* \, .
$$
\end{listeimarge}
\end{defn}

\begin{remarksqed}
		\begin{listeisansmarge}
\item Ainsi, les topos forment une 2-catégorie.

La catégorie des morphismes d'un topos ${\mathcal E}'$ vers un topos ${\mathcal E}$ est notée
$$
{\rm Geom} ({\mathcal E}',{\mathcal E}) \, .
$$

\item Si $X$ et $X'$ sont deux espaces topologiques, toute application continue
$$
f : X' \longrightarrow X
$$
définit un morphisme de topos
$$
(f^* , f_*) : {\mathcal E}_{X'} \longrightarrow {\mathcal E}_X \, .
$$
Réciproquement, si $X$ est un espace topologique ``sobre'', tout morphisme de topos
$$
{\mathcal E}_{X'} \longrightarrow {\mathcal E}_X
$$
provient d'une unique application continue $X' \to X$.

\item En particulier, tout point $x$ d'un espace topologique $X$ définit un morphisme de topos
$$
(x^* , x_*) : {\rm Ens} \longrightarrow {\mathcal E}_X
$$
puisque ${\rm Ens}$ est le topos des faisceaux de l'espace à un point $\{\bullet\}$.

Il y a réciproque si $X$ est sobre.

\item Pour cette raison, on appelle ``points'' d'un topos ${\mathcal E}$ les morphismes de topos
$$
{\rm Ens} \longrightarrow {\mathcal E} \, .
$$
Ils forment la catégorie ${\rm Geom} ({\rm Ens} , {\mathcal E}) = {\rm pt} ({\mathcal E})$.

Les morphismes ${\mathcal E}' \to {\mathcal E}$ sont aussi appelés ``points généralisés'' de ${\mathcal E}$ paramétrés par ${\mathcal E}'$. 
\end{listeisansmarge}
\end{remarksqed}

\subsection{Les topos comme invariants universels}\label{ssec1.1.2}

On a déjà rappelé que l'on peut associer à tout topos ${\mathcal E}$ la catégorie de ses points
$$
{\rm pt} ({\mathcal E}) = {\rm Geom} ({\rm Ens} , {\mathcal E})
$$
ou, plus généralement, la catégorie de ses points généralisés paramétrés par un topos ${\mathcal E}'$
$$
{\rm Geom} ({\mathcal E}' , {\mathcal E}) \, .
$$
Cette construction est ``invariante'' au sens que toute équivalence de topos
$$
{\mathcal E}_1 \xrightarrow{ \ \sim \ } {\mathcal E}_2
$$
induit une équivalence de catégories
$$
{\rm pt} ({\mathcal E}_1) \xrightarrow{ \ \sim \ } {\rm pt} ({\mathcal E}_2)
$$
ou, plus généralement,
$$
{\rm Geom} ({\mathcal E}' , {\mathcal E}_1) \xrightarrow{ \ \sim \ } {\rm Geom} ({\mathcal E}',{\mathcal E}_2) \, .
$$
Elle est même covariante au sens que tout morphisme de topos
$$
{\mathcal E}_1 \longrightarrow {\mathcal E}_2
$$
induit un foncteur
$$
{\rm pt} ({\mathcal E}_1) \longrightarrow {\rm pt} ({\mathcal E}_2)
$$
ou, plus généralement,
$$
{\rm Geom} ({\mathcal E}' , {\mathcal E}_1) \longrightarrow {\rm Geom} ({\mathcal E}',{\mathcal E}_2) \, .
$$

On rappellera de même plus loin que la notion de sous-espace d'un espace topologique se généralise en une notion de sous-topos de n'importe quel topos.

Les sous-topos d'un topos ${\mathcal E}$ forment un ensemble ordonné par l'inclusion
$$
{\rm ST} ({\mathcal E})
$$
dans lequel les réunions et intersections arbitraires sont toujours bien définies, ainsi que les différences de paires de sous-topos.

La construction
$$
{\mathcal E} \longmapsto {\rm ST} ({\mathcal E})
$$
est ``invariante'' au sens que toute équivalence de topos
$$
{\mathcal E}_1 \xrightarrow{ \ \sim \ } {\mathcal E}_2
$$
induit une bijection qui respecte la relation d'ordre $\supseteq$
$$
{\rm ST} ({\mathcal E}_1) \xrightarrow{ \ \sim \ } {\rm ST} ({\mathcal E}_2) \, .
$$
Elle est même à la fois covariante et contravariante au sens que tout morphisme de topos
$$
f : {\mathcal E}_1 \longrightarrow {\mathcal E}_2
$$
définit une opération d'image directe des sous-topos
$$
f_* : {\rm ST} ({\mathcal E}_1) \longrightarrow {\rm ST} ({\mathcal E}_2)
$$
et une opération d'image réciproque adjointe à gauche de la précédente
$$
f^{-1} : {\rm ST} ({\mathcal E}_2) \longrightarrow {\rm ST} ({\mathcal E}_1) \, .
$$

La construction qui associe à tout espace topologique ses espaces ou modules de cohomologie se généralise en une construction qui associe à tout topos ${\mathcal E}$ ses modules de cohomologie
$$
H^i ({\mathcal E},A) \, , \quad i \in {\mathbb N} \, ,
$$
à coefficients dans n'importe quel anneau $A$.

Les modules de cohomologie des topos sont des invariants puisque toute équivalence de topos
$$
{\mathcal E}_1 \xrightarrow{ \ \sim \ } {\mathcal E}_2
$$
induit des isomorphismes de modules
$$
H^i ({\mathcal E}_1, A) \xrightarrow{ \ \sim \ } H^i ({\mathcal E}_2 , A) \, .
$$
Les invariants cohomologiques sont contravariants au sens que tout morphisme de topos
$$
f : {\mathcal E}_1 \longrightarrow {\mathcal E}_2
$$
induit des morphismes de modules en sens inverse
$$
f^* : H^i ({\mathcal E}_2, A) \longrightarrow H^i ({\mathcal E}_1 , A) \, , \quad i \in {\mathbb N} \, .
$$

En fait, les topos ont d'abord été introduits par Grothendieck comme formant le cadre naturel le plus vaste dans lequel les invariants cohomologiques restent définis et conservent leurs structures essentielles. En effet, il avait d'abord dégagé dans son article \cite{Tohoku} les conditions formelles permettant de développer un formalisme cohomologique dans des catégories, et montré que ces conditions étaient vérifiées d'un côté par les catégories de faisceaux linéaires sur des espaces topologiques, et d'un autre côté par les catégories de préfaisceaux linéaires sur des petites catégories. Le pas suivant a consisté à unifier ces deux familles d'exemples en introduisant les notions de site et de faisceau sur un site, et à observer que les catégories de faisceaux linéaires sur des sites vérifient toujours les conditions de l'article \cite{Tohoku}.

\medskip

De même, la définition des groupes fondamentaux de Poincaré des espaces topologiques $X$ connexes et localement connexes par arcs en les points $x \in X$ se généralise dans le cadre des topos. Pour cela, Grothendieck a remplacé dans \cite{SGA1} la définition classique en termes de classes d'homotopie de lacets par la définition plus abstraite mais équivalente comme groupes de symétries des foncteurs fibres en les points $x$ des catégories de revêtements ``étales'', c'est-à-dire de revêtements localement triviaux pour un choix de topologie.

Cette définition permet d'associer à tout topos ${\mathcal E}$ ``localement connexe'' et à tout point $x$ de ${\mathcal E}$ un ``groupe fondamental'' $\pi_1 ({\mathcal E},x)$ qui a une structure de groupe topologique pro-discret. Si ${\mathcal E}$ est ``connexe'', le foncteur fibre en un point $x$ définit une équivalence de la catégorie des revêtements étales de l'objet terminal de ${\mathcal E}$ vers celle des ensembles discrets munis d'une action continue de $\pi_1 ({\mathcal E},x)$.

Les propriétés d'être ``localement connexe'' ou ``connexe'' sont invariantes par équivalence de topos. Elles généralisent les notions topologiques classiques au sens qu'un espace topologique $X$ est ``localement connexe'' ou ``connexe'' si et seulement si son topos associé ${\mathcal E}_X$ est ``localement connexe'' ou ``connexe''.

Les groupes fondamentaux des topos localement connexes constituent des invariants covariants dans la mesure où tout morphisme de tels topos
$$
f : {\mathcal E}' \longrightarrow {\mathcal E}
$$
induit des morphismes de groupes topologiques
$$
\pi_1 ({\mathcal E}' , x') \longrightarrow \pi_1 ({\mathcal E},f(x'))
$$
en les points $x'$ de ${\mathcal E}'$.

La construction des groupes d'homotopie supérieure $\pi_i (X,x)$, $i \geq 2$, se généralise également aux topos localement connexes dont ils constituent encore des invariants covariants.

\medskip

On peut proposer une définition générale de la notion d'invariant des topos.

\begin{defn}\label{defn1.1.4}
	\begin{listeimarge}
\item Une propriété invariante des topos est une propriété qui se formule en des termes de la théorie des catégories qui sont bien définis pour tous les topos, et qui est respectée par les équivalences de topos
$$
{\mathcal E}_1 \xrightarrow{ \ \sim \ } {\mathcal E}_2 \, .
$$
\item Un invariant des topos est une construction catégorique qui associe à chaque topos ${\mathcal E}$ (ou à chaque topos ${\mathcal E}$ de la classe des topos qui possèdent une certaine propriété invariante) un modèle
$$
H({\mathcal E})
$$
d'une certaine théorie ${\mathbb H}$, de telle sorte que les équivalences de topos
$$
{\mathcal E}_1 \xrightarrow{ \ \sim \ } {\mathcal E}_2
$$
induisent des isomorphismes de modèles
$$
H({\mathcal E}_1) \xrightarrow{ \ \sim \ } H({\mathcal E}_2) \, .
$$
\item Un invariant des topos est dit covariant [resp. contravariant] si tout morphisme de topos (ou de topos éléments de la classe où il est bien défini)
$$
f : {\mathcal E}_1 \longrightarrow {\mathcal E}_2
$$
définit un morphisme de modèles
$$
f_* : H({\mathcal E}_1) \longrightarrow H({\mathcal E}_2) \qquad \mbox{[resp.} \ \ f^* : H({\mathcal E}_2) \longrightarrow H({\mathcal E}_1) \, \mbox{],}
$$
de manière compatible avec la composition des morphismes.
\end{listeimarge}
\end{defn}

\begin{remarksqed}
		\begin{listeisansmarge}
\item Par extension, une propriété des sites $({\mathcal C},J)$ ou de tout autre type de présentation ${\mathcal P}$ permettant de construire un topos ${\mathcal E}_{\mathcal P}$ peut être dite ``invariante'' si elle ne dépend que des topos associés $\widehat{\mathcal C}_J$ ou ${\mathcal E}_{\mathcal P}$ et est respectée par les équivalences de topos
$$
\widehat{\mathcal C}_J \xrightarrow{ \ \sim \ } \widehat{\mathcal D}_K \quad \mbox{ou} \quad {\mathcal E}_{{\mathcal P}_1} \xrightarrow{ \ \sim \ } {\mathcal E}_{{\mathcal P}_2} \, .
$$
\item De même, une construction invariante sur des sites $({\mathcal C},J)$ ou d'autres types de présentation ${\mathcal P}$ consiste à leur associer des modèles $H({\mathcal C},J)$ ou $H({\mathcal P})$ d'une certaine théorie ${\mathbb H}$, de telle sorte que toute équivalence des topos associés
$$
\widehat{\mathcal C}_J \xrightarrow{ \ \sim \ } \widehat{\mathcal D}_K \quad \mbox{ou} \quad {\mathcal E}_{{\mathcal P}_1} \xrightarrow{ \ \sim \ } {\mathcal E}_{{\mathcal P}_2}
$$
induise des isomorphismes de modèles
$$
H({\mathcal C},J) \xrightarrow{ \ \sim \ } H({\mathcal D} , K) \quad \mbox{ou} \quad H({\mathcal P}_1) \xrightarrow{ \ \sim \ } H({\mathcal P}_2) \, .
$$
\item Une telle construction invariante est dite covariante [resp. contravariante] si tout morphisme des topos associés
$$
\widehat{\mathcal C}_J \longrightarrow \widehat{\mathcal D}_K \quad \mbox{ou} \quad {\mathcal E}_{{\mathcal P}_1} \longrightarrow {\mathcal E}_{{\mathcal P}_2}
$$
induit un morphisme de modèles
$$
H({\mathcal C},J) \longrightarrow H({\mathcal D} , K) \quad \mbox{ou} \quad H({\mathcal P}_1) \longrightarrow H({\mathcal P}_2)
$$
[resp.
$$
H({\mathcal D},K) \longrightarrow H({\mathcal C} , J) \quad \mbox{ou} \quad H({\mathcal P}_2) \longrightarrow H({\mathcal P}_1) \, \mbox{],}
$$
de manière compatible avec la composition des morphismes de topos. 
\end{listeisansmarge}
\end{remarksqed}

La généralisation des espaces topologiques aux topos de la notion de point, de propriétés telles que la connexité ou la connexité locale, ou des invariants cohomologiques et homotopiques, demande chaque fois une reformulation complète très abstraite dans le langage des catégories et des foncteurs.

En sens inverse, les propriétés et les constructions invariantes générales des topos se reformulent souvent en une multitude de termes concrets très différents dans différents types de présentation des topos par des sites $({\mathcal C},J)$ ou par d'autres données ${\mathcal P}$.

C'est ce qui permet de développer la théorie et la technique des ``topos comme ponts'' introduite dans \cite{Bridges} (voir \cite{Memoire} pour une présentation résumée de multiples applications de cette technique).

\medskip

Une première classe très large d'équivalences entre des topos associés à des sites est fournie par le ``lemme de comparaison'' de Grothendieck:

\begin{lem}\label{lem1.1.5}
Soit un site $({\mathcal C},J)$.

Soit une sous-catégorie pleine ${\mathcal C}'$ de ${\mathcal C}$ qui est ``dense'' au sens que tout objet $X$ de ${\mathcal C}$ admet un $J$-recouvrement par des morphismes $X_i \to X$, $i \in I$, dont les sources $X_i$ sont des objets de ${\mathcal C}'$.

Soit $J'$ la topologie de ${\mathcal C}'$ induite par $J$ au sens qu'une famille de morphismes de ${\mathcal C}'$ est $J'$-couvrante si et seulement si elle est $J$-couvrante en tant que famille de morphismes de ${\mathcal C}$.

Alors le foncteur de restriction à ${\mathcal C}'$ des préfaisceaux sur ${\mathcal C}$
$$
\widehat{\mathcal C} \longrightarrow \widehat{\mathcal C}'
$$
définit une équivalence de topos
$$
\widehat{\mathcal C}_J \xrightarrow{ \ \sim \ } \widehat{{\mathcal C}'}_{\!\!J'} \, .
$$
\end{lem}

\begin{refqed}

Voir le § III.4 de \cite{SGA4tome1}. 

\end{refqed}

Dans le cas des espaces topologiques, ce lemme devient:

\begin{corqed}\label{cor1.1.6}
Soient un espace topologique $X$ et la catégorie ${\mathcal C}_X$ de ses ouverts.

Soit une famille ${\mathcal C}$ d'ouverts de $X$ qui est dense au sens que tout ouvert de $X$ est réunion d'éléments de ${\mathcal C}$.

Alors, si ${\mathcal C}$ est considérée comme une sous-catégorie pleine de ${\mathcal C}_X$ munie de la topologie $J$ induite par la topologie canonique $J_X$ de ${\mathcal C}_X$, le foncteur de restriction des préfaisceaux induit une équivalence de topos
$$
\widehat{({\mathcal C}_X)}_{J_X} \xrightarrow{ \ \sim \ } \widehat{\mathcal C}_J \, .
$$
En particulier, les topos $\widehat{\mathcal C}_J$ ont tous les mêmes points et les mêmes invariants.
\end{corqed}

\begin{exqed}
L'ensemble ${\mathbb R}$ des nombres réels peut être introduit comme ensemble des points des topos de la forme
$$
\widehat{({\mathcal C}_Q)}_{J_Q}
$$
où

\smallskip

\noindent $\left\{\begin{matrix}
\bullet &\mbox{$Q$ est n'importe quelle partie dense de ${\mathbb Q}$,} \hfill \\
\bullet &\mbox{${\mathcal C}_Q$ est la catégorie des intervalles ouverts $]m,M[$ dont les bornes $m < M$} \hfill \\
&\mbox{sont éléments de $Q$, reliés par les inclusions entre intervalles,} \hfill \\
\bullet &\mbox{une famille de sous-intervalles $]m_i,M_i[ \ \subseteq \ ]m,M[$, $i \in I$,} \hfill \\
&\mbox{est couvrante pour la topologie $J_Q$ si et seulement si,} \hfill \\
&\mbox{pour tous éléments $m',M' \in Q$ vérifiant $m < m' <M' < M$,} \hfill \\
&\mbox{elle contient une sous-famille finie $]m_{i_k} , M_{i_k}[$, $1 \leq k \leq n$,} \hfill \\
&\mbox{telle que $m_{i_1} < m'$, $M' < m_{i_n}$ \quad et \quad $m_{i_{k+1}} < M_{i_k} < M_{i_{k+1}}$, $\forall \, k \in \{1,\cdots , n-1\}$.} \hfill
\end{matrix}\right.$. 
\end{exqed}

\subsection{Les topos comme ``pastiches'' de la catégorie des ensembles}\label{ssec1.1.3}

Grothendieck écrit que les topos sont des ``pastiches'' de la catégorie ${\rm Ens}$ des ensembles car ils vérifient les mêmes propriétés catégoriques constructives:

\begin{thm}[\textbf{SGA 4, tome 1}]\label{thm1.1.7}
Soit ${\mathcal E}$ un topos. Alors:

\smallskip
\begin{enumerate}[label=(\arabic*)]
\item ${\mathcal E}$ est une catégorie localement petite.

\item Les colimites arbitraires sont bien définies dans ${\mathcal E}$. 

Plus précisément, pour tout carquois $D$, le foncteur qui associe à tout objet $E$ de ${\mathcal E}$ le $D$-diagramme constant égal à E
$$
\Delta : {\mathcal E} \longmapsto D\mbox{-{\rm diag}} \, ({\mathcal E})
$$
admet un adjoint à gauche
$$
\begin{matrix}
\underset{D}{\varinjlim} \, : D\mbox{-{\rm diag}} \, ({\mathcal E}) &\longrightarrow &{\mathcal E} \hfill \\
\hfill E_{\bullet} &\longmapsto &\underset{D}{\varinjlim} \ E_{\bullet} \, .
\end{matrix}
$$
\item Les limites arbitraires sont bien définies dans ${\mathcal E}$. 

Plus précisément, pour tout carquois $D$, le foncteur des $D$-diagrammes constants
$$
\Delta : {\mathcal E} \longrightarrow D\mbox{-{\rm diag}} \, ({\mathcal E})
$$
admet un adjoint à droite
$$
\varprojlim_D : D\mbox{-{\rm diag}} \, ({\mathcal E}) \longrightarrow {\mathcal E} \, .
$$
\item Pour tout morphisme $e : E_2 \to E_1$, le foncteur de changement de base
$$
\begin{matrix}
\hfill {\mathcal E} / E_1 &\longrightarrow &{\mathcal E}/E_2 \, , \hfill \\
(E\to E_1) &\longmapsto &(E \times_{E_1} E_2 \to E_2)
\end{matrix}
$$
respecte les colimites arbitraires.

\item Les foncteurs de colimites filtrantes, c'est-à-dire indexées par des petites catégories filtrantes ${\mathcal D}$,
$$
\varinjlim_{\mathcal D} : [{\mathcal D},{\mathcal E}] \longrightarrow {\mathcal E}
$$
respectent les limites finies.

\item Les sommes sont disjointes, au sens que pour toute famille $(E_i)_{i \in I}$ d'objets de ${\mathcal E}$ et pour tous indices $i_1 \ne i_2$ de $I$, le produit fibré
$$
E_{i_1} \times_{\underset{i \in I}{\coprod} E_i} E_{i_2}
$$
s'identifie à l'objet initial $\emptyset$ de ${\mathcal E}$.

\item Un morphisme de ${\mathcal E}$ est un isomorphisme si (et seulement si) il est à la fois un monomorphisme et un épimorphisme.

\item Pour tout objet $E$ de ${\mathcal E}$, il y a correspondance bijective entre les quotients de $E$, c'est-à-dire les classes d'équivalence d'épimorphismes $E \twoheadrightarrow \overline E$, et les relations d'équivalence c'est-à-dire les sous-objets
$$
R \xymatrix{\ar@{^{(}->}[r] & } E \times E
$$
qui contiennent la diagonale $E \hookrightarrow E \times E$, sont invariants par permutation des deux facteurs $E$ et sont stables par composition au sens que $R \times_E R \to E \times E$ se factorise à travers $R \hookrightarrow E \times E$.

De plus, quand $E \to \overline E$ et $R \hookrightarrow E \times E$ se correspondent, le carré
$$
\xymatrix{
R \ar[d] \ar[r] &E \ar[d] \\
E \ar[r] &\overline E
}
$$
est à la fois cartésien et cocartésien.

\item Les sous-objets d'un objet $E$ de ${\mathcal E}$ forment un ensemble.

\item Les objets quotients d'un objet $E$ de ${\mathcal E}$ forment un ensemble.

\item La catégorie ${\mathcal E}$ possède des petites familles ${\mathcal C}$ d'objets qui sont ``séparantes'' au sens que tout objet $E$ de ${\mathcal E}$ admet des recouvrements par des familles globalement épimorphiques
$$
(E_i \longrightarrow E)_{i \in I}
$$
dont les sources $E_i$, $i \in I$, sont des éléments de ${\mathcal C}$.
\end{enumerate}
\end{thm}

\begin{remarks}
	\begin{listeisansmarge}
\item Giraud a montré que, réciproquement, une catégorie ${\mathcal E}$ qui ressemble à ${\rm Ens}$ au sens qu'elle possède les propriétés (0), (1), (2), (3), (5), (7) et (10) est nécessairement un topos. En fait, il suffit même de supposer que, au lieu de (2), la catégorie ${\mathcal E}$ possède la propriété (2') que les limites finies soient bien définies dans ${\mathcal E}$.

\item Si ${\mathcal E}$ est un topos et ${\mathcal C}$ n'importe quelle petite sous-catégorie pleine de ${\mathcal E}$ dont les objets forment une famille séparante de ${\mathcal E}$, le foncteur
$$
\begin{matrix}
{\mathcal E} &\longrightarrow &\widehat{\mathcal C} \, , \hfill \\
E &\longmapsto &\left[\begin{matrix} 
X &\mapsto &{\rm Hom} (X,E) \\
^{\mbox{\rotatebox{-90}{$\in$}}} \\ 
{\rm ob}({\mathcal C})
\end{matrix} \right]
\end{matrix}
$$
est pleinement fidèle et définit une équivalence de ${\mathcal E}$ vers la sous-catégorie $\widehat{\mathcal C}_J$ des $J$-faisceaux, si $J$ désigne la topologie pour laquelle une famille de morphismes de ${\mathcal C}$
$$
(X_i \longrightarrow X)_{i \in I}
$$
est couvrante si et seulement si elle est globalement épimorphique dans ${\mathcal E}$.
\end{listeisansmarge}
\end{remarks}

\noindent {\bf Principes de démonstration:}

\smallskip

Toutes les propriétés considérées sont  invariantes par équivalence de catégories. Il suffit donc de considérer le cas où
$$
{\mathcal E} = \widehat{\mathcal C}_J
$$
est la catégorie des $J$-faisceaux sur une petite catégorie ${\mathcal C}$ munie d'une topologie $J$.

La catégorie $\widehat{\mathcal C}$ est localement petite, donc aussi sa sous-catégorie pleine $\widehat{\mathcal C}_J$, ce qui prouve la propriété~(0).

Dans la catégorie ${\rm Ens}$ des ensembles, les foncteurs de colimites et de limites arbitraires sont définis par des formules explicites, qui permettent de vérifier que ${\rm Ens}$ possède les propriétés (1) à (9).

Celles-ci s'étendent à la catégorie de préfaisceaux $\widehat{\mathcal C}$ grâce au lemme suivant:

\begin{lemqed}\label{lem1.1.8}
Soit ${\mathcal C}$ une catégorie essentiellement petite.

Soient $D$ un carquois et $P_{\bullet}$ un $D$-diagramme de préfaisceaux sur ${\mathcal C}$.

Alors les préfaisceaux sur ${\mathcal C}$
$$
X \longmapsto \varinjlim_D \, P_{\bullet} (X) \qquad \mbox{et} \qquad X \longmapsto \varprojlim_D \, P_{\bullet} (X)
$$
constituent respectivement une colimite et une limite du $D$-diagramme $P_{\bullet}$ dans la catégorie $\widehat{\mathcal C}$ des préfaisceaux. 
\end{lemqed} 

La propriété (10) du théorème est vérifiée par la catégorie $\widehat{\mathcal C}$ car, d'après le lemme de Yoneda rappelé ci-dessous, le foncteur de Yoneda
$$
y : {\mathcal C} \longrightarrow \widehat{\mathcal C}
$$
transforme la famille des objets de ${\mathcal C}$ en une famille séparante de $\widehat{\mathcal C}$.

\begin{lemqed}[\textbf{Yoneda}]\label{lem1.1.9}  
Soit ${\mathcal C}$ une catégorie localement petite munie de son foncteur de Yoneda
$$
\begin{matrix}
y : {\mathcal C} &\longrightarrow &\widehat{\mathcal C} = [{\mathcal C}^{\rm op} , {\rm Ens}] \, , \\
\hfill X &\longmapsto &{\rm Hom} (\bullet,X) \, . \hfill
\end{matrix}
$$

Alors, pour tout objet $P$ de $\widehat{\mathcal C}$, on a une bijection naturelle
$$
P(X) \xrightarrow{ \ \sim \ } {\rm Hom} (y(X),P) \, .
$$
\end{lemqed}

Enfin, les propriétés (1) à (10) passent de $\widehat{\mathcal C}$ à $\widehat{\mathcal C}_J$ grâce à la proposition suivante qui rappelle la structure essentielle reliant préfaisceaux et faisceaux:

\begin{prop}\label{prop1.1.10}

Soit ${\mathcal C}$ une catégorie essentiellement petite munie d'une topologie $J$.

Alors:
\begin{listeimarge}
\item Le foncteur pleinement fidèle de plongement des faisceaux dans les préfaisceaux
$$
j_* : \widehat{\mathcal C}_J \longrightarrow \widehat{\mathcal C}
$$
possède un adjoint à gauche entièrement explicitable
$$
j^* : \widehat{\mathcal C} \longrightarrow \widehat{\mathcal C}_J
$$
appelé le foncteur de faisceautisation.

\item La transformation naturelle canonique
$$
j^* \circ j_*  \longrightarrow {\rm id}_{\widehat{\mathcal C}_J}
$$
est un isomorphisme.

\item Le foncteur de faisceautisation $j^*$ respecte non seulement les colimites arbitraires mais aussi les limites finies.

\item La colimite dans $\widehat{\mathcal C}_J$ d'un diagramme de faisceaux se calcule en appliquant le foncteur $j^*$ à sa colimite calculée dans $\widehat{\mathcal C}$ par évaluation en chaque objet $X$ de ${\mathcal C}$.

\item La limite dans $\widehat{\mathcal C}$ d'un diagramme de faisceaux est encore un faisceau, et il est limite du diagramme dans $\widehat{\mathcal C}_J$.

Autrement dit, les limites de diagrammes de faisceaux se calculent comme celles des diagrammes de préfaisceaux, par évaluation en chaque objet $X$ de ${\mathcal C}$.
\end{listeimarge}
\end{prop}

\begin{remarksqed}
		\begin{listeisansmarge}
\item Il existe plusieurs manières de construire le foncteur de faisceautisation $j^*$ par des formules explicites.

L'une consiste à introduire le foncteur
$$
\begin{matrix}
\widehat{\mathcal C} &\longrightarrow &\widehat{\mathcal C} \, , \\
P &\longmapsto &P^+
\end{matrix}
$$
où, pour tout objet $X$ de ${\mathcal C}$,
$$
P^+(X) = \varinjlim_{C \in J(X)} {\rm Hom} (C,P) \, ,
$$
en notant $J(X)$ l'ensemble ordonné filtrant des cribles $J$-couvrants $C$ de $X$, c'est-à-dire des sous-objets dans $\widehat C$
$$
C \xymatrix{\ar@{^{(}->}[r] & } y(X) = {\rm Hom} (\bullet , X)
$$
qui contiennent des familles $J$-couvrantes
$$
(X_i \longrightarrow X)_{i \in I} \, .
$$
Alors le foncteur $j^*$ est construit comme le foncteur composé
$$
P \longmapsto j^* P = (P^+)^+ \, .
$$
\item Le foncteur de faisceautisation $j^*$ se construit aussi comme le foncteur composé
$$
P \longmapsto (P_s)^+
$$
du foncteur $P \mapsto P^+$ et du foncteur de $J$-séparation
$$
P \longmapsto P_s \, .
$$
Un préfaisceau $P$ sur ${\mathcal C}$ est dit $J$-séparé si, pour toute famille $J$-couvrante
$$
\left( X_i \xrightarrow{ \ x_i \ } X \right)_{i \in I}
$$
et toutes sections $p_1 , p_2 \in P(X)$, on a
$$
p_1 = p_2
$$
si $x_i^* (p_1) = x_i^* (p_2)$ dans $P(X_i)$ pour tout $i \in I$.

Le foncteur de plongement dans $\widehat{\mathcal C}_J$ de la sous-catégorie pleine des préfaisceaux $J$-séparés admet un adjoint à gauche, le foncteur de $J$-séparation
$$
P \longmapsto P_s
$$
où, pour tout objet $X$ de ${\mathcal C}$, $P_s (X)$ est l'ensemble des classes d'équivalence d'éléments de $P(X)$ pour la relation de co{\"\i}ncidence locale.

Deux sections $p_1 , p_2 \in P(X)$ co{\"\i}ncident localement s'il existe une famille $J$-couvrante $\left( X_i \xrightarrow{ \ x_i \ } X \right)_{i \in I}$, telle que
$$
x_i^* (p_1) = x_i^* (p_2) \, , \quad \forall \, i \in I \, .
$$
\item Concrètement, pour tout préfaisceau $P$ sur ${\mathcal C}$, les sections en un objet $X$ de ${\mathcal C}$ de
$$
j^* P = (P_s)^+
$$
sont les classes d'équivalence locale des familles de sections locales de $P$ sur $X$ qui sont localement cohérentes.

Une famille de sections locales de $P$ sur $X$ est une famille de sections de $P$
$$
p_i \in P(X_i) \, , \quad i \in I \, ,
$$
en les sources $X_i$ d'une famille $J$-couvrante de morphismes $\left( X_i \xrightarrow{ \ x_i \ } X \right)_{i \in I}$.

Une telle famille est dite localement cohérente si, pour tous indices $i,j \in I$ et tout carré commutatif de ${\mathcal C}$
$$
\xymatrix{
U \ar[d] \ar[r] &X_i \ar[d]^{x_i} \\
X_j \ar[r]_{x_j} &X
}
$$
il existe une famille $J$-couvrante $\left( U_k \xrightarrow{ \ u_k \ } U \right)_{k \in K}$ telle que $p_i \in P(X_i)$ et $p_j \in P(X_j)$ aient même image dans chaque $P(U_k)$, $k \in K$.

Deux familles localement cohérentes de sections locales
$$
\begin{matrix}
\hfill p_i &\in &P(X_i)_{i \in I} \, , \hfill \\
p'_{i'} &\in &P(X'_{i'})_{i' \in I'} \, ,
\end{matrix}
$$
définies sur deux recouvrements $\left( X_i \xrightarrow{ \ x_i \ } X \right)_{i \in I}$ et $\left( X'_{i'} \xrightarrow{ \ x'_{i'} \ } X \right)_{i' \in I'}$ sont dites localement équivalentes si, pour tous indices $i \in I$, $i' \in I'$, et tout carré commutatif de ${\mathcal C}$
$$
\xymatrix{
U \ar[d] \ar[r] &X_i \ar[d]^{x_i} \\
X'_{i'} \ar[r]_{x'_{i'}} &X
}
$$
il existe un recouvrement $\left( U_k \xrightarrow{ \ u_k \ } U \right)_{k \in K}$ tel que $p_i \in P(X_i)$ et $p'_{i'} \in P(X'_{i'})$ aient même image dans chaque $P(U_k)$, $k \in K$. 
\end{listeisansmarge}
\end{remarksqed}

Le composé du foncteur de Yoneda et du foncteur de faisceautisation joue un rôle important:

\begin{defn}\label{defn1.1.11}	
Pour tout site $({\mathcal C},J)$, on appelle ``foncteur canonique'' de ${\mathcal C}$ dans le topos $\widehat{\mathcal C}_J$ le foncteur composé
$$
\ell = j^* \circ y : {\mathcal C} \xymatrix{\ar@{^{(}->}[r] & } \widehat{\mathcal C} \longrightarrow \widehat{\mathcal C}_J \, .
$$
\end{defn}

\begin{remarksqed}
\noindent (i) Si ${\mathcal C}$ est une petite catégorie, la famille des images $\ell (X)$ par $\ell$ de ses objets $X$ est une famille ``séparante'' dans le topos $\widehat{\mathcal C}_J$.

\noindent (ii) On montre que le foncteur $\ell : {\mathcal C} \to \widehat{\mathcal C}_J$ est pleinement fidèle si et seulement si, pour tout objet $X$ de ${\mathcal C}$, le préfaisceau
$$
y(X) = {\rm Hom} (\bullet , X) \qquad \mbox{est un $J$-faisceau.}
$$
On dit dans ce cas que la topologie $J$ est sous-canonique. 
\end{remarksqed}

Non seulement les topos possèdent des foncteurs de limites et colimites arbitraires qui ont les mêmes propriétés que ceux de la catégorie des ensembles, mais les constructions d'ordre supérieur $y$ sont bien définies comme dans le cas des ensembles:

\begin{thm}\label{thm1.1.12}
		
Soit ${\mathcal E}$ un topos.
\begin{listeimarge}
\item Le foncteur contravariant qui associe à tout objet $E$ de ${\mathcal E}$ l'ensemble de ses sous-objets est représentable par un objet $\Omega$ de ${\mathcal E}$ appelé le ``classificateur des sous-objets''.

\item Pour tout objet $E_1$ de ${\mathcal E}$, le foncteur de produit avec $E_1$ dans ${\mathcal E}$
$$
E_1 \times \bullet : {\mathcal E} \longrightarrow {\mathcal E}
$$
admet un adjoint à droite
$$
E \longmapsto E^{E_1} \, .
$$
\item Pour tout morphisme $e : E_2 \to E_1$ de ${\mathcal E}$, le foncteur de changement de base entre catégories relatives
$$
\bullet \times_{E_1} E_2 : {\mathcal E}/E_1 \longrightarrow {\mathcal E}/E_2 \, ,
$$
admet un adjoint à droite
$$
\begin{matrix}
\hfill {\mathcal E}/E_2 &\longrightarrow &{\mathcal E}/E_1 \hfill \\
(E \to E_2) &\longmapsto &(\Pi_e \, E \to E_1) \, .
\end{matrix}
$$
\end{listeimarge}
\end{thm}

\begin{remarks}
	\begin{listeisansmarge}
\item Le classificateur des sous-objets $\Omega$ est muni d'un monomorphisme de source l'objet terminal $1$ de ${\mathcal E}$
$$
1 \xymatrix{\ar@{^{(}->}[r] & } \Omega
$$
tel que, pour tout objet $E$ de ${\mathcal E}$, l'application
$$
\left( E \xrightarrow{ \ \chi \ } \Omega \right) \longmapsto \left( 1 \times_{\Omega} E \xymatrix{\ar@{^{(}->}[r] & } E \right)
$$
définisse une bijection de ${\rm Hom} (E,\Omega)$ sur l'ensemble des sous-objets de ${\mathcal E}$.

Dans le cas où ${\mathcal E} = {\rm Ens}$, l'objet $\Omega$ est l'ensemble à deux éléments $\{0,1\}$ muni de l'injection $\{1\} \hookrightarrow \{0,1\}$. Pour tout sous-ensemble d'un ensemble $E$, l'application
$$
\chi : E \longrightarrow \{0,1\} = \Omega
$$
qui lui correspond est sa fonction caractéristique.

Si ${\mathcal E} = \widehat{\mathcal C}$ est un topos de préfaisceaux, $\Omega$ est le préfaisceau qui associe à tout objet $X$ de ${\mathcal C}$ l'ensemble des sous-objets
$$
C \xymatrix{\ar@{^{(}->}[r] & } y(X)
$$
c'est-à-dire l'ensemble $\Omega (X)$ des cribles $C$ de l'objet $X$.

Si ${\mathcal E} = \widehat{\mathcal C}_J$ est le topos des faisceaux sur un site $({\mathcal C},J)$ muni du foncteur canonique $\ell = j^* \circ y : {\mathcal C} \to \widehat{\mathcal C}_J$, $\Omega$ est le faisceau qui associe à tout objet $X$ de ${\mathcal C}$ l'ensemble des sous-objets de $\ell (X)$ dans $\widehat{\mathcal C}_J$.

Il s'identifie à l'ensemble des sous-objets
$$
C\xymatrix{\ar@{^{(}->}[r] & } y(X)
$$
tels que
$$
C = j_* \, j^* C \times_{j_* j^* y(X)} y(X) \, ,
$$
c'est-à-dire à l'ensemble des cribles $C$ de $X$ qui sont $J$-fermés au sens qu'un morphisme
$$
X' \xrightarrow{ \ x \ } X
$$
est dans $C$ s'il existe une famille $J$-couvrante $\left( X'_i \xrightarrow{ \, x_i \, } X' \right)_{i \in I}$ telle que
$$
(x \circ x_i : X'_i \to X' \to X) \in C \, , \quad \forall \, i \in I \, .
$$
\item Si ${\mathcal E} = \widehat{\mathcal C}_J$ est un topos de faisceaux muni du foncteur canonique $\ell = j^* \circ y : {\mathcal C} \to \widehat{\mathcal C}_J$, les exponentielles sont encore notées
$$
E^{E_1} = {\mathcal H}om (E_1,E) \, .
$$
Le faisceau $E^{E_1} = {\mathcal H}om (E_1,E)$ associe en effet à tout objet $X$ de ${\mathcal C}$ l'ensemble
$$
\begin{matrix}
{\mathcal H}om (E_1,E)(X) &= &{\rm Hom}_{\widehat{\mathcal C}_J} (E_1 \times \ell (X) , E) \hfill \\
&= &{\rm Hom}_{\widehat{\mathcal C}} (E_1 \times y(X),E) \hfill \\
&= &{\rm Hom}_{\widehat{\mathcal C}/y(X)} (E_1 \times y(X), E \times y(X)) \hfill
\end{matrix}
$$
qui est l'ensemble des morphismes
$$
E_1 \bigl\vert_{{\mathcal C}/X} \longrightarrow E_2 \bigl\vert_{{\mathcal C}/X} 
$$
entre les restrictions des faisceaux $E_1$ et $E_2$ à la catégorie relative ${\mathcal C}/X$.

\item Les foncteurs
$$
{\prod}_e : {\mathcal E} / E_2 \longrightarrow {\mathcal E}/E_1
$$ 
associés aux morphismes $e : E_2 \to E_1$ des topos ${\mathcal E}$ sont appelés les ``produits dépendants''.

En effet, si ${\mathcal E} = {\rm Ens}$ est le topos des ensembles et $A$ est un objet de ce topos, la catégorie relative
$$
{\rm Ens}/A
$$
s'identifie à celle des familles d'ensembles $E_a$ indexés par les éléments $a \in A$, et le foncteur $\prod_p$ associé à l'unique application $p : A \to \{\bullet\}$ est
$$
\begin{matrix}
(E_a)_{a \in A} &\longmapsto &\underset{a \in A}{\prod} \, E_a \, , \\
{\rm Ens} / A &\longrightarrow &{\rm Ens} \, . \hfill
\end{matrix}
$$
\end{listeisansmarge}
\end{remarks}

\begin{principedemo}
	
Ce théorème est conséquence de la proposition suivante qui montre que les topos, qui sont construits par complétion de petites catégories, sont effectivement complets:
\end{principedemo}

\begin{prop}\label{prop1.1.13}

Considérons un préfaisceau sur un topos ${\mathcal E}$
$$
P : {\mathcal E}^{\rm op} \longrightarrow {\rm Ens} \, .
$$
Alors:
\begin{listeimarge}
\item Les propriétés suivantes sont équivalentes:
\begin{enumerate}[label=(\arabic*)]
\item Le préfaisceau $P$ est un faisceau pour la topologie canonique de ${\mathcal E}$.
\item Il transforme les colimites arbitraires dans ${\mathcal E}$
$$
\varinjlim_D E_{\bullet}
$$
en limites
$$
\varprojlim_D P(E_{\bullet}) \, .
$$
\item[(3)] Il est représentable par un objet $E$ de ${\mathcal E}$, c'est-à-dire isomorphe à un préfaisceau de la forme
$$
{\rm Hom} (\bullet , E) \, .
$$
\end{enumerate}

\item Si ${\mathcal E} = \widehat{\mathcal C}_J$ est le topos des faisceaux sur un site $({\mathcal C},J)$, avec le foncteur canonique $\ell : {\mathcal C} \to \widehat{\mathcal C}_J$, et que
$$
P : {\mathcal E}^{\rm op} \longrightarrow {\rm Ens}
$$
est un faisceau pour la topologie canonique de ${\mathcal E} = \widehat{\mathcal C}_J$, alors le préfaisceau
$$
\begin{matrix}
{\mathcal C}^{\rm op} &\longrightarrow &{\rm Ens} \, , \hfill \\
\hfill X &\longmapsto &P(\ell (X))
\end{matrix}
$$
est un $J$-faisceau et il représente le préfaisceau $P$ sur ${\mathcal E}$.
\end{listeimarge}
\end{prop}

\begin{principedemosansqed}
	Elle repose sur le lemme suivant:
\end{principedemosansqed}

\begin{lem}\label{lem1.1.14}

Soit une petite catégorie ${\mathcal C}$ munie d'une topologie $J$ et du foncteur canonique $\ell = j^* \circ y : {\mathcal C} \to \widehat{\mathcal C}_J$.

Alors tout objet $F$ de $\widehat{\mathcal C}_J$ s'écrit comme la colimite
$$
F = \varinjlim_{(X,x) \in {\mathcal C}/F} \ell (X)
$$
où ${\mathcal C}/F$ désigne la catégorie des éléments de $F$ c'est-à-dire des paires $(X,x)$ constituées d'un objet $X$ de ${\mathcal C}$ et d'un élément $x \in F(X) = {\rm Hom} (\ell (X),F)$.
\end{lem}

\begin{demo}
Le cas où $J$ est la topologie discrète résulte du lemme de Yoneda.

Le cas général s'en déduit puisque $j^* \circ j_* \xrightarrow{ \, \sim \, } {\rm id}_{\widehat{\mathcal C}_J}$ est un isomorphisme et que $j^* : \widehat{\mathcal C} \to \widehat{\mathcal C}_J$ respecte les colimites. 
\end{demo}

\subsection{Les topos comme lieux d'expression sémantique des théories}\label{ssec1.1.4}

Rappelons les types d'éléments de vocabulaire des théories du premier ordre, et la manière dont ils servent à nommer des éléments de structure dans des topos:

\begin{defn}\label{defn1.1.15}
	\begin{listeimarge}
\item Un vocabulaire du premier ordre, appelé ``signature'', consiste en

\smallskip

$\left\{\begin{matrix}
\bullet &\mbox{une liste de ``noms d'objets'', appelés les ``sortes'',} \hfill \\
\bullet &\mbox{une liste de ``noms d'opérations'', appelés les ``symboles de fonctions'',} \hfill \\
&f : A_1 \cdots A_n \longrightarrow B \\
&\mbox{qui vont d'une famille finie (éventuellement vide) de sortes $A_1 \cdots A_n$ à une sorte $B$,} \hfill \\
\bullet &\mbox{une liste de ``noms de relations'', appelés les ``symboles de relations'',} \hfill \\
&\xymatrix{R \ \ \ar@{>->}[r] &A_1 \cdots A_n} \\
&\mbox{reliant une famille finie (éventuellement vide) de sortes $A_1 \cdots A_n$.} \hfill
\end{matrix}\right.$

%\medskip

\item Etant donné un tel vocabulaire $\Sigma$, une structure de type $\Sigma$ ou $\Sigma$-structure $M$ dans un topos ${\mathcal E}$ consiste en

\smallskip

$\left\{\begin{matrix}
\bullet &\mbox{des objets $M\!A$ de ${\mathcal E}$ indexés par les sortes $A$ de $\Sigma$,} \hfill \\
\bullet &\mbox{des morphismes $M\!f : M\!A_1 \times \cdots \times M\!A_n \to M\!B$ de ${\mathcal E}$} \hfill \\
&\mbox{indexés par les symboles de fonctions $f : A_1 \cdots A_n \to B$ de $\Sigma$} \hfill \\
&\mbox{(et $M\!f : 1 \to M\!B$ défini sur l'objet terminal $1$ de ${\mathcal E}$ si $f : \, \to B$ a pour source la famille vide),} \hfill \\
\bullet &\mbox{des sous-objets $M\!R \hookrightarrow M\!A_1 \times \cdots \times M\!A_n$ de ${\mathcal E}$} \hfill \\
&\mbox{indexés par les symboles de relation $R \rightarrowtail A_1 \cdots A_n$} \hfill \\
&\mbox{(et $M\!R \hookrightarrow 1$ si $R \rightarrowtail$ \quad a pour contexte la famille vide).} \hfill
\end{matrix}\right.$

%\smallskip

\item Un morphisme $u$ entre deux structures $M$ et $N$ de type $\Sigma$ dans un topos ${\mathcal E}$ consiste en une famille de morphismes de ${\mathcal E}$
$$
u_A : M\!A \longrightarrow N\!A
$$
indexés par les sortes $A$ de $\Sigma$, telle que:
\begin{enumerate}
\item[$\bullet$] pour tout symbole de fonction $f : A_1 \cdots A_n \to B$ de $\Sigma$, le carré
$$
\xymatrix{
M\!A_1 \times \cdots \times M\!A_n \ar[d]_{u_{A_1} \times \cdots \times u_{A_n}} \ar[r]^-{M\!f} &M\!B \ar[d]^{u_B} \\
N\!A_1 \times \cdots \times N\!A_n \ar[r]^-{N\!f} &N\!B
}
$$
est commutatif,
\item[$\bullet$] pour tout symbole de relation $R \rightarrowtail A_1 \cdots A_n$ de $\Sigma$, 

\noindent le morphisme $u_{A_1} \times \cdots \times u_{A_n} : M\!A_1 \times \cdots \times M\!A_n \to N\!A_1 \times \cdots \times N\!A_n$ 

\noindent envoie le sous-objet \quad $M\!R \hookrightarrow M\!A_1 \times \cdots \times M\!A_n$ 

\noindent dans le sous-objet \quad \ \  $N\!R \hookrightarrow N\!A_1 \times \cdots \times N\!A_n$.
\end{enumerate}
\end{listeimarge}
\end{defn}

On observe:

\begin{lem}\label{lem1.1.16}

Soit $\Sigma$ une signature du premier ordre. Alors:
	\begin{listeimarge}
\item Pour tout topos ${\mathcal E}$, les $\Sigma$-structures dans ${\mathcal E}$ et leurs morphismes forment une catégorie localement petite que l'on peut noter
$$
\Sigma\mbox{{\rm -str}} \, ({\mathcal E}) \qquad \mbox{ou} \qquad \Sigma\mbox{{\rm -mod}} \, ({\mathcal E}) \, .
$$
Elle possède des limites arbitraires et des colimites filtrantes arbitraires.

\item Pour tout morphisme de topos $f : {\mathcal E}' \to {\mathcal E}$, ses deux composantes
$$
\left( {\mathcal E} \xrightarrow{ \ f^* \ } {\mathcal E}' , {\mathcal E}' \xrightarrow{ \ f_* \ } {\mathcal E} \right)
$$
définissent deux foncteurs adjoints
$$
\left( f^* : \Sigma\mbox{{\rm -str}} \, ({\mathcal E}) \longrightarrow \Sigma\mbox{{\rm -str}} \, ({\mathcal E}'), \ f_* : \Sigma\mbox{{\rm -str}} \, ({\mathcal E}') \longrightarrow \Sigma\mbox{{\rm -str}} \, ({\mathcal E}) \right).
$$
Le foncteur $f_*$ respecte les limites arbitraires, et son adjoint à gauche $f^*$ respecte les limites finies et les colimites filtrantes arbitraires.
\end{listeimarge}
\end{lem}

\begin{esquissedemo}
		\begin{listeisansmarge}
\item résulte de ce que, dans un topos ${\mathcal E}$, les foncteurs de colimites filtrantes, tout comme les foncteurs de limites, respectent les limites finies, donc les produits finis et les sous-objets.

\item résulte de ce que le foncteur $f_* : {\mathcal E}' \to {\mathcal E}$ respecte les limites arbitraires tandis que son adjoint à gauche $f^* : {\mathcal E} \to {\mathcal E}'$ respecte les colimites arbitraires et les limites finies.
\end{listeisansmarge}
\end{esquissedemo}

Une théorie géométrique du premier ordre ${\mathbb T}$ consiste en une signature $\Sigma$ et une famille (éventuellement vide) d'axiomes qui ont la forme de ``séquents'', c'est-à-dire d'implications
$$
\varphi (\vec x) \vdash \psi (\vec x)
$$
reliant des ``formules géométriques'' en une famille finie $\vec x = (x_1^{A_1} , \cdots , x_n^{A_n})$ de variables formelles $x_i^{A_i}$, $1 \leq i \leq n$, associées chacune à une sorte $A_i$ de $\Sigma$.

Pour toute structure $M$ de type $\Sigma$ dans un topos ${\mathcal E}$, les ``formules géométriques'' $\varphi (x_1^{A_1} , \cdots , x_n^{A_n})$, $\psi (x_1^{A_1} , \cdots , x_n^{A_n})$ écrites dans le vocabulaire de $\Sigma$ définissent des sous-objets
$$
\begin{matrix}
M \varphi &\xymatrix{\ar@{^{(}->}[r] & } &M\!A_1 \times \cdots \times M\!A_n \, ,  \\
M \psi &\xymatrix{\ar@{^{(}->}[r] & } &M\!A_1 \times \cdots \times M\!A_n \, .
\end{matrix}
$$
On dit que la $\Sigma$-structure $M$ satisfait un séquent
$$
\varphi (\vec x) \vdash \psi (\vec x)
$$
si les sous-objets définis par $\varphi$ et $\psi$ satisfont la relation
$$
M \varphi \subseteq M\psi \, .
$$
On dit que $M$ est un modèle de ${\mathbb T}$ si $M$ satisfait tous ses axiomes.

\medskip

Voici comment est définie la notion de formule géométrique dans un signature donnée:

\begin{defn}\label{defn1.1.17}
Soit $\Sigma$ une signature.

	\begin{listeimarge}

\item Un ``terme'' de $\Sigma$ à valeurs dans une sorte $B$ est une expression formelle de la forme
$$
f(x_1^{A_1} , \cdots , x_n^{A_n})
$$
qui peut être 

\smallskip

$\left\{\begin{matrix}
\bullet &\mbox{réduite à une seule variable $x^B$,} \hfill \\
\bullet &\mbox{déduite d'un symbole de fonction de $\Sigma$} \hfill \\
&f : A_1 \cdots A_n \longrightarrow B \\
&\mbox{en l'appliquant à des variables formelles $x_1^{A_1} , \cdots , x_n^{A_n}$ associées} \hfill \\
&\mbox{aux sortes $A_1 , \cdots , A_n$ composant la source de $f$,} \hfill \\
\bullet &\mbox{déduite d'un terme déjà défini} \hfill \\
&g(y_1^{B_1} , \cdots , y_m^{B_m}) \\
&\mbox{en remplaçant l'une de ses variables} \hfill \\
&y_i^{B_i} \\
&\mbox{par un symbole de fonction} \hfill \\
&h : C_1 \cdots C_k \longrightarrow B_i \\
&\mbox{appliqué à des variables formelles $z_1^{C_1} , \cdots , z_k^{C_k}$.} \hfill
\end{matrix}\right.$

%\medskip

\item Une ``formule atomique'' de $\Sigma$
$$
\varphi (x_1^{A_1} , \cdots , x_n^{A_n})
$$
est une expression formelle en des variables $x_1^{A_1} , \cdots , x_n^{A_n}$ associées à des sortes $A_1 , \cdots , A_n$, qui peut être

\smallskip

$\left\{\begin{matrix}
\bullet &\mbox{de la forme} \hfill \\
&R(x_1^{A_1} , \cdots , x_n^{A_n}) \\
&\mbox{pour un symbole de relation $R \rightarrowtail A_1 \cdots A_n$ de $\Sigma$,} \hfill \\
\bullet &\mbox{de la forme} \hfill \\
&x_1^A = x_2^A \\
&\mbox{pour une sorte $A$ de $\Sigma$,} \hfill \\
\bullet &\mbox{déduite d'une formule atomique déjà définie} \hfill \\
&\psi (y_1^{B_1} , \cdots , y_m^{B_m}) \\
&\mbox{en remplaçant une variable $y_i^{B_i}$ par un terme à valeurs dans $B_i$} \hfill \\
&h(z_1^{C_1} , \cdots , z_k^{C_k}).
\end{matrix}\right.$

%\medskip

\item Une ``formule géométrique'' de $\Sigma$
$$
\varphi (x_1^{A_1} , \cdots , x_n^{A_n})
$$
est une expression formelle en des variables $x_1^{A_1} , \cdots , x_n^{A_n}$ associées à des sortes $A_1 , \cdots , A_n$, qui peut être

\smallskip

$\left\{\begin{matrix}
\bullet &\mbox{une ``formule atomique'',} \hfill \\
\bullet &\mbox{la ``formule du vrai'' notée} \hfill \\
&\top (x_1^{A_1} , \cdots , x_n^{A_n}), \\
\bullet &\mbox{une conjonction finie de formules géométriques en les mêmes variables déjà définies} \hfill \\
&\varphi_1 (x_1^{A_1} , \cdots , x_n^{A_n}) \wedge \cdots \wedge \varphi_{\ell} (x_1^{A_1} , \cdots , x_n^{A_n}) \, , \\
\bullet &\mbox{déduite d'une formule géométrique déjà définie} \hfill \\
&\psi (x_1^{A_1} , \cdots , x_n^{A_n} ; y_1^{B_1} , \cdots , y_m^{B_m}) = \psi (\vec x ; \vec y) \\
&\mbox{en appliquant à une partie de ses variable le quantificateur existentiel} \hfill \\
&\exists \, (y_1^{B_1} , \cdots , y_m^{B_m}) = \exists \, \vec y \\
&\mbox{pour définir} \hfill \\
&\varphi (\vec x) = (\exists \, \vec y) \, \psi (\vec x ; \vec y) \, , \\
\bullet &\mbox{la ``formule du faux'' notée} \hfill \\
&\perp (x_1^{A_1} , \cdots , x_n^{A_n}) \, , \\
\bullet &\mbox{une disjonction arbitraire de formules géométriques en les mêmes variables déjà définies} \hfill \\
&\underset{i \in I}{\bigvee} \, \varphi_i (x_1^{A_1} , \cdots , x_n^{A_n}) \, .
\end{matrix}\right.$

%\smallskip

\item Un ``séquent géométrique'' de $\Sigma$ est une implication
$$
\varphi (x_1^{A_1} , \cdots , x_n^{A_n}) \vdash \psi (x_1^{A_1} , \cdots , x_n^{A_n})
$$
reliant deux formules géométriques en les mêmes variables $x_1^{A_1} , \cdots , x_n^{A_n}$.

%\smallskip

\item Une ``théorie géométrique'' ${\mathbb T}$ de signature $\Sigma$ consiste en cette signature complétée par une famille de séquents géométriques de $\Sigma$
$$
\varphi_i (\vec x_i) \vdash \psi_i (\vec x_i) \, , \quad i \in I \, ,
$$
appelés les axiomes de ${\mathbb T}$.
\end{listeimarge}
\end{defn}

\begin{remarksqed}
		\begin{listeisansmarge}
\item Une formule géométrique d'une signature $\Sigma$ est dite ``de Horn'' si sa construction à partir des formules atomiques ne fait appel qu'aux symboles $\top$ et $\wedge$.

\noindent Elle est dite ``régulière'' si sa construction ne fait appel qu'aux symboles $\top$, $\wedge$ et $\exists$.

\noindent Elle est dite ``cohérente'' si sa construction ne fait appel qu'aux symboles $\top$, $\wedge$, $\exists$, $\perp$ et au symbole de disjonction finie $\vee$.

\item Un séquent géométrique est dit ``de Horn'' [resp. ``régulier'', resp. ``cohérent''] s'il relie deux formules ``de Horn'' [resp. ``régulières'', resp. ``cohérentes''].

\item Une théorie géométrique est dite ``de Horn'' [resp. ``régulière'', resp. ``cohérente''] si tous ses axiomes sont ``de Horn'' [resp. ``réguliers'', resp. ``cohérents''].

\item Une théorie géométrique est dite algébrique si sa signature $\Sigma$ consiste seulement en des sortes et des symboles de fonctions, sans symboles de relations, et que ses axiomes ont tous la forme
$$
\top (\vec x) \vdash f(\vec x) = g(\vec x)
$$
où $f(\vec x)$ et $g(\vec x)$ sont des termes en une famille finie de variables $\vec x = (x_1^{A_1} , \cdots , x_n^{A_n})$. 
\end{listeisansmarge}
\end{remarksqed}

Voici d'autre part comment les termes, les formules géométriques et les théories géométriques s'interprètent dans les topos:

\begin{defn}\label{defn1.1.18}
Soient $\Sigma$ une signature et $M$ une structure de type $\Sigma$ dans un topos ${\mathcal E}$.
	\begin{listeimarge}
\item Tout terme de $\Sigma$ à valeurs dans une sorte $B$
$$
f(x_1^{A_1} , \cdots , x_n^{A_n})
$$
s'interprète dans $M$ comme un morphisme de ${\mathcal E}$
$$
M\!f : M\!A_1 \times \cdots \times M\!A_n \longrightarrow M\!B
$$
en suivant les règles suivantes:

\smallskip

$\left\{\begin{matrix}
\bullet &\mbox{S'il est une variable $x^B$, il s'interprète comme} \hfill \\
&{\rm id} : M\!B \longrightarrow M\!B \, . \\
\bullet &\mbox{S'il est associé à un symbole de fonction $f : A_1 \cdots A_n \to B$,} \hfill \\
&\mbox{il s'interprète comme l'élément correspondant de la structure $M$} \hfill \\
&M\!f : M\!A_1 \times \cdots \times M\!A_n \longrightarrow M\!B \, . \\
\bullet &\mbox{S'il est déduit d'un terme déjà défini et interprété} \hfill \\
&\quad g(y_1^{B_1} , \cdots , y_m^{B_m}) \\
&\mbox{en remplaçant l'une de ses variables $y_i^{B_i}$ par un symbole de fonction} \hfill \\
&\quad h : C_1 \cdots C_k \longrightarrow y_i^{B_i} \\
&\mbox{appliqué à des variables $z_1^{C_1} , \cdots , z_k^{C_k}$, son interprétation} \hfill \\
&M\!f \\
&\mbox{est le morphisme composé du morphisme d'interprétation de $g$} \hfill \\
&M\!g : M\!B_1 \times \cdots \times M\!B_m \longrightarrow M\!B \\
&\mbox{et du morphisme produit du morphisme} \hfill \\
&M\!h : M\!C_1 \times \cdots \times M\!C_k \longrightarrow M\!B_i \\
&\mbox{et des morphismes} \hfill \\
&{\rm id} : M\!B_{i'} \longrightarrow M\!B_{i'} \, , \quad i' \in \{1,2,\cdots , m\} - \{i\} \, .
\end{matrix}\right.$

%\medskip

\item Toute formule atomique de $\Sigma$
$$
\varphi (x_1^{A_1} , \cdots , x_n^{A_n})
$$
s'interprète dans $M$ comme un sous-objet
$$
M\varphi \xymatrix{\ar@{^{(}->}[r] & } M\!A_1 \times \cdots \times M\!A_n
$$
en suivant les règles suivantes:

\smallskip

$\left\{\begin{matrix}
\bullet &\mbox{Si elle a la forme $R(x_1^{A_1} , \cdots , x_n^{A_n})$ pour un symbole de relation $R$,} \hfill \\
&\mbox{elle s'interprète comme l'élément de structure} \hfill \\
&M\!R \xymatrix{\ar@{^{(}->}[r] & } M\!A_1 \times \cdots \times M\!A_n \, . \\
\bullet &\mbox{Si elle a la forme $x_1^A = x_2^A$, elle s'interprète comme le sous-objet diagonal} \hfill \\
&M\!A  \xymatrix{\ar@{^{(}->}[r] & } M\!A \times M\!A \, . \\
\bullet &\mbox{Si elle se déduit d'une formule atomique déjà définie et interprétée} \hfill \\
&\psi (y_1^{B_1} , \cdots , y_m^{B_m}) \\
&\mbox{en remplaçant une variable $y_i^{B_i}$ par un terme à valeurs dans $B_i$} \hfill \\
&h(z_1^{C_1} , \cdots , z_k^{C_k}) \, , \\
&\mbox{elle s'interprète comme l'image réciproque du sous-objet} \hfill \\
&M\psi  \xymatrix{\ar@{^{(}->}[r] & } M\!B_1 \times \cdots \times M\!B_m \\
&\mbox{par le morphisme produit du morphisme} \hfill \\
&M\!h : M\!C_1 \times \cdots \times M\!C_k \longrightarrow M\!B_i \\
&\mbox{et des morphismes} \hfill \\
&{\rm id} : M\!B_{i'} \longrightarrow M\!B_{i'} \, , \quad i' \in \{1,\cdots , m\}-\{i\} \, .
\end{matrix}\right.$

%\medskip

\item Toute formule géométrique de $\Sigma$
$$
\varphi (x_1^{A_1} , \cdots , x_n^{A_n})
$$
s'interprète dans $M$ comme un sous-objet
$$
M\varphi \xymatrix{\ar@{^{(}->}[r] & } M\!A_1 \times \cdots \times M\!A_n
$$
en suivant les règles suivantes:

\smallskip

$\left\{\begin{matrix}
\bullet &\mbox{Si $\varphi$ est une formule atomique, son interprétation est celle de (ii).} \hfill \\
\bullet &\mbox{Si $\varphi$ est la formule du vrai $\top (x_1^{A_1} , \cdots , x_n^{A_n})$, son interprétation est le sous-objet maximal} \hfill \\
&M\varphi = M\!A_1 \times \cdots \times M\!A_n \, . \\
\bullet &\mbox{Si $\varphi$ s'écrit comme une conjonction finie de formules séjà interprétées,} \hfill \\
&\varphi = \varphi_1 \wedge \cdots \wedge \varphi_{\ell} \, , \\
&\mbox{son interprétation est l'intersection des sous-objets correspondants} \hfill \\
&M\varphi = M\varphi_1 \cap \cdots \cap M\varphi_{\ell} \, . \\
\bullet &\mbox{Si $\varphi (x_1^{A_1} , \cdots , x_n^{A_n})$ se déduit d'une formule déjà interprétée} \hfill \\
&\psi (x_1^{A_1} , \cdots , x_n^{A_n} ; y_1^{B_1} , \cdots , y_m^{B_m}) \\
&\mbox{par application du quantificateur existentiel $\exists \, (y_1^{B_1} , \cdots , y_m^{B_m})$,} \hfill \\
&\mbox{son interprétation $M\varphi$ est l'image du sous-objet} \hfill \\
&M\psi \xymatrix{\ar@{^{(}->}[r] & } M\!A_1 \times \cdots \times M\!A_n \times M\!B_1 \times \cdots \times M\!B_m \\
&\mbox{par le morphisme de projection} \hfill \\
&M\!A_1 \times \cdots \times M\!A_n \times M\!B_1 \times \cdots \times M\!B_m \longrightarrow M\!A_1 \times \cdots \times M\!A_n \, . \\
\bullet &\mbox{Si $\varphi$ est la formule du faux $\perp (x_1^{A_1} , \cdots , x_n^{A_n})$, son interprétation est le sous-objet minimal} \hfill \\
&M\varphi = \emptyset \xymatrix{\ar@{^{(}->}[r] & } M\!A_1 \times \cdots \times M\!A_n \, . \\
\bullet &\mbox{Si $\varphi$ s'écrit comme une disjonction de formules déjà interprétées} \hfill \\
&\varphi = \underset{i \in I}{\bigvee} \, \varphi_i \, , \\
&\mbox{sont interprétation est la réunion des sous-objets correspondants} \hfill \\
&M \varphi = \underset{i \in I}{\bigcup} \, M \varphi_i \, .
\end{matrix}\right.$

%\medskip

\item La structure $M$ de type $\Sigma$ est appelée un modèle d'une théorie géométrique ${\mathbb T}$ de signature $\Sigma$ définie par une famille de séquents géométriques
$$
\varphi_i (\vec x_i) \vdash \psi_i (\vec x_i) \, , \quad i \in I \, ,
$$
si $M$ satisfait tous ces axiomes au sens que, pour tout $i \in I$ avec $\vec x_i = (x_1^{A_1} , \cdots , x_n^{A_n})$, les interprétations
$$
\begin{matrix}
M \varphi_i &\xymatrix{\ar@{^{(}->}[r] & } &M\!A_1 \times \cdots \times M\!A_n \, , \\
M \psi_i &\xymatrix{\ar@{^{(}->}[r] & } &M\!A_1 \times \cdots \times M\!A_n \hfill
\end{matrix}
$$
satisfont la relation d'inclusion
$$
M\varphi_i \subseteq M\psi_i \, .
$$
\end{listeimarge}
\end{defn}

\begin{remarksqed}
\begin{listeisansmarge}
\item Les sous-objets d'un objet $E$ d'un topos ${\mathcal E}$ forment un ensemble ordonné par la relation d'inclusion et dans lequel les réunions arbitraires et les intersections arbitraires sont toujours bien définies.

\item Tout morphisme d'un topos ${\mathcal E}$
$$
e : E' \longrightarrow E
$$
se factorise canoniquement en un épimorphisme suivi d'un monomorphisme
$$
E' \xymatrix{\ar@{->>}[r] & } {\rm Im} (e) \xymatrix{\ar@{^{(}->}[r] & } E \, .
$$
Le sous-objet ${\rm Im} (e)$, qui est appelé l'image de $e$, est le plus petit sous-objet de $E$ à travers lequel $e$ se factorise. Il est aussi l'objet quotient de $E'$ par la relation d'équivalence
$$
E' \times_E E' \, .
$$
\item Pour tout morphisme $e : E' \to E$ d'un topos ${\mathcal E}$, associer à tout sous-objet $E'_1 \hookrightarrow E'$ l'image
$$
e_* \, E'_1 = {\rm Im} (E'_1 \xymatrix{\ar@{^{(}->}[r] & } E' \longrightarrow E)
$$
définit une application d'image directe des sous-objets par $e$
$$
e_* : \{\mbox{sous-objets de $E'$}\} \longrightarrow \{\mbox{sous-objets de $E$}\} \, .
$$
Elle respecte la relation $\supseteq$ et est adjointe à droite de l'application d'image réciproque des sous-objets
$$
\begin{matrix}
e^{-1} : \{\mbox{sous-objets de $E$}\} &\longrightarrow &\{\mbox{sous-objets de $E'$}\} \, , \\
\hfill (E_1 \hookrightarrow E) &\longmapsto &(E_1 \times_E E' \hookrightarrow E') \, . \hfill
\end{matrix}
$$
\end{listeisansmarge}
\end{remarksqed}

On observe:

\begin{prop}\label{prop1.1.19}

Soit ${\mathbb T}$ une théorie géométrique du premier ordre de signature $\Sigma$. Alors:

\begin{listeimarge}

\item Pour tout topos ${\mathcal E}$, les modèles de ${\mathbb T}$ dans ${\mathcal E}$ forment une sous-catégorie pleine
$$
{\mathbb T}\mbox{{\rm -mod}} \, ({\mathcal E}) \quad \mbox{de} \quad \Sigma\mbox{{\rm -str}} \, ({\mathcal E}) = \Sigma\mbox{{\rm -mod}} \, ({\mathcal E}) \, .
$$
Elle est localement petite et possède des colimites filtrantes arbitraires qui se calculent comme dans $\Sigma\mbox{{\rm -str}} \, ({\mathcal E})$.

\item Pour tout morphisme de topos $f : {\mathcal E}' \to {\mathcal E}$, le foncteur induit par sa composante d'image réciproque $f^* : {\mathcal E} \to {\mathcal E}'$
$$
f^* : \Sigma\mbox{{\rm -str}} \, ({\mathcal E}) \longrightarrow \Sigma\mbox{{\rm -str}} \, ({\mathcal E}')
$$
se restreint en un foncteur
$$
f^* : {\mathbb T}\mbox{{\rm -mod}} \, ({\mathcal E}) \longrightarrow {\mathbb T}\mbox{{\rm -mod}} \, ({\mathcal E}')
$$
qui respecte les colimites filtrantes.

De plus, le foncteur
$$
f^* : {\mathcal E} \longrightarrow {\mathcal E}'
$$
respecte les interprétations dans les modèles $M$ de ${\mathcal E}$ de toutes les formules géométriques $\varphi (x_1^{A_1} , \cdots , x_n^{A_n})$ de $\Sigma$, au sens que
$$
f^* (M\varphi) = (f^* M)\,\varphi
$$
comme sous-objets de $f^* (M\!A_1 \times \cdots \times M\!A_n) = f^* M\!A_1 \times \cdots \times f^* M\!A_n$.
\end{listeimarge}
\end{prop}

\begin{esquissedemo}
A partir des éléments d'une structure $M$ de type $\Sigma$ dans un topos ${\mathcal E}$, les interprétations $M\varphi$ dans $M$ des formules géométriques $\varphi$ de $\Sigma$ se construisent exclusivement en termes de limites finies et de colimites.

Alors (i) résulte de ce que, dans un topos, les foncteurs de colimites filtrantes respectent non seulement les colimites mais aussi les limites finies.

Et (ii) résulte de ce que, par définition des morphismes de topos $f : {\mathcal E}' \to {\mathcal E}$, leur composante d'image réciproque $f^* : {\mathcal E} \to {\mathcal E}'$ respecte les limites finies et les colimites. 
\end{esquissedemo}

Etant donnée une théorie géométrique du premier ordre ${\mathbb T}$ définie dans une signature $\Sigma$ par une famille d'axiomes
$$
\varphi_i (\vec x_i) \vdash \psi_i (\vec x_i) \, , \quad i \in I \, ,
$$
on dit qu'un séquent géométrique de $\Sigma$
$$
\varphi (\vec x) \vdash \psi (\vec x)
$$
est ``démontrable'' dans ${\mathbb T}$ s'il se déduit de ses axiomes par les ``règles d'inférence'' générales de la logique géométrique.

On ne rappellera pas ici ces règles d'inférence, c'est-à-dire de déduction logique, commune à toutes les théories.

Elles sont telles que soit vérifié le théorème suivant, qui les justifie:

\begin{thm}\label{thm1.1.20}

Soit ${\mathbb T}$ une théorie géométrique du premier ordre de signature $\Sigma$, et soit un séquent géométrique de $\Sigma$
$$
\varphi (x_1^{A_1} , \cdots , x_n^{A_n}) \vdash \psi (x_1^{A_1} , \cdots , x_n^{A_n}) \, .
$$
Alors ce séquent est démontrable dans ${\mathbb T}$ si et seulement si, pour tout modèle $M$ de ${\mathbb T}$ dans un topos ${\mathcal E}$, ce modèle satisfait le séquent $\varphi \vdash \psi$ au sens que les interprétations
$$
\begin{matrix}
M \varphi &\xymatrix{\ar@{^{(}->}[r] & } &M\!A_1 \times \cdots \times M\!A_n \, , \\
M \psi &\xymatrix{\ar@{^{(}->}[r] & } &M\!A_1 \times \cdots \times M\!A_n \hfill
\end{matrix}
$$
satisfont la relation d'inclusion
$$
M\varphi \subseteq M\psi \, .
$$
\end{thm}

\begin{indicationdemo}
La nécessité de la condition résulte de ce que les règles d'inférence de la logique géométrique correspondent à des propriétés catégoriques générales de tous les topos.

Sa suffisance résulte de l'existence et des propriétés du ``modèle universel'' $U_{\mathbb T}$ d'une théorie géométrique ${\mathbb T}$ dans son ``topos classifiant'' ${\mathcal E}_{\mathbb T}$, qui seront rappelées au paragraphe suivant. 
\end{indicationdemo}

\subsection{Les topos comme incarnations géométriques de la sémantique des théories}\label{ssec1.1.5}

Rappelons d'abord que pour tout petit site $({\mathcal C} , J)$ muni du foncteur canonique
$$
\ell = j^* \circ y : {\mathcal C} \xymatrix{\ar@{^{(}->}[r] & } \widehat{\mathcal C} \longrightarrow \widehat{\mathcal C}_J \, ,
$$
tout point généralisé du topos $\widehat{\mathcal C}_J$ paramétré par un topos ${\mathcal E}$
$$
f = (f^* , f_*) : {\mathcal E} \longrightarrow \widehat{\mathcal C}_J
$$
est caractérisé, à unique isomorphisme près, par le foncteur composé
$$
\rho = f^* \circ \ell : {\mathcal C} \longrightarrow \widehat{\mathcal C}_J \longrightarrow {\mathcal E} \, .
$$
De plus, les foncteurs $\rho : {\mathcal C} \to {\mathcal E}$ qui proviennent d'un morphisme de topos ${\mathcal E} \to \widehat{\mathcal C}_J$ peuvent être vus comme les modèles d'un théorie géométrique du premier ordre ${\mathbb T}_{{\mathcal C},J}$ dont la signature $\Sigma$ est définie par la catégorie~${\mathcal C}$.

Pour ce faire, on a d'abord besoin du lemme suivant:

\begin{lem}\label{lem1.1.21}

Soit ${\mathcal C}$ une petite catégorie munie du foncteur de Yoneda $y : {\mathcal C} \hookrightarrow \widehat{\mathcal C}$.

Soit ${\mathcal E}$ un topos ou, plus généralement, une catégorie localement petite qui possède des foncteurs de colimites arbitraires.

Alors tout foncteur
$$
\rho : {\mathcal C} \longrightarrow {\mathcal E}
$$
se prolonge en un foncteur, unique à unique isomorphisme près,
$$
\widehat\rho : \widehat{\mathcal C} \longrightarrow {\mathcal E}
$$
tel que

\smallskip

$\left\{\begin{matrix}
\bullet &y \circ \widehat\rho = \rho \, , \hfill \\
\bullet &\mbox{$\widehat\rho$ respecte les colimites arbitraires.} \hfill 
\end{matrix} \right.$
\end{lem}

\begin{principedemo}
D'après le lemme \ref{lem1.1.14}, tout objet $P$ de $\widehat{\mathcal C}$ s'écrit
$$
P = \varinjlim_{(X,x) \in {\mathcal C}/P} y(X) \, .
$$
Si $\widehat\rho$ existe, il doit donc transformer $P$ en un objet canoniquement isomorphe à
$$
\varinjlim_{(X,x) \in {\mathcal C}/P} \rho(X) \, .
$$
Réciproquement, on vérifie que cette formule définit bien un foncteur $\widehat\rho : \widehat{\mathcal C} \to {\mathcal E}$ qui respecte les colimites. 
\end{principedemo}

\vspace{-0.5cm}

Ce lemme permet de poser:

\begin{defn}\label{defn1.1.22}

Un foncteur d'une petite catégorie ${\mathcal C}$ dans un topos ${\mathcal E}$
$$
\rho : {\mathcal C} \longrightarrow {\mathcal E}
$$
est appelé ``plat'' si son unique prolongement respectant les colimites
$$
\widehat\rho : \widehat{\mathcal C} \longrightarrow {\mathcal E}
$$
respecte aussi les limites finies.
\end{defn}

\begin{remarkqed}
On montre que si ${\mathcal C}$ possède des limites finies arbitraires, un foncteur
$$
\rho : {\mathcal C} \longrightarrow {\mathcal E}
$$
est plat si et seulement si il respecte les limites finies. 
\end{remarkqed}

%\vspace{-0.6cm}

On peut alors énoncer:

%\vspace{-0.2cm}

\begin{thm}[\textbf{Grothendieck, Diaconescu}]\label{thm1.1.23}  
\nopagebreak
Soit ${\mathcal C}$ une petite catégorie munie d'une topologie $J$ et du foncteur canonique associé $\ell : {\mathcal C} \to \widehat{\mathcal C}_J$.

Alors, pour tout topos ${\mathcal E}$, le foncteur \nopagebreak
$$
\begin{matrix}
\hfill \mbox{\rm Geom} \, ({\mathcal E} , \widehat{\mathcal C}_J) &\longrightarrow &[{\mathcal C},{\mathcal E}] \, , \hfill \\
\left( {\mathcal E} \xrightarrow{ \ f=(f^*,f_*) \ } \widehat{\mathcal C}_J \right) &\longmapsto &\left( {\mathcal C} \xrightarrow{ \ f^* \circ \, \ell \ } {\mathcal E} \right) 
\end{matrix}
$$
définit une équivalence de la catégorie des points ${\mathcal E}$-paramétrés du topos $\widehat{\mathcal C}_J$ vers la sous-catégorie pleine de
$$
[{\mathcal C},{\mathcal E}]
$$
constitué des foncteurs
$$\nopagebreak
\rho : {\mathcal C} \longrightarrow {\mathcal E}
$$
qui sont

%\smallskip

$\left\{\begin{matrix}
\bullet &\mbox{plats,} \hfill \\
\bullet &\mbox{$J$-continus au sens qu'ils transforment toute famille $J$-couvrante de ${\mathcal C}$} \hfill \\
&\mbox{$\left( X_i \xrightarrow{ \ x_i \ } X \right)_{i \in I}$ en une famille globalement épimorphique de ${\mathcal E}$.} \hfill
\end{matrix} \right.$
\end{thm}

\begin{remark}
Le cas particulier de ce théorème où ${\mathcal C}$ possède des limites finies est dû à Grothendieck.

Sa forme générale est due à Diaconescu.
\end{remark}

\begin{esquissedemo}
Si $(f^* , f_*) : {\mathcal E} \to \widehat{\mathcal C}_J$ est un morphisme de topos, le composé
$$
f^* \circ j^* : \widehat{\mathcal C} \longrightarrow \widehat{\mathcal C}_J \longrightarrow {\mathcal E}
$$
respecte les colimites arbitraires et les limites finies. D'autre part, il prolonge $\rho = f^* \circ \ell = f^* \circ j^* \circ y$.

Réciproquement, si $\rho : {\mathcal C} \to {\mathcal E}$ est un foncteur plat, le prolongement
$$
\widehat\rho : \widehat{\mathcal C} \to {\mathcal E}
$$
respecte les colimites arbitraires et les limites finies. Il se factorise en
$$
\widehat{\mathcal C} \xrightarrow{ \ j^* \ } \widehat{\mathcal C}_J \longrightarrow {\mathcal E}
$$
si et seulement si le foncteur $\rho : {\mathcal C} \to {\mathcal E}$ est $J$-continu.
\end{esquissedemo}

On déduit de ce théorème:

\begin{cor}\label{cor1.1.24}

Soit ${\mathcal C}$ une petite catégorie.

Soit $\Sigma_{\mathcal C}$ la signature qui consiste en

\smallskip

$\left\{\begin{matrix}
\bullet &\mbox{des sortes qui sont les objets $X$ de ${\mathcal C}$,} \hfill \\
\bullet &\mbox{des symboles de fonctions qui sont les morphismes $x : X' \to X$ de ${\mathcal C}$,} \hfill \\
\bullet &\mbox{aucun symbole de relation.} \hfill
\end{matrix} \right.$

\smallskip

Alors toute topologie $J$ de ${\mathcal C}$ définit une théorie géométrique ${\mathbb T}_{{\mathcal C},J}$ de signature $\Sigma_{\mathcal C}$ telle que, pour tout topos ${\mathcal E}$, le foncteur
$$
\left( f : {\mathcal E} \longrightarrow \widehat{\mathcal C}_J \right) \longmapsto \left(f^* \circ \ell : {\mathcal C} \longrightarrow {\mathcal E} \right)
$$
définisse une équivalence de la catégorie des points généralisés
$$
\mbox{{\rm Geom}} \, ({\mathcal E} , \widehat{\mathcal C}_J)
$$
vers la catégorie
$$
{\mathbb T}_{{\mathcal C},J}\mbox{{\rm -mod}} \, ({\mathcal E})
$$
des structures de type $\Sigma_{\mathcal C}$ dans ${\mathcal E}$ qui sont des modèles de la théorie ${\mathbb T}_{{\mathcal C},J}$.
\end{cor}

\begin{esquissedemo}
Une structure de type $\Sigma_{\mathcal C}$ dans un topos ${\mathcal E}$ consiste en une double famille

\smallskip

$\left\{\begin{matrix}
\bullet &\mbox{d'objets $\rho(X)$ de ${\mathcal E}$ indexés par les objets $X$ de ${\mathcal C}$,} \hfill \\
\bullet &\mbox{de morphismes $\rho(x) : \rho(X') \to \rho(X)$ de ${\mathcal E}$ indexés par les morphismes $x : X' \to X$ de ${\mathcal C}$.} \hfill 
\end{matrix} \right.$

\smallskip

Demander qu'une telle structure soit 

\smallskip

$\left\{\begin{matrix}
\bullet &\mbox{un foncteur ${\mathcal C} \to {\mathcal E}$,} \hfill \\
\bullet &\mbox{plate,} \hfill \\
\bullet &\mbox{$J$-continue,} \hfill
\end{matrix} \right.$

\smallskip

\noindent équivaut à demander qu'elle satisfasse une certaine liste de séquents géométriques de $\Sigma_{\mathcal C}$.

Ces séquents sont les axiomes de la théorie géométrique ${\mathbb T}_{{\mathcal C},J}$.
\end{esquissedemo}

On rappelle que, réciproquement, pour toute théorie géométrique ${\mathbb T}$, ses modèles dans n'importe quel topos ${\mathcal E}$ peuvent être interprétés comme les foncteurs
$$
\rho : {\mathcal C}_{\mathbb T} \longrightarrow {\mathcal E}
$$
sur une certaine catégorie ${\mathcal C}_{\mathbb T}$, dite la ``catégorie syntactique'' de ${\mathbb T}$, qui sont plats et continus pour une certaine topologie $J_{\mathbb T}$, dite la ``topologie syntactique'' de ${\mathcal C}_{\mathbb T}$.

Voici la définition de ces catégories et de leur topologie:

\begin{defn}\label{defn1.1.25}

Soit ${\mathbb T}$ une théorie géométrique du premier ordre de signature $\Sigma$.

On appelle ``site syntactique'' de ${\mathbb T}$ le site $({\mathcal C}_{\mathbb T} , J_{\mathbb T})$ suivant:
\begin{listeimarge}
\item Les objets de ${\mathcal C}_{\mathbb T}$ sont les formules géométriques de $\Sigma$
$$
\varphi (x_1^{A_1} , \cdots , x_n^{A_n}) \, ,
$$
considérées à substitution près de leurs variables formelles par des variables de mêmes sortes.

\item Les morphismes de ${\mathcal C}_{\mathbb T}$
$$
\varphi (\vec x) \longrightarrow \psi (\vec y)
$$
sont les formules en la réunion disjointe des deux familles de variables
$$
\theta (\vec x , \vec y)
$$
qui sont ``démontrablement fonctionnelles'' au sens que les séquents

\smallskip

$\left\{\begin{matrix}
\theta (\vec x , \vec y) \vdash \varphi (\vec x) \wedge \psi (\vec y) \, , \hfill \\
\varphi (\vec x) \vdash (\exists \, \vec y) \, \theta (\vec x , \vec y) \, , \hfill \\
\theta (\vec x , \vec y) \wedge \theta (\vec x , \vec y\,') \vdash \vec y = \vec y\,' \hfill
\end{matrix}\right.$

\smallskip

\noindent sont démontrables dans la théorie ${\mathbb T}$.

\item La loi de composition associe à deux formules démontrablement fonctionnelles
$$
\varphi (\vec x) \xrightarrow{ \ \theta (\vec x , \vec y) \ } \psi (\vec y) \xrightarrow{ \ \theta' (\vec y , \vec z) \ } \chi (\vec z) 
$$
la formule démontrablement fonctionnelle
$$
\varphi (\vec x) \xrightarrow{ \ (\exists \, \vec y) (\theta (\vec x , \vec y) \wedge \theta' (\vec y , \vec z) \ } \chi (\vec z) \, .
$$
\item La topologie $J_{\mathbb T}$ est celle pour laquelle une famille de morphismes
$$
\varphi_i (\vec x_i) \xrightarrow{ \ \theta_i (\vec x_i , \vec y) \ } \psi (\vec y) \, , \quad i \in I \, ,
$$
est couvrante si et seulement si le séquent
$$
\psi (\vec y) \vdash \bigvee_{i \in I} (\exists \, \vec x_i) \, \theta_i (\vec x_i , \vec z)
$$
est démontrable dans la théorie ${\mathbb T}$.
\end{listeimarge}
\end{defn}

\begin{remarkqed}
Si ${\mathbb T}$ est une théorie cohérente [resp. régulière], on définit de même la ``catégorie syntactique cohérente [resp. régulière]'' de ${\mathbb T}$
$$
{\mathcal C}_{\mathbb T}^{\rm coh} \quad \mbox{[resp. ${\mathcal C}_{\mathbb T}^{\rm reg}$]}
$$
en retenant dans la définition des objets et des morphismes les seules formules cohérentes [reps. régulières] de $\Sigma$.

La topologie syntactique $J_{\mathbb T}^{\rm coh}$ [resp. $J_{\mathbb T}^{\rm reg}$] de cette catégorie est définie en déclarant couvrantes les familles de morphismes qui contiennent une sous-famille finie
$$
\varphi_i (\vec x_i) \xrightarrow{ \ \theta_i (\vec x_i , \vec y) \ } \psi (\vec y) \, , \quad 1 \leq i \leq n \, ,
$$
[resp. un morphisme
$$
\varphi (\vec x) \xrightarrow{ \ \theta (\vec x , \vec y) \ } \psi (\vec y) \ \mbox{]}
$$
dont le séquent cohérent [resp. régulier] associé
$$
\psi (\vec y) \vdash (\exists \, \vec x_1) \, \theta_1 (\vec x_1 , \vec y) \vee \cdots \vee (\exists \, \vec x_n) \, \theta_n (\vec x_n , \vec y)
$$
[resp.
$$
\psi (\vec y) \vdash (\exists \, \vec x) \, \theta (\vec x , \vec y) \ \mbox{]}
$$
est démontrable. 
\end{remarkqed}

Voici les propriétés de base des catégories syntactiques:

\begin{lem}\label{lem1.1.26}
Soit ${\mathbb T}$ une théorie géométrique. Alors:
\begin{listeimarge}
\item La catégorie syntactique ${\mathcal C}_{\mathbb T}$ est essentiellement petite.

\item Elle possède des limites finies arbitraires.

\item Pour tout objet $\psi (\vec x)$ de ${\mathcal C}_{\mathbb T}$, ses sous-objets sont les formules
$$
\varphi (\vec x) \qquad \mbox{en les mêmes variables $\vec x$,}
$$
telle que le séquent
$$
\varphi (\vec x) \vdash \psi (\vec x)
$$
soit démontrable dans ${\mathbb T}$.
\end{listeimarge}
\end{lem}

\begin{remark}
Si ${\mathbb T}$ est une théorie cohérente [resp. régulière], sa catégorie syntactique cohérente ${\mathcal C}_{\mathbb T}^{\rm coh}$ [resp. régulière ${\mathcal C}_{\mathbb T}^{\rm reg}$] est petite.

Elle a des limites finies arbitraires, et les sous-objets d'un objet $\psi (\vec x)$ sont les formules cohérentes [resp. régulières] $\varphi (\vec x)$ telles que le séquent $\varphi (\vec x) \vdash \psi (\vec x)$ soit démontrable.
\end{remark}

\begin{esquissedemosansqed}
	\begin{listeisansmarge}
\item[(ii)] Le produit de formules $\varphi_1 (\vec x_1) , \cdots , \varphi_{n} (\vec x_n)$ est la formule
$$
\varphi (\vec x_1) \wedge \cdots \wedge \varphi_n (\vec x_n)
$$
en la famille de variables réunion disjointe des $\vec x_i$, $1 \leq i \leq n$.

L'égalisateur de deux morphismes
$$
\xymatrix{ \varphi (\vec x) \dar[rr]^{^{^{\mbox{\scriptsize $\theta_1 (\vec x , \vec y)$}}}}_{\theta_2 (\vec x , \vec y)} && \psi (\vec y) }
$$
est le sous-objet de $\varphi (\vec x)$ consistant en la formule en $\vec x$
$$
(\exists \, \vec y) \, (\theta_1 (\vec x , \vec y) \wedge \theta_2 (\vec x , \vec y)) \, .
$$
\item[(iii)] Un morphisme de ${\mathcal C}_{\mathbb T}$
$$
\varphi (\vec x) \xrightarrow{ \ \theta (\vec x , \vec y) \ } \psi (\vec y)
$$
est un monomorphisme si et seulement si le séquent
$$
\theta (\vec x , \vec y) \wedge \theta (\vec x\,' , \vec y) \vdash \vec x = \vec x\,'
$$
est démontrable dans ${\mathbb T}$.

Alors $\theta (\vec x , \vec y)$ définit un isomorphisme de $\varphi (\vec x)$ vers la formule en $\vec y$
$$
(\exists \, \vec x) \, \theta (\vec x , \vec y) \, , \quad \mbox{laquelle implique $\psi (\vec y)$.}
$$
\item[(i)] résulte du lemme suivant:
\end{listeisansmarge}
\end{esquissedemosansqed}

\begin{lem}\label{lem1.1.27}

Soit $\Sigma$ une signature de ${\mathbb T}$.

Alors toute formule géométrique [resp. cohérente, resp. régulière] de $\Sigma$ est équivalente (selon les règles d'inférence de la logique géométrique) à une formule de la forme
$$
\varphi (\vec x) = \bigvee_{i \in I} (\exists \, \vec x_i) \, \varphi_i (\vec x_i , \vec x)
$$
[resp. \qquad $\varphi (\vec x) = (\exists \, \vec x_1) \, \varphi_1 (\vec x_1 , \vec x) \vee \cdots \vee (\exists \, \vec x_n) \, \varphi_n (\vec x_n , \vec x)$,

\noindent [resp. \qquad $\varphi (\vec x) = (\exists \, \vec x_1) \, \varphi_1 (\vec x_1 , \vec x)$ \ ]

\noindent où chaque $\varphi_i (\vec x_i , \vec x)$ est une formule de Horn, c'est-à-dire une conjonction finie de formules atomiques en $(\vec x_i , \vec x)$.
\end{lem}

\begin{indicationdemo}
Cela résulte de ce que les règles d'inférence de la logique géométrique permettent en particulier d'échanger les ordres d'application des symboles $\bigvee$, $\exists$ and $\wedge$. 
\end{indicationdemo}

On a:

\begin{prop}\label{prop1.1.28}

Soit ${\mathbb T}$ une théorie géométrique du premier ordre de signature $\Sigma$.

Soit ${\mathcal E}$ un topos. 

Alors:
\begin{listeimarge}
\item Pour tout modèle $M$ de ${\mathbb T}$ dans ${\mathcal E}$, associer

\smallskip

$\left\{\begin{matrix}
\bullet &\mbox{à tout objet $\varphi (\vec x)$ de ${\mathcal C}_{\mathbb T}$ son interprétation $M\varphi$ dans $M$,} \hfill \\
\bullet &\mbox{à tout morphisme $\varphi (\vec x) \xrightarrow{ \ \theta (\vec x, \vec y) \ } \psi (\vec y)$ de ${\mathcal C}_{\mathbb T}$ le morphisme de ${\mathcal E}$} \hfill \\
&M\varphi \longrightarrow M\psi \\
&\mbox{dont le graphe est l'interprétation $M\theta \hookrightarrow M\varphi \times M\psi$,} \hfill
\end{matrix}\right.$

\medskip

\noindent définit un foncteur
$$
\rho_M : {\mathcal C}_{\mathbb T} \longrightarrow {\mathcal E}
$$
qui

\smallskip

$\left\{\begin{matrix}
\bullet &\mbox{respecte les limites finies,} \hfill \\
\bullet &\mbox{est $J_{\mathbb T}$-continu.} \hfill
\end{matrix}\right.$

%\medskip

\item Réciproquement, tout foncteur vérifiant ces deux propriétés
$$
\rho : {\mathcal C}_{\mathbb T} \longrightarrow {\mathcal E}
$$
définit un modèle $M$ de ${\mathbb T}$ dans ${\mathcal E}$ en posant

\smallskip

$\left\{\begin{matrix}
\bullet &\mbox{$M\!A = \rho (\top (x^A))$ pour toute sorte $A$ de $\Sigma$,} \hfill \\
\bullet &M\!f = \rho \left(\top (x_1^{A_1}) \times \cdots \times \top (x_n^{A_n}) \xrightarrow{ \ y^B = f(x_1^{A_1} , \cdots , x_n^{A_n}) \ } \top (y^B)\right) \hfill \\
&\mbox{pour tout symbole de fonction $f : A_1 \cdots A_n \to B$ de $\Sigma$,} \hfill \\
\bullet &\mbox{$M\!R = \rho (R (x_1^{A_1} , \cdots , x_n^{A_n}))$ pour tout symbole de relation $R \rightarrowtail A_1 \cdots A_n$ de $\Sigma$.} \hfill
\end{matrix}\right.$

%\medskip

\item Ces deux constructions définissent deux équivalences réciproques l'une de l'autre de la catégorie des modèles de ${\mathbb T}$ dans ${\mathcal E}$
$$
{\mathbb T}\mbox{\rm -mod} \, ({\mathcal E})
$$
vers la catégorie des foncteurs
$$
\rho : {\mathcal C}_{\mathbb T} \longrightarrow {\mathcal E}
$$
qui respectent les limites finies et sont $J_{\mathbb T}$-continus.
\end{listeimarge}
\end{prop}

\begin{remark}
Si ${\mathbb T}$ est une théorie cohérente [resp. régulière], le même énoncé s'applique à la catégorie ${\mathcal C}_{\mathbb T}^{\rm coh}$ [resp. ${\mathcal C}_{\mathbb T}^{\rm reg}$] munie de la topologie $J_{\mathbb T}^{\rm coh}$ [resp. $J_{\mathbb T}^{\rm reg}$].
\end{remark}

\noindent {\bf Indication de démonstration:}

\smallskip

Cela résulte des règles d'inférence de la logique géométrique, compte tenu de la manière dont les symboles logiques $\top$, $\wedge$, $\exists$, $\bigvee$, $\perp$ sont interprétés dans les topos. 
\hfill $\Box$

\medskip

En combinant cette proposition et le théorème \ref{thm1.1.23}, on obtient:

\begin{thm}\label{thm1.1.29}

Soit ${\mathbb T}$ une théorie géométrique du premier ordre de signature $\Sigma$.

Soit ${\mathcal E}_{\mathbb T}$ le topos des faisceaux sur le site $({\mathcal C}_{\mathbb T} , J_{\mathbb T})$.

Soit $U_{\mathbb T}$ le modèle de ${\mathbb T}$ dans ${\mathcal E}_{\mathbb T}$ qui correspond au foncteur canonique
$$
\ell : {\mathcal C}_{\mathbb T} \longrightarrow {\mathcal E}_{\mathbb T} \, .
$$
Alors, pour tout topos ${\mathcal E}$, le foncteur
$$
\begin{matrix}
\mbox{\rm Geom} \, ({\mathcal E},{\mathcal E}_{\mathbb T}) &\longrightarrow &{\mathbb T}\mbox{\rm -mod} \, ({\mathcal E}) \, , \\
\hfill \left({\mathcal E} \xrightarrow{ \, f \, } {\mathcal E}_{\mathbb T}\right) &\longmapsto &f^* U_{\mathbb T} \hfill
\end{matrix}
$$
est une équivalence de catégories.
\end{thm}

\begin{remarksqed}{\tiny\color{white}Grothendieck}\vspace{-0.5cm}
	\begin{listeisansmarge}
		\item Les premiers cas particuliers de ce théorème apparaissent dans la thèse \cite{ToposAnneles} de Monique Hakim sous la direction de A. Grothendieck.

Il a été énoncé et démontré en toute généralité dans le livre \cite{CategoricalLogic} de M. Makkai et G. Reyes, à la suite de travaux de plusieurs autres logiciens catégoriciens, en particulier W. Lawvere et A. Joyal.
\item La propriété du théorème caractérise à équivalence près le topos ${\mathcal E}_{\mathbb T}$, appelé le ``topos classifiant'' de ${\mathbb T}$, muni du modèle $U_{\mathbb T}$, appelé le ``modèle universel'' de ${\mathbb T}$.

\item Si ${\mathbb T}$ est une théorie cohérente [resp. régulière], le topos ${\mathcal E}_{\mathbb T}$ peut aussi être construit comme le topos des faisceaux sur le site $({\mathcal C}_{\mathbb T}^{\rm coh} , J_{\mathbb T}^{\rm coh})$ [resp. $({\mathcal C}_{\mathbb T}^{\rm reg} , J_{\mathbb T}^{\rm reg})$].
\end{listeisansmarge}
\end{remarksqed}

Voici des propriétés de base des foncteurs canoniques reliant les catégories syntactiques aux topos classifiants:

\begin{prop}\label{prop1.1.30}

Soit ${\mathbb T}$ une théorie géométrique du premier ordre de signature $\Sigma$.

Alors:
\begin{listeimarge}
\item Le foncteur canonique
$$
\ell : {\mathcal C}_{\mathbb T} \longrightarrow {\mathcal E}_{\mathbb T}
$$
est pleinement fidèle.

\item Pour tout objet $\psi (\vec x)$ de ${\mathcal C}_{\mathbb T}$, le foncteur $\ell$ induit une bijection de l'ensemble ordonné des sous-objets de $\psi (\vec x)$ dans ${\mathcal C}_{\mathbb T}$
$$
\{\varphi (\vec x) \mid \varphi \vdash \psi \quad \mbox{est démontrable dans ${\mathbb T}$}\}
$$
sur l'ensemble ordonné des sous-objets de
$$
\ell (\varphi (\vec x)) = U_{\mathbb T} \ \varphi \quad \mbox{dans} \ {\mathcal E}_{\mathbb T} \, .
$$

\item En particulier, un séquent géométrique de $\Sigma$
$$
\varphi (\vec x) \vdash \psi (\vec x)
$$
est démontrable dans ${\mathbb T}$ si et seulement si il est vérifié par le modèle $U_{\mathbb T}$ de ${\mathbb T}$ dans ${\mathcal E}_{\mathbb T}$, au sens que
$$
U_{\mathbb T} \, \varphi \subseteq U_{\mathbb T} \, \psi
$$
en tant que sous-objets de $U_{\mathbb T} \, \top (\vec x)$.
\end{listeimarge}
\end{prop}

\begin{remark}
Si ${\mathbb T}$ est cohérente [resp. régulière], le foncteur canonique
$$
\ell : {\mathcal C}_{\mathbb T}^{\rm coh} \longrightarrow {\mathcal E}_{\mathbb T} \qquad \mbox{[resp. \ ${\mathcal C}_{\mathbb T}^{\rm reg} \longrightarrow {\mathcal E}_{\mathbb T}$ \, ]}
$$
possède les mêmes propriétés.
\end{remark}

\begin{indicationdemo}
\begin{listeisansmarge}	
\item résulte de ce que, pour tout objet $\varphi (\vec x)$ de ${\mathcal C}_{\mathbb T}$, le préfaisceau qu'il représente
$$
y (\varphi (\vec x)) = \mbox{\rm Hom} (\bullet , \varphi (\vec x))
$$
est un faisceau pour la topologie $J_{\mathbb T}$.

\item Tout sous-objet d'un objet de ${\mathcal E}_{\mathbb T}$ de la forme $\ell(\varphi (\vec x))$ est nécessairement la réunion des images d'une famille de morphismes qui sont les transformés par $\ell$ de morphismes de ${\mathcal C}_{\mathbb T}$
$$
\varphi_i (\vec x_i) \xrightarrow{ \ \theta_i (\vec x_i , \vec x) \ } \varphi (\vec x) \, , \quad i \in I \, .
$$
Donc ce sous-objet est le transformé par $\ell$ du sous-objet de $\varphi (\vec x)$ dans ${\mathcal C}_{\mathbb T}$
$$
\bigvee_{i \in I} (\exists \, \vec x_i) \, \theta_i (\vec x_i , \vec x) \, .
$$
\item Un séquent géométrique $\varphi (\vec x) \vdash \psi (\vec x)$ est démontrable dans ${\mathbb T}$ si et seulement si  le monomorphisme de~${\mathcal C}_{\mathbb T}$
$$
(\varphi \wedge \psi) (\vec x) \xymatrix{\ar@{^{(}->}[r] & } \varphi (\vec x) \quad \mbox{est un isomorphisme,}
$$
ce qui revient d'après (i) à demander que son transformé par $\ell$
$$
U_{\mathbb T} \, \varphi \cap U_{\mathbb T} \, \psi \xymatrix{\ar@{^{(}->}[r] & } U_{\mathbb T} (\varphi) \quad \mbox{soit un isomorphisme.}
$$
\end{listeisansmarge}
\end{indicationdemo}

Les catégories syntactiques et les topos classifiants des théories géométriques incarnent la syntaxe et la sémantique des théories, ainsi que leurs relations. On a en particulier la double définition naturelle:

\begin{defn}\label{defn1.1.31}

Deux théories géométriques ${\mathbb T}_1$ et ${\mathbb T}_2$ sont dites syntactiquement équivalentes si leurs catégories syntactiques ${\mathcal C}_{{\mathbb T}_1}$ et ${\mathcal C}_{{\mathbb T}_2}$ sont équivalentes.

Elles sont dites sémantiquement équivalentes si leurs topos classifiants sont équivalents.
\end{defn}

\begin{remarkqed}{\tiny\color{white}Grothendieck} \vspace{-0.4cm}

Si deux théories sont syntactiquement équivalentes, elles sont a fortiori sémantiquement équivalentes.

Cela résulte de ce que la topologie syntactique $J_{\mathbb T}$ sur la catégorie syntactique ${\mathcal C}_{\mathbb T}$ d'une théorie géométrique ${\mathbb T}$ est déterminée par la structure catégorique de ${\mathcal C}_{\mathbb T}$.

En effet, une famille de morphismes vers un objet $\varphi (\vec x)$ de ${\mathcal C}_{\mathbb T}$
$$
\varphi_i (\vec x_i) \xrightarrow{ \ \theta_i (\vec x_i , \vec x) \ } \varphi (\vec x)
$$
est couvrante pour la topologie $J_{\mathbb T}$ si et seulement si $\varphi (\vec x)$ ne contient aucun sous-objet strict à travers lequel se factorisent tous les morphismes de cette famille. 
\end{remarkqed}

Terminons ce paragraphe de rappels sur les topos classifiants par la notion de théorie cartésienne et la représentation simple du topos classifiant d'une telle théorie:

\begin{defn}\label{defn1.1.32}

Soit ${\mathbb T}$ une théorie géométrique du premier ordre de signature $\Sigma$.

\begin{listeimarge}

\item On appelle formules ${\mathbb T}$-cartésiennes les formules régulières de $\Sigma$ qui s'écrivent sous la forme
$$
(\exists \, \vec y) \, \varphi (\vec x , \vec y)
$$
où $\varphi (\vec x , \vec y)$ est une formule de Horn (c'est-à-dire une conjonction finie de formules atomiques) et le séquent régulier
$$
\varphi (\vec x , \vec y) \wedge \varphi (\vec x , \vec y\,') \vdash \vec y = \vec y\,'
$$
est démontrable dans ${\mathbb T}$.

\item La théorie ${\mathbb T}$ est dite cartésienne si elle peut être définie par une famille d'axiomes qui sont des séquents réguliers de $\Sigma$
$$
\varphi (\vec x) \vdash \psi (\vec x)
$$
dont les composants $\varphi (\vec x) , \psi (\vec x)$ sont des formules ${\mathbb T}$-cartésiennes. 
\end{listeimarge}
\end{defn}

\begin{exsqed}
\begin{listeisansmarge}	
\item Toutes les théories de Horn, en particulier les théories sans axiome (réduites à leur seule signature) et les théories algébriques, sont cartésiennes.

\item La théorie des catégories n'est pas de Horn, mais elle est cartésienne. 
\end{listeisansmarge}
\end{exsqed}

Les théories cartésiennes ont non seulement des catégories syntactiques géométriques, cohérentes et régulières mais aussi des catégories syntactiques cartésiennes:

\begin{defn}\label{defn1.1.33}

Soit ${\mathbb T}$ une théorie cartésienne.

On appelle catégorie syntactique cartésienne de ${\mathbb T}$, notée
$$
{\mathcal C}_{\mathbb T}^{\rm car} \, ,
$$
la sous-catégorie de sa catégorie syntactique régulière ${\mathcal C}_{\mathbb T}^{\rm reg}$ constituée des objets $\varphi (\vec x)$ et des morphismes
$$
\theta (\vec x , \vec y) : \varphi (\vec x) \longrightarrow \psi (\vec y)
$$
constitués de formules $\varphi (\vec x), \theta (\vec x , \vec y) , \psi (\vec y)$ qui sont ${\mathbb T}$-cartésiennes.
\end{defn}

\begin{remarks}
	\begin{listeisansmarge}
\item Comme sous-catégorie de ${\mathcal C}_{\mathbb T}^{\rm reg}$, ${\mathcal C}_{\mathbb T}^{\rm car}$ est une petite catégorie.

\item Elle a des limites finies arbitraires, qui sont respectées par le foncteur de plongement ${\mathcal C}_{\mathbb T}^{\rm car} \hookrightarrow {\mathcal C}_{\mathbb T}^{\rm reg}$.

\item Les sous-objets d'un objet $\psi (\vec x)$ de ${\mathcal C}_{\mathbb T}^{\rm car}$ sont les formules ${\mathbb T}$-cartésiennes $\varphi (\vec x)$ telles que le séquent
$$
\varphi (\vec x) \vdash \psi (\vec x)
$$
soit démontrable dans ${\mathbb T}$.
\end{listeisansmarge}
\end{remarks}

\begin{justifdefi}
${\mathcal C}_{\mathbb T}^{\rm car}$ est bien définie comme sous-catégorie de ${\mathcal C}_{\mathbb T}^{\rm reg}$ car, pour toute paire de morphismes de ${\mathcal C}_{\mathbb T}^{\rm reg}$,
$$
\varphi (\vec x) \xrightarrow{ \ \theta (\vec x , \vec y ) \ } \psi (\vec y) \xrightarrow{ \ \theta' (\vec y , \vec z ) \ } \chi (\vec z) \, ,
$$
la formule qui définit leur composé
$$
(\exists \, \vec y) \left(\theta (\vec x , \vec y ) \wedge \theta' (\vec y , \vec z )\right)
$$
est ${\mathbb T}$-cartésienne si les deux formules $\theta (\vec x , \vec y )$ et $\theta' (\vec y , \vec z )$ sont ${\mathbb T}$-cartésiennes.
\end{justifdefi}

Les théories cartésiennes vérifient une forme simplifiée de la proposition \ref{prop1.1.28}:

\begin{prop}\label{prop1.1.34}

Soit ${\mathbb T}$ une théorie cartésienne de signature $\Sigma$.

Soit ${\mathcal E}$ un topos.

Alors:
\begin{listeimarge}
\item La restriction à ${\mathcal C}_{\mathbb T}^{\rm car}$ de la construction de la proposition I.1.28~(i) associe à tout modèle $M$ de ${\mathbb T}$ dans ${\mathcal E}$ un foncteur
$$
\rho : {\mathcal C}_{\mathbb T}^{\rm car} \longrightarrow {\mathcal E}
$$
qui respecte les limites finies.

\item Réciproquement, la construction de la proposition \ref{prop1.1.28}~(ii) associe à tout foncteur respectant les limites finies
$$
\rho : {\mathcal C}_{\mathbb T}^{\rm car} \longrightarrow {\mathcal E}
$$
un modèle $M$ de ${\mathbb T}$ dans ${\mathcal E}$.

\item Ces deux foncteurs définissent deux équivalences réciproques l'une de l'autre de la catégorie des modèles de ${\mathbb T}$ dans ${\mathcal E}$
$$
{\mathbb T}\mbox{\rm -mod} \, ({\mathcal E})
$$
vers la catégorie des foncteurs respectant les limites finies
$$
\rho : {\mathcal C}_{\mathbb T}^{\rm car} \longrightarrow {\mathcal E} \, .
$$
\end{listeimarge}
\end{prop}

\begin{indicationdemo}
La proposition résulte de ce que
\begin{enumerate}
\item[$\bullet$] les formules qui servent à la construction de la proposition \ref{prop1.1.28}~(ii) sont atomiques, donc définissent des objets, des morphismes ou des sous-objets de ${\mathcal C}_{\mathbb T}^{\rm car}$,
\item[$\bullet$] la théorie ${\mathbb T}$ est définie par des axiomes qui sont des séquents reliant des formules ${\mathbb T}$-cartésiennes.
\end{enumerate}{\tiny\color{white}Grothendieck} \vspace{-0.5cm}
\end{indicationdemo}

En combinant cette proposition et le théorème \ref{thm1.1.23}, on obtient:

\begin{thm}\label{thm1.1.35}

Soit ${\mathbb T}$ une théorie cartésienne de signature $\Sigma$.

Soit ${\mathcal E}_{\mathbb T} = \widehat{{\mathcal C}_{\mathbb T}^{\rm car}}$ le topos des préfaisceaux sur la petite catégorie ${\mathcal C}_{\mathbb T}^{\rm car}$.

Soit $U_{\mathbb T}$ le modèle de ${\mathbb T}$ dans ${\mathcal E}_{\mathbb T}$ qui correspond au foncteur de Yoneda
$$
y : {\mathcal C}_{\mathbb T}^{\rm car} \xymatrix{\ar@{^{(}->}[r] & } \widehat{{\mathcal C}_{\mathbb T}^{\rm car}} = {\mathcal E}_{\mathbb T} \, .
$$
Alors, pour tout topos ${\mathcal E}$, le foncteur
$$
\begin{matrix}
\mbox{\rm Geom} \, ({\mathcal E} , {\mathcal E}_{\mathbb T}) &\longrightarrow &{\mathbb T} \mbox{\rm -mod} \, ({\mathcal E}) \, , \\
\hfill ({\mathcal E} \xrightarrow{ \ f \ } {\mathcal E}_{\mathbb T} ) &\longmapsto &f^* \, U_{\mathbb T} \hfill
\end{matrix}
$$
est une équivalence de catégories.
\end{thm}

\begin{remarkqed}
Une théorie géométrique ${\mathbb T}$ est dite ``de type préfaisceau'' si son topos classifiant ${\mathcal E}_{\mathbb T}$ peut se représenter comme le topos des préfaisceaux
$$
{\mathcal E}_{\mathbb T} \cong \widehat{\mathcal C}
$$
sur une certaine petite catégorie ${\mathcal C}$.

Le théorème dit que toute théorie cartésienne est de type préfaisceau. Il existe cependant des théories non cartésiennes qui sont de type préfaisceau. 
\end{remarkqed}

\section{La notion de sous-topos et ses différentes expressions}\label{sec1.2}

\subsection{La notion catégorique de sous-topos}\label{ssec1.2.1}

Toute application continue entre deux espaces topologiques
$$
f : Y \longrightarrow X
$$
définit un morphisme entre les topos de faisceaux sur $Y$ et $X$
$$
{\mathcal E}_Y \xrightarrow{ \ (f^* , f_*) \ } {\mathcal E}_X
$$
dont la composante d'image directe $f_*$ transforme les faisceaux sur $Y$
$$
F : \underset{\mbox{\footnotesize ouvert de $Y$}}{V} \longmapsto \underset{\mbox{\footnotesize ensemble}}{F(V)}
$$
en les faisceaux sur $X$
$$
f_* F : \underset{\mbox{\footnotesize ouvert de $X$}}{U} \longmapsto F(f^{-1} (U)) \, .
$$
Un morphisme entre deux faisceaux $F_1 , F_2$ sur $Y$ est une famille compatible d'applications indexées par les ouverts $V$ de $Y$
$$
F_1(V) \longrightarrow F_2(V)
$$
tandis qu'un morphisme entre leurs images directes $f_* F_1 \to f_* F_2$ est une famille compatible d'applications indexées par les ouverts $U$ de $X$
$$
F_1 (f^{-1} (U)) \longrightarrow F_2 (f^{-1} (U)) \, .
$$
Si $f : Y \to X$ est le plongement d'un sous-espace $Y$ de $X$, les ouverts de $Y$ sont exactement les traces $f^{-1} (U)$ des ouverts $U$ de $X$, d'où il résulte que le foncteur
$$
f_* : {\mathcal E}_Y \longrightarrow {\mathcal E}_X
$$
est pleinement fidèle.

Cette remarque amène à poser la définition suivante:

\begin{defn}\label{defn1.2.1}

Un morphisme de topos
$$
f = (f^* , f_*) : {\mathcal E}' \longrightarrow {\mathcal E}
$$
est appelé un ``plongement'' s'il satisfait les deux propriétés équivalentes suivantes:

\begin{enumerate}[label=(\arabic*)]

\item Sa composante d'image directe
$$
f_* : {\mathcal E}' \longrightarrow {\mathcal E}
$$
est un foncteur pleinement fidèle.

\item Le morphisme d'adjonction
$$
f^* \circ f_* \longrightarrow {\rm id}_{{\mathcal E}'}
$$
est un isomorphisme de foncteurs.
\end{enumerate}
\end{defn}

\begin{justifdefisansqed}
L'équivalence de ces deux propriétés est un cas particulier du lemme général suivant qui est souvent utile:
\end{justifdefisansqed}
\begin{lem}\label{lem1.2.2}

Pour toute paire de foncteurs adjoints entre deux catégories localement petites
$$
\left( {\mathcal C} \xrightarrow{ \ F \ } {\mathcal D} \, , \ {\mathcal D} \xrightarrow{ \ G \ } {\mathcal C} \right),
$$
le foncteur adjoint à gauche $F$ [resp. adjoint à droite $G$] est pleinement fidèle si et seulement si le morphisme d'adjonction
$$
{\rm id}_{\mathcal C} \longrightarrow G \circ F \qquad \mbox{[resp. \ $F \circ G \longrightarrow {\rm id}_{\mathcal D}$ \, ]}
$$
est un isomorphisme de foncteurs.
\end{lem}

\begin{esquissedemo}
Il suffit de traiter le premier cas.

Il résulte du lemme de Yoneda puisque, par définition de la propriété d'adjonction, le foncteur
$$
\begin{matrix}
{\mathcal C}^{\rm op} \times{\mathcal C} &\longrightarrow &\mbox{\rm Ens} \, , \hfill \\
(X',X) &\longmapsto &\mbox{\rm Hom} (X' , G \circ F (X))
\end{matrix}
$$
est isomorphe au foncteur
$$
(X',X) \longmapsto \mbox{\rm Hom} (F(X') , F (X)) \, .
$$
\end{esquissedemo}

On déduit immédiatement de la définition des plongements de topos:

\begin{lem}\label{lem1.2.3}

Considérons deux morphismes de topos
$$
{\mathcal E}'' \xrightarrow{ \ g \ } {\mathcal E}' \xrightarrow{ \ f \ } {\mathcal E} \, .
$$
Alors, si $f$ est un plongement, le composé
$$
f \circ g : {\mathcal E}'' \longrightarrow {\mathcal E}'
$$
est un plongement si et seulement si l'autre facteur
$$
g : {\mathcal E}'' \longrightarrow {\mathcal E}'
$$
est un plongement.
\end{lem}

\begin{demo}
Il suffit de considérer pour tous objets $E_1$ et $E_2$ de ${\mathcal E}''$ les deux applications
$$
\begin{matrix}
&g_* : \mbox{\rm Hom} (E_1,E_2) \longrightarrow \mbox{\rm Hom} (g_* E_1 , g_* E_2) \, , \hfill \\
&f_* : \mbox{\rm Hom} (g_* E_1 , g_* E_2) \longrightarrow \mbox{\rm Hom} (f_* g_* E_1 , f_* g_* E_2) \, ,
\end{matrix}
$$
et de remarquer que, si $f_*$ est une bijection, la composée
$$
f_* \circ g_* : \mbox{\rm Hom} (E_1,E_2) \longrightarrow \mbox{\rm Hom} (f_* g_* E_1 , f_* g_* E_2)
$$
est une bijection si et seulement si $g_*$ est une bijection.
\end{demo}

Les sous-topos sont d'abord définis comme des classes d'équivalence de plongements dans un topos:

\begin{defn}\label{defn1.2.4}

Soit un topos ${\mathcal E}$.

\begin{listeimarge}

\item Deux plongements de topos
$$
e_1 : {\mathcal E}_1 \xymatrix{\ar@{^{(}->}[r] & } {\mathcal E} \qquad \mbox{et} \qquad e_2 : {\mathcal E}_2 \xymatrix{\ar@{^{(}->}[r] & } {\mathcal E} 
$$
sont dits équivalents s'il existe une équivalence de topos
$$
e : {\mathcal E}_1 \xrightarrow{ \ \sim \ } {\mathcal E}_2
$$
telle que le triangle de morphismes de topos
$$
\xymatrix{
{\mathcal E}_1 \ar[rr]^{\sim} \ar@{_{(}->}[rd] &&{\mathcal E}_2 \ar@{_{(}->}[ld] \\
&{\mathcal E}
}
$$
soit commutatif à isomorphisme près.

\item On appelle sous-topos de ${\mathcal E}$ les classes d'équivalence de plongements de topos  ${\mathcal E}' \hookrightarrow {\mathcal E}$.% \hfill $\Box$
\end{listeimarge}
\end{defn}

Tous sous-topos ainsi défini comme une classe d'équivalence possède un représentant privilégié:

\begin{prop}\label{prop1.2.5}

Soit un topos ${\mathcal E}$.

\smallskip

Pour tout plongement de topos $f : {\mathcal E'} \hookrightarrow {\mathcal E}$, notons $\mbox{\rm Im} \, (f)$ la sous-catégorie pleine de ${\mathcal E}$ constituée des objets qui sont isomorphes à l'image par $f_*$ d'un objet de ${\mathcal E}'$.

Alors:

\begin{listeimarge}

\item La catégorie $\mbox{\rm Im} \, (f)$ est un topos.

\item Le foncteur pleinement fidèle de plongement
$$
j_* : \mbox{\rm Im} \, (f) \xymatrix{\ar@{^{(}->}[r] & } {\mathcal E}
$$
admet un adjoint à gauche
$$
j^* : {\mathcal E} \longrightarrow \mbox{\rm Im} \, (f)
$$
qui respecte non seulement les colimites mais aussi les limites finies.

\item Un objet $E$ de ${\mathcal E}$ est dans $\mbox{\rm Im} \, (f)$ si et seulement si le morphisme
$$
E \longrightarrow j_* \circ j^* E
$$
est un isomorphisme.

\item Le plongement $f : {\mathcal E}' \hookrightarrow {\mathcal E}$ se factorise en une équivalence
$$
{\mathcal E}' \xrightarrow{ \ \sim \ } \mbox{\rm Im} \, (f)
$$
suivie du plongement $\mbox{\rm Im} \, (f) \hookrightarrow {\mathcal E}$.

\item La sous-catégorie pleine $\mbox{\rm Im} \, (f)$ ne dépend que de la classe d'équivalence du plongement $f : {\mathcal E}' \hookrightarrow {\mathcal E}$.
\end{listeimarge}
\end{prop}

\begin{esquissedemo}
Les foncteurs
$$
f_* : {\mathcal E}' \longrightarrow {\mathcal E}
$$
et
$$
f^* : {\mathcal E} \longrightarrow {\mathcal E}'
$$
se restreignent en deux foncteurs
$$
g_* :  {\mathcal E}' \longrightarrow \mbox{\rm Im} \, (f)
$$
et
$$ 
g^* : \mbox{\rm Im} \, (f) \longrightarrow {\mathcal E}'
$$
qui sont les deux composantes d'une équivalence de catégories.

\smallskip

\noindent (i) en résulte puisque ${\mathcal E}'$ est un topos.

\noindent (ii) en résulte également en prenant pour adjoint $j^* : {\mathcal E} \to \mbox{\rm Im} \, (f)$ de $j_* : \mbox{\rm Im} \, (f) \hookrightarrow {\mathcal E}$ le composé
$$
g_* \circ f^* : {\mathcal E} \longrightarrow {\mathcal E}' \xrightarrow{ \ \sim \ } \mbox{\rm Im} \, ({\mathcal E}) \, .
$$
(iii) Pour tout objet $E$ de ${\mathcal E}$, le morphisme
$$
E \longrightarrow j_* \circ j^* E \quad \mbox{est un isomorphisme}
$$
si et seulement si le morphisme
$$
E \longrightarrow f_* \circ f^* E \quad \mbox{est un isomorphisme}.
$$
Comme $f^* \circ f_* \to {\rm id}_{{\mathcal E}'}$ est un isomorphisme, c'est équivalent à demander que $E$ soit isomorphe à l'image par $f_*$ d'un objet de ${\mathcal E}'$.

\noindent (iv) et (v) sont évidents.
\end{esquissedemo}

Cette proposition amène à formuler la définition équivalente suivante de la notion de sous-topos:

\begin{defn}\label{defn1.2.6}

Un sous-topos d'un topos ${\mathcal E}$ est une sous-catégorie pleine
$$
{\mathcal E}' \xymatrix{\ar@{^{(}->}[r] & } {\mathcal E}
$$
telle que:

\smallskip

\noindent $\left\{\begin{matrix}
(0) &\mbox{${\mathcal E}'$ est un topos.} \hfill \\
(1) &\mbox{Le foncteur de plongement $j_* : {\mathcal E}' \hookrightarrow {\mathcal E}$ admet un adjoint à gauche $j^* : {\mathcal E} \to {\mathcal E}'$.} \hfill \\
(2) &\mbox{Cet adjoint à gauche $j^* : {\mathcal E} \to {\mathcal E}'$ respecte non seulement les colimites mais aussi les limites finies.} \hfill \\
(3) &\mbox{Un objet $E$ de ${\mathcal E}$ est dans la sous-catégorie pleine ${\mathcal E}'$ si et seulement si le morphisme canonique} \hfill \\
&E \longrightarrow j_* \circ j^* E \qquad \mbox{est un isomorphisme}.
\end{matrix}\right.
$
\end{defn}

\begin{exs}
\begin{listeisansmarge}
\item Pour tout site $({\mathcal C},J)$, $\widehat{\mathcal C}_J$ est un sous-topos de $\widehat{\mathcal C}$. De plus, pour toute topologie $K$ sur ${\mathcal C}$ qui raffine $J$, $\widehat{\mathcal C}_K$ est un sous-topos de $\widehat{\mathcal C}_J$.

\item Considérons un espace topologique $X$. Son topos associé ${\mathcal E}_X$ est la catégorie des faisceaux sur la petite catégorie ${\mathcal C}_X$ des ouverts de $X$, munie de la topologie $J_X$ qui consiste en la notion usuelle de recouvrement des ouverts.

Alors, pour tout sous-espace $Y$ de $X$, l'image du morphisme de plongement
$$
{\mathcal E}_Y \xymatrix{\ar@{^{(}->}[r] & } {\mathcal E}_X
$$
est le sous-topos de ${\mathcal E}_X$ défini par la topologie $J_Y$ de ${\mathcal C}_X$ pour laquelle une famille d'inclusions d'ouverts de $X$
$$
(U_i \subseteq U)_{i \in I}
$$
est couvrante si et seulement si
$$
\bigcup_{i \in I} (U_i \cap Y) = U \cap Y \, .
$$
\item Dans la situation de (ii), la catégorie ${\mathcal C}_X$ des ouverts de $X$ admet en général davantage de topologies que celles associées à des sous-espaces $Y$ de $X$, et donc ${\mathcal E}_X$ admet davantage de sous-topos que les images des ${\mathcal E}_X$.

Par exemple, ${\mathcal C}_X$ peut être munie de la ``topologie de la densité'' pour laquelle une famille d'inclusions d'ouverts de $X$
$$
(U_i \subseteq U)_{i \in I}
$$
est couvrante si et seulement la réunion
$$
\bigcup_{i \in I} U_i \quad \mbox{est un ouvert dense de $U$.}
$$
\end{listeisansmarge}
\end{exs}

\begin{remarkqed}
On verra au chapitre \ref{chap2} comment la notion de sous-topos ainsi définie et celle de topologie de Grothendieck apparaissent conjointement dans le cadre d'une dualité engendrée par une relation naturelle entre préfaisceaux et cribles sur une catégorie. 

Il en résultera en particulier que la propriété (0) se déduit des propriétés (1), (2) et (3). 
\end{remarkqed}

\subsection{L'expression des sous-topos en termes des topologies de Grothendieck}\label{ssec1.2.2}

On rappelle la proposition suivante:

\begin{prop}\label{prop1.2.7}

Considérons un morphisme de topos
$$
f = (f^* , f_*) : {\mathcal E}' \longrightarrow {\mathcal E} = \widehat{\mathcal C}_J
$$
vers le topos des faisceaux sur un site $({\mathcal C},J)$.

Disons qu'une famille de morphismes de la catégorie sous-jacente ${\mathcal C}$
$$
(X_i \longrightarrow X)_{i \in I}
$$
est $J_f$-couvrante si sa transformée par le foncteur
$$
f^* \circ \ell : {\mathcal C} \xrightarrow{ \ \ell \ } \widehat{\mathcal C}_J = {\mathcal E} \xrightarrow{ \ f^* \ } {\mathcal E}'
$$
composé du foncteur canonique $\ell : {\mathcal C} \to \widehat{\mathcal C}_J = {\mathcal E}$ et du foncteur d'image réciproque
$$
f^* : {\mathcal E} \longrightarrow {\mathcal E}' \, ,
$$
est une famille globalement épimorphique du topos ${\mathcal E}'$.

Alors:
\begin{listeimarge}
\item La notion de famille $J_f$-couvrante de morphismes de ${\mathcal C}$ définit une topologie de ${\mathcal C}$ qui raffine la topologie $J$.

\item La topologie $J_f$ de ${\mathcal C}$ définit un sous-topos
$$
\widehat{\mathcal C}_{J_f} \xymatrix{\ar@{^{(}->}[r] & } \widehat{\mathcal C}_J = {\mathcal E}
$$
de ${\mathcal E}$.

\item Le morphisme de topos
$$
f = (f^* , f_*) : {\mathcal E}' \longrightarrow \widehat{\mathcal C}_J
$$
se factorise en
$$
{\mathcal E}' \xrightarrow{ \ \overline f' \ } \widehat{\mathcal C}_{J_f} \xymatrix{\ar@{^{(}->}[r] & } \widehat{\mathcal C}_J
$$
où le morphisme de topos
$$
\overline f = (\overline f^* , \overline f_*) : {\mathcal E}' \longrightarrow \widehat{\mathcal C}_{J_f}
$$
est ``surjectif'' au sens où sa composante d'image réciproque
$$
\overline f^* : \widehat{\mathcal C}_{J_f} \longrightarrow {\mathcal E}'
$$
est un foncteur fidèle.

\item La factorisation d'un morphisme de topos
$$
f : {\mathcal E}' \longrightarrow {\mathcal E}
$$
en un morphisme surjectif $f'$ suivi d'un plongement
$$
{\mathcal E}' \xrightarrow{ \ f' \ } {\mathcal E}'' \xymatrix{\ar@{^{(}->}[r] & } {\mathcal E}
$$
est unique à équivalence près.

En particulier, tout morphisme de topos qui est à la fois un plongement et un morphisme surjectif est une équivalence.

\item La factorisation canonique d'un morphisme de topos
$$
f : {\mathcal E}' \longrightarrow {\mathcal E}
$$
en un morphisme surjectif suivi d'un plongement définit un sous-topos de ${\mathcal E}$ que l'on peut appeler l'image de $f$ et noter
$$
{\rm Im} \, (f) \xymatrix{\ar@{^{(}->}[r] & } {\mathcal E} \, .
$$
\end{listeimarge}
\end{prop}

\begin{esquissedemo}
\begin{listeisansmarge}
\item La notion de famille $J_f$-couvrante de morphismes de ${\mathcal C}$ respecte les propriétés (0), (1), (2) et (3) de la définition \ref{defn1.1.1} des topologies, car la notion de famille globalement épimorphique de ${\mathcal E}'$ vérifie ces propriétés.

Pour (0), (1) et (3), cette implication résulte simplement de ce que le foncteur
$$
f^* \circ \ell : {\mathcal C} \longrightarrow \widehat{\mathcal C}_J \longrightarrow {\mathcal E}'
$$
respecte les compositions et les morphismes d'identité.

Pour (2), l'implication résulte de ce que le foncteur composé de $f^* : \widehat{\mathcal C}_J \to {\mathcal E}'$ avec le foncteur de faisceautisation $j^* : \widehat{\mathcal C} \to \widehat{\mathcal C}_J$
$$
f^* \circ j^* : \widehat{\mathcal C} \longrightarrow \widehat{\mathcal C}_J \longrightarrow {\mathcal E}'
$$
respecte les limites finies, donc transforme les produits fibrés de $\widehat{\mathcal C}$ de la forme
$$
y(X') \times_{y(X)} y(X_i)
$$
en les produits fibrés
$$
f^* \circ \ell (X') \times_{f^* \circ \ell (X)} f^* \circ \ell (X_i) \, ,
$$
de ce que tout objet de $\widehat{\mathcal C}$, en particulier les $y (X') \times_{y(X)} y(X_i)$, s'écrit comme une colimite d'objets de la forme $y(X'_{i'})$, et de ce que le foncteur $f^* \circ j^* : \widehat{\mathcal C} \to {\mathcal E}'$ respecte les colimites.

Ainsi, $J_f$ est une topologie.

Enfin, la topologie $J_f$ contient la topologie $J$, car si une famille de morphismes de ${\mathcal C}$
$$
X_i \longrightarrow X \, , \quad i \in I \, ,
$$
est $J$-couvrante, la famille de ses transformés par $\ell$
$$
\ell (X_i) \longrightarrow \ell (X) \, , \quad i \in I \, ,
$$
est globalement épimorphique dans $\widehat{\mathcal C}_J$ et car le foncteur
$$
f^* : \widehat{\mathcal C}_J \longrightarrow {\mathcal E}'
$$
respecte les colimites, donc les familles globalement épimorphiques.

\item résulte de (i) puisque toute topologie de ${\mathcal C}$ qui raffine $J$ définit un sous-topos de $\widehat{\mathcal C}_J$.

\item Le foncteur composé
$$
{\mathcal C} \xrightarrow{ \ \ell \ } \widehat{\mathcal C}_J \xrightarrow{ \ f^* \ } {\mathcal E}'
$$
est plat et $J_f$-continu donc définit un morphisme de topos
$$
\overline f = (\overline f^* , \overline f_*) : {\mathcal E}' \longrightarrow \widehat{\mathcal C}_{J_f} \, .
$$
Son composé avec le plongement
$$
\widehat{\mathcal C}_{J_f} \xymatrix{\ar@{^{(}->}[r] & } \widehat{\mathcal C}_J
$$
s'identifie avec le morphisme
$$
f = (f^* , f_*) : {\mathcal E}' \longrightarrow \widehat{\mathcal C}_J
$$
puisque le composé
$$
{\mathcal C} \longrightarrow \widehat{\mathcal C}_{J_f} \xrightarrow{ \ \overline f^* \ } {\mathcal E}'
$$
se confond avec le composé ${\mathcal C} \xrightarrow{ \ \ell \ } \widehat{\mathcal C}_J \xrightarrow{ \ f^* \ } {\mathcal E}'$.

Pour montrer que le foncteur
$$
\overline f^* : \widehat{\mathcal C}_{J_f} \longrightarrow {\mathcal E}'
$$
est fidèle, on commence par remarquer qu'un morphisme de $\widehat{\mathcal C}_{J_f}$
$$
F_1 \longrightarrow F_2
$$
est un épimorphisme si et seulement si son transformé par $\overline f^*$
$$
\overline f^* \, F_1 \longrightarrow \overline f^* \, F_2
$$
est un épimorphisme de ${\mathcal E}'$.

Considérons alors deux morphismes $\widehat{\mathcal C}_{J_f}$ de même source et même but
$$
\xymatrix{ F_1 \dar[r]^{^{^{\mbox{\scriptsize $u$}}}}_{v} & F_2 } \, .
$$
On a
$$
\overline f^* (u) = \overline f^* (v)
$$
si et seulement si le monomorphisme de ${\mathcal E}'$ défini par $\overline f^* (u) \times \overline f^* (v)$
$$
\overline f^* \, F_1 \times_{\overline f^* F_2 \times \overline f^* F_2} \overline f^* \, F_2 \xymatrix{\ar@{^{(}->}[r] & } \overline f^* \, F_1
$$
est un épimorphisme.

Comme $\overline f^*$ respecte les limites finies, cela revient à demander que le monomorphisme de $\widehat{\mathcal C}_{J_f}$ défini par $u \times v$
$$
F_1 \times_{F_2 \times F_2} F_2 \xymatrix{\ar@{^{(}->}[r] & } F_1
$$
soit un épimorphisme, c'est-à-dire que $u=v$.

Ainsi, on a $\overline f^* (u) = \overline f^* (v)$ si et seulement si $u=v$. Cela signifie que le foncteur $\overline f^* : \widehat{\mathcal C}_{J_f} \to {\mathcal E}'$ est fidèle.

\item Considérons un morphisme de topos
$$
f : {\mathcal E}' \longrightarrow {\mathcal E} = \widehat{\mathcal C}_J 
$$
factorisé comme un composé de morphismes de topos
$$
{\mathcal E}' \xrightarrow{ \ f' \ } {\mathcal E}'' \xrightarrow{ \ f'' \ } {\mathcal E} = \widehat{\mathcal C}_J \, .
$$
Dans un topos, un morphisme $E_1 \to E_2$ est un épimorphisme si et seulement si les deux morphismes canoniques
$$
E_2 \textstyle\bigsqcup E_2 \rightrightarrows E_2 \bigsqcup_{E_1} E_2
$$
sont égaux.

Il en résulte que si $f'$ est un morphisme surjectif, c'est-à-dire si le foncteur
$$
f'^* = {\mathcal E}'' \longrightarrow {\mathcal E}'
$$
est fidèle, alors une famille de morphismes de ${\mathcal E}''$
$$
(E_i \longrightarrow E)_{i \in I}
$$
est globalement épimorphique si et seulement sa transformée par $f'^*$
$$
(f'^* E_i \longrightarrow f'^* E)_{i \in I}
$$
est globalement épimorphique.

Cela implique que les deux morphismes de topos
$$
f : {\mathcal E}' \longrightarrow \widehat{\mathcal C}_J
$$
et
$$
f'' : {\mathcal E}'' \longrightarrow \widehat{\mathcal C}_J
$$
définissent la même topologie $J_f$ de ${\mathcal C}$.

Alors $f''$ se factorise en
$$
{\mathcal E}'' \xrightarrow{ \ \overline f'' \ } \widehat{\mathcal C}_{J_f} \xymatrix{\ar@{^{(}->}[r] & } \widehat{\mathcal C}_J
$$
où $\overline f''$ est un morphisme surjectif de topos.

Si $f'' : {\mathcal E}'' \to \widehat{\mathcal C}_J$ est un plongement, sa factorisation
$$
\overline f'' : {\mathcal E}'' \longrightarrow \widehat{\mathcal C}_{J_f}
$$
est aussi un plongement.

Pour tout objet $F$ de $\widehat{\mathcal C}_{J_f}$, le morphisme canonique
$$
F \longrightarrow \overline f''_* \circ \overline f''^* \, F
$$
est transformé par le foncteur fidèle $\overline f''^*$ en un isomorphisme, donc il est un épimorphisme.

Puis le morphisme diagonal
$$
F \longrightarrow F \times_{\overline f_* \circ \overline f''^* F} F
$$
est transformé par $\overline f''^*$ en un isomorphisme, donc il est un épimorphisme, ce qui signifie que le morphisme canonique
$$
F \longrightarrow \overline f''_* \circ \overline f''^* \, F
$$
est toujours un isomorphisme, et donc que
$$
f'' : {\mathcal E}'' \longrightarrow \widehat{\mathcal C}_J
$$
est une équivalence de topos.

\item résulte des points précédents.
\end{listeisansmarge}
\end{esquissedemo}

Cette proposition implique le théorème suivant, qui exprime la notion de sous-topos en termes de topologies de Grothendieck:

\begin{thm}\label{thm1.2.8}

Soit un topos ${\mathcal E}$ écrit comme le topos des faisceaux sur un site $({\mathcal C} , J)$.

Alors:

\begin{listeimarge}

\item Associer à toute topologie $J'$ de ${\mathcal C}$ raffinant $J$ le sous-topos
$$
\widehat{\mathcal C}_{J'} \xymatrix{\ar@{^{(}->}[r] & } \widehat{\mathcal C}_J = {\mathcal E}
$$
définit une bijection de l'ensemble ordonné des topologies de ${\mathcal C}$ vers la classe des sous-topos de ${\mathcal E}$.

Cette bijection renverse la relation d'ordre.

\item La bijection en sens inverse associe à tout sous-topos
$$
f : {\mathcal E}' \xymatrix{\ar@{^{(}->}[r] & } {\mathcal E} = \widehat{\mathcal C}_J
$$
la topologie $J_f$ de ${\mathcal C}$ définie en appelant $J_f$-couvrantes les familles de morphismes de ${\mathcal C}$
$$
(X_i \longrightarrow X)_{i \in I}
$$
dont les transformées par le foncteur
$$
{\mathcal C} \longrightarrow \widehat{\mathcal C}_{J_f} \xrightarrow{ \ f^* \ } {\mathcal E}'
$$
sont globalement épimorphiques dans ${\mathcal E}'$.
\end{listeimarge}
\end{thm}

\begin{remarksqed}
\begin{listeisansmarge}
\item Ce théorème implique en particulier que, pour tout topos ${\mathcal E}$, ses sous-topos forment un ensemble.

On le note ${\rm ST} ({\mathcal E})$.

\item La relation d'ordre de ${\rm ST} ({\mathcal E})$ s'exprime de deux manières différentes suivant que l'on considère l'une ou l'autre des deux définitions de la notion de sous-topos.

Si deux sous-topos de ${\mathcal E}$ sont représentés par deux plongements de topos
$$
{\mathcal E}_1 \xymatrix{\ar@{^{(}->}[r] & }  {\mathcal E} \quad \mbox{et} \quad {\mathcal E}_2  \xymatrix{\ar@{^{(}->}[r] & } {\mathcal E} \, ,
$$
on a ${\mathcal E}_2 \geq {\mathcal E}_1$ si et seulement si le plongement ${\mathcal E}_1 \hookrightarrow {\mathcal E}$ se factorise à isomorphisme près en un composé
$$
{\mathcal E}_1 \xymatrix{\ar@{^{(}->}[r] & } {\mathcal E}_2 \xymatrix{\ar@{^{(}->}[r] & } {\mathcal E}
$$
d'un plongement ${\mathcal E}_1 \hookrightarrow {\mathcal E}_2$ avec le plongement ${\mathcal E}_2 \hookrightarrow {\mathcal E}$.

Si par contre deux tels sous-topos sont vus comme deux sous-catégories pleines ${\mathcal E}_1$ et ${\mathcal E}_2$ de ${\mathcal E}$ qui satisfont les conditions de la définition \ref{defn1.2.6}, on a
$$
{\mathcal E}_2 \geq {\mathcal E}_1
$$
si et seulement si tout objet de ${\mathcal E}$ qui est élément de ${\mathcal E}_1$ est aussi élément de ${\mathcal E}_2$.

\item Il résulte en particulier du théorème que si $({\mathcal C},J)$ et $({\mathcal D},K)$ sont deux sites reliés par une équivalence de topos
$$
\widehat{\mathcal C}_J \xrightarrow{ \ \sim \ } \widehat{\mathcal D}_K \, ,
$$
alors l'ensemble ordonné des topologies de ${\mathcal C}$ raffinant $J$ s'identifie à celui des topologies de ${\mathcal D}$ raffinant $K$.

Autrement dit, ces ensembles ordonnés constituent des invariants topossiques des sites.

\item On donnera au chapitre II une autre démonstration de la bijection du théorème, fondée sur une relation de dualité entre cribles et préfaisceaux sur une catégorie.

Cette autre démonstration aura l'avantage de faire apparaître naturellement la notion de sous-topos et celle de topologie de Grothendieck, sans que leurs définitions soient posées a priori.
\end{listeisansmarge} 
\end{remarksqed}

\subsection{L'expression des sous-topos en termes de théories quotients}\label{ssec1.2.3}

\medskip

On rappelle la définition suivante:

\begin{defn}\label{defn1.2.9}

Soient ${\mathbb T}_1$ et ${\mathbb T}_2$ deux théories géométriques du premier ordre de même signature $\Sigma$.
\begin{listeimarge}
\item On dit que ${\mathbb T}_1$ est une théorie quotient de ${\mathbb T}_2$ si tous les axiomes de ${\mathbb T}_2$ sont démontrables dans ${\mathbb T}_1$, autrement dit si tout séquent géométrique de $\Sigma$
$$
\varphi (\vec x) \vdash \psi (\vec x)
$$
qui est démontrable dans ${\mathbb T}_2$ est démontrable dans ${\mathbb T}_1$.

\item On dit que ${\mathbb T}_1$ et ${\mathbb T}_2$ sont équivalentes si chacune est une théorie quotient de l'autre, autrement dit si n'importe quel séquent géométrique de $\Sigma$
$$
\varphi (\vec x) \vdash \psi (\vec x)
$$
est démontrable dans l'une si et seulement si il est démontrable dans l'autre.
\end{listeimarge}
\end{defn}

\begin{remarkqed}
Si ${\mathbb T}_1$ est une théorie quotient de ${\mathbb T}_2$, alors pour tout topos ${\mathcal E}$, la catégorie des modèles
$$
{\mathbb T}_1\mbox{-mod} \, ({\mathcal E})
$$
s'identifie à une sous-catégorie pleine de
$$
{\mathbb T}_2\mbox{-mod} \, ({\mathcal E}) \, .
$$
Les plongements
$$
{\mathbb T}_1\mbox{-mod} \, ({\mathcal E}_1) \xymatrix{\ar@{^{(}->}[r] & }
{\mathbb T}_2\mbox{-mod} \, ({\mathcal E})
$$
sont nécessairement induits par un morphisme de topos reliant les topos classifiants de ${\mathbb T}_1$ et ${\mathbb T}_2$
$$
{\mathcal E}_{{\mathbb T}_1} \longrightarrow {\mathcal E}_{{\mathbb T}_2} \, .
$$
Ce morphisme correspond au modèle de ${\mathbb T}_2$ dans ${\mathcal E}_{{\mathbb T}_1}$ qui est l'image du modèle universel $U_{{\mathbb T}_1}$ de ${\mathbb T}_1$ dans ${\mathcal E}_{{\mathbb T}_1}$ par le plongement
$$
{\mathbb T}_1\mbox{-mod} \, ({\mathcal E}_{{\mathbb T}_1}) \xymatrix{\ar@{^{(}->}[r] & }
{\mathbb T}_2\mbox{-mod} \, ({\mathcal E}_{{\mathbb T}_1}) \, .
$$
\end{remarkqed}

Le théorème suivant exprime la notion de sous-topos en termes de la logique géométrique du premier ordre:

\begin{thm}\label{thm1.2.10}

Soient ${\mathbb T}$ une théorie géométrique du premier ordre, et ${\mathcal E}_{\mathbb T}$ son topos classifiant.

Alors:

\begin{listeimarge}

\item Pour toute théorie quotient ${\mathbb T}'$ de ${\mathbb T}$, le morphisme de topos naturel entre les topos classifiants
$$
{\mathcal E}_{{\mathbb T}'} \longrightarrow {\mathcal E}_{\mathbb T}
$$
est un plongement.

Autrement dit, il identifie ${\mathcal E}_{{\mathbb T}'}$ à un sous-topos de ${\mathcal E}_{\mathbb T}$.

\item Pour toutes théories quotients ${\mathbb T}_1$ et ${\mathbb T}_2$ de ${\mathbb T}$, on a 
$$
{\mathcal E}_{{\mathbb T}_2} \geq {\mathcal E}_{{\mathbb T}_1} \quad \mbox{comme sous-topos de ${\mathbb T}$}
$$
si et seulement si ${\mathbb T}_1$ est un quotient de ${\mathbb T}_2$.

En particulier, on a
$$
{\mathcal E}_{{\mathbb T}_1} = {\mathcal E}_{{\mathbb T}_2} \quad \mbox{comme sous-topos de ${\mathbb T}$}
$$
si et seulement si les théories quotients ${\mathbb T}_1$ et ${\mathbb T}_2$ sont équivalentes.

\item L'application
$$
{\mathbb T}' \longmapsto ({\mathcal E}_{{\mathbb T}'} \hookrightarrow {\mathcal E}_{\mathbb T})
$$
définit une bijection de l'ensemble des classes d'équivalence de théories quotients ${\mathbb T}'$ de ${\mathbb T}$ sur l'ensemble
$$
{\rm ST} \, ({\mathcal E}_{\mathbb T})
$$
des sous-topos de ${\mathcal E}_{\mathbb T}$.
\end{listeimarge}
\end{thm}

\begin{remarks}
\begin{listeisansmarge}
\item Ce théorème est présenté dans le paragraphe 3.2 du livre \cite{TST} qui en donne deux démonstrations. Il est tiré de la thèse de l'auteur.

\item Il résulte de ce théorème que si ${\mathbb T}_1$ et ${\mathbb T}_2$ sont deux théories géométriques du premier ordre reliées par une équivalence sémantique
$$
{\mathcal E}_{{\mathbb T}_1} \xrightarrow{ \ \sim \ } {\mathcal E}_{{\mathbb T}_2} \, ,
$$
celle-ci définit une bijection de l'ensemble des classes de théories quotients de ${\mathbb T}_1$ sur celui de ${\mathbb T}_2$.

Cette bijection respecte les relations d'ordre définies par la propriété de quotient d'une théorie par rapport à une autre.
\end{listeisansmarge}
\end{remarks}

\begin{esquissedemosansqed}
Le topos classifiant ${\mathcal E}_{\mathbb T}$ de ${\mathbb T}$ est construit comme le topos des faisceaux sur la catégorie syntactique ${\mathcal C}_{\mathbb T}$ de ${\mathbb T}$ munie de la topologie $J_{\mathbb T}$.

Le théorème se démontre en construisant explicitement

\smallskip

$
\left\{\begin{matrix}
\bullet &\mbox{une application qui associe à toute théorie quotient ${\mathbb T}'$ de ${\mathbb T}$ une topologie $J_{{\mathbb T}'}$ sur ${\mathcal C}_{\mathbb T}$ qui raffine $J_{\mathbb T}$,} \hfill \\
\bullet &\mbox{une application qui associe à toute topologie $J$ sur ${\mathcal C}_{\mathbb T}$ qui raffine $J_{\mathbb T}$ une théorie quotient ${\mathbb T}_J$ de ${\mathbb T}$,} \hfill \\
\end{matrix}\right.
$

\smallskip

\noindent et en montrant que

\smallskip

$
\left\{\begin{matrix}
\bullet &\mbox{deux théories quotients ${\mathbb T}_1$ et ${\mathbb T}_2$ de ${\mathbb T}$ définissent la même topologie $J_{{\mathbb T}_1} = J_{{\mathbb T}_2}$ de ${\mathcal C}_{\mathbb T}$} \hfill \\
&\mbox{si elles sont équivalentes,} \hfill \\
\bullet &\mbox{si ${\mathbb T}'$ est une théorie quotient de ${\mathbb T}$, la théorie quotient ${\mathbb T}_J$ associée à la topologie $J = J_{{\mathbb T}'}$} \hfill \\
&\mbox{définie par ${\mathbb T}'$ est équivalente à ${\mathbb T}'$,} \hfill \\
\bullet &\mbox{si $J$ est une topologie de ${\mathcal C}_{\mathbb T}$ qui raffine $J_{\mathbb T}$, la topologie $J_{{\mathbb T}'}$ définie par la théorie quotient} \hfill \\
&\mbox{${\mathbb T}' = {\mathbb T}_J$ de ${\mathbb T}$ associée à $J$ se confond avec $J$.} \hfill
\end{matrix}\right.
$

\smallskip

En effet, le théorème est alors conséquence de ce que, d'après le théorème \ref{thm1.2.8}, l'ensemble ordonné des sous-topos de ${\mathcal E}_{\mathbb T}$
$$
({\rm ST} \, ({\mathcal E}_{\mathbb T}) \, , \ \geq)
$$
s'identifie à l'ensemble des topologies de ${\mathcal C}_{\mathbb T}$ qui raffinent $J_{\mathbb T}$, muni de la relation d'inclusion $\subseteq$.

Pour associer à une théorie quotient ${\mathbb T}'$ de ${\mathbb T}$ une topologie $J_{{\mathbb T}'}$ raffinant $J_{\mathbb T}$ de la catégorie syntactique ${\mathcal C}_{\mathbb T}$, on utilise la notion de topologie engendrée par une collection de familles de morphismes:
\end{esquissedemosansqed}

\begin{defn}\label{defn1.2.11}

Soit ${\mathcal C}$ une catégorie petite ou essentiellement petite.

Soit ${\mathcal X}$ une collection de familles $c$ de morphismes de ${\mathcal C}$ de même but
$$
c = (X_i \xrightarrow{ \ x_i \ } X)_{i \in I_c} \, .
$$
On appelle topologie de ${\mathcal C}$ engendrée par ${\mathcal X}$ la plus petite topologie $J$ pour laquelle chaque famille
$$
c = (X_i \xrightarrow{ \ x_i \ } X)_{i \in I_c} \, , \quad c \in {\mathcal X} \, ,
$$
est couvrante.
\end{defn}

\begin{justifdefi}
Il suffit de remarquer que toute intersection de topologies d'une catégorie ${\mathcal C}$ est encore une topologie ${\mathcal C}$.

En effet, les quatre axiomes (0), (1), (2), (3) de la définition \ref{defn1.1.1} qui définissent la notion de topologie sont respectés par les intersections arbitraires. 
\end{justifdefi}

On rappelle que les objets de la catégorie syntactique ${\mathcal C}_{\mathbb T}$ d'une théorie géométrique du premier ordre ${\mathbb T}$ de signature $\Sigma$ sont les formules géométriques de $\Sigma$
$$
\varphi (\vec x)
$$
considérées à substitution près des variables $\vec x = (x_1^{A_1} , \cdots , x_n^{A_n})$ par des variables de mêmes sortes $A_1 , \cdots , A_n$.

Ses morphismes
$$
\varphi (\vec x) \longrightarrow \psi (\vec y)
$$
sont les formules géométriques en la réunion disjointe des variables
$$
\theta (\vec x , \vec y)
$$
qui sont ``démontrablement fonctionnelles'' dans la théorie ${\mathbb T}$.

Les deux applications en sens inverse qui réalisent la bijection entre l'ensemble des classes d'équivalence de théories quotients ${\mathbb T}'$ de ${\mathbb T}$ et celui des topologies $J$ sur ${\mathcal C}_{\mathbb T}$ raffinant $J_{\mathbb T}$ sont définies par les formules explicites suivantes:

\begin{defn}\label{defn1.2.12}

Soit ${\mathcal C}_{\mathbb T}$ la catégorie syntactique d'une théorie géométrique du premier ordre ${\mathbb T}$ de signature $\Sigma$.

Alors:

\begin{listeimarge}

\item Pour toute théorie quotient ${\mathbb T}'$ de ${\mathbb T}$ définie par une famille d'axiomes de la forme
$$
\varphi_i (\vec x_i) \vdash \psi_i (\vec x_i) \, , \quad i \in I \, ,
$$
notons $J_{{\mathbb T}'}$ la plus petite topologie de ${\mathcal C}_{\mathbb T}$ pour laquelle chaque monomorphisme de ${\mathcal C}_{\mathbb T}$
$$
(\varphi_i \wedge \psi_i) (\vec x_i) \xymatrix{\ar@{^{(}->}[r] & } \varphi_i (\vec x_i) \, , \quad i \in I \, ,
$$
est couvrant.

\item Pour toute topologie $J'$ de ${\mathcal C}_{\mathbb T}$ engendrée sur $J_{\mathbb T}$ par une collection ${\mathcal X}$ de familles de morphismes
$$
c = \left( \varphi_i (\vec x_i) \xrightarrow{ \ \theta_i (\vec x_i , \vec y_c) \ } \psi_c (\vec y_c)\right)_{i \in I_c} \, ,
$$
notons ${\mathbb T}_{\mathcal X}$ la théorie quotient de ${\mathbb T}$ définie en adjoignant aux axiomes de ${\mathbb T}$ la famille d'axiomes
$$
\psi_c (\vec y_c) \vdash \bigvee_{i \in I_c} (\exists \, \vec x_i) \, \theta_i (\vec x_i , \vec y_c) \, , \quad c \in {\mathcal X} \, .
$$
\end{listeimarge}
\end{defn} %\hfill $\Box$

Alors le théorème \ref{thm1.2.10} résulte de la proposition suivante:

\begin{prop}\label{prop1.2.13}

Dans la situation de la définition \ref{defn1.2.12} ci-dessus, on a:
\begin{listeimarge}
\item Pour toutes théories quotients ${\mathbb T}_1$ et ${\mathbb T}_2$ de ${\mathbb T}$ définies par deux familles d'axiomes, on a
$$
J_{{\mathbb T}_1} = J_{{\mathbb T}_2}
$$
si les axiomes de chacune se démontrent à partir des axiomes de l'autre en utilisant les règles d'inférence de la logique du premier ordre.

\item Pour toute topologie $J'$ de ${\mathcal C}_{\mathbb T}$ engendrée sur $J_{\mathbb T}$ par une collection ${\mathcal X}_1$ de générateurs ou par une autre ${\mathcal X}_2$, les deux théories quotients de ${\mathbb T}$ associées
$$
{\mathbb T}_{{\mathcal X}_1} \quad \mbox{et} \quad {\mathbb T}_{{\mathcal X}_2}
$$
sont équivalentes.

\item Tout séquent géométrique de $\Sigma$
$$
\varphi (\vec x) \vdash \psi (\vec x)
$$
est équivalent au séquent
$$
\varphi (\vec x) \vdash (\varphi \wedge \psi)(\vec x) \, .
$$

\item Pour toute collection ${\mathcal X}$ de familles de morphismes de ${\mathcal C}_{\mathbb T}$
$$
c = \left( \varphi_i (\vec x_i) \xrightarrow{ \ \theta_i (\vec x_i , \vec y_c) \ } \psi_c (\vec y_c)\right)_{i \in I_c} \, , \quad c \in {\mathcal X} \, ,
$$
la topologie $J'$ de ${\mathcal C}_{\mathbb T}$ engendrée sur $J_{\mathbb T}$ par cette collection ${\mathcal X}$ est aussi la plus petite topologie raffinant $J_{\mathbb T}$ pour laquelle les monomorphismes
$$
\bigvee_{i \in I_c} (\exists \, \vec x_i) \, \theta_i (\vec x_i , \vec y_c) \xymatrix{\ar@{^{(}->}[r] & } \psi_c (\vec y_c) \, , \quad c \in {\mathcal X} \, ,
$$
sont couvrants.
\end{listeimarge}
\end{prop}

\begin{indicationdemo}
\noindent (i) se démontre en vérifiant que toutes les règles d'inférence de la logique géométrique du premier ordre qui permettent de démontrer pas à pas des propriétés
$$
\varphi (\vec x) \vdash \psi (\vec x)
$$
à partir d'une famille d'axiomes
$$
\varphi_i (\vec x_i) \vdash \psi_i (\vec x_i) \, , \quad i \in I \, ,
$$
s'interprètent en termes du processus d'engendrement d'une topologie de ${\mathcal C}_{\mathbb T}$ à partir d'une collection de générateurs.

\noindent (ii) se démontre en vérifiant que les quatre axiomes (0), (1), (2), (3) de la définition \ref{defn1.1.1} des topologies s'interprètent en termes des règles d'inférence de la logique du premier ordre.

\noindent (iii) résulte de la règle d'inférence qui définit le sens du symbole de conjonction $\wedge$.

\noindent (iv) résulte de la définition de la topologie syntactique $J_{\mathbb T}$ de ${\mathcal C}_{\mathbb T}$ et des axiomes (0) et (3) de la définition \ref{defn1.1.1} de la notion de topologie. \end{indicationdemo}

\subsection{La démontrabilité comme un problème d'engendrement de topologies}\label{ssec1.2.4}

\medskip

Combiner le théorème \ref{thm1.2.10} et le théorème \ref{thm1.2.8} permet d'interpréter tous les problèmes de démontrabilité en logique géométrique du premier ordre comme des problèmes d'engendrements de topologies de Grothendieck à partir de générateurs constitués de familles de morphismes de même but dans des petites catégories sous-jacentes.

Commençons par introduire un procédé général de transformation de séquents géométriques en des familles de morphismes:

\begin{defn}\label{defn1.2.14}

Soit ${\mathbb T}$ une théorie géométrique du premier ordre de signature $\Sigma$.

Soit une présentation de son topos classifiant ${\mathcal E}_{\mathbb T}$ comme topos des faisceaux sur un site $({\mathcal C},J)$
$$
\widehat{\mathcal C}_J \xrightarrow{ \ \sim \ } {\mathcal E}_{\mathbb T}
$$
si bien que le modèle universel $U_{\mathbb T}$ de ${\mathbb T}$ dans ${\mathcal E}_{\mathbb T}$ s'interprète comme à valeurs dans les $J$-faisceaux sur ${\mathcal C}$.

Pour n'importe quel séquent géométrique de $\Sigma$
$$
\varphi (\vec x) \vdash \psi (\vec x)
$$
considérons l'interprétation dans $U_{\mathbb T}$ du séquent
$$
\varphi \wedge \psi (\vec x) \vdash \varphi (\vec x)
$$
comme un monomorphisme de $\widehat{\mathcal C}_J$
$$
U_{\mathbb T} (\varphi \wedge \psi) \xymatrix{\ar@{^{(}->}[r] & } U_{\mathbb T} \, \varphi \, .
$$
Choisissons une famille globalement épimorphique de morphismes de $\widehat{\mathcal C}_J$
$$
\ell (X_b) \longrightarrow U_{\mathbb T} \, \varphi \, , \quad b \in B_{\varphi}
$$
dont les sources $\ell (X_b)$ sont les images d'objets $X_b$ de ${\mathcal C}$ par le foncteur canonique
$$
\ell = j^* \circ y : {\mathcal C} \longrightarrow \widehat{\mathcal C} \longrightarrow \widehat{\mathcal C}_J \, .
$$
Puis choisissons pour tout indice $b$ une famille de morphismes de ${\mathcal C}$
$$
X_s \xrightarrow{ \ x_s \ } X_b \, , \quad s \in S_{\varphi , b , \psi}
$$
dont les transformés par $\ell$
$$
\ell (X_s) \xrightarrow{ \ \ell (x_s) \ } \ell (X_b)
$$
se factorisent à travers le sous-objet
$$
U_{\mathbb T} (\varphi \wedge \psi) \times_{U_{\mathbb T} \varphi} \ell (X_b)
$$
et forment une famille globalement épimorphique sur celui-ci. Notons alors
$$
{\mathcal X}_{\varphi , \psi}
$$
la famille indexée par les $b \in B_{\varphi}$ des familles de morphismes
$$
(X_s \xrightarrow{ \ x_s \ } X_b)_{s \in S_{\varphi , b , \psi}} 
$$
de but $X_b$ dans ${\mathcal C}$.
\end{defn}

\begin{remarksqed}
\begin{listeisansmarge}
\item Notant $j_* : \widehat{\mathcal C}_J \hookrightarrow \widehat{\mathcal C}$ le foncteur de plongement adjoint à droite du foncteur de faisceautisation $j^* : \widehat{\mathcal C} \to \widehat{\mathcal C}_J$, choisir des morphismes de $\widehat{\mathcal C}_J$ de la forme
$$
\ell (X) \longrightarrow F
$$
équivaut à choisir des morphismes de $\widehat{\mathcal C}$
$$
y(X) \longrightarrow j_* F
$$
c'est-à-dire d'après le lemme de Yoneda des éléments des ensembles de sections
$$
F(X) \, .
$$
Si $F$ est lui-même de la forme
$$
\ell (Y) \quad \mbox{pour un objet $Y$ de ${\mathcal C}$,}
$$
il existe pour tout morphisme $\ell (X) \to \ell (Y)$ identifié à un élément de $\ell (Y)(X)$ un $J$-recouvrement de l'objet $X$ par des morphismes de ${\mathcal C}$
$$
X_s \longrightarrow X
$$
tels que les éléments induits dans les $\ell (Y)(X_s)$ se relèvent en des éléments des $y(Y)(X_s) = {\rm Hom} (X_s , Y)$ c'est-à-dire tels que les morphismes composés
$$
\ell (X_s) \longrightarrow \ell (X) \longrightarrow \ell (Y)
$$
soient les images par $\ell$ de morphismes de ${\mathcal C}$
$$
X_s \longrightarrow Y \, .
$$
C'est pourquoi les choix demandés dans la définition sont possibles.

\item Si ${\mathcal C}$ est une petite sous-catégorie pleine de ${\mathcal E}_{\mathbb T}$ stable par produits fibrés et qui contient les interprétations dans ${\mathcal E}_{\mathbb T}$
$$
U_{\mathbb T} (\varphi \wedge \psi) \quad \mbox{et} \quad U_{\mathbb T} \, \varphi
$$
des formules $\varphi \wedge \psi$ et $\varphi$, on peut choisir pour ${\mathcal X}_{\varphi , \psi}$ la famille à un élément constitué de l'unique morphisme de ${\mathcal C}$
$$
U_{\mathbb T} (\varphi \wedge \psi) \xymatrix{\ar@{^{(}->}[r] & } U_{\mathbb T} \, \varphi \, .
$$
(iii) La condition de (ii) est vérifiée en particulier si ${\mathcal C}$ est une petite sous-catégorie pleine de ${\mathcal E}_{\mathbb T}$ stable par produits fibrés et qui contient la catégorie syntactique ${\mathcal C}_{\mathbb T}$ de ${\mathbb T}$.

Si ${\mathbb T}$ est une théorie cohérente [resp. régulière, resp. cartésienne] et si les formules $\varphi$ et $\psi$ sont cohérentes [resp. régulières, resp. ${\mathbb T}$-cartésiennes], il suffit que la sous-catégorie pleine ${\mathcal C}$ de ${\mathcal E}_{\mathbb T}$ contienne la catégorie syntactique cohérente [resp. régulière, resp. cartésienne] de ${\mathbb T}$. 
\end{listeisansmarge}
\end{remarksqed}

Nous pouvons maintenant énoncer la conséquence suivante des théorèmes \ref{thm1.2.10} et \ref{thm1.2.8}:

\begin{cor}\label{cor1.2.15}

Soit ${\mathbb T}$ une théorie géométrique du premier ordre de signature $\Sigma$.

Soit une présentation de son topos classifiant ${\mathcal E}_{\mathbb T}$ comme topos des faisceaux sur un site $({\mathcal C},J)$
$$
\widehat{\mathcal C}_J \xrightarrow{ \ \sim \ } {\mathcal E}_{\mathbb T} \, .
$$
Considérant des séquents géométriques de $\Sigma$
$$
\varphi_i (\vec x_i) \vdash \psi_i (\vec x_i) \, , \quad i \in I \, ,
$$
et
$$
\varphi \vdash \psi \, ,
$$
associons à ces séquents des familles de familles de morphismes de ${\mathcal C}$ de même but
$$
{\mathcal X}_{\varphi_i , \psi_i} \, , \quad i \in I \, ,
$$
et
$$
{\mathcal X}_{\varphi , \psi}
$$
comme dans la définition \ref{defn1.2.14}.

Alors le séquent géométrique de $\Sigma$
$$
\varphi (\vec x) \vdash \psi (\vec x)
$$
est démontrable à partir des axiomes de ${\mathbb T}$ et des séquents géométriques
$$
\varphi_i (\vec x_i) \vdash \psi_i (\vec x_i) \, , \quad i \in I \, ,
$$
si et seulement si chaque famille de morphismes de ${\mathcal C}$
$$
c \in {\mathcal X}_{\varphi , \psi}
$$
est couvrante pour la topologie $J'$ de ${\mathcal C}$ engendrée par $J$ et par les éléments des familles
$$
{\mathcal X}_{\varphi_i , \psi_i} \, , \quad i \in I \, .
$$
\end{cor}

\begin{remarks}
\begin{listeisansmarge}
\item C'est ainsi que tout problème de démontrabilité dans le contexte d'une théorie géométrique du premier ordre peut être réinterprété, d'une infinité de manières différentes, en des problèmes d'engendrement de topologies sur des petites catégories à partir de familles de générateurs constitués de familles de morphismes de même but.

La possibilité de telles traductions est d'autant plus remarquable et intéressante que les notions de catégorie, de topologie sur une catégorie (avec les quatre axiomes qui la définissent) et d'engendrement d'une topologie à partir de générateurs se définissent en quelques lignes, alors que les notions de théorie géométrique du premier ordre et de démontrabilité de séquents géométriques (qui demande d'énoncer la quinzaine des règles d'inférence de la logique géométrique) demandent de nombreuses pages de présentation et d'explications.

\item Dans le cas où ${\mathbb T}$ est une catégorie ``de type préfaisceau'', on peut en particulier choisir une petite catégorie ${\mathcal C}$ et une équivalence de topos
$$
\widehat{\mathcal C} \xrightarrow{ \ \sim \ } {\mathcal E}_{\mathbb T} \, .
$$
Alors un séquent géométrique
$$
\varphi (\vec x) \vdash \psi (\vec x)
$$
est démontrable à partir de la famille de séquents géométriques
$$
\varphi_i (\vec x_i) \vdash \psi_i (\vec x_i) \, , \quad i \in I \, ,
$$
si et seulement si toute famille de morphismes de ${\mathcal C}$ de même but qui est élément de la famille
$$
{\mathcal X}_{\varphi , \psi}
$$
est couvrante pour la topologie de ${\mathcal C}$ qui est engendrée par les éléments des familles
$$
{\mathcal X}_{\varphi_i , \psi_i} \, , \quad i \in I \, .
$$
\item La condition de (ii) est vérifiée si ${\mathbb T}$ est une théorie cartésienne. Dans ce cas, on dispose d'une équivalence
$$
\widehat{\mathcal C} \longrightarrow {\mathcal E}_{\mathbb T}
$$
si ${\mathcal C} = {\mathcal C}_{\mathbb T}^{\rm car}$ désigne la catégorie syntactique cartésienne de ${\mathbb T}$.

\item La condition de (iii) et donc aussi celle de (ii) avec
$$
{\mathcal C} = {\mathcal C}_{\mathbb T}^{\rm car}
$$
sont satisfaites si ${\mathbb T}$ est une théorie sans axiome, réduite à sa signature $\Sigma$.
\end{listeisansmarge}
\end{remarks}

\begin{esquissedemo}
Considérons donc une présentation par un site $({\mathcal C} , J)$ du topos classifiant ${\mathcal E}_{\mathbb T}$ d'une théorie géométrique ${\mathbb T}$
$$
\widehat{\mathcal C}_J \xrightarrow{ \ \sim \ } {\mathcal E}_{\mathbb T} \, .
$$
Se donner une topologie $J'$ de ${\mathcal C}$ raffinant $J$ équivaut à se donner un sous-topos
$$
{\mathcal E}' \xymatrix{\ar@{^{(}->}[r] & } {\mathcal E}_{\mathbb T}
$$
et donc encore à se donner une théorie quotient ${\mathbb T}'$ de ${\mathbb T}$, à équivalence près, avec alors un carré commutatif à isomorphisme près:
$$
\xymatrix{
\widehat{\mathcal C}_{J'} \ar[d]^{\wr} \ar@{^{(}->}[rr]^{(j^* , j_*)} &&\widehat{\mathcal C}_J \ar[d]^{\wr} \\
{\mathcal E}_{{\mathbb T}'} \ar@{^{(}->}[rr] &&{\mathcal E}_{\mathbb T}
}
$$

De plus, le modèle universel $U_{{\mathbb T}'}$ de ${\mathbb T}'$ dans $\widehat{\mathcal C}_{J'}$ est le transformé du modèle universel $U_{\mathbb T}$ de ${\mathbb T}$ dans $\widehat{\mathcal C}_J$ par le foncteur de $J'$-faisceautisation
$$
j^* : \widehat{\mathcal C}_J \longrightarrow \widehat{\mathcal C}_{J'} \, .
$$
Alors un séquent géométrique de la signature $\Sigma$ de ${\mathbb T}$
$$
\varphi (\vec x) \vdash \psi (\vec x)
$$
est démontrable dans la théorie quotient ${\mathbb T}'$ si et seulement si le monomorphisme de $\widehat{\mathcal C}_{J'}$
$$
U_{{\mathbb T}'} (\varphi \wedge \psi) \xymatrix{\ar@{^{(}->}[r] & } U_{{\mathbb T}'} \varphi
$$
est un épimorphisme, et donc si et seulement si le monomorphisme de $\widehat{\mathcal C}_J$
$$
U_{\mathbb T} (\varphi \wedge \psi) \xymatrix{\ar@{^{(}->}[r] & } U_{\mathbb T} \, \varphi
$$
est transformé en un épimorphisme par le foncteur de $J'$-faisceautisation
$$
j^* : \widehat{\mathcal C}_J \longrightarrow \widehat{\mathcal C}_{J'} \, .
$$
Cela revient à demander que la collection
$$
{\mathcal X}_{\varphi , \psi}
$$
associée au séquent géométrique $\varphi \vdash \psi$ par le procédé de la définition \ref{defn1.2.14} soit constituée de familles de morphismes de même but
$$
(X_s \xrightarrow{ \ x_s \ } X_b)_{s \in S_{\varphi , b , \psi}} \, , \quad b \in B_{\varphi} \, ,
$$
qui sont couvrantes pour la topologie $J'$.

La conclusion du corollaire s'en déduit.
\end{esquissedemo}

Les traductions du corollaire signifient que le problème général d'engendrement de topologies sur des petites catégories -- et même comme on verra au chapitre III sur des petites catégories entièrement explicitables -- a la même portée que le problème général de démontrabilité en logique géométrique du premier ordre.

Il est indécidable en général.

\section{Les opérations géométriques sur les sous-topos}\label{sec1.3}

\subsection{Les opérations internes d'intersection et de réunion de sous-topos}\label{ssec1.3.1}

\medskip

Dans tout topos, l'ensemble ordonné de ses sous-topos possède des opérations bien définies de réunion et d'intersection de n'importe quelle famille d'éléments, finie ou infinie.

\begin{lem}\label{lem1.3.1}

Soit un topos ${\mathcal E}$.

Soit $({\mathcal E}_i \hookrightarrow {\mathcal E})_{i \in I}$ une famille d'éléments de l'ensemble ${\rm ST} ({\mathcal E})$ des sous-topos de ${\mathcal E}$ ordonné par la relation~$\geq$.

Alors : 

\begin{listeimarge}

\item Il existe un unique élément
$$
\bigvee_{i \in I} {\mathcal E}_i \xymatrix{\ar@{^{(}->}[r] & } {\mathcal E} \qquad \mbox{de} \quad {\rm ST} ({\mathcal E}) \, ,
$$
caractérisé par la propriété que, pour tout sous-topos ${\mathcal E}' \hookrightarrow {\mathcal E}$, on a
$$
{\mathcal E}' \geq \bigvee_{i \in I} {\mathcal E}_i
$$
si et seulement si
$$
{\mathcal E}' \geq {\mathcal E}_i \, , \quad \forall \, i \in I \, .
$$
\item Il existe un unique élément
$$
\bigwedge_{i \in I} {\mathcal E}_i \xymatrix{\ar@{^{(}->}[r] & } {\mathcal E} \qquad \mbox{de} \quad {\rm ST} ({\mathcal E}) \, ,
$$
caractérisé par la propriété que, pour tout sous-topos ${\mathcal E}' \hookrightarrow {\mathcal E}$, on a
$$
\bigwedge_{i \in I} {\mathcal E}_i \geq {\mathcal E}'
$$
si et seulement si
$$
{\mathcal E}_i \geq {\mathcal E}' \, , \quad \forall \, i \in I \, .
$$
\end{listeimarge}
\end{lem}

\begin{remarks}
\begin{listeisansmarge}
\item Si ${\mathcal E}$ est le topos $\widehat{\mathcal C}_J$ des faisceaux sur un site $({\mathcal C},J)$, et les sous-topos ${\mathcal E}_i \hookrightarrow {\mathcal E}$ correspondent à des topologies $J_i \supseteq J$ sur ${\mathcal C}$, alors on a:

\medskip

$\left\{\begin{matrix}
\bullet &\mbox{Le sous-topos $\underset{i \in I}{\bigvee} \, {\mathcal E}_i$ est associé à la topologie} \hfill \\
&\underset{i \in I}{\bigcap} \, J_i \\
&\mbox{pour laquelle une famille de morphismes de ${\mathcal C}$ de même but est couvrante si et seulement si} \hfill \\
&\mbox{elle est couvrante pour chaque topologie $J_i$, $i \in I$.} \hfill \\
\bullet &\mbox{Le sous-topos $\underset{i \in I}{\bigwedge} \, {\mathcal E}_i$ est associé à la topologie de ${\mathcal C}$ engendrée par les topologies $J_i$, $i \in I$,} \hfill \\
&\mbox{c'est-à-dire la plus petite topologie pour laquelle toutes les familles de morphismes $J_i$-couvrantes} \hfill \\
&\mbox{pour tous les indices $i \in I$ sont couvrantes.} \hfill
\end{matrix}\right.$

\medskip

Si chaque topologie $J_i$, $i \in I$, est engendrée par une collection ${\mathcal X}_i$ de familles de morphismes de ${\mathcal C}$ de même but, alors la topologie qui correspond à l'intersection $\underset{i \in I}{\bigwedge} \, {\mathcal E}_i$ est aussi engendrée par la réunion des collections ${\mathcal X}_i$, $i \in I$.

%\smallskip

\item Si ${\mathcal E}$ est le topos classifiant ${\mathcal E}_{\mathbb T}$ d'une théorie géométrique du premier ordre ${\mathbb T}$, et les sous-topos ${\mathcal E}_i \hookrightarrow {\mathcal E}$ correspondent à des théories quotients ${\mathbb T}_i$ de ${\mathbb T}$, $i \in I$, définies chacune par une famille d'axiomes supplémentaires ${\mathcal A}_i$, alors on a:

\medskip

$\left\{\begin{matrix}
\bullet &\mbox{Le sous-topos $\underset{i \in I}{\bigwedge} \, {\mathcal E}_i$ est associé à la théorie quotient de ${\mathbb T}$ définie en adjoignant aux axiomes de ${\mathbb T}$} \hfill \\
&\mbox{ les axiomes éléments de la réunion $\underset{i \in I}{\bigcup} \, {\mathcal A}_i$ de toutes les familles ${\mathcal A}_i$.} \hfill \\
\bullet &\mbox{Le sous-topos $\underset{i \in I}{\bigvee} \, {\mathcal E}_i$ est associé à la théorie quotient de ${\mathbb T}$ pour laquelle un séquent géométrique} \hfill \\
&\varphi (\vec x) \vdash \psi (\vec x) \\
&\mbox{est démontrable si et seulement si il est démontrable dans chaque théorie ${\mathbb T}_i$, $i \in I$, c'est-à-dire} \hfill \\
&\mbox{à partir de chaque famille d'axiomes ${\mathcal A}_i$ adjointe à la liste des axiomes de ${\mathbb T}$.} \hfill
\end{matrix}\right.$
\end{listeisansmarge}
\end{remarks}

\begin{demo}

Elle résulte de la correspondance bijective entre sous-topos d'un topos ${\mathcal E}$ présenté sous la forme ${\mathcal E} \cong \widehat{\mathcal C}_J$ et topologies de ${\mathcal C}$ raffinant $J$, qui fait l'objet de la première remarque.
\begin{listeisansmarge}
\item résulte de ce que toute intersection de topologies est une topologie.

\item résulte de ce que toute collection de familles de morphismes de ${\mathcal C}$ de même but engendre une topologie, qui est l'intersection de toutes les topologies qui les contiennent. 
\end{listeisansmarge}
\end{demo}

Rappelons déjà:

\begin{prop}\label{prop1.3.2}

Dans l'ensemble ordonné ${\rm ST}({\mathcal E})$ des sous-topos d'un topos ${\mathcal E}$, les opérations de réunion et d'intersection sont liées par les relations suivantes:
\begin{listeimarge}
\item L'opération ${\mathcal E}' \vee \bullet$ de réunion avec un sous-topos ${\mathcal E}'$ respecte les intersections arbitraires au sens que
$$
{\mathcal E}' \vee \left( \bigwedge_{i \in I} {\mathcal E}_i \right) = \bigwedge_{i \in I} ({\mathcal E}' \vee {\mathcal E}_i)
$$
pour toute famille $({\mathcal E}_i)_{i \in I}$ de sous-topos.

\item L'opération ${\mathcal E}' \wedge \bullet$ d'intersection avec un sous-topos ${\mathcal E}'$ respecte les réunions finies au sens que
$$
{\mathcal E}' \wedge ({\mathcal E}_1 \vee \cdots \vee {\mathcal E}_n) = ({\mathcal E}' \wedge {\mathcal E}_1) \vee \cdots \vee ({\mathcal E}' \wedge {\mathcal E}_n)
$$
pour toute famille finie de sous-topos $({\mathcal E}_i)_{1 \leq i \leq n}$.
\end{listeimarge}
\end{prop}

\begin{pourdemo}
	\begin{listeisansmarge}
\item est conséquence formelle de ce que le foncteur
\begin{eqnarray}
{\rm ST} ({\mathcal E}) &\longrightarrow &{\rm ST}({\mathcal E}) \nonumber \\
{\mathcal E}_1 &\longmapsto &{\mathcal E}' \vee {\mathcal E}_1 \nonumber
\end{eqnarray}
admet un adjoint à droite
$$
{\mathcal E}_2 \longmapsto {\mathcal E}_2 \backslash {\mathcal E}' \, ,
$$
caractérisé par l'équivalence
$$
{\mathcal E}' \vee {\mathcal E}_1 \geq {\mathcal E}_2 \Longleftrightarrow {\mathcal E}_1 \geq {\mathcal E}_2 \backslash {\mathcal E}' \, ,
$$
comme il sera rappelé au paragraphe suivant.

\item est conséquence formelle de (i).

Il suffit en effet de montrer que
$$
{\mathcal E}' \wedge ({\mathcal E}_1 \vee {\mathcal E}_2) = ({\mathcal E}' \wedge {\mathcal E}_1) \vee ({\mathcal E}' \wedge {\mathcal E}_2) \, .
$$
Or, on a d'après (i)
\begin{eqnarray}
({\mathcal E}' \wedge {\mathcal E}_1) \vee ({\mathcal E}' \wedge {\mathcal E}_2) &= &({\mathcal E}' \vee {\mathcal E}') \wedge ({\mathcal E}_1 \vee {\mathcal E}') \wedge ({\mathcal E}' \vee {\mathcal E}_2) \wedge ({\mathcal E}_1 \vee {\mathcal E}_2) \nonumber \\
&= &{\mathcal E}' \wedge ({\mathcal E}_1 \vee {\mathcal E}_2) \, . \nonumber
\end{eqnarray}
\end{listeisansmarge}
\end{pourdemo}

\subsection{L'opération interne de différence de deux sous-topos}\label{ssec1.3.2}

Rappelons le théorème d'existence et de caractérisation de l'opération de différence des paires de sous-topos d'un topos donné:

\begin{thm}\label{thm1.3.3}

Soit un topos ${\mathcal E}$.

Pour tout sous-topos ${\mathcal E}_1 \hookrightarrow {\mathcal E}$, le foncteur de réunion des sous-topos de ${\mathcal E}$ avec ${\mathcal E}_1$
\begin{eqnarray}
{\rm ST} ({\mathcal E}) &\longrightarrow &{\rm ST} ({\mathcal E}) \, , \nonumber \\
{\mathcal E}' &\longmapsto &{\mathcal E}_1 \vee {\mathcal E}' \nonumber
\end{eqnarray}
possède un adjoint à droite
\begin{eqnarray}
{\rm ST} ({\mathcal E}) &\longrightarrow &{\rm ST} ({\mathcal E}) \, , \nonumber \\
{\mathcal E}_2 &\longmapsto &{\mathcal E}_2 \backslash {\mathcal E}_1 \nonumber
\end{eqnarray}
caractérisé par la propriété que, pour tout sous-topos ${\mathcal E}'$ de ${\mathcal E}$, on a une équivalence
$$
{\mathcal E}_1 \vee {\mathcal E}' \geq {\mathcal E}_2 \Longleftrightarrow {\mathcal E}' \geq {\mathcal E}_2 \backslash {\mathcal E}_1 \, .
$$
\end{thm}

\begin{demosansqed}

Le topos ${\mathcal E}$ est équivalent au topos $\widehat{\mathcal C}_J$ des faisceaux sur un site $({\mathcal C},J)$.

Compte-tenu du théorème \ref{thm1.2.8}, on est ramené à l'énoncé suivant qui donne de plus une formule de calcul en termes de topologies:
\end{demosansqed}

\begin{prop}\label{prop1.3.4}

Soit une catégorie essentiellement petite ${\mathcal C}$.

\begin{listeimarge}

\item Pour toute topologie $J_1$ sur ${\mathcal C}$, le foncteur d'intersection des topologies $J'$ de ${\mathcal C}$ avec $J_1$
$$
J' \longmapsto J_1 \cap J'
$$
possède un adjoint à gauche
$$
J_2 \longmapsto (J_1 \Rightarrow J_2)
$$
caractérisé par la propriété que, pour toute topologie $J'$ de ${\mathcal C}$, on a l'équivalence
$$
J_1 \cap J' \subseteq J_2 \Longleftrightarrow J' \subseteq (J_1 \Rightarrow J_2) \, .
$$
\item De plus, une famille de morphismes de ${\mathcal C}$ de même but $X$
$$
\Bigl( X_i \xrightarrow{ \ x_i \ } X \Bigl)_{i \in I}
$$
est couvrante pour la topologie $(J_1 \Rightarrow J_2)$ si et seulement si, pour tout morphisme $X' \xrightarrow{ \, x \, } X$ de but $X$ sur ${\mathcal C}$, n'importe quelle famille $J_1$-couvrante de but $X'$
$$
\Bigl( X'_k \xrightarrow{ \ x'_k \ } X' \Bigl)_{k \in K}
$$
est aussi $J_2$-couvrante si elle satisfait la condition:

\smallskip

$
\left\{\begin{matrix}
\mbox{Tout morphisme $X'' \to X'$ qui s'inscrit dans un carré commutatif} \hfill \\
\xymatrix{
X'' \ar[d] \ar[r] &X_i \ar[d]^{x_i} \\
X' \ar[r]^x &X
} \\
\mbox{se factorise en au moins un $X'' \longrightarrow X'_k \xrightarrow{ \, x'_k \, } X'$ avec $k \in K$.} \hfill
\end{matrix}\right.
$
\end{listeimarge}
\end{prop}

\begin{remarks}
\begin{listeisansmarge}
\item Si les topologies $J_1$ et $J_2$ contiennent toutes deux une topologie $J$, on a nécessairement
$$
J \subseteq (J_1 \Rightarrow J_2)
$$
puisque
$$
J_1 \cap J = J \subseteq J_2 \, .
$$
\item Il en résulte que si deux sous-topos
$$
{\mathcal E}_1 \quad \mbox{et} \quad {\mathcal E}_2 \quad \mbox{de} \quad {\mathcal E} \cong \widehat{\mathcal C}_J
$$
sont définis par deux topologies $J_1$ et $J_2$ contenant $J$ de ${\mathcal C}$, alors leur différence
$$
{\mathcal E}_2 \backslash {\mathcal E}_1
$$
est le sous-topos défini par la topologie contenant $J$ de ${\mathcal C}$
$$
(J_1 \Rightarrow J_2) \, .
$$
\end{listeisansmarge}
\end{remarks}

\begin{esquissedemopropsansqed}
Disons qu'une famille de morphismes de ${\mathcal C}$ de même but
$$
\Bigl( X_i \xrightarrow{ \ x_i \ } X \Bigl)_{i \in I}
$$
est $J_0$-couvrante si et seulement si elle satisfait la propriété de (ii).

Il faut montrer que $J_0$ est une topologie, qu'elle satisfait la condition
$$
J_1 \cap J_0 \subseteq J_2
$$
et que toute topologie $J$ de ${\mathcal C}$ qui satisfait la condition
$$
J_1 \cap J \subseteq J_2
$$
est contenue dans $J_0$.

Montrons d'abord que $J_0$ est une topologie, c'est-à-dire satisfait les axiomes (0), (1), (2) et (3) de la définition \ref{defn1.1.1}.

L'axiome (0) est vérifié car

\smallskip

$
\left\{\begin{matrix}
\bullet &\mbox{toute famille de morphismes de ${\mathcal C}$ de même but $X$ qui contient une famille $J_0$-couvrante} \hfill \\
&\mbox{est également $J_0$-couvrante} \hfill \\
\bullet &\mbox{la propriété pour une famille $\bigl( X_i \xrightarrow{ \ x_i \ } X \bigl)_{i \in I}$ d'être $J_0$-couvrante ne dépend que du crible} \hfill \\
&\mbox{engendré par cette famille.} \hfill
\end{matrix}\right.
$

\smallskip

\noindent L'axiome (1) est vérifié car le crible maximal de n'importe quel objet $X'$ de ${\mathcal C}$ est couvrant pour n'importe quelle topologie. L'axiome (2) est vérifié car la propriété qui définit la notion de famille $J_0$-couvrante $\bigl( X_i \xrightarrow{ \ x_i \ } X \bigl)_{i \in I}$ d'un objet $X$ comporte une quantification sur tous les morphismes
$$
X' \xrightarrow{ \ x \ } X
$$
et qu'on dispose du préfaisceau, c'est-à-dire du foncteur contravariant
$$
\Omega : {\mathcal C}^{\rm op} \longrightarrow {\rm Ens}
$$
qui associe

\smallskip

$
\left\{\begin{matrix}
\bullet &\mbox{à tout objet $X$ de ${\mathcal C}$ l'ensemble $\Omega (X)$ des cribles $C$ sur $X$,} \hfill \\
\bullet &\mbox{à tout morphisme $X' \xrightarrow{ \, x \, } X$ l'application} \hfill \\
&x^* : \Omega (X) \longrightarrow \Omega (X') \\
&\mbox{qui transforme tout crible $C$ engendré par une famille} \hfill \\
&\Bigl( X_i \xrightarrow{ \ x_i \ } X \Bigl)_{i \in I} \\
&\mbox{en le crible sur $X'$ constitué des morphismes} \hfill \\
&X'' \longrightarrow X' \\
&\mbox{qui s'inscrivent dans au moins un carré commutatif de la forme:} \hfill \\
&\xymatrix{
X'' \ar[d] \ar[r] &X_i \ar[d]^{x_i} \\
X' \ar[r]^x &X
}
\end{matrix}\right.
$

\medskip

Pour démontrer que $J_0$ satisfait l'axiome (3) de composition des familles couvrantes, on a besoin du lemme suivant:
\end{esquissedemopropsansqed}

\begin{lem}\label{lem1.3.5}

Soit un site $({\mathcal C},J)$.

Pour  tout crible $C$ sur un objet $X$ de ${\mathcal C}$, notons $\overline C$ la collection des morphismes de but $X$
$$
X' \xrightarrow{ \ x \ } X
$$
qui sont $J$-localement dans $C$ au sens qu'existe une famille $J$-couvrante de morphismes de but $X'$
$$
\Bigl( X'_i \xrightarrow{ \ x_i \ } X' \Bigl)_{i \in I}
$$
telle que tous les morphismes composés
$$
x \circ x_i : X'_i \longrightarrow X' \longrightarrow X
$$
soient éléments de $C$.

Alors:

\begin{listeimarge}

\item Pour tout crible $C$ sur un objet $X$, $\overline C$ est un crible sur $X$ qui contient $C$.

\item On a toujours $\overline{\overline C} = \overline C$.

\item Pour tout morphisme $X' \xrightarrow{ \, x \, } X$ de ${\mathcal C}$ et tout crible $C$ sur $X$, on a
$$
\overline{x^* (C)} = x^* (\overline C) \, .
$$
\end{listeimarge}
\end{lem}

\begin{remarks}
	\begin{listeisansmarge}
\item On dit qu'un crible $C$ est $J$-fermé si $\overline C = C$.

\item Il résulte de (i) et (ii) que, pour tout crible $C,\overline C$ est le plus petit crible $J$-fermé qui contient $C$. On l'appelle la $J$-fermeture de $C$.

\item Un crible $C$ sur un objet $X$ est $J$-couvrant si et seulement si sa fermeture $\overline C$ est le crible maximal.
\end{listeisansmarge}
\end{remarks}

\begin{esquissedemolem}
\begin{listeisansmarge}
\item résulte des axiomes (2) et (1) des topologies.

\item résulte de l'axiome (3) des topologies.

\item résulte de la formule de définition de $\overline C$ à partir de $C$.
\end{listeisansmarge} 
\end{esquissedemolem}

\noindent {\bf Fin de l'esquisse de démonstration de la proposition \ref{prop1.3.4}:}

\smallskip

Pour vérifier que $J_0$ est une topologie, il reste à montrer que $J_0$ satisfait l'axiome (3) des topologies.

Considérons donc une famille $J_0$-couvrante d'un objet $X$
$$
\Bigl( X_i \xrightarrow{ \ x_i \ } X \Bigl)_{i \in I}
$$
et, pour chaque indice $i \in I$, une famille $J_0$-couvrante de $X_i$
$$
\Bigl( X_{i,j} \xrightarrow{ \ x_{i,j} \ } X_i \Bigl)_{j \in I_i} \, , \quad i \in I \, .
$$

Il faut montrer que la famille des composés
$$
x_i \circ x_{i,j} : X_{i,j} \longrightarrow X_i \longrightarrow X \, , \quad j \in I_i \, , \quad i \in I \, ,
$$
est encore $J_0$-couvrante.

Comme on sait déjà que $J_0$ satisfait l'axiome (2), il suffit de prouver que si une famille de morphismes
$$
\Bigl( X'_k \xrightarrow{ \ x'_k \ } X \Bigl)_{k \in K}
$$
est $J_1$-couvrante et engendre un crible $C$ qui contient tous les morphismes
$$
x_i \circ x_{i,j} : X_{i,j} \longrightarrow X_i \longrightarrow X \, , \quad j \in I_i \, , \quad i \in I \, ,
$$
alors cette famille est aussi $J_2$-couvrante.

Notons $\overline C$ la $J_2$-fermeture du crible $C$ engendré par la famille $(X'_k \xrightarrow{ \, x'_k \, } X)_{k \in K}$.

Pour tout indice $i \in I$, le crible $x_i^* C$ sur $X_i$ est $J_1$-couvrant et contient les morphismes
$$
X_{i,j} \xrightarrow{ \ x_{i,j} \ } X_i \, , \quad j \in I_i \, .
$$
Comme ceux-ci composent une famille $J_0$-couvrante de $X_i$, il résulte de la définition de $J_0$ que le crible $x_i^* C$ est $J_2$-couvrant, autrement dit que $x_i^* \overline C$ est le crible maximal de $X_i$ c'est-à-dire
$$
x_i \in \overline C \, .
$$

Ainsi, le crible $\overline C$ contient la famille $J_0$-couvrante de $X$
$$
\Bigl( X_i \xrightarrow{ \ x_i \ } X \Bigl)_{i \in I}
$$
et il est $J_1$-couvrant puisqu'il contient $C$.

Donc il est $J_2$-couvrant.

Comme il est $J_2$-fermé, il est le crible maximal sur $X$, ce qui signifie comme voulu que le crible $C$ est $J_2$-couvrant.

On sait maintenant que $J_0$ est une topologie.

Il résulte de la définition de $J_0$ que si une famille de morphismes de ${\mathcal C}$ de même but
$$
\Bigl( X_i \xrightarrow{ \ x_i \ } X \Bigl)_{i \in I}
$$
est à la fois $J_0$-couvrante et $J_1$-couvrante, alors elle est aussi $J_2$-couvrante.

Autrement dit, on a
$$
J_1 \cap J_0 \subseteq J_2 \, .
$$

Enfin, considérons une topologie $J$ de ${\mathcal C}$ telle que
$$
J_1 \cap J \subseteq J_2 \, .
$$

Considérons une famille $J$-couvrante d'un objet $X$
$$
\Bigl( X_i \xrightarrow{ \ x_i \ } X \Bigl)_{i \in I}
$$
et le crible $C$ qu'elle engendre.

Puis considérons un morphisme
$$
X' \xrightarrow{ \ x \ } X
$$
et une famille $J_1$-couvrante de $X'$
$$
\Bigl( X'_k \xrightarrow{ \ x'_k \ } X' \Bigl)_{k \in K}
$$
dont le crible engendré $C'$ contient $x^* C$.

Comme $x^* C$ est $J$-couvrant, $C'$ est aussi $J$-couvrant et donc $J_2$-couvrant puisque $J_1 \cap J \subseteq J_2$.

Cela prouve comme voulu que $J \subseteq J_0$. \hfill $\Box$

\medskip

La formule de la proposition \ref{prop1.3.4} (ii) permet de calculer la différence de deux sous-topos ${\mathcal E}_1$ et ${\mathcal E}_2$ d'un topos ${\mathcal E}$ lorsque ${\mathcal E}$ est présenté par un site $({\mathcal C},J)$ et que les sous-topos ${\mathcal E}_1$ et ${\mathcal E}_2$ correspondent à deux topologies $J_1$ et $J_2$ de ${\mathcal C}$ qui sont explicites.

Lorsque les topologies $J_1$ et $J_2$ sont présentées comme engendrées par des générateurs, le calcul doit passer par les formules d'engendrement des topologies par des générateurs.

De telles formules seront présentées au chapitre III.

Elles s'appliqueront en particulier au cas où ${\mathcal E}$ est présenté comme le topos classifiant d'une théorie géométrique du premier ordre ${\mathbb T}$ et les sous-topos ${\mathcal E}_1$ et ${\mathcal E}_2$ de ${\mathcal E}$ correspondent à deux théories quotients ${\mathbb T}_1$ et ${\mathbb T}_2$ de ${\mathbb T}$ définies par des familles d'axiomes supplémentaires.

\subsection{L'opération d'image directe des sous-topos}\label{ssec1.3.3}

On rappelle que tout morphisme de topos définit une application d'image directe des sous-topos:

\begin{defn}\label{defn1.3.6}

Pour tout morphisme de topos $f : {\mathcal E}' \to {\mathcal E}$, on appelle ``application d'image directe'' définie par $f$ et on note
$$
f_* : {\rm ST} ({\mathcal E}') \longrightarrow {\rm ST} ({\mathcal E})
$$
l'application qui transforme tout sous-topos
$$
{\mathcal E}'_1 \xymatrix{\ar@{^{(}->}[r] & } {\mathcal E}'
$$
en l'unique sous-topos
$$
f_* \, {\mathcal E}'_1 = {\mathcal E}_1 \xymatrix{\ar@{^{(}->}[r] & } {\mathcal E}
$$
à travers lequel le morphisme composé
$$
{\mathcal E}'_1 \xymatrix{\ar@{^{(}->}[r] & } {\mathcal E}' \xrightarrow{ \ f \ } {\mathcal E}
$$
se factorise en un morphisme surjectif
$$
{\mathcal E}'_1 \xymatrix{\ar@{->>}[r] & } f_* \, {\mathcal E}'_1 = {\mathcal E}_1 \, .
$$
\end{defn}

\begin{remarkqed}

L'existence et l'unicité de la factorisation du composé
$$
{\mathcal E}'_1 \xymatrix{\ar@{^{(}->}[r] & } {\mathcal E}_1 \xrightarrow{ \ f \ } {\mathcal E}
$$
en un morphisme surjectif suivi d'un plongement
$$
{\mathcal E}'_1 \xymatrix{\ar@{->>}[r] & } {\mathcal E}_1 \xymatrix{\ar@{^{(}->}[r] & } {\mathcal E}
$$
ont été rappelées dans la proposition \ref{prop1.2.7}.
\end{remarkqed}

Les propriétés de base des applications d'image directe sont les suivantes:

\begin{prop}\label{prop1.3.7}
\begin{listeimarge}
\item Pour tout morphisme de topos $f : {\mathcal E}' \to {\mathcal E}$, l'application induite
$$
f_* : {\rm ST} ({\mathcal E}') \longrightarrow {\rm ST} ({\mathcal E})
$$
respecte la relation d'ordre $\geq$.

\item Pour tous morphismes de topos ${\mathcal E}'' \xrightarrow{ \, g \, } {\mathcal E}' \xrightarrow{ \, f \, } {\mathcal E}$, on a
$$
(f \circ g)_* = f_* \circ g_* : {\rm ST} ({\mathcal E}'') \longrightarrow {\rm ST} ({\mathcal E}) \, .
$$
\end{listeimarge}
\end{prop}

\begin{remarks}
\begin{listeisansmarge}
\item On rappellera au sous-paragraphe suivant que chaque application
$$
f_* : {\rm ST} ({\mathcal E}') \longrightarrow {\rm ST} ({\mathcal E})
$$
respecte les réunions arbitraires de sous-topos.

\item Les ensembles ordonnés ${\rm ST}({\mathcal E})$ constituent donc des invariants covariants des topos ${\mathcal E}$, au sens qu'ils composent un foncteur covariant de la catégorie des topos vers celle des ensembles ordonnés.
\end{listeisansmarge}
\end{remarks}

\begin{demo}
\begin{listeisansmarge}
\item résulte de ce que, pour tout sous-topos ${\mathcal E}'_1 \hookrightarrow {\mathcal E}'$, son image $f_* \, {\mathcal E}'_1$ est le plus petit sous-topos ${\mathcal E}_1 \hookrightarrow {\mathcal E}$ au travers duquel se factorise le morphisme composé ${\mathcal E}'_1 \hookrightarrow {\mathcal E}' \xrightarrow{ \, f \, } {\mathcal E}$.

%\smallskip

\item Si ${\mathcal E}''_1 \hookrightarrow {\mathcal E}''$ est un sous-topos, $g_* \, {\mathcal E}''_1 = {\mathcal E}'_1 \hookrightarrow {\mathcal E}'$ est son image par $g$, et $f_* \, {\mathcal E}'_1 = {\mathcal E}_1 \hookrightarrow {\mathcal E}$ est l'image de ${\mathcal E}'_1$ par $f$, les morphismes
$$
{\mathcal E}''_1 \xymatrix{\ar@{^{(}->}[r] & } {\mathcal E}'' \xrightarrow{ \ g \ } {\mathcal E}' \quad \mbox{et} \quad {\mathcal E}'_1 \xymatrix{\ar@{^{(}->}[r] & } {\mathcal E}' \xrightarrow{ \ f \ } {\mathcal E}
$$
se factorisent en deux morphismes surjectifs
$$
{\mathcal E}''_1 \xymatrix{\ar@{->>}[r] & } {\mathcal E}'_1 \quad \mbox{et} \quad {\mathcal E}'_1 \xymatrix{\ar@{->>}[r] & } {\mathcal E}_1 \, .
$$
Le composé de ceux-ci est un morphisme surjectif
$$
{\mathcal E}''_1 \longrightarrow {\mathcal E}_1
$$
ce qui, par unicité de la factorisation de ${\mathcal E}''_1 \hookrightarrow {\mathcal E}'' \xrightarrow{ \, f \circ g \, } {\mathcal E}$, assure que
$$
(g \circ f)_* \, {\mathcal E}''_1 = {\mathcal E}_1 \, .
$$
\end{listeisansmarge}
\end{demo}

La proposition \ref{prop1.2.7} qui permet de définir les applications d'images directes de sous-topos fournit en même temps une expression de ces images directes en termes de topologies:

\begin{prop}\label{prop1.3.8}

Considérons un morphisme d'un topos ${\mathcal E}'$ dans le topos $\widehat{\mathcal C}_J$ des faisceaux sur un site $({\mathcal C},J)$
$$
f : {\mathcal E}' \longrightarrow {\mathcal E} = \widehat{\mathcal C}_J
$$
défini par un foncteur plat et $J$-continu
$$
\rho : {\mathcal C} \longrightarrow {\mathcal E}' \, .
$$

On a dans ces conditions:

\begin{listeimarge}

\item Pour tout sous-topos de ${\mathcal E}'$
$$
{\mathcal E}'_1 \xymatrix{\ar@{^{(}->}[rr]^{(j^*, \, j_*)} && } {\mathcal E}' \, ,
$$
son image directe $f_* \, {\mathcal E}'_1 = {\mathcal E}_1$ par $f : {\mathcal E}' \to {\mathcal E} = \widehat{\mathcal C}_J$ est le sous-topos de $\widehat{\mathcal C}_J$ défini par la topologie $J_1$ de ${\mathcal C}$ pour laquelle une famille de morphismes de même but
$$
\Bigl( X_i \xrightarrow{ \ x_i \ } X \Bigl)_{i \in I}
$$
est couvrante si et seulement si sa transformée dans ${\mathcal E}'_1$ par le foncteur
$$
j^* \circ \rho : {\mathcal C} \longrightarrow {\mathcal E}' \xrightarrow{ \ j^* \ } {\mathcal E}'_1
$$
est globalement épimorphique.

\item Supposons que le topos ${\mathcal E}'$ est présenté comme topos $\widehat{\mathcal D}_K$ des faisceaux sur un site $({\mathcal D},K)$ et que le foncteur
$$
\rho : {\mathcal C} \longrightarrow {\mathcal E}' \xrightarrow{ \ \sim \ } \widehat{\mathcal D}_K
$$
se factorise à travers le foncteur canonique $\ell : {\mathcal D} \to \widehat{\mathcal D}_K$ en
$$
{\mathcal C} \xrightarrow{ \ \widetilde \rho \ } {\mathcal D} \xrightarrow{ \ \ell \ } \widehat{\mathcal D}_K \, .
$$
Alors, si ${\mathcal E}'_1 \hookrightarrow {\mathcal E}'$ est un sous-topos défini par une topologie $K_1$ de ${\mathcal D}$, son image directe $f_* \, {\mathcal E}'_1 = {\mathcal E}_1 \hookrightarrow {\mathcal E} = \widehat{\mathcal C}_J$ est définie comme sous-topos de $\widehat{\mathcal C}_J$ par la topologie $J_1$ de ${\mathcal C}$ pour laquelle une famille de morphismes 
$$
\Bigl( X_i \xrightarrow{ \ x_i \ } X \Bigl)_{i \in I}
$$
est couvrante si et seulement si sa transformée par le foncteur $\widetilde\rho : {\mathcal C} \to {\mathcal D}$ est $K_1$-couvrante.
\end{listeimarge}
\end{prop}

\begin{remark}
Dans la situation de (ii), si la topologie $K_1$ de ${\mathcal D}$ qui définit un sous-topos ${\mathcal E}'_1 \hookrightarrow {\mathcal E}_1 \xrightarrow{ \, \sim \, } \widehat{\mathcal D}_K$ est définie par une collection de générateurs, le calcul de la topologie $J_1$ de ${\mathcal C}$ qui définit le sous-topos ${\mathcal E}_1$ de ${\mathcal E} = \widehat{\mathcal C}_J$ image directe de ${\mathcal E}'_1$ par $f : {\mathcal E}' \to {\mathcal E}$ demande de faire appel à une formule qui exprime la topologie engendrée par des générateurs, comme on en verra au chapitre \ref{chap3}.
\end{remark}

\begin{demo}
\begin{listeisansmarge}
\item résulte de la proposition \ref{prop1.2.7}.

\item est une simple traduction de (i) puisqu'une famille de morphismes de ${\mathcal D}$ de même but
$$
\Bigl( X'_i \xrightarrow{ \ x'_i \ } X' \Bigl)_{i \in I}
$$
est couvrante pour une topologie $K_1$ de ${\mathcal D}$ si et seulement si sa transformée par le foncteur canonique
$$
\ell : {\mathcal D} \longrightarrow \widehat{\mathcal D}_{K_1}
$$
est globalement épimorphique. 
\end{listeisansmarge}
\end{demo}

De même, la théorie des topos classifiants et le théorème \ref{thm1.2.10} fournissent une expression en termes logiques des images directes de sous-topos:

\begin{prop}\label{prop1.3.9}

Considérons un morphisme d'un topos ${\mathcal E}'$ dans le topos classifiant ${\mathcal C}_{\mathbb T}$ d'une théorie géométrique du premier ordre ${\mathbb T}$
$$
m : {\mathcal E}' \longrightarrow {\mathcal E} = {\mathcal E}_{\mathbb T}
$$
qui correspond à un modèle $M'$ de la théorie ${\mathbb T}$ dans le topos ${\mathcal E}'$.

On a dans ces conditions:

\begin{listeimarge}

\item Pour tout sous-topos de ${\mathcal E}'$
$$
{\mathcal E}'_1 \xymatrix{\ar@{^{(}->}[rr]^{(j^* , \, j_*)} && } {\mathcal E}' \, ,
$$
son image directe $m_* \, {\mathcal E}'_1 = {\mathcal E}_1 \hookrightarrow {\mathcal E}$ par $m : {\mathcal E}' \to {\mathcal E} = {\mathcal E}_{\mathbb T}$ s'interprète comme le topos classifiant
$$
{\mathcal E}_{{\mathbb T}_1} \xymatrix{\ar@{^{(}->}[r] & } {\mathcal E}_{\mathbb T}
$$
d'une théorie quotient ${\mathbb T}_1$ de ${\mathbb T}$ pour laquelle un séquent géométrique
$$
\varphi (\vec x) \vdash \psi (\vec x)
$$
est démontrable si et seulement si il est vérifié par le modèle 
$$
j^* \, M'
$$
de la théorie ${\mathbb T}$ dans ${\mathcal E}'_1$.

\item Ecrivant le topos classifiant ${\mathcal E}_{\mathbb T}$ comme topos des faisceaux sur la catégorie syntactique géométrique ${\mathcal C}_{\mathbb T}$ munie de la topologie syntactique $J_{\mathbb T}$, supposons que le topos ${\mathcal E}'$ est présenté comme topos $\widehat{\mathcal D}_K$ des faisceaux sur un site $({\mathcal D},K)$ et que le foncteur
$$
m^* : {\mathcal E}_{\mathbb T} \longrightarrow {\mathcal E}' \xrightarrow{ \ \sim \ } \widehat{\mathcal D}_K
$$
s'inscrit dans un carré commutatif:
$$
\xymatrix{
{\mathcal C}_{\mathbb T} \ar[d]_-{\ell} \ar[r]^-M &{\mathcal D} \ar[d]^-{\ell} \\
{\mathcal E}_{\mathbb T} \ar[r]^-{m^*} &\widehat{\mathcal D}_K
}
$$
Alors, si ${\mathcal E}'_1 \hookrightarrow {\mathcal E}'$ est un sous-topos défini par une topologie $K_1$ de ${\mathcal D}$, son image directe $m_* \, {\mathcal E}'_1 = {\mathcal E}_1 \hookrightarrow {\mathcal E} = {\mathcal E}_{\mathbb T}$ s'interprète comme le topos classifiant
$$
{\mathcal E}_{{\mathbb T}_1} \xymatrix{\ar@{^{(}->}[r] & } {\mathcal E}_{\mathbb T}
$$
d'une théorie quotient ${\mathbb T}_1$ de ${\mathbb T}$ pour laquelle un séquent géométrique
$$
\varphi (\vec x) \vdash \psi (\vec x)
$$
est démontrable si et seulement si le monomorphisme de ${\mathcal C}_{\mathbb T}$
$$
(\varphi \wedge \psi)(\vec x) \xymatrix{\ar@{^{(}->}[r] & } \varphi (\vec x)
$$
est transformé par le foncteur $M : {\mathcal C}_{\mathbb T} \to {\mathcal D}$ en un morphisme $K_1$-couvrant.
\end{listeimarge}
\end{prop}

\begin{remarks}
	\begin{listeisansmarge}
\item Si un morphisme $m : {\mathcal E}' \to {\mathcal E}_{\mathbb T}$ correspond à un modèle $M'$ de la théorie ${\mathbb T}$ dans le topos ${\mathcal E}'$, et qu'un sous-topos ${\mathcal E}'_1 \xymatrix{\ar@{^{(}->}[rr]^{(j_* , \, j^*)} && } {\mathcal E}'$ a pour image directe par $m$ le sous-topos ${\mathcal E}_{{\mathbb T}_1} \hookrightarrow {\mathcal E}_{\mathbb T}$ qui classifie une théorie quotient ${\mathbb T}_1$ de ${\mathbb T}$, celle-ci peut être appelée la théorie du modèle $j^* M'$.

En effet, il résulte de (i) que les séquents géométriques du langage de ${\mathbb T}$ qui sont démontrables dans ${\mathbb T}_1$ sont exactement ceux qui sont vérifiés par le modèle $j^* M'$.

\item Cette caractérisation ne permet pas a priori d'extraire de la collection de tous les séquents géométriques démontrables dans ${\mathbb T}_1$ des familles d'axiomes qui suffisent à définir de telles théories ${\mathbb T}_1$.

\item La caractérisation de (ii) demande de savoir distinguer ceux des morphismes de ${\mathcal D}$ qui sont $K_1$-couvrants.

Si la topologie $K_1$ est définie par une collection de générateurs, cela doit passer par une formule exprimant la topologie engendrée par des générateurs.
\end{listeisansmarge}
\end{remarks}

\begin{demo}
\begin{listeisansmarge}
\item L'image directe d'un sous-topos
$$
{\mathcal E}'_1 \xymatrix{\ar@{^{(}->}[rr]^{(j^* , \, j_*)} && } {\mathcal E}'
$$
par un morphisme de topos
$$
m : {\mathcal E}' \longrightarrow {\mathcal E} = {\mathcal E}_{\mathbb T}
$$
correspondant à un modèle $M'$ de ${\mathbb T}$ dans ${\mathcal E}'$, est le plus petit sous-topos
$$
{\mathcal E}_1 \xymatrix{\ar@{^{(}->}[r] & } {\mathcal E}
$$
à travers lequel se factorise le morphisme ${\mathcal E}'_1 \hookrightarrow {\mathcal E}' \xrightarrow{ \, m \, } {\mathcal E}$.

Or, les sous-topos ${\mathcal E}_1 \hookrightarrow {\mathcal E} = {\mathcal E}_{\mathbb T}$ correspondent aux classes d'équivalence syntactique de théories quotients ${\mathbb T}_1$ de ${\mathbb T}$.

Comme le morphisme composé
$$
{\mathcal E}'_1 \xymatrix{\ar@{^{(}->}[rr]^{(j^* , \, j_*)} && } {\mathcal E}' \xrightarrow{ \ m \ } {\mathcal E} = {\mathcal E}_{\mathbb T}
$$
correspond au modèle $j^* M'$ de ${\mathbb T}$ dans ${\mathcal E}'_1$, il se factorise à travers le sous-topos défini par une théorie quotient ${\mathbb T}_1$
$$
{\mathcal E}_{{\mathbb T}_1} \xymatrix{\ar@{^{(}->}[r] & } {\mathcal E}_{\mathbb T}
$$
si et seulement si le modèle $j^* M'$ satisfait tous les axiomes de ${\mathbb T}_1$.

\item se déduit de (i).

En effet, un séquent géométrique écrit dans le langage de ${\mathbb T}$
$$
\varphi (\vec x) \vdash \psi (\vec x)
$$
est vérifié par le modèle $j^* M'$ de ${\mathbb T}$ dans ${\mathcal E}'_1$ si et seulement si le monomorphisme de ${\mathcal C}_{\mathbb T}$
$$
(\varphi \wedge \psi) (\vec x) \xymatrix{\ar@{^{(}->}[r] & } \varphi (\vec x)
$$
est transformé par le foncteur
$$
{\mathcal C}_{\mathbb T} \xrightarrow{ \ \ell \ } {\mathcal E}_{\mathbb T} = {\mathcal E} \xrightarrow{ \, m^* \ } {\mathcal E}' \xrightarrow{ \ j^* \ } {\mathcal E}'_1
$$
en un isomorphisme  de ${\mathcal E}'_1$.

Pour cela, il faut et il suffit qu'il soit transformé en un épimorphisme de ${\mathcal E}'_1$.

Si l'on a des équivalences
\begin{eqnarray}
{\mathcal E}' &\xrightarrow{ \ \sim \ } &\widehat{\mathcal D}_K \, , \nonumber \\
{\mathcal E}'_1 &\xrightarrow{ \ \sim \ } &\widehat{\mathcal D}_{K_1} \, , \nonumber
\end{eqnarray}
et que le foncteur
$$
{\mathcal C}_{\mathbb T} \xrightarrow{ \ \ell \ } {\mathcal E}_{\mathbb T} = {\mathcal E} \xrightarrow{ \, m^* \ } {\mathcal E}' \xrightarrow{ \ \sim \ } \widehat{\mathcal D}_K
$$
se factorise en
$$
{\mathcal C}_{\mathbb T} \xrightarrow{ \ M \ } {\mathcal D} \xrightarrow{ \, \ell \ } \widehat{\mathcal D}_K \, ,
$$
cela revient à demander que le foncteur $M$ transforme le monomorphisme de ${\mathcal C}_{\mathbb T}$
$$
(\varphi \wedge \psi)(\vec x) \xymatrix{\ar@{^{(}->}[r] & } \varphi (\vec x)
$$
en un morphisme $K_1$-couvrant de ${\mathcal D}$. 
\end{listeisansmarge}
\end{demo}

La remarque (i) qui suit l'énoncé de la proposition \ref{prop1.3.9} conduit à poser la définition suivante:

\begin{defn}\label{defn1.3.10}

Soit un modèle $M'$  d'une théorie géométrique du premier ordre ${\mathbb T}$ dans un topos ${\mathcal E}'$.

On appelle description des sous-topos de ${\mathcal E}'$ dans le langage de ${\mathbb T}$ induite par le modèle $M'$ l'application
$$
({\mathcal E}'_1 \xymatrix{\ar@{^{(}->}[r] & } {\mathcal E}') \longmapsto {\mathbb T}_{{\mathcal E}'_1}
$$
qui associe à tout sous-topos ${\mathcal E}'_1$ de ${\mathcal E}'$ la théorie quotient ${\mathbb T}_{{\mathcal E}'_1}$ de ${\mathbb T}$, uniquement déterminée à équivalence près, telle que le sous-topos
$$
{\mathcal E}_{{\mathbb T}_{{\mathcal E}'_1}} \xymatrix{\ar@{^{(}->}[r] & } {\mathcal E}_{\mathbb T}
$$
soit l'image directe du sous-topos ${\mathcal E}'_1 \hookrightarrow {\mathcal E}'$ par le morphisme
$$
m : {\mathcal E}' \longrightarrow {\mathcal E}_{\mathbb T}
$$
qui correspond au modèle $M'$ de ${\mathbb T}$ dans le topos ${\mathcal E}'$.
\end{defn}

\begin{remarkqed}

Une description de choses d'un certain type -- ici les sous-topos d'un topos donné -- dans un certain langage est en général partielle. Il se peut en particulier que deux choses différentes -- ici deux sous-topos -- aient la même description.

Décrire dans un certain langage suppose de donner un sens aux mots de son vocabulaire. Ici, la signification des éléments de vocabulaire d'un langage formel ${\mathbb T}$ est donnée par le choix d'un modèle $M'$ de ${\mathbb T}$ dans le topos ${\mathcal E}'$. 
\end{remarkqed}

\subsection{L'opération d'image réciproque des sous-topos}\label{ssec1.3.4}

\medskip

On rappelle que tout morphisme de topos définit aussi une application d'image réciproque des sous-topos:

\begin{thm}\label{thm1.3.11}

Soit un morphisme de topos $f : {\mathcal E}' \to {\mathcal E}$.

Alors:

\begin{listeimarge}

\item L'application d'image directe des sous-topos par $f$
$$
f_* : ({\rm ST} ({\mathcal E}') , \geq) \longrightarrow ({\rm ST} ({\mathcal E}) , \geq)
$$
admet une adjointe à gauche, appelée l'application d'image réciproque des sous-topos par $f$, notée
$$
f^{-1} : ({\rm ST} ({\mathcal E}) , \geq) \longrightarrow ({\rm ST} ({\mathcal E}') , \geq) \, .
$$
Elle est caractérisée par la propriété que, pour tous sous-topos
$$
{\mathcal E}_1 \xymatrix{\ar@{^{(}->}[r] & } {\mathcal E} \quad \mbox{et} \quad {\mathcal E}'_1 \xymatrix{\ar@{^{(}->}[r] & } {\mathcal E}' \, ,
$$
on a
$$
f^{-1} {\mathcal E}_1 \geq {\mathcal E}'_1
$$
si et seulement si
$$
{\mathcal E}_1 \geq f_* \, {\mathcal E}'_1 \, .
$$
\item L'application d'image directe $f_*$ respecte les réunions arbitraires de sous-topos, tandis que son adjointe l'application d'image réciproque $f^{-1}$ respecte les intersections arbitraires.
\end{listeimarge}
\end{thm}

\begin{remarks}
\begin{listeisansmarge}
\item L'application $f_*$ ne respecte pas en général les intersections de sous-topos, même finies.

\item En revanche, on démontrera au chapitre IV que l'application d'image réciproque $f^{-1}$ respecte toujours les réunions finies de sous-topos.

Elle ne respecte pas en général les réunions infinies.

\item Il résulte de la propriété de caractérisation de l'opération $f^{-1}$ que, si $f$ est un plongement
$$
{\mathcal E}' \xymatrix{\ar@{^{(}->}[r] & } {\mathcal E} \, ,
$$
l'opération d'image réciproque $f^{-1}$ se confond avec l'opérateur d'intersection de sous-topos
$$
({\mathcal E}_1 \xymatrix{\ar@{^{(}->}[r] & } {\mathcal E}) \longmapsto {\mathcal E}' \wedge {\mathcal E}_1 \, .
$$
On a déjà rappelé dans la proposition \ref{prop1.3.2} (ii) qu'un tel opérateur respecte les réunions finies de sous-topos.

\item Pour toute paire de morphismes de topos
$$
{\mathcal E}'' \xrightarrow{ \ g \ } {\mathcal E}' \xrightarrow{ \ f \ } {\mathcal E} \, ,
$$
l'identité de la proposition \ref{prop1.3.7} (ii)
$$
(f \circ g)_* = f_* \circ g_*
$$
induit par passage aux adjoints une identité
$$
(f \circ g)^{-1} = g^{-1} \circ f^{-1} \, .
$$

Cela signifie que les ensembles ordonnés ${\rm ST} ({\mathcal E})$ constituent aussi des invariants contravariants des topos ${\mathcal E}$. Autrement dit, ils composent un foncteur contravariant
\begin{eqnarray}
{\mathcal E} &\longmapsto &{\rm ST} ({\mathcal E}) \, , \nonumber \\
({\mathcal E}' \xrightarrow{ \, f \, } {\mathcal E}) &\longmapsto &(f^{-1} : {\rm ST}({\mathcal E}) \to {\rm ST} ({\mathcal E}')) \nonumber
\end{eqnarray}
de la catégorie des topos vers celle des ensembles ordonnés.
\end{listeisansmarge}
\end{remarks}

\begin{demosansqed}
\begin{listeisansmarge}
\item[(ii)] est conséquence formelle de (i).

\item[(i)] résulte de la proposition suivante, qui explicite comment les opérations d'image réciproque de sous-topos s'expriment en termes de topologies, et même de générateurs des topologies:
\end{listeisansmarge}
\end{demosansqed}

\begin{prop}\label{prop1.3.12}

Considérons un morphisme d'un topos ${\mathcal E}'$ dans le topos $\widehat{\mathcal C}_J$ des faisceaux sur un site $({\mathcal C},J)$
$$
f : {\mathcal E}' \longrightarrow {\mathcal E} : \widehat{\mathcal C}_J
$$
défini par un foncteur plat et $J$-continu
$$
\rho : {\mathcal C} \longrightarrow {\mathcal E}' \, .
$$

On a dans ces conditions:

\begin{listeimarge}

\item Pour tout sous-topos ${\mathcal E}_1 \hookrightarrow {\mathcal E} = \widehat{\mathcal C}_J$ défini par une topologie $J_1$ de ${\mathcal C}$, l'ensemble des sous-topos de ${\mathcal E}'$
$$
{\mathcal E}'_1 \xymatrix{\ar@{^{(}->}[rr]^{(j^* , \, j_*)} && } {\mathcal E}'
$$
tels que le foncteur composé
$$
j^* \circ \rho : {\mathcal C} \longrightarrow {\mathcal E}' \longrightarrow {\mathcal E}'_1
$$
transforme toute famille $J_1$-couvrante de morphismes de ${\mathcal C}$ en une famille globalement épimorphique de ${\mathcal E}'_1$, possède un plus grand élément $f^{-1} {\mathcal E}_1 \hookrightarrow {\mathcal E}'$.

\item Supposons que le topos ${\mathcal E}'$ est présenté comme topos $\widehat{\mathcal D}_K$ des faisceaux sur un site $({\mathcal D},K)$ et que le foncteur
$$
\rho : {\mathcal C} \longrightarrow {\mathcal E}' \xrightarrow{ \ \sim \ } \widehat{\mathcal D}_K
$$
se factorise à travers le foncteur canonique $\ell : {\mathcal D} \to \widehat{\mathcal D}_K$ en
$$
{\mathcal C} \xrightarrow{ \ \widetilde\rho \ } {\mathcal D} \xrightarrow{ \ \ell \ } \widehat{\mathcal D}_K \, .
$$

Alors, pour tous sous-topos ${\mathcal E}_1 \hookrightarrow {\mathcal E} = \widehat{\mathcal C}_J$ et ${\mathcal E}'_1 \hookrightarrow {\mathcal E}' \cong \widehat{\mathcal D}_K$ définis par des topologies $J_1$ et $K'_1$ de ${\mathcal C}$ et ${\mathcal D}$, on a
$$
{\mathcal E}_1 \geq f_* \, {\mathcal E}'_1
$$
si et seulement si
$$
f^{-1} {\mathcal E}_1 \geq {\mathcal E}'_1 \quad \mbox{soit} \quad K_1 \subseteq K'_1 \, ,
$$
où $f^{-1} {\mathcal E}_1 \hookrightarrow {\mathcal E}' \cong \widehat{\mathcal D}_K$ désigne le sous-topos défini par la topologie $K_1$ de ${\mathcal D}$ engendrée sur $K$ par les transformées par
$$
\widetilde\rho : {\mathcal C} \longrightarrow {\mathcal D}
$$
des familles $J_1$-couvrantes de morphismes de ${\mathcal C}$ de même but.

\item Dans la situation de (ii), la topologie $K_1$ de ${\mathcal D}$ qui définit le sous-topos
$$
f^{-1} {\mathcal E}_1 \xymatrix{\ar@{^{(}->}[r] & } {\mathcal E}' \cong \widehat{\mathcal D}_K
$$
image réciproque par $f$ du sous-topos défini par une topologie $J_1$
$$
{\mathcal E}_1 \xymatrix{\ar@{^{(}->}[r] & } {\mathcal E} = \widehat{\mathcal C}_J \, ,
$$
est la topologie engendrée par la réunion de

\smallskip

$\left\{\begin{matrix}
\bullet &\mbox{n'importe quelle famille de générateurs de la topologie $K$ de ${\mathcal D}$,} \hfill \\
\bullet &\mbox{les transformés par le foncteur $\widetilde\rho : {\mathcal C} \to {\mathcal D}$ de n'importe quelle famille de générateurs dans ${\mathcal C}$} \hfill \\
&\mbox{de la topologie $J_1$ sur $J$.} \hfill
\end{matrix}\right.$
\end{listeimarge}
\end{prop}

\begin{demo}
\begin{listeisansmarge}
\item résulte de (ii) puisqu'une famille de morphismes de ${\mathcal D}$
$$
\Bigl( X'_i \xrightarrow{ \ x'_i \ } X' \Bigl)_{i \in I}
$$
est couvrante pour une topologie $K'_1$ de ${\mathcal D}$ si et seulement si le foncteur canonique
$$
{\mathcal D} \longrightarrow \widehat{\mathcal D}_{K'_1}
$$
la transforme en une famille globalement épimorphique.

\item résulte de la proposition \ref{prop1.3.8} (ii).

\item résulte de (i) car si $J_1$ est une topologie de ${\mathcal C}$ engendrée sur $J$ par une famille ${\mathcal X}$ de générateurs et que
$$
{\mathcal E}'_1 \xymatrix{\ar@{^{(}->}[rr]^{(j^* , \, j_*)} && } {\mathcal E}'
$$
est un sous-topos de ${\mathcal E}'$ tel que le foncteur plat et $J$-continu
$$
j^* \circ \rho : {\mathcal C} \longrightarrow {\mathcal E}' \xrightarrow{ \ j^* \ } {\mathcal E}'_1
$$
transforme en familles globalement épimorphiques de ${\mathcal E}'_1$ toutes les familles de morphismes de ${\mathcal C}$ de même but
$$
\Bigl( X_i \xrightarrow{ \ x_i \ } X \Bigl)_{i \in I}
$$
qui sont éléments de la famille génératrice ${\mathcal X}$, alors il transforme toutes les familles $J_1$-couvrantes de ${\mathcal C}$ en familles globalement épimorphiques de ${\mathcal E}'_1$. 
\end{listeisansmarge}
\end{demo}

La théorie des topos classifiants et le théorème \ref{thm1.2.10} permettent de reformuler en termes logiques l'expression des images réciproques de sous-topos:

\begin{cor}\label{cor1.3.13}

Considérons un morphisme d'un topos ${\mathcal E}'$ dans le topos classifiant ${\mathcal E}_{\mathbb T}$ d'une théorie géométrique du premier ordre ${\mathbb T}$
$$
m : {\mathcal E}' \longrightarrow {\mathcal E} = {\mathcal E}_{\mathbb T}
$$
qui correspond à un modèle $M'$ de la théorie ${\mathbb T}$ dans le topos ${\mathcal E}'$.

On a dans ces conditions:

\begin{listeimarge}

\item Pour toute théorie quotient ${\mathbb T}_1$ de ${\mathbb T}$ définie par une famille d'axiomes supplémentaires ajoutés à ceux de ${\mathbb T}$
$$
\bigl( \varphi_i (\vec x_i) \vdash \psi_i (\vec x_i)\bigl)_{i \in I} \, ,
$$
l'image réciproque du sous-topos
$$
{\mathcal E}_{{\mathbb T}_1} \xymatrix{\ar@{^{(}->}[r] & } {\mathcal E}_{\mathbb T} = {\mathcal E}
$$
par le morphisme $m : {\mathcal E}' \to {\mathcal E} = {\mathcal E}_{\mathbb T}$ est le plus grand sous-topos
$$
{\mathcal E}'_1 \xymatrix{\ar@{^{(}->}[rr]^{(j^* , \, j_*)} && } {\mathcal E}'
$$
tel que le modèle $j^* M'$ de ${\mathbb T}$ dans ${\mathcal E}'_1$ satisfasse tous les axiomes
$$
\varphi_i (\vec x_i) \vdash \psi_i (\vec x_i) \, , \quad i \in I \, .
$$
\item Ecrivant le topos classifiant ${\mathcal E}_{\mathbb T}$ comme topos des faisceaux sur la catégorie syntactique géométrique ${\mathcal C}_{\mathbb T}$ munie de la topologie syntactique $J_{\mathbb T}$, supposons que le topos ${\mathcal E}'$ est présenté comme topos $\widehat{\mathcal D}_K$ des faisceaux sur un site $({\mathcal D},K)$ et que le foncteur
$$
m^* : {\mathcal E}_{\mathbb T} \longrightarrow {\mathcal E}' \xrightarrow{ \ \sim \ } \widehat{\mathcal D}_K
$$
s'inscrit dans un carré commutatif:
$$
\xymatrix{
{\mathcal C}_{\mathbb T}\ar[d]_-{\ell} \ar[r]^-M &{\mathcal D} \ar[d]^-{\ell} \\
{\mathcal E}_{\mathbb T} \ar[r]^-{m^*} &\widehat{\mathcal D}_K
}
$$

Alors, pour toute théorie quotient ${\mathbb T}_1$ de ${\mathbb T}$ définie par une famille d'axiomes supplémentaires
$$
\bigl( \varphi_i (\vec x_i) \vdash \psi_i (\vec x_i)\bigl)_{i \in I} \, ,
$$
l'image réciproque par $m : {\mathcal E}' \to {\mathcal E}_{\mathbb T}$ du sous-topos
$$
{\mathcal E}_{{\mathbb T}_1} \xymatrix{\ar@{^{(}->}[r] & } {\mathcal E}_{\mathbb T}
$$
est le sous-topos
$$
m^{-1} {\mathcal E}_{{\mathbb T}_1} \xymatrix{\ar@{^{(}->}[r] & } {\mathcal E}' \xrightarrow{ \ \sim \ } \widehat{\mathcal D}_K
$$
défini par la plus petite topologie $K_1$ de ${\mathcal D}$ qui contient $K$ et pour laquelle les transformés par le foncteur
$$
M : {\mathcal C}_{\mathbb T} \longrightarrow {\mathcal D}
$$
des monomorphismes de ${\mathcal C}_{\mathbb T}$
$$
(\varphi_i \wedge \psi_i) (\vec x_i) \xymatrix{\ar@{^{(}->}[r] & } \varphi_i (\vec x_i)
$$
sont des morphismes couvrants.
\end{listeimarge}
\end{cor}

\begin{demo}
\begin{listeisansmarge}
\item résulte de (ii).

\item Pour tout sous-topos
$$
{\mathcal E}'_1 \xymatrix{\ar@{^{(}->}[r] & } {\mathcal E}' \xrightarrow{ \ \sim \ } \widehat{\mathcal D}_K
$$
défini par une topologie $K'_1$ de ${\mathcal D}$ contenant $K$, il résulte de la proposition I.3.9 (ii) que son image directe par $m : {\mathcal E}' \to {\mathcal E} = {\mathcal E}_{\mathbb T}$
$$
M_* \, {\mathcal E}'_1 \xymatrix{\ar@{^{(}->}[r] & } {\mathcal E} = {\mathcal E}_{\mathbb T}
$$
est contenue dans le sous-topos 
$$
{\mathcal E}_{{\mathbb T}_1} \xymatrix{\ar@{^{(}->}[r] & } {\mathcal E}_{\mathbb T}
$$
si et seulement si, pour tout élément de la famille d'axiomes supplémentaires
$$
\bigl( \varphi_i (\vec x_i) \vdash \psi_i (\vec x_i)\bigl)_{i \in I} 
$$
qui définit la théorie quotient ${\mathbb T}_1$ de ${\mathbb T}$, le foncteur
$$
\widetilde\rho : {\mathcal C}_{\mathbb T} \longrightarrow {\mathcal D}
$$
transforme chacun des monomorphismes
$$
(\varphi_i \wedge \psi_i) (\vec x_i) \xymatrix{\ar@{^{(}->}[r] & } \varphi_i (\vec x_i) \, , \quad i \in I \, ,
$$
en un morphisme $K'_1$-couvrant de ${\mathcal D}$.

Cela montre que la topologie $K_1$ de ${\mathcal D}$ qui définit le sous-topos
$$
m^{-1} {\mathcal E}_{{\mathbb T}_1} \xymatrix{\ar@{^{(}->}[r] & } {\mathcal E}' \xrightarrow{ \ \sim \ } \widehat{\mathcal D}_K
$$
est la plus petite topologie contenant $K$ qui possède cette propriété.
\end{listeisansmarge} 
\end{demo}

\subsection{Correspondances topossiques et actions sur les sous-topos}\label{ssec1.3.5}

On a la définition naturelle:

\begin{defn}\label{defn1.3.14}
\begin{listeimarge}
\item On appelle correspondance d'un topos ${\mathcal E}'$ dans un topos ${\mathcal E}$ la donnée d'un topos ${\mathcal M}$ muni de deux morphismes:
$$
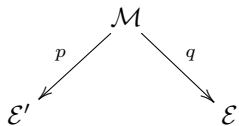
\xymatrix{
&{\mathcal M} \ar[dl]_p \ar[dr]^q \\
{\mathcal E}' &&{\mathcal E}
}
$$

\item On appelle chaîne topossique d'un topos ${\mathcal E}'$ dans un topos ${\mathcal E}$ un diagramme de morphismes de topos de la forme:
$$
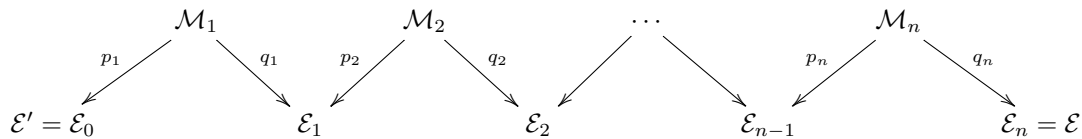
\xymatrix{
&{\mathcal M}_1 \ar[dl]_{p_1} \ar[dr]^{q_1} &&{\mathcal M}_2 \ar[dl]_{p_2} \ar[dr]^{q_2} &&\cdots \ar[dl] \ar[dr] &&{\mathcal M}_n \ar[dl]_{p_n} \ar[dr]^{q_n} \\
{\mathcal E}' = {\mathcal E}_0 &&{\mathcal E}_1 &&{\mathcal E}_2 &&{\mathcal E}_{n-1} &&{\mathcal E}_n = {\mathcal E}
}
$$

\end{listeimarge}
\end{defn}

\begin{remarksqed}
\begin{listeisansmarge}
\item Si les topos ${\mathcal E}'$ et ${\mathcal E}$ sont présentés comme les topos classifiants de théories géométriques du premier ordre ${\mathbb T}'$ et ${\mathbb T}$, une correspondance de ${\mathcal E}' = {\mathcal E}_{{\mathbb T}'}$ dans ${\mathcal E} = {\mathcal E}_{\mathbb T}$ consiste en un topos ${\mathcal M}$ muni à la fois d'un modèle de ${\mathbb T}'$ et d'un modèle de ${\mathbb T}$.

De telles correspondances peuvent être construites à partir d'une théorie géométrique du premier ordre ${\mathbb T}'$ en

\smallskip

$\left\{\begin{matrix}
\bullet &\mbox{choisissant un modèle $M'$ de ${\mathbb T}'$ dans un topos ${\mathcal E}'$,} \hfill \\
\bullet &\mbox{réalisant dans ${\mathcal E}'$ des constructions à partir du modèle $M'$,} \hfill \\
\bullet &\mbox{montrant que les constructions réalisées sont munies de nouvelles structures qui les font apparaître} \hfill \\
&\mbox{comme un modèle d'une nouvelle théorie géométrique du premier ordre ${\mathbb T}$.} \hfill
\end{matrix} \right.$

\smallskip

Si les constructions réalisées dans le topos ${\mathcal M}$ à partir du modèle $M'$ font appel à des procédés ``d'ordre supérieur'' comme de recourir aux foncteurs exponentiels ${\mathcal H}om (\bullet , \bullet)$, on peut considérer qu'on aura ``monté dans l'abstraction'' en passant du langage formel ${\mathbb T}'$ au langage formel ${\mathbb T}$.

\item Plus généralement, construire une chaîne topossique reliant les topos classifiants ${\mathcal E}' = {\mathcal E}_{{\mathbb T}'}$ et ${\mathcal E} = {\mathcal E}_{\mathbb T}$ de deux théories géométriques ${\mathbb T}'$ et ${\mathbb T}$ peut consister à élaborer une suite de théories géométriques ${\mathbb T}' = {\mathbb T}_0 , {\mathbb T}_1 , {\mathbb T}_2 , \cdots , {\mathbb T}_{n-1} , {\mathbb T}_n = {\mathbb T}$ reliées par des topos ${\mathcal M}_i$, $1 \leq i \leq n$, qui sont munis chacun d'un modèle de ${\mathbb T}_{i-1}$ et d'un modèle de ${\mathbb T}_i$.
\end{listeisansmarge} 
\end{remarksqed}

Les opérations d'images réciproques et d'images directes des sous-topos permettent de définir les actions sur les sous-topos de correspondances ou plus généralement de chaînes topossiques:

\begin{defn}\label{defn1.3.15}
\begin{listeimarge}
\item L'action sur les sous-topos d'une correspondance d'un topos ${\mathcal E}'$ dans un topos ${\mathcal E}$
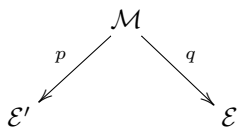
$$
\xymatrix{
&{\mathcal M} \ar[dl]_p \ar[dr]^q \\
{\mathcal E}' &&{\mathcal E}
}
$$

est l'application composée
$$
q_* \circ p^{-1} : {\rm ST} ({\mathcal E}') \longrightarrow {\rm ST} ({\mathcal M}) \longrightarrow {\rm ST} ({\mathcal E}) \, .
$$
\item Plus généralement, l'action sur les sous-topos d'une chaîne topossique d'un topos ${\mathcal E}'$ dans un topos ${\mathcal E}$ 
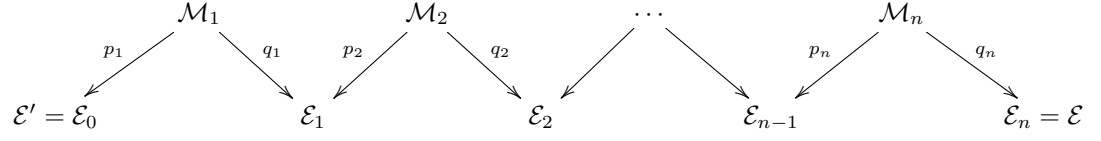
$$
\xymatrix{
&{\mathcal M}_1 \ar[dl]_{p_1} \ar[dr]^{q_1} &&{\mathcal M}_2 \ar[dl]_{p_2} \ar[dr]^{q_2} &&\cdots \ar[dl] \ar[dr] &&{\mathcal M}_n \ar[dl]_{p_n} \ar[dr]^{q_n} \\
{\mathcal E}' = {\mathcal E}_0 &&{\mathcal E}_1 &&{\mathcal E}_2 &&{\mathcal E}_{n-1} &&{\mathcal E}_n = {\mathcal E}
}
$$

est l'application composée
$$
(q_n)_* \circ p_n^{-1} \circ \cdots \circ (q_2)_* \circ p_2^{-1} \circ (q_1)_* \circ p_1^{-1} :  {\rm ST} ({\mathcal E}') \longrightarrow {\rm ST} ({\mathcal E}) \, .
$$
\end{listeimarge}
\end{defn}

\begin{remarksqed}
\begin{listeisansmarge}
\item Les actions sur les sous-topos de correspondances et plus généralement de chaînes topossiques respectent les relations d'ordres.

\item Ces actions respectent toujours les réunions finies de sous-topos, puisque les applications d'images directes respectent les réunions arbitraires et que celles d'images réciproques respectent les réunions finies.

\item En revanche, elles ne respectent pas en général les intersections, non plus que les réunions infinies. 
\end{listeisansmarge}
\end{remarksqed}

\subsection{Morphismes localement connexes et images directes extraordinaires de sous-topos}\label{ssec1.3.6}

\medskip

On a rappelé que tout morphisme de topos $f : {\mathcal E}' \to {\mathcal E}$ définit une paire de foncteurs adjoints $(f^{-1} , f_*)$ entre les ensembles ordonnés des sous-topos de ${\mathcal E}'$ et ${\mathcal E}$
$$
 \xymatrix{ ({\rm ST} ({\mathcal E}') , \geq) \ar@<-3pt>[rr]_{f_*} && ({\rm ST} ({\mathcal E}) , \geq) \ar@<-3pt>[ll]_{f^{-1}} } \, .
 $$
 
 En revanche, le foncteur $f^{-1}$ n'admet pas en général de foncteur adjoint à gauche:
 
 Nous montrerons cependant au chapitre IV:
 
\begin{thm}\label{thm1.3.16}
 
Soit un morphisme de topos
 $$
 f : {\mathcal E}' \longrightarrow {\mathcal E}
 $$
 qui est ``localement connexe''.
 
 Alors l'opération d'image réciproque des sous-topos par $f$
 $$
 f^{-1} : {\rm ST}({\mathcal E}) \longrightarrow {\rm ST} ({\mathcal E}')
 $$
 admet non seulement un adjoint à droite $f_*$ mais aussi un adjoint à gauche, que l'on peut appeler ``opération d'image directe extraordinaire''
 $$
 f_! : {\rm ST}({\mathcal E}') \longrightarrow {\rm ST} ({\mathcal E})
 $$
caractérisée par la propriété que, pour tous sous-topos
$$
{\mathcal E}'_1 \xymatrix{\ar@{^{(}->}[r] \  &}  {\mathcal E}' \quad \mbox{et} \quad {\mathcal E}_1 \xymatrix{\ar@{^{(}->}[r] \  &} {\mathcal E} \, ,
$$
on a
$$
f_! ({\mathcal E}'_1) \geq {\mathcal E}_1
$$
si et seulement si
$$
{\mathcal E}'_1 \geq f^{-1} ({\mathcal E}_1) \, .
$$
\end{thm}

\begin{remarksqed}
\begin{listeisansmarge}
\item L'existence d'un tel adjoint à gauche $f_!$ de $f^{-1}$ signifie que l'opération $f^{-1}$ respecte les réunions arbitraires de sous-topos.

\item Son adjointe l'opération d'image directe extraordinaire $f_!$ respecte les intersections arbitraires de sous-topos.
\end{listeisansmarge} 
\end{remarksqed}

Afin de rappeler la notion de morphisme ``localement connexe'' de topos, on a d'abord besoin de rappeler:

\begin{lem}\label{lem1.3.17}
\begin{listeimarge}
\item Pour tout objet $E$ d'un topos ${\mathcal E}$, la catégorie relative
$$
{\mathcal E}/E
$$
est un topos.

\item Pour tout morphisme $e : E_2 \to E_1$ d'un topos ${\mathcal E}$, le foncteur de produit fibré
\begin{eqnarray}
{\mathcal E}/E_1 &\longrightarrow &{\mathcal E}/E_1 \, , \nonumber \\
(E_3 \to E_1) &\longmapsto &(E_2 \times_{E_1} E_3 \to E_1) \nonumber
\end{eqnarray}
admet un adjoint à droite noté
$$
(E \to E_1) \longmapsto {\mathcal H}om_{E_1}(E_2 , E) \, .
$$
\item Pour tout morphisme $e : E_2 \to E_1$ d'un topos ${\mathcal E}$, le foncteur de changement de base
$$
\begin{matrix}
e^* &: &\hfill {\mathcal E}/E_1 &\longrightarrow &{\mathcal E}/E_2 \, , \hfill \\
&&(E \to E_1) &\longmapsto &(E \times_{E_1} E_2 \to E_2)  
\end{matrix}
$$
admet un adjoint à gauche qui est
$$
\begin{matrix}
e_! &: &\hfill {\mathcal E}/E_2 &\longrightarrow &{\mathcal E}/E_1 \, , \hfill \\
&&(E_3 \to E_2) &\longmapsto &(E_3 \to E_2 \to E_1)  
\end{matrix}
$$
et un adjoint à droite noté
$$
\textstyle\prod_e : {\mathcal E}/E_2 \longrightarrow {\mathcal E}/E_1 \, .
$$
\item Pour tout morphisme de topos ${\mathcal E}' \xrightarrow{ \ (f^* , \, f_*) \ } {\mathcal E}$, tout objet $E_1$ de ${\mathcal E}$ et tous objets $(E_2 \to E_1)$ et $(E \to E_1)$ de ${\mathcal E}/E_1$, on a dans ${\mathcal E}'$ un morphisme canonique
$$
f^* \, {\mathcal H}om_{E_1} (E_2 , E) \longrightarrow {\mathcal H}om_{f^* E_1} (f^* E_2 , f^* E) \, .
$$
\end{listeimarge}
\end{lem}

\begin{remarks}
\begin{listeisansmarge}
\item Si $E$ est l'objet terminal $1$ de ${\mathcal E}$, ${\mathcal E}/E$ s'identifie à ${\mathcal E}$ et le foncteur ${\mathcal H}om_E (\bullet , \bullet)$ est simplement noté ${\mathcal H}om (\bullet , \bullet)$.

\item De même, si $e$ est l'unique morphisme d'un objet $E$ de ${\mathcal E}$ vers l'objet terminal $1$ de ${\mathcal E}$, le foncteur
$$
\textstyle\prod_e : {\mathcal E}/E \longrightarrow {\mathcal E} \, ,
$$
adjoint à droite de
\begin{eqnarray}
{\mathcal E} &\longrightarrow &{\mathcal E}/E \, , \nonumber \\
E_1 &\longmapsto &(E_1 \times E \to E) \, , \nonumber
\end{eqnarray}
est simplement noté
$$
\textstyle\prod_E : {\mathcal E}/E \longrightarrow {\mathcal E} \, .
$$
\end{listeisansmarge}
\end{remarks}

\begin{esquissedemo}
	\begin{listeisansmarge}
\item Si le topos ${\mathcal E}$ est présenté comme topos $\widehat{\mathcal C}_J$ des faisceaux sur un site $({\mathcal C},J)$, alors la catégorie relative
$$
{\mathcal E}/E
$$
est équivalente au topos des faisceaux sur le site
$$
({\mathcal C} / E , J)
$$
où

\smallskip

$\left\{\begin{matrix}
\bullet &\mbox{${\mathcal C}/E$ désigne la ``catégorie des éléments'' de $E$, qui a pour objets les paires} \hfill \\
&\mbox{$(X,x)$ constituées d'un objet $X$ de ${\mathcal C}$ et d'un élément $x\in E(X)$,} \hfill \\
&\mbox{et pour morphismes} \hfill \\
&(X_2 , x_2) \longrightarrow (X_1 , x_1) \\
{ \ } \\
&\mbox{les morphismes $u : X_2 \to X_1$ de ${\mathcal C}$ tels que} \hfill \\
{ \ } \\
&E(u) (x_1) = x_2 \, , \\
{ \ } \\
\bullet &\mbox{$J$ est la topologie de ${\mathcal C}/E$ pour laquelle une famille de morphismes} \hfill \\
{ \ } \\
&(X_i , x_i) \longrightarrow (X,x) \, , \quad i \in I \, , \\
{ \ } \\
&\mbox{est couvrante si et seulement si la famille sous-jacente de morphismes de ${\mathcal C}$} \hfill \\
{ \ } \\
&(X_i \longrightarrow X)_{i \in I} \\
&\mbox{est $J$-couvrante.} \hfill
\end{matrix}\right.$

%\medskip

\item se ramène d'après (i) au cas où $E_1$ est l'objet terminal de ${\mathcal E}$. 

Ce cas a été rappelé dans le théorème \ref{thm1.1.12} (iii).

%\medskip 

\item L'existence de l'adjoint à droite $\prod_e$ a été rappelée dans le théorème \ref{thm1.1.12} (iii).

Le foncteur
$$
\begin{matrix}
e_! &: &\hfill {\mathcal E}/E_2 &\longrightarrow &{\mathcal E}/E_1 \, , \hfill \\
&&(E_3 \to E_2) &\longmapsto &(E_3 \to E_2 \to E_1)  
\end{matrix}
$$
est adjoint à gauche du foncteur
$$
e^* : (E \to E_1) \longmapsto (E \times_{E_1} E_2 \to E_2)
$$
puisqu'on a par définition des produits fibrés
$$
{\rm Hom}_{E_1} (E_3 , E) = {\rm Hom}_{E_2} (E_3 , E \times_{E_1} E_2) \, .
$$
\item Pour tout objet $E_3 \to E_1$ de ${\mathcal E}/E_1$, tout morphisme de ${\mathcal E}/E_1$
$$
E_3 \times_{E_1} E_2 \longrightarrow E
$$
induit un morphisme de ${\mathcal E}' / f^* E_1$
$$
f^* E_3 \times_{f^* E_1} f^* E_2 \longrightarrow f^* E \, .
$$
Autrement dit, tout morphisme de ${\mathcal E}/E_1$
$$
E_3 \longrightarrow {\mathcal H}om_{E_1} (E_2 , E)
$$
induit un morphisme de ${\mathcal E}' / f^* E_1$
$$
f^* E_3 \longrightarrow {\mathcal H}om_{f^* E_1} (f^* E_2 , f^* E_3)
$$
ou, ce qui revient au même, un morphisme de ${\mathcal E}/E_1$
$$
E_3 \longrightarrow f_* {\mathcal H}om_{f^* E_1} (f^* E_2, f^* E_3) \, .
$$
Cela définit un morphisme naturel de ${\mathcal E}/E_1$
$$
{\mathcal H}om_{E_1} (E_2 , E) \longrightarrow f_* {\mathcal H}om_{f^* E_1} (f^* E_2 , f^* E)
$$
ou, ce qui revient au même, un morphisme de ${\mathcal E}' / f^* E_1$
$$
f^* {\mathcal H}om_{E_1} (E_2 , E) \longrightarrow {\mathcal H}om_{f^* E_1} (f^* E_2 , f^* E) \, .
$$
\end{listeisansmarge}
\end{esquissedemo}

Nous pouvons maintenant rappeler les notions de morphisme de topos ``essentiel'' ou ``localement connexe'':

\begin{defn}\label{defn1.3.18}
 \begin{listeimarge}
\item Un morphisme de topos ${\mathcal E}' \xrightarrow{ \ (f^* , \, f_*) \ } {\mathcal E}$ est dit ``essentiel'' si sa composante d'image inverse
$$
f^* : {\mathcal E} \longrightarrow {\mathcal E}'
$$
admet aussi un adjoint à gauche
$$
f_! : {\mathcal E}' \longrightarrow {\mathcal E} \, .
$$
\item Un morphisme de topos est dit ``localement connexe'' si c'est un morphisme essentiel
$$
{\mathcal E}' \xrightarrow{ \ (f_! , \, f^* , \, f_*) \ } {\mathcal E}
$$
et qu'il possède les trois propriétés équivalentes suivantes:

%\medskip

$\left\lmoustache\begin{matrix}
(1) &\mbox{Pour tous morphismes $E_2 \to E_1 \leftarrow E$ de ${\mathcal E}$, le morphisme canonique de ${\mathcal E}'$} \hfill \\
{ \ } \\
&f^* {\mathcal H}om_{E_1} (E_2,E) \longrightarrow {\mathcal H}om_{f^* E_1} (f^* E_2 , f^* E) \\
{ \ } \\
&\mbox{est un isomorphisme.} \hfill \\
%\end{matrix}\right.$

%$\left\rmoustache\begin{matrix}
(2) &\mbox{Pour tout morphisme $e : E_2 \to E_1$ de ${\mathcal E}$, le carré} \hfill \\
{ \ } \\
&\xymatrix{
{\mathcal E}' / f^* E_2 \ar[d]_{f_!} &&{\mathcal E}' / f^* E_1 \ar[ll]_{\bullet \times_{f^* E_1} f^* E_2} \ar[d]^{f_!} \\
{\mathcal E}/E_2 &&{\mathcal E}/E_1 \ar[ll]_{\bullet \times_{E_1} E_2}
} \\
{ \ } \\
&\mbox{est commutatif à isomorphisme canonique près.} \hfill \\
{ \ } \\
\end{matrix}\right.$

$\left\rmoustache\begin{matrix}
(3) &\mbox{Pour tout morphisme $e : E_2 \to E_1$ de ${\mathcal E}$, le carré} \hfill \\
{ \ } \\
&\xymatrix{
{\mathcal E}' / f^* E_2 \ar[rr]^{\prod_{f^* e}} &&{\mathcal E}' / f^* E_1 \\
{\mathcal E}/E_2 \ar[u]^{f^*} \ar[rr]_{\prod_e} &&{\mathcal E}/E_1 \ar[u]_{f^*}
} \\
{ \ } \\
&\mbox{est commutatif à isomorphisme canonique près.} \hfill 
\end{matrix}\right.$
\end{listeimarge}
\end{defn}

\begin{remarks}
\begin{listeisansmarge}[label=(\roman*)]
\item Si ${\mathcal E} = {\rm Ens}$ est le ``topos ponctuel'', il suffit qu'un morphisme
$$
f = (f^* , f_*) : {\mathcal E}' \longrightarrow {\mathcal E} = {\rm Ens}
$$
soit ``essentiel'' pour qu'il soit ``localement connexe''.

En effet, tout objet $E_1$ de ${\rm Ens}$ s'écrit
$$
E_1 = \coprod_{i \in E_1} \{\bullet\} = \coprod_{i \in E_1} 1_{{\mathcal E}}
$$
et l'on a donc aussi
$$
f^* E_1 = \coprod_{i \in E_1} 1_{{\mathcal E}'} \, .
$$
Pour un objet $E'$ de ${\mathcal E}'$, se donner un morphisme
$$
E' \longrightarrow f^* E_1
$$
équivaut à se donner une décomposition
$$
E' = \coprod_{i \in E_1} E'_i \, .
$$
De même, se donner une application $E_2 \to E_1$ équivaut à se donner un ensemble $E_2$ décomposé en
$$
E_2 = \coprod_{i \in E_1} E_{2,i} \, .
$$
Celle-ci induit une décomposition
$$
E' \times_{f^* E_1} f^* E_2 = \coprod_{i \in E_1} \ \coprod_{j \in E_{2,i}} E'_i \, .
$$
La conclusion résulte de ce que le foncteur
$$
f_! : {\mathcal E}' \longrightarrow {\mathcal E} = {\rm Ens} \, ,
$$
qui est un adjoint à gauche, respecte toutes le colimites.

\item Le qualificatif ``localement connexe'' a été choisi du fait que, pour tout espace topologique $X$, le morphisme de topos qu'il définit
$$
p : {\mathcal E}_X \longrightarrow {\mathcal E}_{\{\bullet\}} = {\rm Ens}
$$
est ``localement connexe'' si et seulement si l'espace $X$ est ``localement connexe'' au sens que tout ouvert est somme disjointe d'ouverts connexes.

En effet, si $X$ est un espace topologique localement connexe, le topos associé ${\mathcal E}_X$ s'écrit comme le topos $\widehat{\mathcal C}_J$ des faisceaux sur la catégorie ${\mathcal C}$ des ouverts connexes de $X$, munie de la topologie $J$ des recouvrements. Pour tout ouvert connexe $U$ de $X$ vu comme un objet de ${\mathcal E}_X$, et pour tout ensemble $E$, les morphismes de ${\mathcal E}_X$
$$
U \longrightarrow p^* E = \coprod_{i \in E} p^* \{\bullet\} = \coprod_{i \in E} 1_{{\mathcal E}_X}
$$
sont les applications continues
$$
U \longrightarrow \coprod_{i \in E} X
$$
qui rendent commutatif le triangle:
$$
\xymatrix{
U \ar[rr] \ar[rd] &&\underset{i \in E}{\coprod} \, X \ar[ld] \\
&X
}
$$
Comme $U$ est connexe, l'ensemble de ces applications s'identifie à $E$. On en déduit que le foncteur
$$
\begin{matrix}
p^* &: &{\rm Ens} &\longrightarrow &{\mathcal E}_X \, , \hfill \\
&&\hfill E &\longmapsto &\underset{i \in E}{\coprod} \, X
\end{matrix}
$$
admet pour adjoint à gauche le foncteur
$$
\begin{matrix}
p_! &: &{\mathcal E}_X = \widehat{\mathcal C}_J &\longrightarrow &{\rm Ens} \, , \hfill \\
&&\hfill F &\longmapsto &\underset{(X,x) \in {\mathcal C}/F}{\varprojlim} \, \{\bullet\}
\end{matrix}
$$
qui associe à tout faisceau $F$ sur $({\mathcal C},J)$ la limite du foncteur constant
$$
(X,x) \longmapsto \{\bullet \}
$$
sur la catégorie ${\mathcal C}/F$ des éléments de $F$.

Ainsi, le morphisme de topos $p : {\mathcal E}_X \to {\rm Ens}$ est essentiel et, d'après la remarque (i), il est localement connexe.

Réciproquement, supposons que le foncteur
$$
\begin{matrix}
p^* &: &{\rm Ens} &\longrightarrow &{\mathcal E}_X \, , \hfill \\
&&\hfill E &\longmapsto &\underset{i \in E}{\coprod} \, X
\end{matrix}
$$
admet un adjoint à gauche
$$
p_! : {\mathcal E}_X \longrightarrow {\rm Ens} \, .
$$
Pour tout ouvert $U$ de $X$, le morphisme canonique de ${\mathcal E}_X$
$$
U \longrightarrow p^* p_! U = \coprod_{i \in p_! U} p^* \{\bullet\} = \coprod_{i \in p_! U} X
$$
induit une décomposition en somme disjointes d'ouverts
$$
U = \coprod_{i \in p_! U} U_i
$$
dont chaque composante $U_i$, $i \in I$, vérifie
$$
p_! U_i = \{\bullet\} \, .
$$
Pour tout ensemble $E$, l'ensemble des applications continues
$$
U_i \longrightarrow \coprod_{j \in E} X
$$
qui rendent commutatif le triangle
$$
\xymatrix{
U_i \ar[rd] \ar[rr] &&\underset{j \in E}{\coprod} \, X \ar[ld] \\
&X
}
$$
s'identifie à $E$.

Cela signifie que chacune des composantes $U_i$ de $U$, $i \in p_! U$, est connexe.
\end{listeisansmarge}
\end{remarks}

\begin{justifdefi}
Il faut vérifier l'équivalence des propriétés (1), (2) et (3).

Pour tous morphismes $E_2 \to E_1 \leftarrow E$ de ${\mathcal E}$ et $E' \to f^* E_1$ de ${\mathcal E}'$, les propriétés d'adjonction se traduisent par des identifications
\begin{eqnarray}
&&{\rm Hom} (E' , f^* ({\mathcal H}om_{E_1} (E_2,E))) \nonumber \\
&= &{\rm Hom} (f_! E' , {\mathcal H}om_{E_1} (E_2,E)) \nonumber \\
&= &{\rm Hom} (f_! E' \times_{E_1} E_2 , E) \nonumber
\end{eqnarray}
et
\begin{eqnarray}
&&{\rm Hom} (E' , {\mathcal H}om_{f^* E_1} (f^* E_2,f^* E)) \nonumber \\
&= &{\rm Hom} (E' \times_{f^* E_1} f^* E_2 , f^* E) \nonumber \\
&= &{\rm Hom} (f_! (E' \times_{f^* E_1} f^* E_2) , E) \, . \nonumber
\end{eqnarray}
Par conséquent, les morphismes naturels
$$
f^* {\mathcal H}om_{E_1} (E_2,E)) \longrightarrow {\mathcal H}om_{f^* E_1} (f^* E_2 , f^* E)
$$
induisent des morphismes naturels
$$
f_! (E' \times_{f^* E_1} f^* E_2) \longrightarrow f_! E' \times_{E_1} E_2 \, ,
$$
et les premiers sont des isomorphismes si et seulement si les seconds sont des isomorphismes.

C'est l'équivalence des propriétés (1) et (2).

L'équivalence de (2) et (3) résulte de ce que les carrés de (2) sont adjoints à gauche des carrés de (3). 
\end{justifdefi}

Le lemme \ref{lem1.3.17} (iii) fournit une grande famille d'exemples de morphismes localement connexes de topos:

\begin{cor}\label{cor1.3.19}
 
Pour tout morphisme $e : E_2 \to E_1$ d'un topos ${\mathcal E}$, la paire de foncteurs adjoints
$$
({\mathcal E} / E_1 \xrightarrow{ \ \bullet \times_{E_1} E_2 \ } {\mathcal E}/E_2 , {\mathcal E}/E_2 \xrightarrow{ \ \prod_e \ } {\mathcal E}/E_1)
$$
définit un morphisme de topos
$$
(e^* , e_*) : {\mathcal E} / E_2 \longrightarrow {\mathcal E}/E_1
$$
qui est localement connexe.
\end{cor}

\begin{demo}

On sait déjà que le foncteur
$$
\bullet \times_{E_1} E_2 : {\mathcal E}/E_1 \longrightarrow {\mathcal E}/E_2
$$
admet pour adjoint à gauche le foncteur
$$
\begin{matrix}
e_! &: &\hfill {\mathcal E}/E_2 &\longrightarrow &{\mathcal E}/E_1 \, , \hfill \\
&&(E_3 \to E_2) &\longmapsto &(E \to E_2 \to E_1) \, .
\end{matrix}
$$
Pour montrer que $(e^* , e_*)$ est localement connexe, vérifions que $e_!$ possède la propriété (2) de la définition I.3.18.

Il suffit de traiter le cas où $e$ est l'unique morphisme
$$
E \longrightarrow 1
$$
d'un objet $E$ de ${\mathcal E}$ vers l'objet terminal $1$.

Pour tous morphismes $E_2 \to E_1$ et $E_3 \to E \times E_1 = e^* E_1$ de ${\mathcal E}$, l'objet
$$
E_3 \times_{E \times E_1} E \times E_2
$$
s'identifie dans ${\mathcal E} / E_2$ à l'objet
$$
E_3 \times_{E_1} E_2 \, .
$$
Cela signifie comme voulu que les foncteurs
$$
e_! (\bullet \times_{e^* E_1} e^* E_2)
$$
et
$$
e_! (\bullet) \times_{E_1} E_2
$$
sont canoniquemet isomorphes. 
\end{demo}

On déduira du théorème \ref{thm1.3.16}:

\begin{cor}\label{cor1.3.20}
 
Pour tout morphisme de topos 
$$
f = (f^* , f_*) : {\mathcal E}' \longrightarrow {\mathcal E} \, ,
$$
l'application induite d'image réciproque des sous-topos
$$
f^{-1} : {\rm ST} ({\mathcal E}) \longrightarrow {\rm ST} ({\mathcal E}')
$$
respecte les réunions finies de sous-topos.
\end{cor}

\noindent {\bf Indication sur la démonstration:}

\smallskip

On rappellera en effet au chapitre IV que tout morphisme de topos
$$
{\mathcal E}' \longrightarrow {\mathcal E}
$$
peut s'écrire comme le composé
$$
{\mathcal E}' \xymatrix{\ar@{^{(}->}[r] & } {\mathcal E}'' \longrightarrow {\mathcal E}
$$
d'un plongement de topos 
$$
{\mathcal E}' \xymatrix{\ar@{^{(}->}[r] & } {\mathcal E}''
$$
et d'un morphisme localement connexe
$$
{\mathcal E}'' \longrightarrow {\mathcal E}' \, .
$$
La conclusion en résultera d'après le théorème \ref{thm1.3.16} et la proposition \ref{prop1.3.2} (ii). 
\hfill $\Box$

	\chapter{Sous-topos et topologies de Grothendieck}\label{chap2}

\section{Une notion générale de correspondance galoisienne}\label{sec2.1}

\subsection{Les équivalences induites par des paires de foncteurs adjoints}\label{ssec2.1.1}

\medskip

Commençons par rappeler que toute paire de foncteurs adjoints induit une équivalence entre les catégories d'objets fixes de part et d'autre:

\begin{prop}\label{prop2.1.1}

Soit une paire de foncteurs adjoints
$$
\left({\mathcal C} \xrightarrow{ \ F \ } {\mathcal D} \, , \ {\mathcal D} \xrightarrow{ \ G \ } {\mathcal C}\right)
$$
entre deux catégories localement petites ${\mathcal C}$ et ${\mathcal D}$.

Soient ${\mathcal C}'$ et ${\mathcal D}'$ les deux sous-catégories pleines de ${\mathcal C}$ et ${\mathcal D}$ constituées des ``points fixes'', c'est-à-dire des objets $X$ de ${\mathcal C}$ ou $Y$ de ${\mathcal D}$ tels que le morphisme d'adjonction canonique
$$
X \longrightarrow G \circ F (X) \quad \mbox{ou} \quad F \circ G (Y) \longrightarrow Y
$$
soit un isomorphisme.

Alors $F$ et $G$ induisent des équivalences de catégories
$$
\xymatrix{
{\mathcal C}' \ar@<2pt>[r]^{\sim} &{\mathcal D}' \ar@<2pt>[l]^{\sim}
}
$$
réciproques l'une de l'autre au sens que les transformations naturelles
$$
{\rm id}_{{\mathcal C}'} \longrightarrow G \circ F \quad \mbox{et} \quad F \circ G \longrightarrow {\rm id}_{{\mathcal D}'}
$$
sont des isomorphismes.
\end{prop}

\begin{remark}

En particulier, si la transformation naturelle
$$
{\rm id}_{\mathcal C} \longrightarrow G \circ F \quad \mbox{[resp. \ $F \circ G \longrightarrow {\rm id}_{\mathcal D}$ \, ]}
$$
est un isomorphisme, alors le foncteur $F$ [resp. $G$] définit une équivalence de ${\mathcal C}$ [resp. ${\mathcal D}$] sur la sous-catégorie pleine de ${\mathcal D}$ [resp. ${\mathcal C}$] constituée des ``points fixes''.
\end{remark}

Cela complète l'énoncé du lemme \ref{lem1.2.2}.

\medskip

\begin{esquissedemo}
Il suffit de démontrer que si le morphisme canonique
$$
X \longrightarrow G \circ F (X) \quad \mbox{[resp. \ $F \circ G (Y) \longrightarrow Y$ \, ]}
$$
est un isomorphisme, et
$$
Y=F(X) \qquad \mbox{[resp. \ $X=G(Y)$ \, ]},
$$
alors le morphisme canonique
$$
F \circ G(Y) \longrightarrow Y \quad \mbox{[resp. \ $X \longrightarrow G \circ F (X)$ \, ]}
$$
est un isomorphisme.

En effet, le transformé par $F$ [resp. par $G$] de l'isomorphisme
$$
X \xrightarrow{ \ \sim \ } G \circ F(X) \quad \mbox{[resp. \ $F \circ G (Y) \xrightarrow{ \ \sim \ } Y$ \, ]}
$$
est alors un isomorphisme
$$
Y \xrightarrow{ \ \sim \ } F \circ G(Y) \quad \mbox{[resp. \ $G \circ F (X) \xrightarrow{ \ \sim \ } X$ \, ]}
$$
puisque l'on note
$$
Y = F(X) \quad \mbox{[resp. \ $X = G(Y)$ \, ]}.
$$

Or, on sait d'autre part que le composé de celui-ci avec le morphisme d'adjonction
$$
F \circ G(Y) \longrightarrow Y \quad \mbox{[resp. \ $G \circ F (X) \longrightarrow X$ \, ]}
$$
est le morphisme identité de $Y=F(X)$ [resp. de $X=G(Y)$]
$$
{\rm id}_{F(X)} : F(X) \longrightarrow F \circ G \circ F(X) \longrightarrow F(X)
$$
[resp.
$$
{\rm id}_{G(Y)} : G(Y) \longrightarrow G \circ F \circ G(Y) \longrightarrow G(Y) \ ].
$$
D'où le résultat. 
\end{esquissedemo}

\subsection{Le cas particulier des structures ordonnées}\label{ssec2.1.2}

\smallskip

La proposition précédente d'équivalence entre les sous-catégories des points fixes d'une paire de foncteurs adjoints reliant deux catégories s'applique en particulier lorsque ces deux catégories sont deux ensembles ou plus généralement deux classes munis de relations d'ordre.

Dans ce cas, on a de surcroît la propriété que les points fixes sont exactement les points images:

\begin{cor}\label{cor2.1.2}

Soient deux ensembles ordonnés ou plus généralement deux classes $C$ et $D$ munies de relations d'ordre~$\leq$.

Soit une paire d'applications respectant les relations d'ordre
$$
\xymatrix{
(C , \leq) \ar@<2pt>[r]^F & (D,\leq) \ar@<2pt>[l]^G
}
$$
qui sont adjointes au sens que, pour tous éléments $c \in C$, $d \in D$, on a
$$
F(c) \leq d \Longleftrightarrow c \leq G(d) \, .
$$
Alors:

\begin{listeimarge}

\item Si l'on note les sous-ensembles ou sous-classes de points fixes
$$
C' = \{c \in C \mid G \circ F(c) = c \}
$$
et
$$
D' = \{d \in D \mid F \circ G(d) = d \} \, ,
$$
$F$ et $G$ induisent deux bijections inverses l'une de l'autre
$$
\xymatrix{
C' \ar@<2pt>[r]^{F \atop \sim} &D' \, . \ar@<2pt>[l]^{\overset{\sim}{G}}
}
$$

\item Un élément $c \in C$ [resp. $d \in D$] est fixé par $G \circ F$ [resp. par $F \circ G$] si et seulement si il est une image au sens que
$$
\exists \, d \in D \, , \quad c = G(d)
$$
[resp.
$$
\exists \, c \in C \, , \quad d = F(c) \mbox{\ ]}.
$$
\end{listeimarge}
\end{cor}

\begin{remarks}
\begin{listeisansmarge}
\item Pour tout $c \in C$ [resp. $d \in D$], on a
$$
c \leq G \circ F(c) \quad \mbox{[resp. \ $F \circ G (d) \leq d$ \ ].}
$$
\item Pour tout $c \in C$ [resp. $d \in D$] et tout $c' \in C'$ [resp. $d' \in D'$], on a
$$
G \circ F(c) \leq c' \quad \mbox{si} \quad c \leq c'
$$
[resp. 
$$
d' \leq F \circ G (d) \quad \mbox{si} \quad d' \leq d \mbox{ \, ].}
$$

\item Ainsi, pour tout élément $c$ de $C$ [resp. $d$ de $D$],
$$
G \circ F(c) \quad \mbox{[resp. \ $F \circ G (d)$]}
$$
est le plus petit [resp. le plus grand] point fixe qui soit plus grand que $c$ [resp. plus petit que $d$].
\end{listeisansmarge}
\end{remarks}

\begin{demo}
\begin{listeisansmarge}
\item est un cas particulier de la proposition \ref{prop2.1.1}.

En effet, dans une catégorie qui consiste en une relation d'ordre sur l'ensemble ou la classe de ses objets, les seuls isomorphismes sont les égalités.

\item Si un élément $c \in C$ [resp. $d \in D$] est un point fixe, il est évidemment une image.

Réciproquement, s'il est de la forme
$$
c= G(d) \quad \mbox{[resp. \ $d=F(c)$ \, ],}
$$
les morphismes naturels
$$
G(d) \longrightarrow G \circ F \circ G(d) \longrightarrow G(d)
$$
[resp.
$$
F(c) \longrightarrow F \circ G \circ F(c) \longrightarrow F(c) \mbox{ \, ]}
$$
s'écrivent
$$
c \leq G \circ F(c) \leq c
$$
[resp. 
$$
d \leq F \circ G(d) \leq d \mbox{ \, ].}
$$

\noindent Cela signifie que $c$ [resp. $d$] est un point fixe.
\end{listeisansmarge} 
\end{demo}

\subsection{Les paires de foncteurs adjoints définies par des relations}\label{ssec2.1.3}

\smallskip

Une observation élémentaire mais très utile et qui s'applique dans d'innombrables situations est que toute relation entre deux ensembles ou plus généralement deux classes induit une paire de foncteurs adjoints entre les ensembles ou classes ordonnés de leur parties:

\begin{lem}\label{lem2.1.3}

Soient deux ensembles ou classes $T$ et $S$ reliés par une relation arbitraire
$$
R \xymatrix{\ar@{^{(}->}[r] & } T \times S \, .
$$
Alors $R$ définit une paire adjointe d'applications respectant les relations d'ordre
$$
\xymatrix{
({\mathcal P} (T) , \subseteq) \ar@<2pt>[r]^{F_R} & ({\mathcal P} (S) , \supseteq) \ar@<2pt>[l]^{G_R}
}
$$
entre les ensembles ou classes partiellement ordonnés ${\mathcal P} (T)$ et ${\mathcal P} (S)$ des sous-ensembles ou sous-classes de $T$ et $S$
$$
(J \subseteq T) \longmapsto F_R(J) = \{ s \in S \mid (t,s) \in R \, , \ \forall \, t \in J \} \, ,
$$
$$
(I \subseteq S) \longmapsto G_R(I) = \{ t \in T \mid (t,s) \in R \, , \ \forall \, s \in I \} \, .
$$
\end{lem}

\begin{remark}

Il résulte de la propriété d'adjonction de $F_R$ et $G_R$ que pour toute famille de sous-ensembles ou sous-classes
$$
J_k \subseteq T \, , \ k \in K \qquad \mbox{[resp. \ $I_k \subseteq S \, , \ k \in K$ \ ],}
$$
on a
$$
F_R \left(\bigcup_{k \in K} J_k \right) = \bigcap_{k \in K} F_R (J_k) \qquad \mbox{[resp.  $\displaystyle G_R \left( \bigcup_{k \in K} I_k \right) = \bigcap_{k \in K} G_R (I_k)$ \ ].}
$$
\end{remark}

\begin{demo}

Il est évident sur les définitions des applications $F_R$ et $G_R$ que
$$
J_1 \subseteq J_2 \Longrightarrow F_R(J_1) \supseteq F_R(J_2) \, ,
$$
$$
I_1 \supseteq I_2 \Longrightarrow G_R(I_1) \subseteq G_R(I_2) \, .
$$
Ainsi, ces deux applications respectent les relations d'ordre $\subseteq$ de ${\mathcal P}(T)$ et $\supseteq$ de ${\mathcal P}(S)$. Elles sont adjointes car, pour toutes parties
$$
J \subseteq T \quad \mbox{et} \quad I \subseteq S \, ,
$$
on a les équivalences
$$
\begin{matrix}
&F_R (J) \supseteq I \hfill \\
\Longleftrightarrow &(t,s) \in R \, , \ \forall \, t \in J \, , \ \forall \, s \in I \hfill \\
\Longleftrightarrow &J \subseteq G_R (I) \, . \hfill
\end{matrix}
$$
\end{demo}

\subsection{Correspondances galoisiennes induites et procédés d'engendrement}\label{ssec2.1.4}

\smallskip

En combinant le corollaire \ref{cor2.1.2} et le lemme \ref{lem2.1.3}, on obtient que toute relation entre deux ensembles ou deux classes définit une correspondance dite ``galoisienne'' entre parties qui sont des points fixes:

\begin{cor}\label{cor2.1.4}

Soit une relation entre deux ensembles ou deux classes $T$ et $S$
$$
R \xymatrix{\ar@{^{(}->}[r] & } T \times S \, .
$$
Considérons la paire adjointe d'applications respectant les relations d'ordre que définit cette relation
$$
\xymatrix{
({\mathcal P} (T) , \subseteq) \ar@<2pt>[r]^{F_R} & ({\mathcal P} (S) , \supseteq) \, . \ar@<2pt>[l]^{G_R}
}
$$

Alors:

\begin{listeimarge}

\item Les applications $F_R$ et $G_R$ induisent deux bijections inverses l'une de l'autre entre les points fixes de ${\mathcal P} (T)$ et ceux de ${\mathcal P} (S)$
$$
\xymatrix{
\{J \subseteq T \mid G_R \circ F_R (J) = J\} \ar@<2pt>[r]^{\sim} & I \subseteq S \mid F_R \circ G_R (I) = I \} \, .\ar@<2pt>[l]^{\sim}
}
$$

\item Pour tout
$$
J \subseteq T \qquad \mbox{[resp. $I \subseteq S$ \ ]},
$$
on a
$$
G_R \circ F_R (J) = J \qquad \mbox{[resp. $F_R \circ G_R (I) = I$ \ ]}
$$
si et seulement si $J$ [resp. $I$] est une image au sens qu'existe $I \subseteq S$ [resp. $J \subseteq T$] tel que
$$
J = G_R(I) \qquad \mbox{[resp. $I = F_R(J)$ \ ]}.
$$

\item Pour tout
$$
J \subseteq T \qquad \mbox{[resp. $I \subseteq S$ \ ]},
$$
on a l'inclusion
$$
J \subseteq G_R \circ F_R(J) \qquad \mbox{[resp. $F_R \circ G_R (I) \supseteq I$ \ ]}
$$
et, réciproquement, pour tout point fixe $J' \subseteq T$ [resp. $I' \subseteq S$], on a l'implication
$$
J \subseteq J' \Longrightarrow G_R \circ F_R(J) \subseteq J'
$$
[resp. \
$$
I' \supseteq I \Longrightarrow I' \supseteq F_R \circ G_R(I) \mbox{\, ].}
$$
\end{listeimarge}
\end{cor}

\begin{remarks}
\begin{listeisansmarge}
\item Ce procédé général d'engendrement de correspondances bijectives avait déjà été identifié et utilisé à de multiples reprises par Birkhoff, comme exposé dans la section 3.2 de \cite{Galois}.

\item On parle de ``correspondances galoisiennes" car ce procédé généralise la correspondance de Galois.

Considérons en effet une extension finie de corps
$$
K \xymatrix{\ar@{^{(}->}[r] & } L
$$
et le groupe $G$ des automorphismes de $L$ qui fixent les points de $K$.

Puis considérons la relation
$$
R \xymatrix{\ar@{^{(}->}[r] & } L \times G
$$
définie par la formule
$$
R = \{(x , \sigma) \in L \times G \mid \sigma (x) = x \} \, .
$$
Celle-ci définit une paire adjointe d'applications respectant les relations d'ordre
$$
\xymatrix{
({\mathcal P} (L) , \subseteq) \ar@<2pt>[r]^{F_R} & ({\mathcal P} (G) , \supseteq) \, , \ar@<2pt>[l]^{G_R}
}
$$
et ces applications induisent deux bijections inverses l'une de l'autre entre les points fixes de ${\mathcal P} (L)$ et ceux de ${\mathcal P} (G)$.

Or, la théorie de Galois montre que les points fixes de ${\mathcal P} (G)$ sont les sous-groupes
$$
H \subseteq G \, ,
$$
et que les points fixes de ${\mathcal P} (L)$ sont les corps intermédiaires
$$
K \subseteq K_1 \subseteq L
$$
qui contiennent le sous-corps des points fixes
$$
K_0 = G_R (G) = \{x \in L \mid \sigma (x) = x \, , \ \forall \, x \in L \} \, .
$$

\item Pour tout
$$
J \subseteq T \qquad \mbox{[resp. \ $I \subseteq S$ \, ],}
$$
son image
$$
G_R \circ F_R (J) \qquad \mbox{[resp. \ $F_R \circ G_R(I)$ \, ]}
$$
peut être appelée le point fixe engendré par $J$ [resp. $I$].

En effet, il résulte de (iii) que cette image est le plus petit point fixe qui contienne $J$ [resp. $I$].
\end{listeisansmarge}
\end{remarks}

\begin{demo}
\begin{listeisansmarge}
\item[(i) et (ii)] sont des cas particuliers des points (i) et (ii) du corollaire \ref{cor2.1.2} compte tenu du lemme \ref{lem2.1.3}.

\item[(iii)] est un cas particulier de la remarque (iii) qui suit l'énoncé du corollaire \ref{cor2.1.2}. 
\end{listeisansmarge}
\end{demo}

\section{La dualité des cribles et des préfaisceaux}\label{sec2.2}

\subsection{Définition d'une relation entre cribles et préfaisceaux}\label{ssec2.2.1}

\medskip

Nous allons introduire une relation naturelle entre les cribles d'une catégorie essentiellement petite et les préfaisceaux sur cette catégorie, puis appliquer à cette relation les résultats généraux rappelés au paragraphe précédent.

Commençons par définir la relation:

\begin{defn}\label{defn2.2.1}

Soit une catégorie essentiellement petite ${\mathcal C}$, munie du foncteur de Yoneda de plongement pleinement fidèle dans la catégorie des préfaisceaux
$$
y : {\mathcal C} \xymatrix{\ar@{^{(}->}[r] & } \widehat{\mathcal C} = [{\mathcal C}^{\rm op} , {\rm Ens}] \, .
$$
Soit $T$ la classe des cribles de ${\mathcal C}$, constitués d'un objet $X$ de ${\mathcal C}$ et d'un sous-préfaisceau
$$
C \xymatrix{\ar@{^{(}->}[r] & } y(X) \, .
$$
Soit $S$ la classe des préfaisceaux sur ${\mathcal C}$, c'est-à-dire des objets de $\widehat{\mathcal C}$.

On appellera ``dualité des cribles et des préfaisceaux'' la relation
$$
R \xymatrix{\ar@{^{(}->}[r] & } T \times S
$$
constituée des paires d'éléments
$$
(C \xymatrix{\ar@{^{(}->}[r] & } y(X),P)
$$
telles que, pour tout morphisme de ${\mathcal C}$
$$
X' \xrightarrow{ \ x \ } X \, ,
$$
l'application de restriction au sous-préfaisceau $C \times_{y(X)} y(X')$
$$
P(X') = {\rm Hom} (y(X'),P) \longrightarrow {\rm Hom} (C \times_{y(X)} y(X') , P)
$$
est une bijection.
\end{defn}

\begin{remarksqed}
\begin{listeisansmarge}
\item  Dans la pratique, les cribles $C \hookrightarrow y(X)$ sont présentés par des familles de morphismes de but $X$
$$
(X_i \xrightarrow{ \ x_i \ } X)_{i \in I}
$$
qui sont ``génératrices'' au sens que les éléments de $C$ sont les morphismes de but $X$ qui se factorisent à travers l'un au moins des morphismes $X_i \xrightarrow{ \, x_i \, } X$ de la famille.

\item Le produit fibré $C \times_{y(X)} y(X')$ est l'image réciproque
$$
x^* C \xymatrix{\ar@{^{(}->}[r] & } y(X')
$$
par le morphisme $x : X' \to X$ du crible $C \hookrightarrow y(X)$. Il est constitué des morphismes $X'' \to X'$ dont le composé avec $X' \xrightarrow{ \, x \, } X$ est élément de $C$.

\item Si le crible $x^* C = C \times_{y(X)} y(X')$ est engendré par une famille de morphismes
$$
(X'_i \xrightarrow{ \ x'_i \ } X')_{i \in I'} \, ,
$$
se donner un morphisme vers un préfaisceau $P$
$$
C \times_{y(X)} y(X') \longrightarrow P
$$
consiste concrètement à se donner une famille de sections
$$
p_i \in P(X'_i) \, , \quad i \in I' \, ,
$$
telle que, pour tous indices $i_1 , i_2 \in I'$ et tout carré commutatif de ${\mathcal C}$
$$
\xymatrix{
X'' \ar[d] \ar[r] &X'_{i_1} \ar[d]^-{x'_{i_1}} \\
X'_{i_2} \ar[r]^-{x'_{i_2}} &X'
}
$$
les deux éléments
$$
p'_{i_1} \in P(X'_{i_1}) \quad \mbox{et} \quad p'_{i_2} \in P(X'_{i_2})
$$
aient même image dans $P(X'')$.

\item Ainsi, la paire $(C \hookrightarrow y(X),P)$ est élément de la relation $R \hookrightarrow T \times S$ si et seulement si pour tout morphisme
$$
X' \xrightarrow{ \ x \ } X \, ,
$$
et pour n'importe quelle famille génératrice $(X'_i \xrightarrow{ \ x'_i \ } X')_{i \in I'}$ de $x^* C = C \times_{y(X)} y(X')$, toute famille de sections comme dans (iii)
$$
p'_i \in P(X'_i) \, , \quad i \in I' \, ,
$$
provient d'une unique section $p' \in P(X')$. 
\end{listeisansmarge}
\end{remarksqed}

\subsection{La dualité induite des topologies et des sous-topos}\label{ssec2.2.2}

\smallskip

On considère toujours une catégorie essentiellement petite ${\mathcal C}$, la classe $T$ des cribles $C \hookrightarrow y(X)$ de ${\mathcal C}$, la classe $S$ des préfaisceaux $P$ sur ${\mathcal C}$ et la relation
$$
R \xymatrix{\ar@{^{(}->}[r] & } T \times S
$$
de la définition \ref{defn2.2.1}.

On sait d'après le lemme \ref{lem2.1.3} que cette relation $R$ définit une paire adjointe d'applications respectant les relations d'ordre
$$
\xymatrix{ ({\mathcal P} (T) , \subseteq) \ar@<2pt>[r]^{F_R} & ({\mathcal P} (S) , \supseteq), \ar@<2pt>[l]^{G_R} }
$$
et on sait d'après le corollaire \ref{cor2.1.4} que ces applications induisent deux bijections inverses l'une de l'autre entre les points fixes de ${\mathcal P} (T)$ et ceux de ${\mathcal P} (S)$.

Nous allons démontrer le théorème suivant qui fait apparaître naturellement la notion de topologie de Grothendieck en tant que point fixe de l'application $G_R \circ F_R$ et la notion de sous-topos ${\mathcal E} \hookrightarrow \widehat{\mathcal C}$ du topos de préfaisceaux $\widehat{\mathcal C}$ comme point fixe de l'application $F_R \circ G_R$:

\begin{thm}\label{thm2.2.2}

Considérons la relation  de dualité $R$ entre la classe $T$ des cribles d'une catégorie essentiellement petite ${\mathcal C}$ et la classe $S$ des préfaisceaux sur ${\mathcal C}$, avec la paire adjointe qu'elle définit
$$
\xymatrix{ ({\mathcal P} (T) , \subseteq) \ar@<2pt>[r]^{F_R} & ({\mathcal P} (S) , \supseteq). \ar@<2pt>[l]^{G_R} }
$$
Alors:

\begin{listeimarge}

\item Une sous-classe
$$
J \subseteq T = \{\mbox{cribles $C \hookrightarrow y(X)$ des objets $X$ de ${\mathcal C}$}\}
$$
est un point fixe de l'application $G_R \circ F_R$ si et seulement si elle est une topologie de Grothendieck au sens qu'elle satisfait les trois axiomes suivants:

\medskip

$
\left\{\begin{matrix}
(1) &\mbox{Pour tout objet $X$ de ${\mathcal C}$, le crible maximal} \hfill \\
%{ \ } \\
&C=y(X) \quad \mbox{est élément de $J$.} \\
%{ \ } \\
(2) &\mbox{Pour tout élément $C \hookrightarrow y(X)$ de $J$, son image réciproque} \hfill \\
%{ \ } \\
&x^* C = C \times_{y(X)} y(X') \hookrightarrow y(X') \quad \mbox{par tout morphisme} \quad X' \xrightarrow{ \ x \ } X \\
&\mbox{est encore élément de $C$.} \hfill \\
%{ \ } \\
(3) &\mbox{Si un crible $C \hookrightarrow y(X)$ engendré par une famille $(X_i \xrightarrow{ \, x_i \, } X)_{i \in I}$ est élément de $J$,} \hfill \\
&\mbox{tout crible $C' \hookrightarrow y(X)$ de $X$ tel que} \hfill \\
%{ \ } \\
&(x_i^* C \hookrightarrow y(X_i)) \in J \, , \quad \forall \, i \in I \, , \\
&\mbox{est encore élément de $J$.} \hfill
\end{matrix} \right.
$

%\medskip

\item Une sous-classe
$$
I \subseteq S = \{\mbox{préfaisceaux $P$ sur ${\mathcal C}$}\}
$$
est un point fixe de l'application $F_R \circ G_R$ si et seulement si la sous-catégorie pleine
$$
{\mathcal E}_I \xymatrix{\ar@{^{(}->}[r] & } \widehat{\mathcal C}
$$
sur les objets $P$ de $\widehat{\mathcal C}$ qui sont éléments de $I$ est un sous-topos du topos de préfaisceaux $\widehat{\mathcal C}$ au sens qu'elle satisfait les trois axiomes suivants:

\medskip

$\left\{\begin{matrix}
(1) &\mbox{Le foncteur pleinement fidèle de plongement} \hfill \\
&j_* : {\mathcal E}_I \xymatrix{\ar@{^{(}->}[r] & } \widehat{\mathcal C} \\
&\mbox{admet un adjoint à gauche} \hfill \\
&j^* : \widehat{\mathcal C} \longrightarrow {\mathcal E}_I \, . \\
(2) &\mbox{Cet adjoint à gauche respecte non seulement les colimites mais aussi les limites finies.} \hfill \\
(3) &\mbox{Un objet $P$ de $\widehat{\mathcal C}$ est un objet de la sous-catégorie pleine ${\mathcal E}_I$, c'est-à-dire est un élément de $I$,} \hfill \\ &\mbox{si et seulement si le morphisme canonique} \hfill \\
&P \longrightarrow j_* \circ j^* P \\
&\mbox{est un isomorphisme.} \hfill
\end{matrix}\right.
$
\end{listeimarge}
\end{thm}

\begin{demo}

Elle sera donnée dans les trois sous-paragraphes suivants c), d) et e). 
\end{demo}

Ce théorème permet d'appliquer à la dualité des cribles et des préfaisceaux les conclusions générales du corollaire \ref{cor2.1.4}.

On obtient ainsi les deux corollaires suivants:

\begin{cor}\label{cor2.2.3}

Considérons la paire adjointe
$$
\xymatrix{ ({\mathcal P} (T) , \subseteq) \ar@<2pt>[r]^{F_R} & ({\mathcal P} (S) , \supseteq) \ar@<2pt>[l]^{G_R} }
$$
définie par la relation de dualité $R$ entre la classe $T$ des cribles d'une catégorie essentiellement petite ${\mathcal C}$ et la classe $S$ des préfaisceaux sur ${\mathcal C}$.

Alors les deux applications $F_R$ et $G_R$ induisent deux bijections inverses l'une de l'autre
$$
\xymatrix{ \left\{ {{\mbox{topologies} \atop \mbox{$J$}} \atop \mbox{de ${\mathcal C}$}} \right\} \ar@<2pt>[rr]^-{\sim} && \left\{ {{\mbox{sous-topos} \atop \mbox{${\mathcal E} \hookrightarrow \widehat{\mathcal C}$}} \atop \mbox{du topos $\widehat{\mathcal C}$}} \right\}. \ar@<2pt>[ll]^-{\sim} }
$$
\end{cor}

\begin{remarks}
\begin{listeisansmarge}
\item Ici, le mot ``topologie'' a le sens d'une sous-classe $J$ de la classe $T$ des cribles de ${\mathcal C}$ qui satisfait les trois axiomes (1), (2), (3) du théorème \ref{thm2.2.2} (i).

\item De même, le mot ``sous-topos'' a ici le sens d'une sous-catégorie pleine
$$
{\mathcal E} \xymatrix{\ar@{^{(}->}[r] & } \widehat{\mathcal C}
$$
qui satisfait les trois axiomes du théorème \ref{thm2.2.2} (ii).

Il résulte du corollaire \ref{cor2.2.3} que toute telle sous-catégorie pleine vérifiant ces axiomes
$$
{\mathcal E} \xymatrix{\ar@{^{(}->}[r] & } \widehat{\mathcal C}
$$
est nécessairement de la forme
$$
{\mathcal E} = \widehat{\mathcal C}_J
$$
pour une topologie $J$ de ${\mathcal C}$.

Par conséquent, elle est un topos.

Cette conséquence avait été annoncée dans la remarque qui suit la définition \ref{defn1.2.6}.
\end{listeisansmarge}
\end{remarks}

\begin{demo}

C'est un cas particulier du corollaire \ref{cor2.1.4} (i) compte tenu de l'identification des points fixes réalisée dans le théorème \ref{thm2.2.2}. 
\end{demo}

\begin{cor}\label{cor2.2.4}

Considérons toujours la paire adjointe
$$
\xymatrix{ ({\mathcal P} (T) , \subseteq) \ar@<2pt>[r]^{F_R} & ({\mathcal P} (S) , \supseteq) \ar@<2pt>[l]^{G_R} }
$$
définie par la relation de dualité $R$ entre les classes $T$ et $S$ des cribles et des préfaisceaux sur une catégorie essentiellement petite ${\mathcal C}$.

Alors:

\begin{listeimarge}

\item Pour toute collection $J$ de cribles sur ${\mathcal C}$, l'image
$$
\overline J = G_R \circ F_R(J)
$$
est la topologie de ${\mathcal C}$ engendrée par $J$, c'est-à-dire la plus petite topologie de ${\mathcal C}$ qui contient $J$, et la sous-classe
$$
F_R(J) = F_R(\overline J)
$$
est constituée des préfaisceaux $P$ sur ${\mathcal C}$ qui sont des faisceaux pour la topologie $\overline J$.

\item Pour toute collection $I$ de préfaisceaux sur ${\mathcal C}$, l'image
$$
\overline I = F_R \circ G_R(I)
$$
est la sous-classe des préfaisceaux qui sont des objets du plus petit sous-topos de $\widehat{\mathcal C}$ qui contient les éléments de $I$, et ce sous-topos est défini par la topologie
$$
G_R (I) = G_R (\overline I) \, .
$$
\end{listeimarge}
\end{cor}

\begin{remark}
En particulier, pour toute famille de sous-topos
$$
{\mathcal E}_k = \widehat{\mathcal C}_{J_k} \xymatrix{\ar@{^{(}->}[r] & } \widehat{\mathcal C} \, , \quad k \in K \, ,
$$
définis par des topologies $J_k$ de ${\mathcal C}$, leur intersection
$$
\bigwedge_{k \in K} {\mathcal E}_k \xymatrix{\ar@{^{(}->}[r] & } \widehat{\mathcal C}
$$
a pour objets les éléments de la sous-classe
$$
F_R \left( \bigcup_{k \in K} J_k \right) = \bigcap_{k \in K} F_R (J_k) \, ,
$$
c'est-à-dire les préfaisceaux sur ${\mathcal C}$ qui sont des faisceaux pour chacune des topologies $J_k$, $k \in K$. Autrement dit, la classe des objets du sous-topos
$$
\bigwedge_{k \in K} {\mathcal E}_k \xymatrix{\ar@{^{(}->}[r] & } \widehat{\mathcal C}
$$
est l'intersection des sous-classes des objets des sous-topos ${\mathcal E}_k \hookrightarrow \widehat{\mathcal C}$, $k \in K$.
\end{remark}

\begin{demo}

Les égalités
$$
F_R(J) = F_R \circ G_R \circ F_R(J) = F_R(\overline J)
$$
et
$$
G_R(I) = G_R \circ F_R \circ G_R(I) = G_R(\overline I)
$$
sont vérifiées du fait que, d'après le corollaire \ref{cor2.1.4} (ii), toute image $F_R(J)$ [resp. $G_R(I)$] est un point fixe de l'application $F_R \circ G_R$ [resp. $G_R \circ F_R$].

Le reste des assertions (i) et (ii) est un cas particulier du corollaire \ref{cor2.1.4} (iii). 
\end{demo}

\subsection{Apparition de la notion de topologie de Grothendieck}\label{ssec2.2.3}

\smallskip

La première étape de la démonstration du théorème \ref{thm2.2.2} est la proposition suivante qui fait apparaître les trois axiomes de définition de la notion de topologie:

\begin{prop}\label{prop2.2.5}

Considérons l'application
\begin{eqnarray}
G_R : ({\mathcal P} (S), \supseteq) &\longrightarrow &({\mathcal P} (T), \subseteq) \, , \nonumber \\
(I \subseteq S) &\longmapsto &G_R(I) = \{C \hookrightarrow y(X) \mid (C \hookrightarrow y(X) , P) \in R \, , \ \forall \, P \in I \} \nonumber
\end{eqnarray}
définie par la relation de dualité $R$ entre la classe $S$ des préfaisceaux sur une catégorie essentiellement petite ${\mathcal C}$ et la classe $T$ des cribles de ${\mathcal C}$. Alors, pour toute collection $I$ de préfaisceaux, la sous-classe
$$
J = G_R(I)
$$
est une topologie au sens que

\medskip

$\left\{\begin{matrix}
(1) &\mbox{elle comprend les cribles maximaux,} \hfill \\
(2) &\mbox{elle est respectée par les applications $x^*$ d'image réciproque des cribles} \hfill \\
&\mbox{par les morphismes $x:X' \to X$ de ${\mathcal C}$,} \hfill \\
(3) &\mbox{elle est locale au sens qu'un crible $C$ d'un objet $X$ est dans $J$} \hfill \\
&\mbox{s'il en est ainsi de ses images réciproques $x_k^* C$par des morphismes} \hfill \\
&\mbox{$x_k : X_k \to X$ qui engendrent un crible élément de $J$.} \hfill
\end{matrix}\right.
$
\end{prop}

\begin{remarks}
\begin{listeisansmarge}
\item Ainsi, toute collection $I$ de préfaisceaux définit une topologie
$$
J = G_R(I) \, .
$$
C'est la topologie la plus fine pour laquelle tous les préfaisceaux $P$ éléments de $I$ sont des faisceaux.

\item On peut prendre en particulier pour $I$ la classe de tous les préfaisceaux représentables de ${\mathcal C}$ c'est-à-dire de ceux de la forme
$$
y(X) = {\rm Hom} (\bullet , X)
$$
associés aux objets $X$ de ${\mathcal C}$.

La topologie correspondante
$$
J = G_R(I)
$$
est appelée la ``topologie canonique'' de ${\mathcal C}$.

C'est la topologie la plus fine pour laquelle tous les préfaisceaux de la forme
$$
y(X) = {\rm Hom} (\bullet , X)
$$
sont des faisceaux.
\end{listeisansmarge}
\end{remarks}

\begin{demo}

Toute intersection de topologies est encore une topologie.

Par conséquent, il suffit de vérifier la proposition lorsque la classe $I$ est réduite à un unique préfaisceau $P$ et donc
$$
G_R(I) = J
$$
consiste en les cribles
$$
C \xymatrix{\ar@{^{(}->}[r] & } y(X)
$$
tels que, pour tout morphisme $X' \xrightarrow{ \, x \, } X$ de ${\mathcal C}$, l'application de restriction
$$
P(X') = {\rm Hom} (y(X'),P) \longrightarrow {\rm Hom} (C \times_{y(X)} y(X'),P)
$$
est bijective.

Si $C = y(X)$, on a $C \times_{y(X)} y(X') = y(X')$ et donc la condition (1) est vérifiée.

La condition (2) l'est également puisque la formation des images réciproques de cribles
$$
x^* : (C \hookrightarrow y(X)) \longmapsto (C \times_{y(X)} y(X') \hookrightarrow y(X'))
$$
est compatible avec la composition des morphismes de ${\mathcal C}$
$$
X'' \xrightarrow{ \ x' \ } X' \xrightarrow{ \ x \ } X
$$
au sens que
$$
x'^* (x^* C) = (x \circ x')^* C \, .
$$

Pour la condition (3) enfin, considérons un crible
$$
C \xymatrix{\ar@{^{(}->}[r] & } y(X)
$$
et une famille de morphismes
$$
x_k : X_k \longrightarrow X \, , \ k \in K \, ,
$$
engendrant un crible $C'$ de $X$ qui soit élément de $J$ et telle que tous les cribles
$$
x_k^* C \xymatrix{\ar@{^{(}->}[r] & } y(X_k)
$$
soient éléments de $J$.

Montrons que tout morphisme de $\widehat{\mathcal C}$
$$
C \longrightarrow P
$$
provient d'un unique élément de $P(X)$.

En effet, ce morphisme induit des morphismes composés
$$
x_k^* C = C \times_{y(X)} y(X_k) \longrightarrow C \longrightarrow P
$$
qui proviennent chacun d'éléments uniques
$$
p_k \in P(X_k) \, .
$$

De plus, pour tous indices $k_1 , k_2$ et tout carré commutatif de ${\mathcal C}$
$$
\xymatrix{
X' \ar[d] \ar[r] &X_{k_1} \ar[d]^{x_{k_1}} \\
X_{k_2} \ar[r]^{x_{k_2}} &X
}
$$
les sections $p_{k_1} \in P (X_{k_1})$ et $p_{k_2} \in P(X_{k_2})$ induisent la même section dans $P(X')$ qui correspond au morphisme
$$
C \times_{y(X)} y(X')  \longrightarrow C \longrightarrow P \, .
$$

Ainsi, le crible $C' \hookrightarrow y(X)$ engendré par les morphismes $X_k \xrightarrow{ \, x_k \, } X$ se trouve muni d'un morphisme
$$
C' \longrightarrow P \, .
$$
Comme ce crible est par hypothèse élément de $J$, le morphisme $C' \to P$ provient d'un unique élément
$$
p \in P(X) \, .
$$
Cet élément relève le morphisme $C \to P$.

Ainsi l'application
$$
P(X) = {\rm Hom} (y(X),P) \longrightarrow {\rm Hom} (C,P)
$$
est bijective.

De même, pour tout morphisme $X' \xrightarrow{ \, x \, }X$, l'application
$$
P(X') = {\rm Hom} (y(X'),P) \longrightarrow {\rm Hom} (C \times_{y(X)} y(X'),P)
$$
est bijective.

En effet, on sait déjà que $J$ vérifie l'axiome (2) et donc que le crible $x^* C'$ est élément de $J$.

Le résultat précédent s'applique alors au crible
$$
x^* C = C \times_{y(X)} y(X') \longrightarrow y(X')
$$
puisque, pour tout morphisme $X'' \xrightarrow{ \, x' \, }X'$ qui est élément de $C'$, le crible
$$
x'^* (x^* C) = (x \circ x')^* C \xymatrix{\ar@{^{(}->}[r] & } y(X'')
$$
est élément de $J$.

Cela montre que $J$ satisfait l'axiome (3) de localité en plus des axiomes (1) et (2) déjà vérifiés.

Autrement dit, $J$ est une topologie. 
\end{demo}

\subsection{Apparition de la notion de sous-topos}\label{ssec2.2.4}

La seconde étape de la démonstration du théorème \ref{thm2.2.2} est le corollaire suivant qui fait apparaître les trois axiomes de définition de la notion de sous-topos d'un topos de préfaisceaux:

\begin{cor}\label{cor2.2.6}

Considérons l'application
$$
F_R : ({\mathcal P}(T),\subseteq) \longrightarrow ({\mathcal P}(S),\supseteq)
$$
définie par la relation de dualité $R$ entre la classe $T$ des cribles d'une catégorie essentiellement petite ${\mathcal C}$ et la classe $S$ de ses préfaisceaux.

Alors, pour toute collection $J$ de cribles, la sous-catégorie pleine
$$
{\mathcal E}_I \xymatrix{\ar@{^{(}->}[r] & } \widehat{\mathcal C}
$$
dont les objets sont les éléments de la classe
$$
I = F_R (J) \, ,
$$
est un sous-topos au sens que

\medskip

$\left\{\begin{matrix}
(1) &\mbox{le foncteur de plongement $j_* : {\mathcal E}_I \hookrightarrow \widehat{\mathcal C}$ admet un adjoint à gauche $j^*$,} \hfill \\
(2) &\mbox{celui-ci respecte non seulement les colimites mais aussi les limites finies,} \hfill \\
(3) &\mbox{un objet $P$ de $\widehat{\mathcal C}$ est dans ${\mathcal E}_I$ si et seulement si le morphisme canonique} \hfill \\
&P \longrightarrow j_* \circ j^{*P} \\
&\mbox{est un isomorphisme.} \hfill
\end{matrix} \right.$
\end{cor}

\begin{demo}

L'image $I = F_R(J)$ est un point fixe de l'application composée
$$
F_R \circ G_R
$$
de $F_R$ avec son adjointe à droite $G_R$.

Autrement dit, on a aussi
$$
I = F_R (\overline J)
$$
si l'on note $\overline J = G_R \circ F_R(J)$.

Or on sait déjà d'après la proposition II.2.5 que
$$
\overline J = G_R \circ F_R(J)
$$
est une topologie.

Par conséquent, il n'y a pas de restriction à supposer que $J$ est une topologie et donc que
$$
{\mathcal E}_I = \widehat{\mathcal C}_J \xymatrix{\ar@{^{(}->}[r] & } \widehat{\mathcal C}
$$
est la sous-catégorie pleine des faisceaux pour la topologie $J$.

Alors le corollaire résulte de l'existence et des propriétés du foncteur de faisceautisation
$$
j^* : \widehat{\mathcal C} \longrightarrow \widehat{\mathcal C}_J
$$
adjoint à gauche du foncteur de plongement
$$
j_* : \widehat{\mathcal C}_J \xymatrix{\ar@{^{(}->}[r] & } \widehat{\mathcal C} \, ,
$$
telles qu'elles ont été rappelées dans la proposition \ref{prop1.1.10}. 
\end{demo}

\subsection{Vérification de ce que les topologies et les sous-topos sont des points fixes}\label{ssec2.2.5}

On a déjà montré que les sous-classes de cribles qui sont des points fixes de la relation de dualité entre cribles et préfaisceaux sont nécessairement des topologies de Grothendieck.

Réciproquement, les topologies sont toujours des points fixes:

\begin{lem}\label{lem2.2.7}

Considérons la paire adjointe d'applications respectant les relations d'ordre
$$
\xymatrix{ ({\mathcal P} (T) , \subseteq) \ar@<2pt>[r]^{F_R} & ({\mathcal P} (S) , \supseteq) \ar@<2pt>[l]^{G_R} }
$$
qui est définie par la relation de dualité $R$ entre la classe $T$ des cribles d'une catégorie essentiellement petite ${\mathcal C}$ et la classe $S$ des préfaisceaux sur ${\mathcal C}$.

Alors toute sous-classe $J \subseteq T$ de cribles qui est une topologie est un point fixe de la relation de dualité $R$ au sens que
$$
G_R \circ F_R (J) = J \, .
$$
\end{lem}

\begin{demo}

Partons donc d'une topologie $J$ de ${\mathcal C}$.

Alors $I = F_R(J)$ est la classe des objets de $\widehat{\mathcal C}$ qui sont dans la sous-catégorie pleine
$$
{\mathcal E}_I = \widehat{\mathcal C}_J \xymatrix{\ar@{^{(}->}[r] & } \widehat{\mathcal C} \, ,
$$
laquelle est munie du foncteur de faisceautisation
$$
j^* : \widehat{\mathcal C} \longrightarrow \widehat{\mathcal C}_J
$$
adjoint à gauche du foncteur de plongement
$$
j_* : \widehat{\mathcal C}_J \xymatrix{\ar@{^{(}->}[r] & } \widehat{\mathcal C} \, .
$$
Un crible
$$
C \xymatrix{\ar@{^{(}->}[r] & } y(X)
$$
est élément de la topologie $J$ si et seulement si le monomorphisme induit de $\widehat{\mathcal C}_J$
$$
j^* C \xymatrix{\ar@{^{(}->}[r] & } j^* y(X)
$$
est un épimorphisme et donc un isomorphisme.

Cela revient à demander que pour tout objet $P$ de $\widehat{\mathcal C}$ qui est un $J$-faisceau, c'est-à-dire un élément de $I$, l'application de restriction
$$
P(X) = {\rm Hom} (y(X) , P) \longrightarrow {\rm Hom} (C,P)
$$
est une bijection.

Comme la topologie $J$ est stable par les applications d'images réciproques par les morphismes $X' \xrightarrow{ \, x \, } X$ de ${\mathcal C}$
$$
x^* : (C \hookrightarrow y(X)) \longmapsto (x^* C = C \times_{y(X)} y(X') \hookrightarrow y(X')) \, ,
$$
cela montre comme annoncé que
$$
G_R (I) = J \, .
$$
\end{demo}

De même, on a déjà montré que les sous-classes de préfaisceaux qui sont des points fixes de la relation de dualité entre préfaisceaux et cribles sont nécessairement les classes des objets de sous-topos du topos des préfaisceaux.

Réciproquement, les classes d'objets des sous-topos sont toujours des points fixes:

\begin{prop}\label{prop2.2.8}

Considérons la paire adjointe d'applications respectant les relations d'ordre
$$
\xymatrix{ ({\mathcal P} (T) , \subseteq) \ar@<2pt>[r]^{F_R} & ({\mathcal P} (S) , \supseteq) \ar@<2pt>[l]^{G_R} }
$$
qui est définie par la relation de dualité $R$ entre la classe $S$ des préfaisceaux sur une catégorie essentiellement petite ${\mathcal C}$ et la classe $T$ des cribles de ${\mathcal C}$.

Soient une sous-classe $I \subseteq S$ de préfaisceaux sur ${\mathcal C}$ et la sous-catégorie pleine qu'elle définit
$$
{\mathcal E}_I \xymatrix{\ar@{^{(}->}[r] & } \widehat{\mathcal C} \, .
$$
Alors, si ${\mathcal E}_I \hookrightarrow \widehat{\mathcal C}$ est un sous-topos, la sous-classe $I$ de ses objets est un point fixe de la relation de dualité au sens que
$$
F_R \circ G_R (I) = I \, .
$$
\end{prop}

\begin{demo}

Par hypothèse, le foncteur de plongement
$$
j_* : {\mathcal E}_I \xymatrix{\ar@{^{(}->}[r] & } \widehat{\mathcal C}
$$
admet un adjoint à gauche
$$
j^* : \widehat{\mathcal C} \longrightarrow {\mathcal E}_I
$$
qui respecte les limites finies; et un objet $P$ de $\widehat{\mathcal C}$ est un objet de ${\mathcal E}_I$, c'est-à-dire un élément de $I$, si et seulement si le morphisme canonique
$$
P \longrightarrow j_* \circ j^* P
$$
est un isomorphisme. Enfin, comme $j_*$ est pleinement fidèle, la transformation naturelle $j^* \circ j_* \to {\rm Id}_{{\mathcal E}_I}$ est un isomorphisme.

Pour tout crible de ${\mathcal C}$
$$
C \xymatrix{\ar@{^{(}->}[r] & } y(X) \, ,
$$
demander que le morphisme induit de ${\mathcal E}_I$
$$
j^* C \longrightarrow j^* y(X)
$$
soit un isomorphisme équivaut à demander que, pour tout objet $E$ de ${\mathcal E}_I$, c'est-à-dire pour tout élément $E$ de $I$, l'application de restriction
$$
{\rm Hom} (y(X) , j_* E) \longrightarrow {\rm Hom} (C , j_* E)
$$
soit bijective.

De plus, si le morphisme
$$
j^* C \longrightarrow j^* y(X)
$$
est un isomorphisme, alors pour tout morphisme $X' \xrightarrow{ \, x \, } X$ de ${\mathcal C}$, le morphisme
$$
j^* (C \times_{y(X)} y(X')) \longrightarrow j^* y(X')
$$
est aussi un isomorphisme puisque
$$
j^* (C \times_{y(X)} y(X')) = j^* C \times_{j^* y(X)} j^* y(X') \, .
$$

Cela prouve que la sous-classe des cribles de ${\mathcal C}$
$$
C \xymatrix{\ar@{^{(}->}[r] & } y(X)
$$
qui satisfont ces deux conditions équivalentes est la topologie
$$
J = G_R(I) \, .
$$
Il reste à montrer que
$$
I = F_R(J)
$$
c'est-à-dire que, pour tout préfaisceau $P$ qui est un $J$-faisceau, le morphisme canonique de $\widehat{\mathcal C}$
$$
P \longrightarrow j_* \circ j^* P
$$
est un isomorphisme.

Etant donné un préfaisceau $P$ qui est un $J$-faisceau, considérons le plongement diagonal
$$
P \xymatrix{\ar@{^{(}->}[r] & } P \times_{j_* \circ j^* P} P \, .
$$
Pour tout morphisme
$$
y(X) \longrightarrow P \times_{j_* \circ j^* P} P \, ,
$$
c'est-à-dire toute section du préfaisceau $P \times_{j_* \circ j^* P} P$ sur un objet $X$ de ${\mathcal C}$, son produit fibré avec le plongement diagonal
$$
P \xymatrix{\ar@{^{(}->}[r] & } P \times_{j_* \circ j^* P} P
$$
est un crible
$$
C \xymatrix{\ar@{^{(}->}[r] & } y(X)
$$
qui est un élément de $J$ puisque $j^*$ respecte les produits fibrés et que 
$$
j^* P \longrightarrow j^* (P \times_{j_* \circ j^* P} P) = j^* P \times_{j^* P} j^* P
$$
est un isomorphisme, et donc aussi
$$
j^* C \longrightarrow j^* y(X) \, .
$$
Par conséquent, le morphisme
$$
C \longrightarrow P
$$
se relève en un unique morphisme
$$
y(X) \longrightarrow P
$$
c'est-à-dire en une section du préfaisceau $P$ en l'objet $X$.

Cela montre que le morphisme de $\widehat{\mathcal C}$
$$
P \longrightarrow j_* \circ j^* P
$$
est un monomorphisme.

Son produit fibré avec n'importe quel morphisme de la forme
$$
y(X) \longrightarrow j_* \circ j^* P
$$
est un crible
$$
C \xymatrix{\ar@{^{(}->}[r] & } y(X)
$$
qui est un élément de $J$ puisque $j^*$ respecte les produits fibrés et que
$$
j^* P \longrightarrow j^* \circ j_* \circ j^* P
$$
est un isomorphisme.

Par conséquent, le morphisme
$$
C \longrightarrow P
$$
se relève en un unique morphisme
$$
y(X) \longrightarrow P
$$
c'est-à-dire en une section du préfaisceau $P$ en l'objet $X$.

Cela montre que le monomorphisme de $\widehat{\mathcal C}$
$$
P \xymatrix{\ar@{^{(}->}[r] & } j_* \circ j^* P
$$
est un isomorphisme, autrement dit que $P$ est un élément de $I$.

On conclut comme annoncé que
$$
I = F_R (J) \, .
$$
\end{demo}

Le lemme \ref{lem2.2.7} et la proposition \ref{prop2.2.8} terminent la démonstration du théorème \ref{thm2.2.2}.

\section{La dualité des monomorphismes et des objets d'un topos}\label{sec2.3}

\subsection{Définition d'une relation entre monomorphismes et objets}\label{ssec2.3.1}

Nous allons donner maintenant une variante de la dualité du paragraphe précédent entre cribles et préfaisceaux, pour la rendre intrinsèque aux topos et la généraliser au-delà du cadre des topos de préfaisceaux.

Commençons par définir une relation naturelle entre les monomorphismes d'un topos de Grothendieck et ses objets.

\begin{propqed}\label{prop2.3.1}

Soit un topos ${\mathcal E}$.

Soit $T$ la classe des monomorphismes de ${\mathcal E}$
$$
C \xymatrix{\ar@{^{(}->}[r] & } X \, ,
$$
et soit $S$ la classe des objets $E$ de ${\mathcal E}$.

On appellera ``dualité des monomorphismes et des objets'' la relation
$$
R \xymatrix{\ar@{^{(}->}[r] & } T \times S
$$
constituée des paires d'éléments
$$
(C \xymatrix{\ar@{^{(}->}[r] & } X , E)
$$
telles que, pour tout morphisme de ${\mathcal E}$
$$
X' \xrightarrow{ \ x \ } X \, ,
$$
l'application de restriction au sous-objet $C \times_X X' \hookrightarrow X'$
$$
{\rm Hom} (X',E) \longrightarrow {\rm Hom} (C \times_X X' , E)
$$
est une bijection. 
\end{propqed}

Toujours d'après les résultats généraux du lemme \ref{lem2.1.3}, la relation $R$ de dualité entre les classes $T$ et $S$ des monomorphismes et des objets d'un topos définit une paire adjointe d'applications respectant les relations d'ordre
$$
\xymatrix{
({\mathcal P} (T) , \subseteq) \ar@<2pt>[r]^{F_R} & ({\mathcal P} (S) , \supseteq), \ar@<2pt>[l]^{G_R}
}
$$
et on sait d'après le corollaire \ref{cor2.1.4} que ces applications induisent deux bijections inverses l'une de l'autre entre les points fixes de ${\mathcal P} (T)$ et ceux de ${\mathcal P} (S)$.

L'identification des points fixes va faire apparaître une reformulation de la notion de topologie sur un topos en termes de classes de monomorphismes, et de l'autre côté la notion générale de sous-topos d'un topos comme sous-catégorie pleine satisfaisant certaines conditions.

\subsection{La dualité induite des topologies et des sous-topos d'un topos}\label{ssec2.3.2}

\medskip

Nous allons démontrer le théorème suivant qui identifie les points fixes de part et d'autre de la dualité naturelle entre monomorphismes et objets d'un topos de Grothendieck.

\begin{thm}\label{thm2.3.2}

Considérons la relation de dualité $R$ entre la classe $T$ des monomorphismes d'un topos ${\mathcal E}$ et la classe $S$ de ses objets, avec la paire adjointe qu'elle définit
$$
\xymatrix{
({\mathcal P} (T) , \subseteq) \ar@<2pt>[r]^{F_R} & ({\mathcal P} (S) , \supseteq). \ar@<2pt>[l]^{G_R}
}
$$

Alors :

\begin{listeimarge}

\item Une sous-classe
$$
J \subseteq T = \{\mbox{monomorphismes $C \hookrightarrow X$ de ${\mathcal E}$}\}
$$
est un point fixe de l'application $G_R \circ F_R$ si et seulement si elle est une topologie au sens qu'elle satisfait les trois axiomes suivants:

\medskip

\noindent $\left\{\begin{matrix}
(1) &\mbox{Tout isomorphisme $C \xrightarrow{ \ \sim \ }X$ est élément de $J$.} \hfill \\
(2) &\mbox{Pour tout élément $C \hookrightarrow X$ de $J$ et tout morphisme $X' \to X$ de ${\mathcal E}$,} \hfill \\
&\mbox{le monomorphisme $C \times_X X' \hookrightarrow X'$ est encore élément de $J$.} \hfill \\
(3) &\mbox{Pour qu'un monomorphisme $C \hookrightarrow X$ soit élément de $J$, il suffit qu'existe un élément $C' \hookrightarrow X$ de $J$} \hfill \\
&\mbox{et une famille globalement épimorphique $(X_k \to C')_{k \in K}$ tels que tous les produits fibrés} \hfill \\
&C \times_X X_k \xymatrix{\ar@{^{(}->}[r] & } X_k \, , \ k \in K \, , \ \mbox{soient éléments de $J$.}
\end{matrix}\right.
$

%\medskip

\item Une sous-classe
$$
I \subseteq S = \{\mbox{objets $E$ de ${\mathcal E}$}\}
$$
est un point fixe de l'application $F_R \circ G_R$ si et seulement si la sous-catégorie pleine
$$
{\mathcal E}_I \xymatrix{\ar@{^{(}->}[r] & } {\mathcal E}
$$
sur les objets $E$ de ${\mathcal E}$ qui sont éléments de $I$ est un sous-topos du topos ${\mathcal E}$ au sens qu'elle satisfait les trois axiomes suivants:

\medskip

\noindent $\left\{\begin{matrix}
(1) &\mbox{Le foncteur pleinement fidèle de plongement} \hfill \\
&j_* : {\mathcal E}_I \xymatrix{\ar@{^{(}->}[r] & } {\mathcal E} \\
&\mbox{admet un adjoint à gauche} \hfill \\
&j^* : {\mathcal E} \longrightarrow {\mathcal E}_I \, . \\
(2) &\mbox{Cet adjoint à gauche respecte non seulement les colimites mais aussi les limites finies.} \hfill \\
(3) &\mbox{Un objet $E$ de ${\mathcal E}$ est un objet de ${\mathcal E}_I$, c'est-à-dire est élément de $I$,} \hfill \\
&\mbox{si et seulement si le morphisme canonique} \hfill \\
&E \longrightarrow j_* \circ j^* E \\
&\mbox{est un isomorphisme.} \hfill
\end{matrix}\right.
$
\end{listeimarge}
\end{thm}

\begin{remarks}
\begin{listeisansmarge}
\item Si le topos ${\mathcal E}$ est présenté comme topos des faisceaux sur un site $({\mathcal C} , J_{\mathcal E})$
$$
\widehat{\mathcal C}_{J_{\mathcal E}} \xrightarrow{ \ \sim \ } {\mathcal E} \, ,
$$
muni du foncteur canonique
$$
\ell : {\mathcal C} \longrightarrow \widehat{\mathcal C}_{J_{\mathcal E}} \xrightarrow{ \ \sim \ } {\mathcal E} \, ,
$$
se donner une topologie sur ${\mathcal E}$ au sens de la partie (i) du théorème équivaut à se donner une topologie de Grothendieck sur ${\mathcal C}$ qui raffine la topologie $J_{\mathcal E}$.

En effet, toute topologie $J \subseteq T$ de ${\mathcal E}$ définit une notion de famille $J$-couvrante de morphismes de ${\mathcal C}$ de même but
$$
(X_k \xrightarrow{ \ x_k \ } X)_{k \in K}
$$
par la condition que le sous-objet
$$
C \xymatrix{\ar@{^{(}->}[r] & } \ell (X)
$$
défini comme la réunion des images des morphismes
$$
\ell (x_k) : \ell (X_k) \longrightarrow \ell(X) \, , \ k \in K \, ,
$$
soit élément de la classe $J$.

Cette notion raffine celle de famille $J_{\mathcal E}$-couvrante
$$
(X_k \xrightarrow{ \ x_k \ } X)_{k \in K}
$$
caractérisée par la condition que la famille transformée par $\ell$
$$
\ell (x_k) : \ell (X_k) \longrightarrow \ell(X) \, , \quad k \in K \, ,
$$
soit globalement épimorphique.

Réciproquement, si $J$ est une topologie de Grothendieck sur ${\mathcal C}$ qui raffine $J_{\mathcal E}$, elle définit une topologie sur le topos ${\mathcal E}$ en décidant qu'un monomorphisme de ${\mathcal E}$
$$
C \xymatrix{\ar@{^{(}->}[r] & } X
$$
est élément de la topologie induite par $J$ si, pour n'importe quelle famille globalement épimorphique de morphismes de ${\mathcal E}$
$$
\ell (X_k) \longrightarrow X \, , \quad k \in K \, ,
$$
dont les sources sont les images par $\ell$ d'objets $X_k$ de ${\mathcal C}$, tous les cribles
$$
y(X_k) \times_X C \xymatrix{\ar@{^{(}->}[r] & } y(X_k) \, , \quad k \in K \, ,
$$
sont $J$-couvrants.

\item Toujours si le topos ${\mathcal E}$ est présenté comme topos des faisceaux sur un site $({\mathcal C} , J_{\mathcal E})$, se donner un sous-topos de ${\mathcal E}$ équivaut à se donner un sous-topos du topos de préfaisceaux $\widehat{\mathcal C}$ qui est contenu dans le sous-topos
$$
\widehat{\mathcal C}_{J_{\mathcal E}} \xymatrix{\ar@{^{(}->}[r] & } \widehat{\mathcal C} \, .
$$
\end{listeisansmarge}
\end{remarks}

\begin{demo}

Elle sera donnée dans les trois sous-paragraphes suivants c), d) et e). 
\end{demo}

Ce théorème permet d'appliquer à la dualité des monomorphismes et des objets d'un topos les conclusions générales du corollaire \ref{cor1.1.6}.

On obtient ainsi les deux corollaires suivants:

\begin{cor}\label{cor2.3.3}

Considérons la paire adjointe
$$
\xymatrix{ ({\mathcal P} (T) , \subseteq) \ar@<2pt>[r]^{F_R} & ({\mathcal P} (S) , \supseteq) \ar@<2pt>[l]^{G_R} }
$$
définie par la relation de dualité $R$ entre la classe $T$ des monomorphismes d'un topos ${\mathcal E}$ et la classe $S$ de ses objets.

Alors les deux applications $F_R$ et $G_R$ induisent deux bijections inverses l'une de l'autre
$$
\xymatrix{ \left\{ {{\mbox{topologies} \atop \mbox{$J$}} \atop \mbox{de ${\mathcal E}$}} \right\} \ar@<2pt>[rr]^-{\sim} && \left\{ {{\mbox{sous-topos} \atop \mbox{${\mathcal E}_I \hookrightarrow {\mathcal E}$}} \atop \mbox{du topos ${\mathcal E}$}} \right\}. \ar@<2pt>[ll]^-{\sim} }
$$
\end{cor}

\begin{remarks}
	\begin{listeisansmarge}
\item Ici, le mot ``topologie'' a le sens d'une sous-classe $J$ de la classe $T$ des monomorphismes de ${\mathcal E}$ qui satisfait les trois axiomes (1), (2), (3) du théorème \ref{thm2.3.2} (i).

\item De même, le mot ``sous-topos'' a ici le sens d'une sous-catégorie pleine
$$
{\mathcal E} \xymatrix{\ar@{^{(}->}[r] & } \widehat{\mathcal C}
$$
qui satisfait les trois axiomes du théorème \ref{thm2.3.2} (ii).

Il résulte de la remarque (ii) qui suit le théorème \ref{thm2.3.2} et du corollaire \ref{cor2.2.3} que toute telle sous-catégorie pleine vérifiant ces axiomes est un topos.

Cette conséquence avait été annoncée dans la remarque qui suit la définition \ref{defn1.2.6}.

\item Si ${\mathcal E}$ est présenté comme le topos des faisceaux sur un site $({\mathcal C}, J_{\mathcal E})$, il résulte des remarques (i) et (ii) qui suivent le théorème \ref{thm2.3.2} que le corollaire retrouve la correspondance bijective entre topologies de Grothendieck de ${\mathcal C}$ qui raffinent $J_{\mathcal E}$ et sous-topos du topos de faisceaux
$$
\widehat{\mathcal C}_{J_{\mathcal E}} \xrightarrow{ \ \sim \ } {\mathcal E} \, .
$$
\end{listeisansmarge}
\end{remarks}

\begin{demo}

C'est un cas particulier du corollaire \ref{cor2.1.4} (i) compte tenu de l'identification des points fixes réalisée dans le théorème \ref{thm2.3.2}. 
\end{demo}

\begin{cor}\label{cor2.3.4}

Considérons toujours la paire adjointe
$$
\xymatrix{ ({\mathcal P} (T) , \subseteq) \ar@<2pt>[r]^{F_R} & ({\mathcal P} (S) , \supseteq) \ar@<2pt>[l]^{G_R} }
$$
définie par la relation de dualité $R$ entre les classes $T$ et $S$ des monomorphismes et des objets d'un topos ${\mathcal E}$.

Alors:

\begin{listeimarge}

\item Pour toute collection $J$ de monomorphismes de ${\mathcal E}$, l'image
$$
\overline J = G_R \circ F_R (J)
$$
est la topologie de ${\mathcal E}$ engendrée par $J$, c'est-à-dire la plus petite topologie de ${\mathcal E}$ qui contient $J$, et la sous-classe
$$
F_R(J) = F_R (\overline J)
$$
est constituée des objets $E$ de ${\mathcal E}$ qui sont dans le sous-topos correspondant.

\item Pour toute collection $I$ d'objets de ${\mathcal E}$, l'image
$$
\overline I = F_R \circ G_R(I)
$$
est la sous-classe des objets du plus petit sous-topos de ${\mathcal E}$ qui contient les éléments de $I$, et ce sous-topos est défini par la topologie
$$
G_R(I) = G_R (\overline I) \, .
$$
\end{listeimarge}
\end{cor}

\begin{remark}

En particulier, pour toute famille de sous-topos
$$
{\mathcal E}_k \xymatrix{\ar@{^{(}->}[r] & }  {\mathcal E} \, , \quad k \in K \, ,
$$
définis par des topologies $J_k$ de ${\mathcal E}$, leur intersection
$$
\bigwedge_{k \in K} {\mathcal E}_k \xymatrix{\ar@{^{(}->}[r] & } {\mathcal E}
$$
a pour objets les éléments de la sous-classe
$$
F_R \left( \bigcup_{k \in K} J_k \right) = \bigcap_{k \in K} F_R (J_k) \, .
$$
Autrement dit, la classe des objets du sous-topos
$$
\bigwedge_{k \in K} {\mathcal E}_k \xymatrix{\ar@{^{(}->}[r] & } {\mathcal E}
$$
est l'intersection des sous-classes des objets des sous-topos ${\mathcal E}_k \hookrightarrow {\mathcal E}$, $k \in K$.
\end{remark}

\begin{demo}

Les égalités
$$
F_R(J) = F_R \circ G_R \circ F_R (J) = F_R (\overline J)
$$
et
$$
G_R(I) = G_R \circ F_R \circ G_R (I) = G_R (\overline I)
$$
sont vérifiées du fait que, d'après le corollaire \ref{cor2.1.4} (ii), chaque image de l'application $F_R$ [resp. $G_R$] est un point fixe de l'application $F_R \circ G_R$ [resp. $G_R \circ F_R$].

Le reste des assertions (i) et (ii) est un cas particulier du corollaire \ref{cor2.1.4} (iii). 
\end{demo}

\subsection{Apparition de la notion de topologie sur un topos}\label{ssec2.3.3}

\medskip

La première étape de la démonstration du théorème \ref{thm2.3.2} est la proposition suivante:

\begin{prop}\label{prop2.3.5}

Considérons l'application
$$
\begin{matrix}
G_R : &({\mathcal P} (S) , \supseteq) &\longrightarrow &({\mathcal P} (T) , \subseteq) \hfill \\
&\hfill (I \subseteq S) &\longmapsto &G_R(I) = \{ C \hookrightarrow X \mid (C \hookrightarrow X,E) \in R \, , \ \forall \, E \in I \}
\end{matrix}
$$
définie par la relation de dualité $R$ entre la classe $S$ des objets d'un topos ${\mathcal E}$ et la classe $T$ des monomorphismes $C \hookrightarrow X$ de ${\mathcal E}$.

Alors, pour toute collection $I$ d'objets de ${\mathcal E}$, la sous-classe
$$
J = G_R(I)
$$
est une topologie de ${\mathcal E}$ au sens que

\medskip

$\left\{\begin{matrix}
(1) &\mbox{elle comprend les isomorphismes $C \xrightarrow{ \, \sim \, } X$,} \hfill \\
(2) &\mbox{elle est respectée par changement de base par tout morphisme $X' \to X$ de ${\mathcal E}$} \hfill \\
&(C \hookrightarrow X) \longmapsto (C \times_X X' \hookrightarrow X') \, , \\
(3) &\mbox{elle est locale au sens qu'un monomorphisme $C \hookrightarrow X$ est dans $J$} \hfill \\
&\mbox{s'il en est ainsi de ses images réciproques $C \times_X X_k \hookrightarrow X_k$ par} \hfill \\
&\mbox{une famille de morphismes $X_k \to X$, $k \in K$, dont la réunion des images} \hfill \\
&\mbox{est un sous-objet $C' \hookrightarrow X$ élément de $J$.} \hfill
\end{matrix}\right.$
\end{prop}

\begin{remark}

Si $I$ est la classe de tous les objets de ${\mathcal E}$,
$$
J = G_R(I)
$$
est la topologie dont les seuls éléments sont les isomorphismes
$$
C \xrightarrow{ \ \sim \ } X \, .
$$
\end{remark}

\begin{demo}

Toute intersection de topologies est encore une topologie.

Par conséquent, il suffit de vérifier la proposition lorsque la classe $I$ est résuite à un unique objet $E$ de ${\mathcal E}$.

Dans ce cas, $J = G_R(I)$ consiste en les monomorphismes de ${\mathcal E}$
$$
C \xymatrix{\ar@{^{(}->}[r] & } X
$$
tels que, pour tout morphisme $X' \to X$ de ${\mathcal E}$, l'application de restriction au sous-objet $C \times_X X' \hookrightarrow X'$
$$
{\rm Hom} (X',E) \longrightarrow {\rm Hom} (C \times_X X' ,E)
$$
est bijective.

Si $C \xrightarrow{ \ \sim \ } X$ est un isomorphisme, il en est de même des
$$
C \times_X X' \xrightarrow{ \ \sim \ } X'
$$
et donc la condition (1) est vérifiée.

La condition (2) l'est également puisque la formation des produits fibrés $C \times_X X'$ est compatible avec la composition des morphismes de changement de base $X'' \to X' \to X$ au sens que
$$
(C \times_X X') \times_{X'} X'' = C \times_X X'' \, .
$$

Pour la condition (3) enfin, considérons un monomorphisme
$$
C \xymatrix{\ar@{^{(}->}[r] & } X
$$
et une famille de morphismes
$$
x_k : X_k \longrightarrow X \, , \quad k \in K \, ,
$$
dont la réunion des images est un élément de $J$
$$
C' \xymatrix{\ar@{^{(}->}[r] & } X
$$
et telle que chaque monomorphisme
$$
C \times_X X_k \longrightarrow X_k \, , \ k \in K \, ,
$$
soit élément de $J$.

Montrons que tout morphisme de ${\mathcal E}$
$$
C \longrightarrow E
$$
provient par restriction d'un unique morphisme $X \to E$.

Pour tout indice $k \in K$, le morphisme composé
$$
C \times_X X_k \longrightarrow C \longrightarrow E
$$
provient par restriction d'un unique morphisme
$$
X_k \longrightarrow E \, .
$$
De plus, pour tous indices $k_1 , k_2 \in K$, les deux morphismes
$$
X_{k_1} \longrightarrow E \quad \mbox{et} \quad X_{k_2} \longrightarrow E
$$
induisent le même morphisme
$$
X_{k_1} \times_X X_{k_2} \longrightarrow E
$$
puisqu'ils induisent le même morphisme
$$
X_{k_1} \times_X C \times_X X_{k_2} \longrightarrow C \longrightarrow E \, .
$$
Or le sous-objet
$$
C' \xymatrix{\ar@{^{(}->}[r] & } X \, ,
$$
qui est réunion des images des morphismes $X_k \to X$, $k \in K$, s'écrit comme la colimite du diagramme constitué des objets
$$
X_k \, , \quad k \in K \, ,
$$
et
$$
X_{k_1} \times_X X_{k_2} \, , \quad k_1 , k_2 \in K \, ,
$$
reliés par les morphismes de projection
$$
X_{k_1} \longleftarrow X_{k_1} \times_X X_{k_2} \longrightarrow X_{k_2} \, .
$$
On en déduit que la famille compatible de morphismes
$$
X_k \longrightarrow E \, , \quad k \in K \, ,
$$
provient d'un unique morphisme
$$
C' \longrightarrow E \, .
$$
Comme le monomorphisme $C' \hookrightarrow X$ est élément de $J$, ce morphisme
$$
C' \longrightarrow E
$$
provient à son tour d'un unique morphisme
$$
X \longrightarrow E \, .
$$
Celui-ci relève le morphisme $C \to E$.

Ainsi, l'application
$$
{\rm Hom} (X,E) \longrightarrow {\rm Hom} (C,E)
$$
est bijective.

De même, pour tout morphisme $X' \to X$, l'application
$$
{\rm Hom} (X',E) \longrightarrow {\rm Hom} (C \times_X X' ,E)
$$
de restriction au sous-objet $C \times_X X' \hookrightarrow X'$ est bijective.

En effet, on sait déjà que $J$ vérifie l'axiome (2) et donc que le sous-objet
$$
C' \times_X X' \xymatrix{\ar@{^{(}->}[r] & } X' \, ,
$$
qui est la réunion des images des morphismes
$$
X_k \times_X X' \longrightarrow X' \, , \ k \in K \, ,
$$
est élément de $J$.

Le résultat précédent s'applique alors au monomorphisme
$$
C \times_X X' \xymatrix{\ar@{^{(}->}[r] & } X'
$$
puisque, pour tout indice $k \in K$, le monomorphisme
$$
(C \times_X X') \times_{X'} (X_k \times_X X') \xymatrix{\ar@{^{(}->}[r] & } X_k \times_X X'
$$
est élément de $J$.

Cela montre que $J$ satisfait l'axiome (3) de localité en plus des axiomes (1) et (2) déjà vérifiés.

Autrement dit, $J$ est une topologie.
\end{demo}

\subsection{Apparition de la notion de sous-topos}\label{ssec2.3.4}

\smallskip

La seconde étape de la démonstration du théorème \ref{thm2.3.2} est le corollaire suivant:

\begin{cor}\label{cor2.3.6}

Considérons l'application
$$
F_R : ({\mathcal P} (T) , \subseteq) \longrightarrow ({\mathcal P} (S) , \supseteq)
$$
définie par la relation de dualité $R$ entre la classe $T$ des monomorphismes $C \hookrightarrow X$ d'un topos ${\mathcal E}$ et la classe $S$ de ses objets.

Alors, pour toute collection $J$ de monomorphismes de ${\mathcal E}$, la sous-catégorie pleine
$$
{\mathcal E}_I \xymatrix{\ar@{^{(}->}[r] & } {\mathcal E}
$$
dont les objets sont les éléments de la classe
$$
I = F_R(J)
$$
est un sous-topos au sens que

\medskip

$\left\{\begin{matrix}
(1) &\mbox{le foncteur de plongement $j_* : {\mathcal E}_I \hookrightarrow {\mathcal E}$ admet un adjoint à gauche $j^* : {\mathcal E} \to {\mathcal E}_I$,} \hfill \\
(2) &\mbox{celui-ci respecte non seulement les colimites mais aussi les limites finies,} \hfill \\
(3) &\mbox{un objet $E$ de ${\mathcal E}$ est dans ${\mathcal E}_I$ si et seulement si le morphisme canonique} \hfill \\
&E \longrightarrow j_* \circ j^* E \\
&\mbox{est un isomorphisme.} \hfill
\end{matrix}\right.$
\end{cor}

\begin{demo}

L'image $I = F_R(J)$ est un point fixe de l'application composée
$$
F_R \circ G_R
$$
de $F_R$ avec son adjoint à droite $G_R$, si bien que l'on a aussi
$$
I = F_R (\overline J)
$$
si l'on note $\overline J = G_R \circ F_R(J)$.

Or, on sait déjà d'après la proposition \ref{prop2.3.5} que $\overline J$ est une topologie de ${\mathcal E}$.

Par conséquent, il n'y a pas de restriction à supposer que $J$ est une topologie.

Si le topos ${\mathcal E}$ est présenté comme topos des faisceaux sur un site $({\mathcal C}, J_{\mathcal E})$, il résulte de la remarque (i) qui suit l'énoncé du théorème \ref{thm2.3.2} que la topologie $J$ de ${\mathcal E}$ correspond à une topologie de Grothendieck de ${\mathcal C}$ qui raffine $J_{\mathcal E}$ et que l'on peut encore noter $J$.

Alors la classe
$$
I = F_R(J)
$$
a pour éléments les objets de
$$
\widehat{\mathcal C}_{J_{\mathcal E}} \xrightarrow{ \ \sim \ } {\mathcal E}
$$
qui sont des faisceaux pour la topologie $J$.

Ils définissent la sous-catégorie pleine
$$
\widehat{\mathcal C}_J \xrightarrow{ \ \sim \ } {\mathcal E}_I \quad \mbox{de} \quad \widehat{\mathcal C}_{J_{\mathcal E}} \xrightarrow{ \ \sim \ } {\mathcal E} \, .
$$
Alors le foncteur de plongement $j_* : {\mathcal E}_I \hookrightarrow {\mathcal E}$ admet pour adjoint à gauche la restriction du foncteur de faisceautisation $\widehat{\mathcal C} \to \widehat{\mathcal C}_J$.

D'où la conclusion. 
\end{demo}

\subsection{Vérification de ce que les topologies et les sous-topos sont des points fixes}\label{ssec2.3.5}

\smallskip

On a déjà montré que les sous-classes de monomorphismes d'un topos ${\mathcal E}$ qui sont des points fixes de la relation de dualité entre monomorphismes et objets de ${\mathcal E}$ sont nécessairement des topologies.

Réciproquement, les topologies sont des points fixes:

\begin{lem}\label{lem2.3.7}

Considérons la paire adjointe d'applications respectant les relations d'ordre
$$
 \xymatrix{ ({\mathcal P} (T) , \subseteq) \ar@<2pt>[r]^{F_R} & ({\mathcal P} (S) , \supseteq) \ar@<2pt>[l]^{G_R} }
$$
qui est définie par la relation de dualité $R$ entre la classe $T$ des monomorphismes d'un topos ${\mathcal E}$ et la classe $S$ de ses objets.

Alors toute sous-classe $J \subseteq T$ de monomorphismes qui est une topologie est un point fixe de la relation de dualité $R$ au sens que
$$
G_R \circ F_R(J) = J \, .
$$
\end{lem}

\begin{demo}

Partons donc d'une topologie $J$ de ${\mathcal C}$.

Alors la classe $I = F_R(J)$ définit la sous-catégorie pleine
$$
{\mathcal E}_I \xymatrix{\ar@{^{(}->}[r] & } {\mathcal E}
$$
dont le foncteur de plongement $j_*$ admet un adjoint à gauche
$$
j^* : {\mathcal E} \longrightarrow {\mathcal E}_I
$$
qui respecte les limites finies.

Un monomorphisme $C \hookrightarrow X$ de ${\mathcal E}$ est élément de la classe
$$
\overline J = G_R (I) = G_R \circ F_R (J)
$$
si et seulement si, pour tout morphisme $X' \to X$ et tout objet $E$ de ${\mathcal E}_I$, l'application de restriction
$$
\begin{matrix}
{\rm Hom} (X' , j_* E) &\longrightarrow &{\rm Hom} (C \times_X X' , j_* E) \\
\Vert &&\Vert \\
{\rm Hom} (j^* X',E) &&{\rm Hom} (j^* C \times_{j^* X} j^* X' , E)
\end{matrix}
$$
est une bijection.

Cela revient à demander que le morphisme de ${\mathcal E}_I$
$$
j^* C \longrightarrow j^* X
$$
soit un isomorphisme, autrement dit que le morphisme,
$$
C \longrightarrow X
$$
soit couvrant pour la topologie $J$ qui définit le sous-topos
$$
{\mathcal E}_I \xymatrix{\ar@{^{(}->}[r] & } {\mathcal E} \, .
$$
D'où le résultat. 
\end{demo}

De même, on a déjà montré que les sous-classes d'objets d'un topos ${\mathcal E}$ qui sont des points fixes de la relation de dualité entre monomorphismes et objets de ${\mathcal E}$ sont des sous-topos.

Réciproquement, les classes d'objets des sous-topos sont toujours des points fixes.

\begin{prop}\label{prop2.3.8}

Considérons la paire adjointe d'applications respectant les relations d'ordre
$$
 \xymatrix{ ({\mathcal P} (T) , \subseteq) \ar@<2pt>[r]^{F_R} & ({\mathcal P} (S) , \supseteq) \ar@<2pt>[l]^{G_R} }
$$
définie par la relation de dualité $R$ entre la classe $S$ des objets d'un topos ${\mathcal E}$ et la classe $T$ de ses monomorphismes.

Soient une sous-classe $I \subseteq S$ d'objets de ${\mathcal E}$ et la sous-catégorie pleine qu'elle définit
$$
{\mathcal E}_I \xymatrix{\ar@{^{(}->}[r] & } {\mathcal E} \, .
$$

Alors, si ${\mathcal E}_I \hookrightarrow {\mathcal E}$ est un sous-topos, la sous-classe $I$ de ses objets est un point fixe de la relation de dualité au sens que
$$
F_R \circ G_R (I) = I \, .
$$
\end{prop}

\begin{demo}

Par hypothèse, le foncteur de plongement
$$
j_* : {\mathcal E}_I \xymatrix{\ar@{^{(}->}[r] & } {\mathcal E}
$$
admet un adjoint à gauche
$$
j^* : {\mathcal E} \longrightarrow {\mathcal E} _I
$$
qui respecte les limites finies; et un objet $E$ de ${\mathcal E} $ est un objet de ${\mathcal E}_I$, c'est-à-dire un élément de $I$, si et seulement si le morphisme canonique
$$
E \longrightarrow j_* \circ j^* E
$$
est un isomorphisme. Enfin, comme $j_*$ est pleinement fidèle, la transformation naturelle $j^* \circ j_* \to {\rm Id}_{{\mathcal E}_I}$ est un isomorphisme.

Les éléments de la topologie de ${\mathcal E}$
$$
J = G_R(I)
$$
sont les monomorphismes
$$
C \xymatrix{\ar@{^{(}->}[r] & } X
$$
dont le transformé par $j^*$
$$
j^* C \xymatrix{\ar@{^{(}->}[r] & } j^* X
$$
sont des isomorphismes de ${\mathcal E}_I$.

Il reste à montrer que
$$
I = F_R(J)
$$
c'est-à-dire que, pour tout objet $E$ de ${\mathcal E}$ tel que
$$
{\rm Hom} (X,E) = {\rm Hom} (C,E) \quad \mbox{pour tout élément $(C \hookrightarrow X)$ de $J$,}
$$
le morphisme canonique de ${\mathcal E}$
$$
E \longrightarrow j_* \circ j^* E
$$
est un isomorphisme.

Or, le monomorphisme diagonal
$$
E \xymatrix{\ar@{^{(}->}[r] & } E \times_{j_* \circ j^* E} E
$$
est élément de $J$ puisque le foncteur $j^*$ le transforme en un isomorphisme de ${\mathcal E}_I$.

Donc le morphisme ${\rm Id}_E : E \to E$ admet un unique relèvement
$$
E \times_{j_* \circ j^* E} E \longrightarrow E
$$
et les deux projections
$$
E \times_{j_* \circ j^* E} E \rightrightarrows E
$$
coïncident.

Cela signifie que le morphisme
$$
E \longrightarrow j_* \circ j^* E
$$
est un monomorphisme.

Il est alors élément de $J$ puisque son transformé par le foncteur $j^*$ est un isomorphisme.

Pour tout morphisme $X \to j_* \circ j^* E$ de ${\mathcal E}$, le produit fibré
$$
E \times_{j_* \circ j^* E} X \xymatrix{\ar@{^{(}->}[r] & } X
$$
est encore élément de $J$, et donc le morphisme
$$
E \times_{j_* \circ j^* E} X \longrightarrow E
$$
se relève en un unique morphisme
$$
X \longrightarrow E \, .
$$
Cela implique que le monomorphisme
$$
E \xymatrix{\ar@{^{(}->}[r] & } j_* \circ j^* E
$$
est un isomorphisme, comme on voulait. 
\end{demo}

Le lemme \ref{lem2.3.7} et la proposition \ref{prop2.3.8} terminent la démonstration du théorème \ref{thm2.3.2}.

\section{La dualité des cribles et des plongements de préfaisceaux}\label{sec2.4}

\subsection{Définition d'une relation entre cribles et plongements de préfaisceaux}\label{ssec2.4.1}

\medskip

Introduisons encore une relation naturelle, cette fois entre les cribles d'une catégorie essentiellement petite et les plongements de préfaisceaux sur cette catégorie, à laquelle nous appliquerons à nouveau les résultats généraux du paragraphe \ref{sec2.1} sur les correspondances galoisiennes.

Voici la définition de cette relation:

\begin{defn}\label{defn2.4.1}

Soit une catégorie essentiellement petite ${\mathcal C}$, munie du foncteur de Yoneda
$$
y : {\mathcal C} \xymatrix{\ar@{^{(}->}[r] & } \widehat{\mathcal C} = [{\mathcal C}^{\rm op} , {\rm Ens}] \, .
$$
Soit $T$ la classe des cribles $C \hookrightarrow y(X)$ des objets $X$ de ${\mathcal C}$.

Soit $S$ la classe des plongements de $\widehat{\mathcal C}$
$$
Q \xymatrix{\ar@{^{(}->}[r] & } P \, .
$$
On appellera ``dualité des cribles et des sous-préfaisceaux'' de ${\mathcal C}$ la relation
$$
R \xymatrix{\ar@{^{(}->}[r] & } T \times S
$$
constituée des paires d'éléments
$$
(C \xymatrix{\ar@{^{(}->}[r] & } y(X) \, , \ Q \xymatrix{\ar@{^{(}->}[r] & } P)
$$
telles que
\medskip

$\left\{\begin{matrix}
\mbox{pour tout morphisme} &X' \xrightarrow{ \ x \ } X \ \mbox{de ${\mathcal C}$ et tout élément} \quad p \in P(X') \, ,\hfill \\
\mbox{on a} \hfill &p \in Q (X') \subseteq P(X') \hfill \\
\mbox{si tout morphisme} \hfill&(X'' \xrightarrow{ \ x' \ } X') \in x^* C \ \mbox{envoie $p$ dans le sous-ensemble $Q(X'') \subseteq P(X'')$.}
\end{matrix}\right.$
\end{defn}

\begin{remarksqed}
	\begin{listeisansmarge}
\item Par ``plongements'' de préfaisceaux sur ${\mathcal C}$, on entend les monomorphismes de $\widehat{\mathcal C}$
$$
Q \xymatrix{\ar@{^{(}->}[r] & } P
$$
tels que, en tout objet $X$ de ${\mathcal C}$, la composante
$$
Q(X) \xymatrix{\ar@{^{(}->}[r] & } P(X)
$$
soit l'application d'inclusion d'un sous-ensemble $Q(X)$ de l'ensemble $P(X)$.

\item Si le crible $x^* C \hookrightarrow y(X')$ image réciproque du crible $C \hookrightarrow y(X)$ par un morphisme $X' \xrightarrow{ \ x \ } X$ est engendré par une famille de morphismes
$$
x'_k : X'_k \longrightarrow X' \, , \quad k \in K \, ,
$$
la condition de la définition se résume à demander que tout élément
$$
p \in P(X')
$$
soit élément du sous-ensemble
$$
Q(X') \subseteq P(X')
$$
si les morphismes $x'_k : X'_k \to X'$ le transforment en des éléments des sous-ensembles
$$
Q(X'_k) \subseteq P(X'_k) \, , \quad k \in K \, .
$$
\end{listeisansmarge}
\end{remarksqed}

\subsection{La dualité induite des topologies et des propriétés de fermeture}\label{ssec2.4.2}

\medskip

Nous allons démontrer le théorème suivant qui identifie les points de part et d'autre de la relation de dualité de la définition \ref{defn2.4.1} entre les cribles d'une catégorie essentiellement petite et les sous-préfaisceaux sur une telle catégorie:

\begin{thm}\label{thm2.4.2}

Considérons la relation de dualité $R$ entre la classe $T$ des cribles $C \hookrightarrow y(X)$ d'une catégorie essentiellement petite ${\mathcal C}$ et la classe $S$ des plongements $Q \hookrightarrow P$ de préfaisceaux sur ${\mathcal C}$, avec la paire adjointe d'applications respectant les relations d'ordre qu'elle définit
$$
 \xymatrix{ ({\mathcal P} (T) , \subseteq) \ar@<2pt>[r]^{F_R} & ({\mathcal P} (S) , \supseteq). \ar@<2pt>[l]^{G_R} }
$$

Alors:

\begin{listeimarge}

\item Une sous-classe
$$
J \subseteq T = \{\mbox{cribles $C \hookrightarrow y(X)$ des objets $X$ de ${\mathcal C}$}\}
$$
est un point fixe de l'application
$$
G_R \circ F_R : {\mathcal P} (T) \longrightarrow {\mathcal P} (S) \longrightarrow {\mathcal P} (T)
$$
si et seulement si elle est une topologie.

\item Une sous-classe
$$
I \subseteq S = \{\mbox{plongements $Q \hookrightarrow P$ de $\widehat{\mathcal C}$}\}
$$
est un point fixe de l'application
$$
F_R \circ G_R : {\mathcal P} (S) \longrightarrow {\mathcal P} (T) \longrightarrow {\mathcal P} (S)
$$
si et seulement si elle est une ``propriété  de fermeture'' au sens qu'elle satisfait les conditions suivantes:

\medskip

$\left\{\begin{matrix}
(1) &\mbox{Les égalités $Q = P$ sont éléments de $I$.} \hfill \\
(2) &\mbox{Le changement de base pour tout morphisme $P' \to P$ de $\widehat{\mathcal C}$ transforme les sous-préfaisceaux} \hfill \\
&\mbox{$Q \hookrightarrow P$ qui sont éléments de $I$ en des sous-préfaisceaux $Q \times_P P' \hookrightarrow P'$ qui sont éléments de $I$.} \hfill \\
(3) &\mbox{Pour toute famille de sous-préfaisceaux} \hfill \\
&Q_k \xymatrix{\ar@{^{(}->}[r] & } P \, , \ k \in K \, , \ \mbox{qui sont éléments de $I$,} \\
&\mbox{leur intersection} \hfill \\
&\underset{k \in K}{\bigcap} \ Q_k \xymatrix{\ar@{^{(}->}[r] & } P \\
&\mbox{est encore élément de $I$.} \hfill \\
(4) &\mbox{Si, pour tout sous-préfaisceau $Q \hookrightarrow P$, on note $\overline Q \hookrightarrow P$ le plus petit sous-préfaisceau} \hfill \\
&\mbox{qui contient $Q$, on a pour tout morphisme de changement de base $P' \to P$} \hfill \\
&\overline{P' \times_P Q} = P' \times_P \overline Q \, .
\end{matrix}\right.$
\end{listeimarge}
\end{thm}

\begin{remarks}
	\begin{listeisansmarge}
\item Si $J$ est une topologie, ou plus généralement une classe de cribles de ${\mathcal C}$, et
$$
I = F_R(J)
$$
est la ``propriété de fermeture'' correspondante, les sous-préfaisceaux
$$
Q \xymatrix{\ar@{^{(}->}[r] & } P
$$
qui sont éléments de $I$ sont dits $J$-fermés.

De plus, pour tout sous-préfaisceau
$$
Q \xymatrix{\ar@{^{(}->}[r] & } P \, ,
$$
le plus petit sous-préfaisceau $J$-fermé qui le contient, noté
$$
\overline Q \xymatrix{\ar@{^{(}->}[r] & } P \, ,
$$
est appelé sa $J$-fermeture.

\item Dans le cas des préfaisceaux de ${\mathcal C}$ de la forme
$$
P = y(X) \quad \mbox{pour les objets $X$ de ${\mathcal C}$,}
$$
ses sous-préfaisceaux sont des cribles.

On retrouve les notions de cribles $J$-fermés et de $J$-fermeture d'un crible qui avaient été rappelées dans le lemme \ref{lem1.3.5}.

\item Généralisant la caractérisation de ce lemme, si $J$ est une topologie, la $J$-fermeture d'un sous-préfaisceau
$$
Q \xymatrix{\ar@{^{(}->}[r] & } P
$$
est le sous-préfaisceau
$$
\overline Q \xymatrix{\ar@{^{(}->}[r] & } P
$$
qui associe à tout objet $X$ de ${\mathcal C}$ le sous-ensemble
$$
\overline Q(X) \subseteq P(X)
$$
des sections $p \in P(X)$ qui sont $J$-localement dans $Q$ au sens qu'existe une famille $J$-couvrante de morphismes de but $X$
$$
x_k : X_k \longrightarrow X \, , \quad k \in K \, ,
$$
qui transforment $p$ en des éléments des sous-ensembles
$$
Q(X_k) \subseteq P(X_k) \, .
$$
\end{listeisansmarge}
\end{remarks}

\begin{demo}

Elle sera donnée dans les trois sous-paragraphes suivants c), d) et e). 
\end{demo}

Ce théorème permet d'appliquer aussi à la dualité des cribles et des sous-préfaisceaux les conclusions générales du corollaire \ref{cor2.1.4}.

On obtient ainsi les deux corollaires suivants:

\begin{cor}\label{cor2.4.3}

Considérons la paire adjointe
$$
 \xymatrix{ ({\mathcal P} (T) , \subseteq) \ar@<2pt>[r]^{F_R} & ({\mathcal P} (S) , \supseteq) \ar@<2pt>[l]^{G_R} }
$$
définie par la relation de dualité $R$ entre la classe $T$ des cribles d'une catégorie essentiellement petite ${\mathcal C}$ et la classe $S$ des plongements $Q \hookrightarrow P$ de préfaisceaux sur ${\mathcal C}$.

Alors les deux applications $F_R$ et $G_R$ induisent deux bijections inverses l'une de l'autre
$$
\xymatrix{ \left\{ {{\mbox{topologies} \atop \mbox{$J$}} \atop \mbox{de ${\mathcal C}$}} \right\} \ar@<2pt>[rr]^-{\sim} && \left\{ {{\mbox{propriétés de fermeture} \atop \mbox{des sous-préfaisceaux}} \atop \mbox{$Q \hookrightarrow P$}} \right\}. \ar@<2pt>[ll]^-{\sim} }
$$
\end{cor}

\begin{remarks}
\begin{listeisansmarge}
\item Ici, le mot ``topologie'' a toujours le sens d'une sous-classe $J$ de la classe $T$ des cribles de ${\mathcal C}$ qui satisfait les trois axiomes (1), (2), (3) du théorème \ref{thm2.2.2} (ii).

\item De même, l'expression ``propriété de fermeture'' a le sens d'une sous-classe $I$ de la classe $S$ des sous-préfaisceaux $Q \hookrightarrow P$ sur ${\mathcal C}$ qui satisfait les quatre axiomes (1), (2), (3), (4) du théorème \ref{thm2.4.2}.

Une telle ``propriété de fermeture'' détermine une ``opération de fermeture''
$$
(Q \hookrightarrow P) \longmapsto (\overline Q \hookrightarrow P)
$$
qui consiste à associer à tout sous-préfaisceau $Q$ d'un préfaisceau $P$ le plus petit sous-préfaisceau fermé $\overline Q$ de $P$ qui contient $Q$.

Réciproquement, elle est déterminée par cette opération puisqu'un sous-préfaisceau
$$
Q \xymatrix{\ar@{^{(}->}[r] & } P
$$
est fermé si et seulement si
$$
\overline Q = Q \, .
$$
\end{listeisansmarge}
\end{remarks}

\begin{demo}

C'est un cas particulier du corollaire \ref{cor2.1.4} (i) compte tenu de l'identification des points fixes réalisée dans le théorème \ref{thm2.4.2}. 
\end{demo}

\begin{cor}\label{cor2.4.4}

Considérons toujours la paire adjointe
$$
 \xymatrix{ ({\mathcal P} (T) , \subseteq) \ar@<2pt>[r]^{F_R} & ({\mathcal P} (S) , \supseteq) \ar@<2pt>[l]^{G_R} }
$$
définie par la relation de dualité $R$ entre les classe $T$ et $S$ des cribles et des sous-préfaisceaux sur une catégorie essentiellement petite ${\mathcal C}$.

Alors:

\begin{listeimarge}

\item Pour toute collection $J$ de cribles sur ${\mathcal C}$, l'image
$$
\overline J = G_R \circ F_R (J)
$$
est la topologie de ${\mathcal C}$ engendrée par $J$, et la sous-classe
$$
F_R(J) = F_R (\overline J)
$$
est constituée des sous-préfaisceaux sur ${\mathcal C}$
$$
Q \xymatrix{\ar@{^{(}->}[r] & } P
$$
qui sont fermés pour la topologie $\overline J$, définissant l'opération de $\overline J$-fermeture des sous-préfaisceaux
$$
(Q \hookrightarrow P) \longmapsto (\overline Q \hookrightarrow P) \, .
$$
\item Pour toute collection $I$ de sous-préfaisceaux sur ${\mathcal C}$, l'image
$$
\overline I = F_R \circ G_R (I)
$$
est la plus petite propriété de fermeture qui contient les éléments de $I$, et la topologie correspondante
$$
G_R(I) = G_R(\overline I)
$$
est constituée des cribles
$$
C \xymatrix{\ar@{^{(}->}[r] & } y(X)
$$
dont la fermeture est le crible maximal $y(X)$.
\end{listeimarge}
\end{cor}

\begin{remark}
En particulier, pour toute famille de sous-topos
$$
{\mathcal E}_k = \widehat{\mathcal C}_{J_k} \xymatrix{\ar@{^{(}->}[r] & } \widehat{\mathcal C} \, , \quad k \in K \, ,
$$
définis par des topologies $J_k$ de ${\mathcal C}$, leur intersection
$$
\bigwedge_{k \in K} {\mathcal E}_k \xymatrix{\ar@{^{(}->}[r] & } \widehat{\mathcal C}
$$
est définie par la topologie $J$ engendrée par la réunion
$$
\bigcup_{k \in K} J_k \, .
$$
Il lui correspond la propriété de fermeture
$$
F_R \left( \bigcup_{k \in K} J_k \right) = \bigcap_{k \in K} F_R (J_k) \, .
$$
Autrement dit, un sous-préfaisceau
$$
Q \xymatrix{\ar@{^{(}->}[r] & } P
$$
est fermé pour la topologie $J$ engendrée par les $J_k$, $k \in K$, si et seulement si il est fermé pour chacune des topologies $J_k$.
\end{remark}

\begin{demo}

Encore une fois, les égalités
$$
F_R(J) = F_R \circ G_R \circ F_R (J) = F_R (\overline J)
$$
et
$$
G_R(I) = G_R \circ F_R \circ G_R (I) = G_R (\overline I)
$$
sont vérifiées du fait que, d'après le corollaire \ref{cor2.1.4} (ii), chaque image de $F_R$ ou $G_R$ est un point fixe.

Le reste des assertions (i) et (ii) est un cas particulier du corollaire \ref{cor2.1.4} (iii).
\end{demo}

\subsection{Réapparition de la notion de topologie}\label{ssec2.4.3}

\medskip

La première étape de la démonstration du théorème \ref{thm2.4.2} est la proposition suivante qui fait à nouveau apparaître les trois axiomes de définition de la notion de topologie, cette fois à partir de la relation de dualité entre cribles et sous-préfaisceaux:

\begin{prop}\label{prop2.4.5}

Considérons l'application

$$
\begin{matrix}
G_R : &({\mathcal P} (S) , \supseteq) &\longrightarrow &({\mathcal P} (T) , \subseteq) \hfill \\
&\hfill (I \subseteq S) &\longmapsto &G_R(I) = \{ C \hookrightarrow y(X) \mid (C \hookrightarrow (X), Q \hookrightarrow P) \in R \, , \ \forall (Q \hookrightarrow P) \in I \}
\end{matrix}
$$
définie par la relation de dualité $R$ entre la classe des plongements $Q \hookrightarrow P$ de préfaisceaux sur  ${\mathcal C}$ et la classe $T$ des cribles $C \hookrightarrow y(X)$ des objets $X$ de ${\mathcal C}$.

Alors, pour toute collection $I$ de plongements de préfaisceaux, la sous-classe
$$
J = G_R(I)
$$
est une topologie.
\end{prop}

\begin{demo}

Toute intersection de topologies est une topologie.

Par conséquent, il suffit de vérifier la proposition lorsque la classe $I$ est réduite à un unique plongement de préfaisceaux
$$
Q \xymatrix{\ar@{^{(}->}[r] & } P \, .
$$

Alors $J = G_R(I)$ est la classe des cribles d'objets $X$ de ${\mathcal C}$
$$
C \xymatrix{\ar@{^{(}->}[r] & } y(X)
$$
tels que, pour tout morphisme $x : X' \to X$ de ${\mathcal C}$ et pour tout $p \in P(X')$, on a
$$
p \in Q (X') \subseteq P(X')
$$
si tout morphisme $(X'' \xrightarrow{ \, x' \, } X') \in x^* C$ envoie $p$ dans le sous-ensemble
$$
Q(X'') \subseteq P(X'') \, .
$$

Si $C = y(X)$, on a $x^* C = C \times_{y(X)} y(X') = y(X')$ pour tout morphisme $X' \xrightarrow{ \, x \, } X$, et donc la condition (1) du théorème \ref{thm2.2.2} est vérifiée.

La condition (2) l'est également puisque la formation des images réciproques de cribles
$$
x^* : (C \hookrightarrow y(X)) \longmapsto (C \times_{y(X)} y(X') \hookrightarrow y(X'))
$$
est compatible avec la composition des morphismes de ${\mathcal C}$
$$
X'' \xrightarrow{ \ x' \ } X' \xrightarrow{ \ x \ } X
$$
au sens que $x'^* (x^* C) = (x \circ x')^* C$.

Pour la condition (3), considérons deux cribles d'un objet $X$
$$
C \xymatrix{\ar@{^{(}->}[r] & } y(X) \quad \mbox{et} \quad C' \xymatrix{\ar@{^{(}->}[r] & } y(X)
$$
tels que
$$
(C \xymatrix{\ar@{^{(}->}[r] & } y(X)) \in J
$$
et
$$
(x^* C' \xymatrix{\ar@{^{(}->}[r] & } y(X')) \in J \, , \ \forall \, (X' \xrightarrow{ \ x \ } X) \in C \, .
$$
Comme ces conditions sont respectées par changement de base, on est réduit à démontrer qu'un élément
$$
p \in P(X)
$$
est dans le sous-ensemble
$$
Q(X) \subseteq P(X)
$$
si tout morphisme $(X' \xrightarrow{ \ x \ } X) \in C'$ l'envoie dans le sous-ensemble
$$
Q(X') \subseteq P(X') \, .
$$
Or, sous cette hypothèse, le composé de tout morphisme
$$
(X' \xrightarrow{ \ x \ } X) \in C
$$
et de tout morphisme
$$
(X'' \xrightarrow{ \ x' \ } X) \in x^* C'
$$
envoie la section $p \in P(X)$ dans le sous-ensemble
$$
Q(X'') \subseteq P(X'') \, .
$$
Comme les cribles $(x^* C' \hookrightarrow y(X'))$, $(X' \xrightarrow{ \, x \, } X) \in C$, sont supposés éléments de $J$, on en déduit que chaque $(X' \xrightarrow{ \, x \, } X) \in C$ envoie $p$ dans le sous-ensemble
$$
Q(X') \subseteq P(X') \, .
$$
Comme de plus $C$ est supposé élément de $J$, on conclut comme voulu que $p$ est élément du sous-ensemble
$$
Q(X) \subseteq P(X) \, .
$$

Ainsi, $J$ vérifie bien la condition (3) de localité. C'est une topologie. 
\end{demo}

\subsection{Apparition de la notion de propriété de fermeture}\label{ssec2.4.4}

\medskip

La seconde étape de la démonstration du théorème \ref{thm2.4.2} est le corollaire suivant qui fait apparaître les quatre axiomes de définition de la notion de propriété de fermeture des sous-préfaisceaux:

\begin{cor}\label{cor2.4.6}

Considérons l'application
$$
F_R : ({\mathcal P} (T) , \subseteq) \longrightarrow ({\mathcal P} (S) , \supseteq)
$$
définie par la relation de dualité $R$ entre la classe $T$ des cribles d'une catégorie essentiellement petite ${\mathcal C}$ et la classe $S$ des plongements $Q \hookrightarrow P$ de préfaisceaux sur ${\mathcal C}$. Alors, pour toute collection $J$ de cribles, la sous-classe
$$
I = F_R (J) \subseteq S = \{\mbox{plongements de préfaisceaux $Q \hookrightarrow P$}\}
$$
est une propriété de fermeture au sens qu'elle satisfait les conditions suivantes:

\medskip

$\left\{\begin{matrix}
(1) &\mbox{Elle contient les égalités $Q = P$.} \hfill \\
(2) &\mbox{Elle est stable par changement de base.} \hfill \\
(3) &\mbox{Elle est stable par intersections de sous-préfaisceaux d'un préfaisceau.} \hfill \\
(4) &\mbox{Les applications qui associent à chaque sous-préfaisceau $Q$ d'un préfaisceau $P$ le plus petit} \hfill \\
&\mbox{élément $\overline Q \hookrightarrow P$ de $I$ qui le contient sont compatibles avec les changements de base.} \hfill
\end{matrix}\right.$
\end{cor}

\begin{demo}

L'image $I = F_R(J)$ est un point fixe de l'application composée
$$
F_R \circ G_R
$$
de $F_R$ avec son adjointe à droite $G_R$.

Autrement dit, on a aussi
$$
I = F_R(\overline J)
$$
si l'on note $\overline J = G_R \circ F_R (J)$.

Or, on sait d'après la proposition II.4.5 que
$$
\overline J = G_R \circ F_R(J)
$$
est une topologie.

Donc il n'y a pas de restriction à supposer que $J$ et une topologie.

Dans ce cas, $I = F_R(J)$ est constituée des plongements de préfaisceaux
$$
Q \xymatrix{\ar@{^{(}->}[r] & } P
$$
tels qu'une section $p \in P(X)$ est dans le sous-ensemble
$$
Q(X) \subseteq P(X)
$$
s'il existe une famille $J$-couvrante de morphismes de but $X$
$$
x_k : X_k \longrightarrow X \, , \quad k \in K \, ,
$$
qui envoient $p$ dans les sous-ensembles
$$
Q(X_k) \subseteq P(X_k) \, .
$$
Il est évident que la classe $I$ satisfait les conditions (1) et (3). De plus, notant
$$
j^* : \widehat{\mathcal C} \longrightarrow \widehat{\mathcal C}_J
$$
le foncteur de faisceautisation adjoint à gauche du foncteur de plongement
$$
j_* : \widehat{\mathcal C}_J \xymatrix{\ar@{^{(}->}[r] & } \widehat{\mathcal C} \, ,
$$
on note qu'un plongement de préfaisceaux
$$
Q \xymatrix{\ar@{^{(}->}[r] & } P
$$
est élément de $I$ si et seulement si le carré commutatif
$$
\xymatrix{
Q \ \ar[d] \ar@{^{(}->}[r] &P \ar[d] \\
j_* \circ j^* Q \ \ar@{^{(}->}[r] &j_* \circ j^* P
}
$$
est cartésien.

Comme les foncteurs $j_*$ et $j^*$ respectent les produits fibrés, on en déduit que la classe $I$ est respectée par les applications de changement de base
$$
(Q \hookrightarrow P) \longmapsto (Q \times_P P' \hookrightarrow P')
$$
par les morphismes de préfaisceaux $P' \to P$.

C'est la condition (2).

Enfin, la $J$-fermeture de tout sous-préfaisceau
$$
Q \xymatrix{\ar@{^{(}->}[r] & } P
$$
s'écrit
$$
\overline Q = j_* \circ j^* Q \times_{j_* \circ j^* P} P \xymatrix{\ar@{^{(}->}[r] & } P \, .
$$
Pour tout morphisme de changement de base $P' \to P$, on obtient
\begin{eqnarray}
\overline{Q \times_P P'} &= &j_* \circ j^* (Q \times_P P') \times_{j_* \circ j^* P'} P' \nonumber \\
&= &j_* \circ j^* Q \times_{j_* \circ j^* P} P' \nonumber \\
&= &(j_* \circ j^* Q \times_{j_* \circ j^* P} P) \times_P P' \nonumber \\
&= &\overline Q \times_P P' \, . \nonumber
\end{eqnarray}

C'est la condition (4).

Ainsi, $I$ est bien une ``propriété de fermeture''. 
\end{demo}

\subsection{Vérification de ce que les topologies et les propriétés de fermeture sont des points fixes}\label{ssec2.4.5}
%\medskip
\nopagebreak
On a déjà montré que les sous-classes de cribles qui sont des points fixes de la relation de dualité entre cribles et sous-préfaisceaux sont nécessairement des topologies.
\nopagebreak
Réciproquement, les topologies sont toujours des points fixes:

\nopagebreak
\begin{lem}\label{lem2.4.7}
\nopagebreak
Considérons  la paire adjointe d'applications respectant les relations d'ordre
$$
 \xymatrix{ ({\mathcal P} (T) , \subseteq) \ar@<2pt>[r]^{F_R} & ({\mathcal P} (S) , \supseteq) \ar@<2pt>[l]^{G_R} }
$$
qui est définie par la relation de dualité $R$ entre la classe $T$ des cribles d'une catégorie essentiellement petite ${\mathcal C}$ et la classe $S$ des plongements de préfaisceaux sur ${\mathcal C}$.
\nopagebreak
Alors toute sous-classe $J \subseteq T$ de cribles qui est une topologie est un point fixe de la relation de dualité $R$ au sens que
$$
G_R \circ F_R (J) = J \, .
$$
\end{lem}
\nopagebreak
\begin{demo}
\nopagebreak
Considérons une topologie $J$ et la propriété de fermeture
$$
I = F_R (J)
$$
qui lui est associée, avec l'opération de fermeture des sous-préfaisceaux qu'elle définit
$$
(Q \hookrightarrow P) \longmapsto (\overline Q \hookrightarrow P) \, .
$$
Cette opération est caractérisée par la formule que, en tout objet $X$ de ${\mathcal C}$, le sous-ensemble
$$
\overline Q (X) \subseteq P(X)
$$
est celui des éléments $p \in P(X)$ tels qu'existe une famille $J$-couvrante de morphismes de but $X$
$$
x_k : X_k \longrightarrow X \, , \quad k \in K \, ,
$$
qui envoient $p$ dans les sous-ensembles
$$
Q(X_k) \subseteq P(X_k) \, .
$$

Un crible d'un objet $X$ de ${\mathcal C}$
$$
C \xymatrix{\ar@{^{(}->}[r] & }  y(X)
$$

est élément de la topologie $J$ si et seulement si
$$
\overline C = y(X)
$$

c'est-à-dire si et seulement si tout carré commutatif
$$
\xymatrix{
C \ \ar[d] \ar@{^{(}->}[r] &y(X) \ar[d] \\
Q \ \ar@{^{(}->}[r] &P
}
$$

se relève en un morphisme
$$
y(X) \longrightarrow \overline Q \, .
$$

Cela revient à demander que le crible $C \hookrightarrow y(X)$ soit élément de
$$
G_R(I) = G_R \circ F_R(J) \, .
$$

Ainsi, la topologie $J$ est bien un point fixe. 
\end{demo}

\newpage
De même, on a déjà montré que les sous-classes de plongements de préfaisceaux qui sont des points fixes de la relation de dualité entre sous-préfaisceaux et cribles sont nécessairement des ``propriétés de fermeture'' au sens du théorème \ref{thm2.4.2} (ii).
\nopagebreak
Réciproquement, les ``propriétés de fermeture'' sont toujours des points fixes:
\nopagebreak
\begin{prop}\label{prop2.4.8}
\nopagebreak
Considérons la paire adjointe d'applications respectant les relations d'ordre
$$
 \xymatrix{ ({\mathcal P} (T) , \subseteq) \ar@<2pt>[r]^{F_R} & ({\mathcal P} (S) , \supseteq) \ar@<2pt>[l]^{G_R} }
$$
qui est définie par la relation de dualité $R$ entre la classe $S$ des plongements $Q \hookrightarrow P$ de préfaisceaux sur une catégorie essentiellement petite ${\mathcal C}$ et la classe $T$ des cribles sur ${\mathcal C}$.
\nopagebreak
Alors toute sous-classe
$$
I \subseteq S = \{\mbox{sous-préfaisceaux $Q \hookrightarrow P$}\}
$$
qui est une ``propriété de fermeture'' est un point fixe de la relation de dualité au sens que
$$
I = F_R \circ G_R(I) \, .
$$
\end{prop}
\nopagebreak
\begin{demo}
\nopagebreak
Considérons une propriété de fermeture $I$ et l'opération de fermeture qu'elle définit
$$
(Q \hookrightarrow P) \longmapsto (\overline Q \hookrightarrow P) \, .
$$
\nopagebreak
Considérons d'autre part la topologie
$$
J = G_R(I)
$$
\nopagebreak
qui lui est associée.
\nopagebreak
Alors un crible d'un objet $X$ de ${\mathcal C}$
$$
C \xymatrix{\ar@{^{(}->}[r] & } y(X)
$$
\nopagebreak
est élément de $J$ si et seulement si, pour tout morphisme de ${\mathcal C}$
$$
x : X' \longrightarrow X
$$

et tout plongement de préfaisceaux
$$
Q \xymatrix{\ar@{^{(}->}[r] & } P \, ,
$$
tout carré commutatif
$$
\xymatrix{
x^* C \ \ar[d] \ar@{^{(}->}[r] &y(X') \ar[d] \\
Q \ \ar@{^{(}->}[r] &P
}
$$

se relève en un morphisme
$$
y(X') \longrightarrow \overline Q \, .
$$

Cela signifie qu'un tel crible
$$
C \xymatrix{\ar@{^{(}->}[r] & } y(X)
$$
et élément de la topologie 
$$
J = G_R(I)
$$

si et seulement si
$$
\overline C = y(X) \, .
$$
\newpage
Alors un plongement de préfaisceaux
$$
Q \xymatrix{\ar@{^{(}->}[r] & } P 
$$
est élément de
$$
F_R(J) = F_R \circ G_R(I)
$$
si et seulement si, pour tout crible $J$-couvrant d'un objet $X$ de ${\mathcal C}$
$$
C \xymatrix{\ar@{^{(}->}[r] & } y(X) \, ,
$$
tout carré commutatif 
$$
\xymatrix{
C \ \ar[d] \ar@{^{(}->}[r] &y(X) \ar[d] \\
Q \ \ar@{^{(}->}[r] &P
}
$$
se relève en un morphisme
$$
y(X) \longrightarrow Q \, .
$$
Cela revient à demander que
$$
Q = \overline Q \, ,
$$
autrement dit que
$$
Q \xymatrix{\ar@{^{(}->}[r] & } P
$$
soit élément de $I$.

Ainsi, on a
$$
I = F_R(J) = F_R \circ G_R(I) \, .
$$
La propriété de fermeture $I$ est bien un point fixe. 
\end{demo}

Le lemme \ref{lem2.4.7} et la proposition \ref{prop2.4.8} terminent la démonstration du théorème \ref{thm2.4.2}.

	\chapter{Engendrement de topologies et démontrabilité}\label{chap3}

\section{Une formule fermée d'engendrement des topologies}\label{sec3.1}

\subsection{Cribles et précribles couvrants}\label{ssec3.1.1}

Dès la définition \ref{defn1.1.1}, nous avons rappelé la première définition de la notion de topologie sur une catégorie ${\mathcal C}$, en termes de familles de morphismes de ${\mathcal C}$ de même but $X$
$$
\left(X_i \xrightarrow{ \ x_i \ } \ X\right)_{i \in I} \, .
$$
Cette définition a un sens sur n'importe quelle catégorie ${\mathcal C}$.

Dans le cas d'une catégorie ${\mathcal C}$ essentiellement petite, la dualité des cribles des objets $X$ de ${\mathcal C}$, c'est-à-dire des sous-préfaisceaux
$$
C \xymatrix{\ar@{^{(}->}[r] &} y(X) \qquad \mbox{dans} \quad \widehat{\mathcal C} \, ,
$$
et des préfaisceaux sur ${\mathcal C}$, a fait réapparaître la notion de topologie, cette fois en termes de cribles dans la proposition \ref{prop2.2.5}.

Ces deux définitions sont équivalentes du fait du lemme suivant:

\begin{lem}\label{lem3.1.1}

Soit ${\mathcal C}$ une catégorie essentiellement petite.

Disons que deux précribles d'un objet $X$ de ${\mathcal C}$, c'est-à-dire deux familles de morphismes de ${\mathcal C}$ de but $X$
$$
\left(X_i \xrightarrow{ \ x_i \ } \ X \right)_{i \in I} \qquad \mbox{et} \qquad \left(X'_j \xrightarrow{ \ x'_j \ } \ X \right)_{j \in I'} \, ,
$$
sont équivalents si

\medskip

$\left\{\begin{matrix}
\bullet &\mbox{pour tout indice $i \in I$ \quad [resp. $j \in I'$],} \hfill \\
&\mbox{il existe un indice $j \in I'$ \quad [resp. $i \in I$]} \hfill \\
&\mbox{et une factorisation en un triangle commutatif de ${\mathcal C}$} \hfill \\
&\xymatrix{
X_i \ar[rr] \ar[rd]_{x_i} &&X'_j \ar[ld]^{x'_j} \\
&X
} \qquad \xymatrix{
X'_j \ar[rr] \ar[rd]_-{\mbox{[resp.} \qquad x'_j} &&X_i \ar[ld]^-{x_i \quad \mbox{].}} \\
&X
}
\end{matrix}\right.$

\medskip

\noindent Alors se donner un crible d'un objet $X$ de ${\mathcal C}$ équivaut à se donner une classe d'équivalence de précribles de ${\mathcal C}$.
\end{lem}

\begin{remark}
Dans la pratique, un crible d'un objet $X$ d'une catégorie ${\mathcal C}$ est presque toujours donné sous la forme d'un précrible qui l'engendre. Autrement dit, on peut considérer qu'un crible est une classe d'équivalence de précribles, presque toujours représenté dans la pratique par un élément de la classe qu'il constitue.
\end{remark}

\begin{demo}
En effet, tout précrible d'un objet $X$
$$
\left(X_i \xrightarrow{ \ x_i \ } \ X \right)_{i \in I}
$$
engendre le crible 
$$
C \xymatrix{\ar@{^{(}->}[r] & } y(X)
$$
dont la composante en tout objet $X'$ de ${\mathcal C}$ est le sous-ensemble de 
$$
{\rm Hom} (X',X) = y(X)(X')
$$
constitué des morphismes
$$
X' \xrightarrow{ \ x \ } X
$$
tels qu'existent un indice $i \in I$ et une factorisation dans ${\mathcal C}$:
$$
\xymatrix{
X' \ar[rr] \ar[rd]_{x} &&X_i \ar[ld]^{x_i} \\
&X
}
$$

La conclusion résulte de ce que

\medskip

$\left\{\begin{matrix}
\bullet &\mbox{deux précribles d'un objet $X$ engendrent le même crible si et seulement si ils sont équivalents,} \hfill \\
\bullet &\mbox{tout crible d'un objet $X$ est engendré par des précribles de $X$.} \hfill
\end{matrix}\right.$
{\color{white}Grothendieck}
\end{demo}

Sur une catégorie essentiellement petite ${\mathcal C}$, n'importe quelle collection $J$ de cribles d'objets $X$ de ${\mathcal C}$
$$
C \xymatrix{\ar@{^{(}->}[r] & } y(X)
$$
peut être vue comme une ``notion de crible couvrant'': on dit alors qu'un crible est $J$-couvrant s'il est élément de la collection $J$.

De même, n'importe quelle collection de précribles d'objets $X$ de ${\mathcal C}$
$$
\left(X_i \xrightarrow{ \ x_i \ } \ X \right)_{i \in I}
$$
peut être vue comme une ``notion de précrible couvrant'': on dit alors qu'un précrible est couvrant s'il est élément de cette collection.

\newpage
Il est naturel de poser:

\begin{defn}\label{defn3.1.2}

Soit ${\mathcal C}$ une catégorie essentiellement petite.

\begin{listeimarge}

\item On dira qu'une notion de crible couvrant est ``stable'' si

\medskip

$\left\{\begin{matrix}
\bullet &\mbox{tout crible d'un objet $X$ de ${\mathcal C}$ qui contient un crible couvrant (ou qui est le crible maximal)} \hfill \\
&\mbox{est lui-même couvrant,} \hfill \\
\bullet &\mbox{pour tout morphisme $x : X' \to X$ de ${\mathcal C}$, l'image réciproque par $x$} \hfill \\
&\mbox{de tout crible couvrant $C \hookrightarrow y(X)$ de $X$} \hfill \\
&x^* C = C \times_{y(X)} y(X') \xymatrix{\ar@{^{(}->}[r] & } y(X') \\
&\mbox{est un crible couvrant de $X'$.} \hfill
\end{matrix}\right.$

\medskip

\item On dira qu'une notion de précrible couvrant est ``stable '' si

\medskip

$\left\{\begin{matrix}
\bullet &\mbox{pour tout morphisme $x : X' \to X$ de ${\mathcal C}$ et pour tout précrible couvrant de $X$} \hfill \\
&\left(X_i \xrightarrow{ \ x_i \ } \ X \right)_{i \in I} \, , \\
&\mbox{il existe un précrible couvrant de $X'$,} \hfill \\
&\left(X'_j \xrightarrow{ \ x'_j \ } \ X' \right)_{j \in I'} \\
&\mbox{tel que, pour tout indice $j \in I'$, il existe un indice $i \in I$ et un carré commutatif de ${\mathcal C}$ de la forme:} \hfill \\
&\xymatrix{
X'_j \ar[d]_{x'_j} \ar[r] &X_i \ar[d]^{x_i} \\
X' \ar[r]^x &X
}
\end{matrix}\right.$
\end{listeimarge}
\end{defn}

\begin{remarks}
\begin{listeisansmarge}
\item Toute notion de précrible couvrant qui est stable induit une notion de crible couvrant qui est stable: elle est définie en disant qu'un crible d'un objet $X$ est couvrant s'il contient au moins un précrible couvrant (ou s'il est le crible maximal).

%\smallskip

\item En sens inverse, toute notion de crible couvrant qui est stable induit une notion stable de précrible couvrant définie en disant qu'un précrible est couvrant si et seulement si le crible qu'il engendre est couvrant.
\end{listeisansmarge}
\end{remarks}

Cette notion de précrible couvrant associée est non seulement stable mais aussi compatible avec la relation d'équivalence entre précribles. De plus, tout précrible qui contient un précrible couvrant est encore couvrant et tous les précribles qui contiennent un morphisme identité $X \xrightarrow{ \, {\rm id}_X \, } X$ sont couvrants.

\medskip

Il est facile d'engendrer une notion de crible couvrant ou de précrible couvrant qui soit stable à partir d'une famille de cribles ou de précribles:

\begin{lem}\label{lem3.1.3}

Soit ${\mathcal C}$ une catégorie essentiellement petite.

\begin{listeimarge}

\item Considérons une famille de cribles d'objets de ${\mathcal C}$
$$
( C_i \xymatrix{\ar@{^{(}->}[r] & } y(X_i))_{i \in I} \, .
$$
Alors la plus petite notion de crible couvrant qui soit stable et contienne cette famille est constituée des cribles des objets $X$ de ${\mathcal C}$
$$
C \xymatrix{\ar@{^{(}->}[r] & } y(X)
$$
tels qu'il existe un indice $i \in I$ et un morphisme
$$
X \xrightarrow{ \ x \ } X_i
$$
vérifiant
$$
C \supseteq x^* C_i
$$
ou que $C$ soit le crible maximal.

%\smallskip

\item Supposons que pour tout morphisme $x : X' \to X$ de ${\mathcal C}$ et tout précrible de son but $X$
$$
P = \left(X_i \xrightarrow{ \ x_i \ } \ X \right)_{i \in I}
$$
engendrant un crible $C$, on puisse construire un crible de $X'$
$$
x^* P = P' = \left(X'_{i'} \xrightarrow{ \ x'_{i'} \ } \ X' \right)_{i' \in I'}
$$
dont le crible engendré soit l'image réciproque $x^* C$.

Alors, toute famille de précribles d'objets $X_k$ de ${\mathcal C}$
$$
P_k = \left(X_{k,i} \xrightarrow{ \ x_{k,i} \ } \ X_k \right)_{i \in I_k} \, , \quad k \in K \, ,
$$
définit une notion de précrible couvrant qui est stable, en disant qu'un précrible d'un objet $X$ de ${\mathcal C}$ est couvrant s'il est de la forme
$$
x^* P_k
$$
pour un indice $k \in K$ et un morphisme $X \xrightarrow{ \, x \, } X_k$.
\end{listeimarge}
\end{lem}

\begin{remarksqed}
\begin{listeisansmarge}
\item L'hypothèse de (ii) est satisfaite en particulier si la catégorie ${\mathcal C}$ possède des produit fibrés.

En effet, pour tout morphisme $x : X' \to X$ de ${\mathcal C}$ et tout précrible de son but $X$
$$
P = \left(X_i \xrightarrow{ \ x_i \ } \ X \right)_{i \in I} \, ,
$$
l'image réciproque $x^* C$ par $x$ du crible $C$ engendré par ce précrible $P$ est le crible engendré par le précrible de $X'$
$$
x^* P = \left(X_i \times_X X' \longrightarrow X' \right)_{i \in I} \, .
$$
\item Dans la situation de (ii), la notion stable de crible couvrant associée à la notion stable de précrible couvrant ainsi définie, est celle engendrée au sens de (i) par les cribles engendrés par les précribles
$$
x^* P_k
$$
construits à partir des précribles $P_k$ des objets $X_k$, $k \in K$, par action de tous les morphismes $x : X \to X_k$.
\end{listeisansmarge} 
\end{remarksqed}

\subsection{Notions stables de précribles couvrants et propriétés de fermeture associées}\label{ssec3.1.2}

Rappelons que nous avons introduit dans la définition II.4.1 une relation
$$
R \xymatrix{\ar@{^{(}->}[r] & } T \times S
$$
entre la classe $T$ des cribles d'objets $X$
$$
C \xymatrix{\ar@{^{(}->}[r] & } y(X)
$$
d'une catégorie essentiellement petite ${\mathcal C}$, et la classe $S$ des plongements de préfaisceaux sur ${\mathcal C}$
$$
Q \xymatrix{\ar@{^{(}->}[r] & } P \, .
$$
Cette relation $R$ est constituée des paires d'éléments
$$
(C \xymatrix{\ar@{^{(}->}[r] & } y(X) , Q \xymatrix{\ar@{^{(}->}[r] & } P)
$$
telles que, pour tout morphisme $X' \xrightarrow{ \, x \, } X$ de ${\mathcal C}$ de but $X$ et pour tout élément $p \in P(X')$, on a
$$
p \in Q(X') \subseteq P(X')
$$
si tout morphisme $(X'' \xrightarrow{ \, x' \, } X') \in x^* C$ envoie $p$ dans le sous-ensemble
$$
Q(X'') \subseteq P(X'') \, .
$$
D'après le lemme \ref{lem2.1.3}, $R$ définit comme toute relation une paire adjointe d'applications respectant les relation d'ordre
$$
\xymatrix{
({\mathcal P} (T) , \subseteq) \ar@<2pt>[r]^{F_R} & ({\mathcal P} (S) , \supseteq) \ar@<2pt>[l]^{G_R}
}
$$
entre la classe ordonnée ${\mathcal P} (T)$ des sous-classes $J$ de $T$ et la classe ordonnée ${\mathcal P} (S)$ des sous-classes $I$ de $S$.

On a montré que, pour toute sous-classe de plongements de préfaisceaux
$$
I \subseteq S \, ,
$$
son image
$$
G_R (I) = \left\{ (C \hookrightarrow y(X)) \mid (C \hookrightarrow y(X) , Q \hookrightarrow P) \in R \, , \ \forall \, (Q \hookrightarrow P) \in I \right\}
$$
est une topologie, et que, pour toute sous-classe de cribles
$$
J \subseteq T \, ,
$$
son image
$$
F_R (J) = \left\{ (Q \hookrightarrow P) \mid (C \hookrightarrow y(X) , Q \hookrightarrow P) \in R \, , \ \forall \, (C \hookrightarrow y(X)) \in J \right\}
$$
est une ``propriété de fermeture'' au sens du théorème \ref{thm2.4.2} (ii).

Dans le cas où $J \subseteq T$ est une notion de crible couvrant qui est stable, la formule de définition de la propriété de fermeture $F_R(J)$ à partir de $J$ se simplifie:

\begin{lem}\label{lem3.1.4}

Considérons la relation  de dualité $R$ entre la classe $T$ des cribles $C \hookrightarrow y(X)$ d'une catégorie essentiellement petite ${\mathcal C}$ et la classe $S$ des plongements $Q \hookrightarrow P$ de préfaisceaux sur ${\mathcal C}$, avec la paire adjointe d'applications respectant les relations d'ordre qu'elle définit
$$
\xymatrix{
({\mathcal P} (T) , \subseteq) \ar@<2pt>[r]^{F_R} & ({\mathcal P} (S) , \supseteq). \ar@<2pt>[l]^{G_R}
}
$$

Alors:
\begin{listeimarge}
\item Si $J$ est une notion de crible couvrant qui est stable, un plongement de préfaisceaux sur ${\mathcal C}$
$$
Q \xymatrix{\ar@{^{(}->}[r] & } P
$$
est élément de la propriété de fermeture
$$
I = F_R (J)
$$
si et seulement si, pour tout crible $C$ d'un objet $X$ de ${\mathcal C}$ qui est élément de $J$, un élément
$$
p \in P(X)
$$
est dans le sous-ensemble
$$
Q(X) \subseteq P(X)
$$
dès que tout morphisme $( X' \xrightarrow{ \, x \, } X ) \in C$ envoie $p$ dans le sous-ensemble
$$
Q(X') \subseteq P(X') \, .
$$

\item Supposons de plus que la notion stable de crible couvrant $J$ soit définie au sens du lemme III.1.3 (ii) par une notion de précrible couvrant qui est stable.

Sous cette hypothèse, un plongement de préfaisceaux sur ${\mathcal C}$
$$
Q \xymatrix{\ar@{^{(}->}[r] & } P
$$
est élément de la propriété de fermeture
$$
I = F_R(J)
$$
si et seulement si, pour tout précrible couvrant d'un objet $X$ de ${\mathcal C}$
$$
\left( X_i \xrightarrow{ \ x_i \ } X \right)_{i \in I} \, ,
$$
un élément
$$
p \in P(X)
$$
est dans le sous-ensemble
$$
Q(X) \subseteq P(X)
$$
dès que tout morphisme $X_i \xrightarrow{ \, x_i \, } X$, $i \in I$, envoie $p$ dans le sous-ensemble
$$
Q(X_i) \subseteq P(X_i) \, .
$$
\end{listeimarge}
\end{lem}

\begin{demo}
\begin{listeisansmarge}
\item Par définition, si une notion de crible couvrant est stable, elle est respectée par les applications d'images réciproques
$$
C \longmapsto x^* C
$$
associées aux morphismes $x : X' \to X$ de ${\mathcal C}$.

C'est pourquoi la formule de définition de la propriété de fermeture
$$
F_R(J)
$$
se réduit à l'expression de (i) si $J$ est une notion de crible couvrant qui est stable.

%\smallskip

\item Si l'on part d'une notion de précrible couvrant qui est stable, la notion stable de crible couvrant qu'elle définit consiste en les cribles d'objets $X$ de ${\mathcal C}$
$$
C \xymatrix{\ar@{^{(}->}[r] & } y(X)
$$
qui sont maximaux ou contiennent un précrible couvrant
$$
( X_i \xrightarrow{ \ x_i \ } X )_{i \in I} \, .
$$
Un plongement de préfaisceaux de ${\mathcal C}$
$$
Q \xymatrix{\ar@{^{(}->}[r] & } P
$$
vérifie la propriété de (i) relativement aux cribles couvrants d'un objet $X$
$$
C \xymatrix{\ar@{^{(}->}[r] & } y(X)
$$
si et seulement si il satisfait la propriété de (ii) relativement aux précribles couvrants de $X$
$$
\left( X_i \xrightarrow{ \ x_i \ } X \right)_{i \in I} 
$$
ou, ce qui revient au même, relativement aux cribles qu'ils engendrent.
\end{listeisansmarge} 
\end{demo}

\subsection{Propriété de fermeture et formule d'engendrement des topologies}\label{ssec3.1.3}

\medskip

On déduit du corollaire \ref{cor2.4.4} et du lemme \ref{lem3.1.4}:

\begin{thm}\label{thm3.1.5}

Soit une catégorie essentiellement petite ${\mathcal C}$.

Soit une ``notion de crible couvrant'' $J$ qui est stable au sens que:

\smallskip

$\left\{\begin{matrix}
\bullet &\mbox{tout crible qui contient un crible couvrant (ou est maximal) est couvrant,} \hfill \\
\bullet &\mbox{toute image réciproque $x^* C$ par un morphisme $x : X' \to X$ de ${\mathcal C}$ d'un crible couvrant $C$} \hfill \\
&\mbox{de son but $X$ est un crible couvrant de $X'$.} \hfill
\end{matrix}\right.$

\smallskip

Alors un crible d'un objet $X$ de ${\mathcal C}$
$$
C \xymatrix{\ar@{^{(}->}[r] & } y(X)
$$
est élément de la topologie $\overline J$ engendrée par $J$ si et seulement si le crible maximal de $X$ est l'unique crible
$$
C' \xymatrix{\ar@{^{(}->}[r] & } y(X)
$$
qui satisfait les deux conditions suivantes:

\smallskip

$\left\{\begin{matrix}
\bullet &\mbox{Il contient $C$.} \hfill \\
\bullet &\mbox{Il est $J$-fermé au sens qu'un morphisme} \hfill \\
&x : X' \longrightarrow X \\
&\mbox{est élément de $C'$ dès que le crible de $X'$} \hfill \\
&\{ X'' \xrightarrow{ \ x' \ } X' \mid (x \circ x' : X'' \to X) \in C' \} \\
&\mbox{est $J$-couvrant.} \hfill
\end{matrix}\right.$
\end{thm}

\begin{remark}

Cet énoncé rend plus explicite et raffine la proposition 4.1.1 du livre \cite{TST}.

\end{remark}

\begin{demo}
Un crible
$$
C \xymatrix{\ar@{^{(}->}[r] & } y(X)
$$
est couvrant pour la topologie $\overline J$ si et seulement si sa $\overline J$-fermeture est le crible maximal, autrement dit si tout crible $\overline J$-fermé qui le contient est le crible maximal.

Or, d'après le corollaire \ref{cor2.4.4}, un crible est $\overline J$-fermé si et seulement si il est $J$-fermé.

On conclut d'après le lemme \ref{lem3.1.4} (i). 
\end{demo}

Par combinaison avec le lemme \ref{lem3.1.3} (i), on obtient:

\begin{cor}\label{cor3.1.6}

Considérons une famille de cribles
$$
C_k \xymatrix{\ar@{^{(}->}[r] & } y(X_k) \, , \quad k \in K \, ,
$$
d'objets $X_k$ d'une catégorie essentiellement petite ${\mathcal C}$.

Alors un crible d'un objet $X$ de ${\mathcal C}$
$$
C \xymatrix{\ar@{^{(}->}[r] & } y(X)
$$
est élément de la topologie $J$ engendrée par la famille $(C_k \hookrightarrow y(X_k))_{k \in K}$ si et seulement si le crible maximal de $X$ est l'unique crible
$$
C' \xymatrix{\ar@{^{(}->}[r] & } y(X)
$$
qui satisfait les deux conditions suivantes:

\smallskip

$\left\{\begin{matrix}
\bullet &\mbox{Il contient $C$.} \hfill \\
\bullet &\mbox{Il est fermé au sens qu'un morphisme} \hfill \\
&x : X' \longrightarrow X \\
&\mbox{est élément de $C'$ dès qu'existent un indice $k \in K$ et un morphisme $x_k :X' \to X_k$ tels que} \hfill \\
&(x \circ x' : X'' \to X) \in C' \, , \quad  \forall \, (X'' \xrightarrow{ \ x' \ } X') \in x_k^* \, C_k \, . 
\end{matrix}\right.$
\end{cor}

\begin{demo}

C'est une conséquence immédiate du lemme \ref{lem3.1.3} (i) appliqué au théorème \ref{thm3.1.5}. 
\end{demo}

Par combinaison avec le lemme \ref{lem3.1.3} (ii), on obtient encore l'énoncé suivant qui est la forme la plus concrète du théorème \ref{thm3.1.5}:

\begin{cor}\label{cor3.1.7}

Soit une catégorie essentiellement petite ${\mathcal C}$ où, pour tout morphisme $x : X' \to X$ et tout précrible de son but $X$
$$
P = \left( X_i \xrightarrow{ \ x_i \ } X \right)_{i \in I}
$$
engendrant un crible $C$, on sait construire un précrible de $X'$
$$
x^* P = P' = \left( X'_{i'} \xrightarrow{ \ x'_{i'} \ } X' \right)_{i' \in I'}
$$
dont le crible engendré soit l'image réciproque $x^* C$.

Considérons une famille de précribles d'objets $X_k$ de ${\mathcal C}$
$$
P_k = \left( X_{k,i} \xrightarrow{ \ x_{k,i} \ } X_k \right)_{i \in I_k} \, , \quad k \in K \, ,
$$
et la plus petite topologie $J$ de ${\mathcal C}$ pour laquelle ces précribles sont couvrants.

Alors un précrible d'un objet $X$
$$
\left( X_i \xrightarrow{ \ x_i \ } X \right)_{i \in I}
$$
est couvrant pour la topologie $J$ si le crible maximal de $X$ est l'unique crible
$$
C \xymatrix{\ar@{^{(}->}[r] & } y(X)
$$
qui satisfait les deux conditions suivantes:

\smallskip

$\left\{\begin{matrix}
\bullet &\mbox{Il contient tous les morphismes $x_i : X_i \to X$, $i \in I$.} \hfill \\
\bullet &\mbox{Il est fermé au sens qu'un morphisme} \hfill \\
&x : X' \longrightarrow X \\
&\mbox{est élément de $C$ dès qu'existent un indice $k \in K$ et un morphisme $x_k : X' \to X_k$ tels que} \hfill \\
&(x \circ x_{i'} : X'_{i'} \longrightarrow X) \in C \\
&\mbox{pour tout élément $( X'_{i'} \xrightarrow{ \ x_{i'} \ } X' )$ du précrible} \hfill \\
&x^*_k \, P_k = \left( X'_{i'} \xrightarrow{ \ x_{i'} \ } X' \right)_{i' \in I'} \, .
\end{matrix}\right.$
\end{cor}

\begin{remark}

Ce corollaire s'applique en particulier lorsque la catégorie ${\mathcal C}$ possède des produits fibrés.

Alors pour tout précrible d'un objet $X$
$$
P = \left( X_i \xrightarrow{ \ x_i \ } X \right)_{i \in I}
$$
et tout morphisme $x : X' \to X$ de but $X$, on peut prendre pour précrible image réciproque
$$
x^* \, P
$$
le précrible de $X'$ constitué des morphismes
$$
X_i \times_X X' \longrightarrow X' \, , \quad i \in I \, .
$$
\end{remark}

\begin{demo}

C'est une traduction en termes de précribles du corollaire \ref{cor3.1.6}. 
\end{demo}

\subsection{Application aux réunions de topologies}\label{ssec3.1.4}

\medskip

Un cas particulier important du théorème \ref{thm3.1.5} est le corollaire suivant:

\begin{cor}\label{cor3.1.8}

Soit une famille de topologies $J_k$, $k \in K$, sur une catégorie essentiellement petite ${\mathcal C}$.

Soit 
$$
J = \bigvee_{k \in K} J_k
$$
la plus petite topologie sur ${\mathcal C}$ qui contient toutes les topologies $J_k$, $k \in K$.

Alors un crible sur un objet $X$ de ${\mathcal C}$
$$
C \xymatrix{\ar@{^{(}->}[r] & } y(X)
$$
est élément de $J$ si et seulement si le crible maximal est l'unique crible de $X$ qui satisfait les deux conditions suivantes:

\smallskip

$\left\{\begin{matrix}
\bullet &\mbox{Il contient le crible $C$.} \hfill \\
\bullet &\mbox{Il est $J_k$-fermé pour tout $k \in K$.} \hfill 
\end{matrix}\right.$
\end{cor}

\begin{demo}

Associer à tout objet $X$ la famille
$$
J'(X) = \bigcup_{k \in K} J_k(X)
$$
des cribles de $X$ qui sont couvrants pour l'une au moins des topologies $J_k$, $k \in K$, définit sur ${\mathcal C}$ une notion de crible couvrant qui est stable.

Pour conclure d'après le théorème \ref{thm3.1.5}, il suffit de remarquer qu'un crible d'un objet $X$
$$
C \xymatrix{\ar@{^{(}->}[r] & } y(X)
$$
est $J'$-fermé si et seulement si il est fermé pour chacune des topologies $J_k$, $k \in K$. 
\end{demo}

Il résulte de ce corollaire que si une famille de cribles ou de précribles
$$
{\mathcal X}
$$
d'une catégorie essentiellement petite ${\mathcal C}$ est écrite comme une réunion de sous-familles de cribles ou de précribles
$$
{\mathcal X} = \bigcup_{k \in K} {\mathcal X}_k \, ,
$$
alors la topologie $J$ engendrée par ${\mathcal X}$ s'exprime par une formule simple en termes des topologies $J_k$ engendrées par les sous-familles ${\mathcal X}_k$ ou, plus précisément, des notions de cribles $J_k$-fermés qui leur correspondent.

\subsection{Produits finis de topos}\label{ssec3.1.5}

\medskip

On rappelle que, pour toutes catégories ${\mathcal C}_1 , \cdots , {\mathcal C}_n$, la notation
$$
{\mathcal C}_1 \times \cdots \times {\mathcal C}_n
$$
désigne la catégorie dont

\smallskip

$\left\{\begin{matrix}
\bullet &\mbox{les objets sont les suites $(X_1 , X_2 , \cdots , X_n)$ d'objets $X_1$ de ${\mathcal C}_1$, $X_2$ de ${\mathcal C}_2$, $\cdots$, $X_n$ de ${\mathcal C}_n$,} \hfill \\
\bullet &\mbox{les morphismes $(X'_1 , X'_2 , \cdots , X'_n) \to (X_1 , X_2 , \cdots , X_n)$ sont les suites} \hfill \\
&\left( X'_1 \xrightarrow{ \ x_1 \ } X_1 , X'_2 \xrightarrow{ \ x_2 \ } X_2 , \cdots , X'_n \xrightarrow{ \ x_n \ } X_n \right) \\
&\mbox{de morphismes $x_1$ de ${\mathcal C}_1$, $x_2$ de ${\mathcal C}_2$, $\cdots$, $x_n$ de ${\mathcal C}_n$,} \hfill \\
\bullet &\mbox{la loi de composition associe à deux morphismes} \hfill \\
&(X''_1 , X''_2 , \cdots , X''_n) \to (X'_1 , X'_2 , \cdots , X'_n) \to (X_1 , X_2 , \cdots , X_n) \\
&\mbox{la suite des composés de leurs composantes} \hfill \\
&X''_i \to X'_i \to X_i \, , \quad 1 \leq i \leq n \, .
\end{matrix}\right.$

\smallskip

La catégorie ${\mathcal C}_1 \times \cdots \times {\mathcal C}_n$ est un produit des catégories ${\mathcal C}_1 , \cdots , {\mathcal C}_n$ au sens qu'elle est munie de foncteurs de projection

\smallskip

$\left\{
\begin{matrix}
\hfill p_i : {\mathcal C}_1 \times \cdots \times {\mathcal C}_n &\longrightarrow &{\mathcal C}_i \, , \quad 1 \leq i \leq n \, , \\
\hfill (X_1 , \cdots , X_n) &\longmapsto &X_i \, , \hfill \\
\left( X'_1 \xrightarrow{ \ x_1 \ } X_1 , \cdots , X'_n \xrightarrow{ \ x_n \ } X_n \right) &\longmapsto &\left( X'_i \xrightarrow{ \ x_i \ } X_i \right) \, , \hfill
\end{matrix}
\right.$

\smallskip

tels que, pour toute catégorie ${\mathcal C}$, on ait
\smallskip

$\left\{
\begin{matrix}
\bullet &\mbox{toute suite de foncteurs} \hfill \\
&F_1 : {\mathcal C} \to {\mathcal C}_1 , \cdots , F_n : {\mathcal C} \to {\mathcal C}_n \\
&\mbox{provient d'un unique foncteur} \hfill \\
&F : {\mathcal C} \to {\mathcal C}_1 \times \cdots \times {\mathcal C}_n \\
&\mbox{au sens que} \hfill \\
&F_i = p_i \circ F \, , \quad 1 \leq i \leq n \, , \\
\bullet &\mbox{toute suite de transformations naturelles} \hfill \\
&p_i \circ F \xrightarrow{ \ \alpha_i \ } p_i \circ G \, , \quad 1 \leq i \leq n \, , \\
&\mbox{entre les projections de deux foncteurs $F,G : {\mathcal C} \rightrightarrows {\mathcal C}_1 \times \cdots \times {\mathcal C}_n$} \hfill \\
&\mbox{provient d'une unique transformation naturelle} \hfill \\
&\alpha : F \longrightarrow G \\
&\mbox{via les foncteurs de projection $p_i : {\mathcal C}_1 \times \cdots \times {\mathcal C}_n \to {\mathcal C}_i \, , \quad 1 \leq i \leq n$.}
\end{matrix}
\right.$

\medskip

Si ${\mathcal C}_1 , {\mathcal C}_2 , \cdots , {\mathcal C}_n$ sont des catégories essentiellement petites, leur produit ${\mathcal C}_1 \times \cdots \times {\mathcal C}_n$ est une catégorie essentiellement petite et les foncteurs de projection 
$$
p_i : {\mathcal C}_1 \times \cdots \times {\mathcal C}_n \longrightarrow {\mathcal C}_i
$$
induisent des morphismes de topos de préfaisceaux
$$
(p_i^* , (p_i)_*) : \reallywidehat{{\mathcal C}_1 \times \cdots \times {\mathcal C}_n} \longrightarrow \widehat{\mathcal C}_i \, .
$$
On commence par montrer:

\begin{lem}\label{lem3.1.9}

Considérons des catégories essentiellement petites ${\mathcal C}_1 , \cdots , {\mathcal C}_n$ munies de topologies $J_1 , \cdots J_n$.

Alors on a pour tout indice $i$, $1 \leq i \leq n$:

\begin{listeimarge}

\item Il existe une topologie sur ${\mathcal C}_1 \times \cdots \times {\mathcal C}_n$ pour laquelle un crible $C$ d'un objet $(X_1 , \cdots , X_n)$ est couvrant s'il contient une famille de morphismes de la forme
$$
\left( X_1 \xrightarrow{ \, {\rm id} \, } X_1 , \cdots , X_{i-1} \xrightarrow{ \, {\rm id} \, } X_{i-1} , X_{i,k} \xrightarrow{ \,x_k \, } X_i , X_{i+1} \xrightarrow{ \, {\rm id} \, } X_{i+1} , \cdots , X_n \xrightarrow{ \, {\rm id} \, } X_n \right)
$$
où les
$$
X_{i,k} \xrightarrow{ \ x_k \ } X_i \, , \quad k \in K \, ,
$$
constituent une famille $J_i$-couvrante de $X_i$.

%\smallskip

\item Cette topologie définit le sous-topos du topos des préfaisceaux
$$
\reallywidehat{{\mathcal C}_1 \times \cdots \times {\mathcal C}_n}
$$
qui est l'image réciproque du sous-topos
$$
\widehat{({\mathcal C}_1)}_{J_i} \xymatrix{\ar@{^{(}->}[r] & } \widehat{\mathcal C}_i
$$
par le morphisme de topos
$$
(p_i^* , (p_i)_*) : \reallywidehat{{\mathcal C}_1 \times \cdots \times {\mathcal C}_n} \longrightarrow \widehat{\mathcal C}_i \, .
$$
\end{listeimarge}
\end{lem}

\begin{remark}

La topologie de ${\mathcal C}_1 \times \cdots \times {\mathcal C}_n$ ainsi induite par la topologie $J_i$ de ${\mathcal C}_i$ sera encore notée $J_i$.
\end{remark}

\begin{demo}
\begin{listeisansmarge}
\item Il faut vérifier que cette notion de crible couvrant satisfait les trois axiomes des topologies.

Les cribles maximaux sont couvrants puisque le morphisme identité d'un objet $(X_1 , \cdots , X_n)$ de 

${\mathcal C}_1 \times \cdots \times {\mathcal C}_n$ est la suite
$$
\left( X_1 \xrightarrow{ \, {\rm id} \, } X_1 , \cdots , X_i \xrightarrow{ \, {\rm id} \, } X_i , \cdots , X_n \xrightarrow{ \, {\rm id} \, } X_n \right) \, .
$$
D'autre part l'image réciproque $C'$ d'un crible couvrant $C$ par un morphisme
$$
\left( X'_1 \to X_1 , \cdots , X'_i \to X_i , \cdots , X'_n \to X_n \right) 
$$
est encore un crible couvrant.

En effet, $C$ contient par hypothèse les morphismes de la forme
$$
\left( X_1 \xrightarrow{ \, {\rm id} \, } X_1 , \cdots , X_{i-1} \xrightarrow{ \, {\rm id} \, } X_{i-1} , X_{i,k} \xrightarrow{ \,x_k \, } X_i , X_{i+1} \xrightarrow{ \, {\rm id} \, } X_{i+1} , \cdots , X_n \xrightarrow{ \, {\rm id} \, } X_n \right)
$$
pour une famille $J_i$-couvrante de morphismes
$$
x_k : X_{i,k} \longrightarrow X_i \, , \quad k \in K \, .
$$
Alors il existe une famille $J_i$-couvrante de morphismes
$$
x'_{k'} : X'_{i,k'} \longrightarrow X'_i \, , \quad k' \in K' \, ,
$$
dont les composés avec le morphisme
$$
X'_i \longrightarrow X_i
$$
sont éléments du crible engendré par les morphismes
$$
x_k : X_{i,k} \longrightarrow X_i \, , \quad k\in K \, .
$$
Ainsi, l'image réciproque $C'$ du crible $C$ contient les morphismes
$$
\left( X'_1 \xrightarrow{ \, {\rm id} \, } X'_1 , \cdots , X'_{i-1} \xrightarrow{ \, {\rm id} \, } X'_{i-1} , X'_{i,k'} \xrightarrow{ \,x'_{k'} \, } X'_i , X'_{i+1} \xrightarrow{ \, {\rm id} \, } X'_{i+1} , \cdots , X'_n \xrightarrow{ \, {\rm id} \, } X_n \right) \, , \quad k' \in K' \, ,
$$
et donc il est couvrant.

Enfin, la notion de crible couvrant induite par $J_i$ dans ${\mathcal C}_1 \times \cdots \times {\mathcal C}_n$ satisfait l'axiome de transitivité puisque, pour toutes familles $J_i$-couvrantes
$$
X_{i,k} \xrightarrow{ \ x_k \ } X_i \, , \quad k \in K 
$$
et
$$
X_{i,k,k'} \xrightarrow{ \ x_{k,k'} \ } X_k \, , \quad k' \in K_k \, , \ k \in K \, , 
$$
la famille des composés
$$
X_{i,k,k'} \xrightarrow{ \ x_{k,k'} \ } X_{i,k} \xrightarrow{ \ x_k \ } X_i \, , \quad k \in K \, , \ k' \in K_k \, ,
$$
est encore $J_i$-couvrante.

%\smallskip

\item Pour tout objet $X_i$ de ${\mathcal C}_i$, le foncteur
$$
p_i^* : \widehat{\mathcal C}_i \longrightarrow \reallywidehat{{\mathcal C}_1 \times \cdots \times {\mathcal C}_n}
$$
de composition avec $p_i : {\mathcal C}_1 \times \cdots \times {\mathcal C}_n \to {\mathcal C}_i$ transforme le préfaisceau représentable
$$
y(X_i) = {\rm Hom} (\bullet , X_i)
$$
en le préfaisceau
$$
(X'_1 , \cdots , X'_i , \cdots , X'_n) \longmapsto {\rm Hom} (X'_i , X_i) \, .
$$
De plus, pour tout crible $J_i$-couvrant de $X_i$
$$
C_i \xymatrix{\ar@{^{(}->}[r] & } y(X_i)
$$
engendré par une famille de morphismes
$$
X_{i,k} \xrightarrow{ \ x_k \ } X_i \, , \quad k \in K \, ,
$$
et pour tout objet $(X'_1 , \cdots , X'_i , \cdots , X'_n)$ de ${\mathcal C}_1 \times \cdots \times {\mathcal C}_n$ muni d'un morphisme
$$
X'_i \longrightarrow X_i
$$
vu comme un morphisme de $\reallywidehat{{\mathcal C}_1 \times \cdots \times {\mathcal C}_n}$
$$
y(X'_1 , \cdots , X'_i , \cdots , X'_n) \longrightarrow p_i^* (y(X_i)) \, ,
$$
le produit fibré
$$
y(X'_1 , \cdots , X'_i , \cdots , X'_n) \times_{p_i^* (y(X_i))} p_i^* (C_i)
$$
est le crible
$$
C'_i \xymatrix{\ar@{^{(}->}[r] & } y(X'_1 , \cdots , X'_i , \cdots , X'_n)
$$
engendré par les morphismes de la forme
$$
\left( X'_1 \xrightarrow{ \, {\rm id} \, } X'_1 , \cdots , X'_{i-1} \xrightarrow{ \, {\rm id} \, } X'_{i-1} , X''_i \longrightarrow X'_i , X'_{i+1} \xrightarrow{ \, {\rm id} \, } X'_{i+1} , \cdots , X'_n \xrightarrow{ \, {\rm id} \, } X_n \right) 
$$
où les $X''_i \to X'_i$ décrivent le crible de $X'_i$ image réciproque par le morphisme $X'_i \to X_i$ du crible $C_i$ de $X_i$.

Cela montre que le sous-topos de $\reallywidehat{{\mathcal C}_1 \times \cdots \times {\mathcal C}_n}$ image réciproque du sous-topos $\widehat{({\mathcal C}_i)}_{J_i} \hookrightarrow \widehat{\mathcal C}_i$ est celui défini par la topologie de ${\mathcal C}_1 \times \cdots \times {\mathcal C}_n$ introduite dans (i).
\end{listeisansmarge}
\end{demo}

Nous pouvons maintenant énoncer:

\begin{thm}\label{thm3.1.10}

Considérons des catégories essentiellement petites ${\mathcal C}_1 , \cdots , {\mathcal C}_n$ munies de topologies $J_1 , \cdots , J_n$.

Notons $J_1 \vee \cdots \vee J_n$ la plus petite topologie de la catégorie produit
$$
{\mathcal C}_1 \times \cdots \times {\mathcal C}_n
$$
qui contient les topologies $J_1 , \cdots , J_n$ induites par les foncteurs de projection $p_i : {\mathcal C}_1 \times \cdots \times {\mathcal C}_n \to {\mathcal C}_i$, $1 \leq i \leq n$.

Alors:

\begin{listeimarge}

\item Un crible $C$ d'un objet $(X_1 , \cdots , X_n)$ de ${\mathcal C}_1 \times \cdots \times {\mathcal C}_n$ est couvrant pour la topologie $J_1 \vee \cdots \vee J_n$ si et seulement si le crible maximal est l'unique crible $C'$ de $(X_1 , \cdots , X_n)$ qui satisfait les conditions suivantes:

\smallskip

$\left\{\begin{matrix}
\bullet &\mbox{Il contient le crible $C$.} \hfill \\
\bullet &\mbox{Pour tout indice $i$, $1 \leq i \leq n$, le crible $C'$ est $J_i$-fermé au sens qu'un morphisme} \hfill \\
&\left( X'_1 \xrightarrow{ \, x_1 \, } X_1 , \cdots , X'_i \xrightarrow{ \, x_i \, } X_i , \cdots , X'_n \xrightarrow{ \, x_n \, } X_n \right) \\
&\mbox{est élément de $C'$ dès qu'il existe une famille $J_i$-couvrante de morphismes de ${\mathcal C}_i$} \hfill \\
&x_{i,k} : X'_{i,k} \longrightarrow X'_i \, , \quad k \in K \, , \\
&\mbox{telle que tous les morphismes} \hfill \\
&\left( X'_1 \xrightarrow{ \, x_1 \, } X_1 , \cdots , X'_{i-1} \xrightarrow{ \, x_{i-1} \, } X_{i-1} , X'_{i,k} \xrightarrow{ \, x_i \circ x_{i,k} \, } X_i , X'_{i+1} \xrightarrow{ \, x_{i+1} \, } X_{i+1} , \cdots , X'_n \xrightarrow{ \, x_n \, } X_n \right) \, , \ k \in K \, , \\
&\mbox{soient éléments de $C'$.} \hfill
\end{matrix} \right.$

\smallskip

\item Le topos
$$
\reallywidehat{({\mathcal C}_1 \times \cdots \times {\mathcal C}_n)}_{J_1 \vee \cdots \vee J_n}
$$
muni des morphismes de topos induits par les foncteurs de projection $p_i : {\mathcal C}_1 \times \cdots \times {\mathcal C}_n \to {\mathcal C}_i$, $1 \leq i \leq n$,
$$
\reallywidehat{({\mathcal C}_1 \times \cdots \times {\mathcal C}_n)}_{J_1 \vee \cdots \vee J_n} \longrightarrow \widehat{({\mathcal C}_i)}_{J_i}
$$
est un produit des topos
$$
\widehat{({\mathcal C}_1)}_{J_1} , \cdots , \widehat{({\mathcal C}_n)}_{J_n}
$$
dans la $2$-catégorie des topos, au sens que pour tout topos ${\mathcal E}$, le foncteur
$$
{\rm Geom} \left( {\mathcal E} , \reallywidehat{({\mathcal C}_1 \times \cdots \times {\mathcal C}_n)}_{J_1 \vee \cdots \vee J_n} \right) \longrightarrow {\rm Geom} \left( {\mathcal E} , \widehat{({\mathcal C}_1)}_{J_1} \right) \times \cdots \times {\rm Geom} \left( {\mathcal E} ,  \widehat{({\mathcal C}_n)}_{J_n} \right)
$$
de la catégorie des morphismes de topos
$$
{\mathcal E} \longrightarrow \reallywidehat{({\mathcal C}_1 \times \cdots \times {\mathcal C}_n)}_{J_1 \vee \cdots \vee J_n}
$$
dans le produit des catégories de morphismes de topos
$$
{\mathcal E} \longrightarrow \widehat{({\mathcal C}_i)}_{J_i} \, , \quad 1 \leq i \leq n \, ,
$$
est une équivalence de catégories.
\end{listeimarge}
\end{thm}

\begin{demosansqed}
\begin{listeisansmarge}
\item résulte du lemme \ref{lem3.1.9} (i) combiné avec le corollaire \ref{cor3.1.8}.

%\smallskip

\item résulte du lemme \ref{lem3.1.9} (ii) combiné avec le cas particulier suivant du théorème:
\end{listeisansmarge}
\end{demosansqed}

\begin{prop}\label{prop3.1.11}

Considérons des catégories essentiellement petites ${\mathcal C}_1 , \cdots , {\mathcal C}_n$.

Alors le topos de préfaisceaux
$$
\reallywidehat{{\mathcal C}_1 \times \cdots \times {\mathcal C}_n}
$$
muni des morphismes de topos
$$
\reallywidehat{{\mathcal C}_1 \times \cdots \times {\mathcal C}_n} \longrightarrow \widehat{\mathcal C}_i
$$
induits par les foncteurs de projection $p_i : {\mathcal C}_1 \times \cdots \times {\mathcal C}_n \to {\mathcal C}_i$, $1 \leq i \leq n$, est un produit des topos $\widehat{\mathcal C}_i$, $1 \leq i \leq n$, dans la $2$-catégorie des topos.
\end{prop}

\begin{demo}

On doit montrer que les morphismes de topos
$$
\reallywidehat{{\mathcal C}_1 \times \cdots \times {\mathcal C}_n} \longrightarrow \widehat{\mathcal C}_i \, , \quad 1 \leq i \leq n \, ,
$$
induisent, pour tout topos ${\mathcal E}$, une équivalence entre

\smallskip

$\left\{\begin{matrix}
\bullet &\mbox{la catégorie ${\mathcal F}$ des foncteurs} \hfill \\
&\reallywidehat{{\mathcal C}_1 \times \cdots \times {\mathcal C}_n} \longrightarrow {\mathcal E} \\
&\mbox{qui respectent les colimites arbitraires et les limites finies,} \hfill \\
\bullet &\mbox{le produit ${\mathcal F}_1 \times \cdots \times {\mathcal F}_n$ des catégories ${\mathcal F}_i$, $1\leq i \leq n$, des foncteurs} \hfill \\
&\widehat{\mathcal C}_i \longrightarrow {\mathcal E} \\
&\mbox{qui respectent les colimites arbitraires et les limites finies.} \hfill
\end{matrix}\right.$

\smallskip

Les foncteurs
$$
p_i^* : \widehat{\mathcal C}_i \longrightarrow \reallywidehat{{\mathcal C}_1 \times \cdots \times {\mathcal C}_n} \, , \quad 1 \leq i \leq n \, ,
$$
de composition avec les foncteurs de projection
$$
p_i : {\mathcal C}_1 \times \cdots \times {\mathcal C}_n \longrightarrow {\mathcal C}_i
$$
respectent toutes les colimites et limites.

Donc la composition avec chaque $p_i^*$, $1 \leq i \leq n$, définit un foncteur
$$
\begin{matrix}
\hfill {\mathcal F} &\longrightarrow &{\mathcal F}_i \, , \hfill \\
\widehat\rho &\longmapsto &\widehat\rho \circ p_i^* = \widehat\rho_i \, .
\end{matrix}
$$
\end{demo}

On a :

\begin{lem}\label{lem3.1.12}

Soient toujours des catégories essentiellement petites ${\mathcal C}_1 , \cdots , {\mathcal C}_n$.

Considérons un foncteur à valeurs dans un topos ${\mathcal E}$ qui respecte les colimites arbitraires et les limites finies
$$
\widehat\rho : \reallywidehat{{\mathcal C}_1 \times \cdots \times {\mathcal C}_n} \longrightarrow {\mathcal E} \, ,
$$
puis les foncteurs induits
$$
\widehat\rho \circ p_i^* = \widehat\rho_i : \widehat{\mathcal C}_i \longrightarrow {\mathcal E} \, .
$$
Notons
$$
\rho : {\mathcal C}_1 \times \cdots \times {\mathcal C}_n \longrightarrow {\mathcal E}
$$
et
$$
\rho_i : {\mathcal C}_i \longrightarrow {\mathcal E} \, , \quad 1 \leq i \leq n
$$
les composés des $\widehat\rho$ et $\widehat\rho_i$ avec les foncteurs de Yoneda
$$
y : {\mathcal C}_1 \times \cdots \times {\mathcal C}_n \xymatrix{\ar@{^{(}->}[r] & } \reallywidehat{{\mathcal C}_1 \times \cdots \times {\mathcal C}_n} \, ,
$$
$$
y_i : {\mathcal C}_i \xymatrix{\ar@{^{(}->}[r] & } \widehat{\mathcal C}_i \, , \quad 1 \leq i \leq n \, .
$$
Alors le foncteur
$$
\rho : {\mathcal C}_1 \times \cdots \times {\mathcal C}_n \longrightarrow {\mathcal E}
$$
est canoniquement isomorphe au foncteur
$$
\begin{matrix}
\rho_1 \times \cdots \times \rho_n &: &{\mathcal C}_1 \times \cdots \times {\mathcal C}_n &\longrightarrow &{\mathcal E} \, , \hfill  \\
&&(X_1 , \cdots , X_n) &\longmapsto &\rho_1 (X_1) \times \cdots \times \rho_n (X_n) \, .
\end{matrix}
$$
\end{lem}

\begin{demolem}

Par définition de la catégorie produit ${\mathcal C}_1 \times \cdots \times {\mathcal C}_n$, on a pour tous objets $(X_1 , \cdots , X_n)$ et $(X'_1 , \cdots , X'_n)$ la formule
$$
{\rm Hom} \left( (X'_1 , \cdots , X'_n) , (X_1 , \cdots , X_n) \right) = {\rm Hom} (X'_1 , X_1) \times \cdots \times {\rm Hom} (X'_n , X_n) \, .
$$
Cette formule signifie que le foncteur de Yoneda
$$
y : {\mathcal C}_1 \times \cdots \times {\mathcal C}_n \longrightarrow \reallywidehat{{\mathcal C}_1 \times \cdots \times {\mathcal C}_n}
$$
s'identifie au foncteur produit
$$
p_1^* (y_1 \circ p_1) \times \cdots \times p_n^* (y_n \circ p_n) \, .
$$
Comme le foncteur
$$
\widehat\rho : \reallywidehat{{\mathcal C}_1 \times \cdots \times {\mathcal C}_n} \longrightarrow {\mathcal E}
$$
respecte par hypothèse les limites finies, on obtient que son composé avec le foncteur de Yoneda
$$
y : {\mathcal C}_1 \times \cdots \times {\mathcal C}_n \longrightarrow \reallywidehat{{\mathcal C}_1 \times \cdots \times {\mathcal C}_n}
$$
est canoniquement isomorphe au foncteur produit
$$
\begin{matrix}
\widehat\rho_1 (y_1 \circ p_1) \times \cdots \times \widehat\rho_n (y_n \circ p_n) &: &{\mathcal C}_1 \times \cdots \times {\mathcal C}_n &\longrightarrow &{\mathcal E} \, , \hfill  \\
&&(X_1 , \cdots , X_n) &\longmapsto &\rho_1 (X_1) \times \cdots \times \rho_n (X_n) \, .
\end{matrix}
$$
\end{demolem}

\begin{suitedemoprop}
Réciproquement, considérons des foncteurs à valeurs dans un topos ${\mathcal E}$ qui respectent les colimites arbitraires et les limites finies
$$
\widehat\rho_i : \widehat{\mathcal C}_i \longrightarrow {\mathcal E} \, , \quad 1 \leq i \leq n \, .
$$
Considérons leurs composés avec les foncteurs de Yoneda
$$
\rho_i = \widehat\rho_i \circ y_i : {\mathcal C}_i \longrightarrow {\mathcal E} \, , \quad 1 \leq i \leq n \, ,
$$
puis le foncteur produit
$$
\begin{matrix}
\rho &: &{\mathcal C}_1 \times \cdots \times {\mathcal C}_n &\longrightarrow &{\mathcal E} \hfill  \\
&&(X_1 , \cdots , X_n) &\longmapsto &\rho_1 (X_1) \times \cdots \times \rho_n (X_n) \, .
\end{matrix}
$$
Celui-ci s'étend en un foncteur respectant les colimites
$$
\widehat\rho : \reallywidehat{{\mathcal C}_1 \times \cdots \times {\mathcal C}_n} \longrightarrow {\mathcal E} \, .
$$
Montrons que, pour $1 \leq i \leq n$, le composé
$$
\widehat\rho \circ p_i^* : \widehat{\mathcal C}_i \longrightarrow {\mathcal E}
$$
est canoniquement isomorphe au foncteur $\widehat\rho_i$.

Comme tous deux respectent les colimites, il suffit de montrer que leurs restrictions à la catégorie ${\mathcal C}_i$ sont canoniquement isomorphes.

Or, pour tout objet $X_i$ de ${\mathcal C}_i$, le préfaisceau sur ${\mathcal C}_1 \times \cdots \times {\mathcal C}_n$
$$
p_i^* (y_i (X_i))
$$
s'écrit comme une colimite sur la catégorie de ses éléments
$$
\varinjlim_{(X'_1 , \cdots , X'_n) \in {\mathcal C}_1 \times \cdots \times {\mathcal C}_n / p_i^* (y_i(X_i))} y(X'_1 , \cdots , X'_n)
$$
et cette catégorie d'éléments 
$$
{\mathcal C}_1 \times \cdots \times {\mathcal C}_n / p_i^* (y_i(X_i))
$$
s'identifie à la catégorie produit
$$
{\mathcal C}_1 \times \cdots \times {\mathcal C}_{i-1} \times ({\mathcal C}_i / X_i) \times {\mathcal C}_{i+1} \times \cdots \times {\mathcal C}_n \, .
$$
Il en résulte que
$$
\widehat\rho \circ p_i^* \circ y_i(X_i) 
$$
est canoniquement isomorphe à
\begin{eqnarray}
&&\varinjlim_{{\mathcal C}_1} \rho_1 (X'_1) \times \cdots \times \varinjlim_{{\mathcal C}_{i-1}} \rho_{i-1} (X'_i) \times \varinjlim_{{\mathcal C}_i / X_i} \rho_i (X'_i) \times \varinjlim_{{\mathcal C}_{i+1}} \rho_{i+1} (X'_{i+1}) \times \cdots \times \varinjlim_{{\mathcal C}_n} \rho_n (X'_n) \nonumber \\
&= &\widehat\rho_1 (1_{\widehat{\mathcal C}_1}) \times \cdots \times \widehat\rho_{i-1} (1_{\widehat{\mathcal C}_{i-1}}) \times \rho_i (X_i) \times \widehat\rho_{i+1} (1_{\widehat{\mathcal C}_{i+1}}) \times \cdots \times \widehat\rho_n (1_{\widehat{\mathcal C}_n}) \nonumber \\
&= &\rho_i (X_i) \nonumber
\end{eqnarray}
puisque chaque foncteur
$$
\widehat\rho_{i'} : \widehat{\mathcal C}_{i'} \longrightarrow \widehat{\mathcal E}
$$
transforme l'objet terminal $1_{\widehat{\mathcal C}_{i'}}$ de $\widehat{\mathcal C}_{i'}$ en l'objet terminal $1_{\mathcal E}$ du topos ${\mathcal E}$.

On a bien montré que chaque composé
$$
\widehat\rho \circ p_i^* : \widehat{\mathcal C}_i \longrightarrow \reallywidehat{{\mathcal C}_1 \times \cdots \times {\mathcal C}_n} \longrightarrow {\mathcal E}
$$
est canoniquement isomorphe au foncteur $\widehat\rho_i$.

Pour conclure, il reste à montrer que le foncteur respectant les colimites arbitraires
$$
\widehat{\rho} : \reallywidehat{{\mathcal C}_1 \times \cdots \times {\mathcal C}_n} \longrightarrow {\mathcal E}
$$
respecte aussi les limites finies.

Nous allons utiliser pour cela les identités
$$
\widehat{\rho} \circ p_i^* = \widehat\rho_i \, , \quad 1 \leq i \leq n \, ,
$$
et le fait que les foncteurs $\widehat\rho_i$ respectent les limites finies.

Commençons par montrer que pour tout diagramme fini $X_{\bullet}$ de ${\mathcal C} = {\mathcal C}_1 \times \cdots \times {\mathcal C}_n$ indexé par un carquois fini $D$, le morphisme canonique
$$
\widehat\rho \left( \varprojlim_D y(X_{\bullet})\right) \longrightarrow \varprojlim_D \rho (X_{\bullet})
$$
est un isomorphisme.

En effet, notant pour tout indice $i$, $1 \leq i \leq n$
$$
X_{i,\bullet} = p_i (X_{\bullet})
$$
le diagramme fini de ${\mathcal C}_i$ indexé par $D$ qui est déduit de $X_{\bullet}$ par composition avec le foncteur $p_i : {\mathcal C} = {\mathcal C}_1 \times \cdots \times {\mathcal C}_n \to {\mathcal C}_i$, on a
$$
\varprojlim_D y(X_{\bullet}) = p_1^* \left(\varprojlim_D y_1 (X_{1,\bullet}) \right) \times \cdots \times p_n^* \left(\varprojlim_D y_n (X_{n,\bullet}) \right) 
$$
puis
$$
\widehat\rho \left(\varprojlim_D y(X_{\bullet}) \right) = \widehat\rho_1 \left(\varprojlim_D y_1 (X_{1,\bullet}) \right) \times \cdots \times \widehat\rho_n \left(\varprojlim_D y_n (X_{n,\bullet}) \right) \, .
$$
Pour tout objet $E$ de ${\mathcal E}$, se donner un morphisme
$$
E \longrightarrow \widehat\rho \left(\varprojlim_D y(X_{\bullet}) \right)
$$
équivaut donc à se donner une famille de $n$ morphismes
$$
E \longrightarrow \widehat\rho_i \left(\varprojlim_D y_i (X_{i,\bullet}) \right) \, , \quad 1 \leq i \leq n \, ,
$$
ou, ce qui revient au même,
$$
E \longrightarrow \varprojlim_D \rho_i (X_{i,\bullet}) \, , \quad 1 \leq i \leq n \, .
$$
Comme on a pour tout objet $d$ de $D$
$$
\rho (X_d) = \rho_1 (X_{1,d}) \times \cdots \times \rho_n (X_{n,d})
$$
et que, pour toute flèche $d' \to d$ de $D$, le morphisme correspondant
$$
X_{d'} \longrightarrow X_d
$$
consiste en une famille de $n$ morphismes
$$
X_{i,d'} \longrightarrow X_{i,d} \, , \quad 1 \leq i \leq n \, ,
$$
on voit que se donner une famille de morphismes
$$
E \longrightarrow \varprojlim_D \rho_i (X_{i,\bullet}) \, , \quad 1 \leq i \leq n \, ,
$$
équivaut à se donner un morphisme
$$
E \longrightarrow \varprojlim_D \rho (X_{\bullet}) \, .
$$
On a bien montré que le morphisme
$$
\widehat\rho \left(\varprojlim_D y(X_{\bullet}) \right) \longrightarrow \varprojlim_D \rho (X_{\bullet}) 
$$
est un isomorphisme.

Prenant pour $D$ le carquois vide, on obtient en particulier que le foncteur
$$
\widehat\rho : \widehat{\mathcal C} \longrightarrow {\mathcal E}
$$
transforme l'objet terminal de $\widehat{\mathcal C}$ en un objet terminal de ${\mathcal E}$.

Puis considérons deux préfaisceaux $P$ et $Q$ sur ${\mathcal C} = {\mathcal C}_1 \times \cdots \times {\mathcal C}_n$.

Notant ${\mathcal C}/P$ et ${\mathcal C}/Q$ les catégories de leurs éléments, on a
$$
P = \varinjlim_{(X,p) \in {\mathcal C}/P} y(X) \quad \mbox{et} \quad Q = \varinjlim_{(Y,q) \in {\mathcal C}/Q} y(Y) \, .
$$
Or, on a pour tous objets $X$ et $Y$ de ${\mathcal C}$
$$
\widehat\rho \, (y(X) \times y(Y)) = \widehat\rho \, (y(X)) \times \widehat\rho \, (y(Y)) \, .
$$
Comme les produits respectent les colimites dans $\widehat{\mathcal C}$ et dans ${\mathcal E}$ et que le foncteur $\widehat\rho$ respecte les colimites, on obtient
$$
\widehat\rho \, (P \times Q) = \widehat\rho \, (P) \times \widehat\rho \, (Q) \, .
$$

Pour conclure, on doit encore prouver que pour toute paire de morphismes de préfaisceaux sur ${\mathcal C} = {\mathcal C}_1 \times \cdots \times {\mathcal C}_n$
$$
\xymatrix{ P \dar[r]^{^{^{\mbox{\footnotesize$q_1$}}}}_{q_2} & Q } ,
$$
le morphisme canonique
$$
\begin{matrix}
\widehat\rho \left( {\rm eg} \left( {\xymatrix{ P \dar[r]^{^{^{\mbox{\footnotesize$q_1$}}}}_{q_2} & Q }} \right)\right) &\longrightarrow &{\rm eg} \left( {\xymatrix{ \widehat\rho (P) \dar[r] & \widehat\rho (Q) }} \right) \\
\Vert && \Vert \\
\widehat\rho \, (P \times_{Q \times Q} Q) &&\widehat\rho \, (P) \times_{\widehat\rho(Q) \times \widehat\rho (Q)} \widehat\rho \, (Q)
\end{matrix}
$$
est un isomorphisme.

Ecrivant le préfaisceau $P$ comme une colimite
$$
P = \varinjlim_{(X,p) \in {\mathcal C}/P} y(X) \, ,
$$
on se réduit au cas où $P$ est un préfaisceau représentable
$$
P = y(X)
$$
et où donc les morphismes $y(X) = {\xymatrix{ P \dar[r]^{^{^{\mbox{\footnotesize$q_1$}}}}_{q_2} & Q }}$ peuvent être vus comme des sections de $Q$
$$
q_1 , q_2 \in Q(X) \, .
$$
Considérons un objet $E$ de ${\mathcal E}$ muni d'un morphisme
$$
E \longrightarrow {\rm eg} \left({\xymatrix{ \rho \, (X) \dar[r] & \widehat\rho \, (Q) }}\right)
$$
c'est-à-dire d'un morphisme
$$
e : E \longrightarrow \rho (X)
$$
tel que
$$
\widehat\rho \, (q_1) \circ e = \widehat\rho \, (q_2) \circ e : E \longrightarrow \widehat\rho \, (Q) \, .
$$
Ecrivant le préfaisceau $Q$ comme une colimite
$$
Q = \varinjlim_{(Y,q) \in {\mathcal C}/Q} y(Y) \, ,
$$
on a aussi
$$
\widehat\rho \, (Q) = \varinjlim_{(Y,q) \in {\mathcal C}/Q} \rho (Y) \, .
$$
L'égalité
$$
\widehat\rho \, (q_1) \circ e = \widehat\rho \, (q_2) \circ e
$$
signifie qu'existent une famille globalement épimorphique de morphismes de ${\mathcal E}$
$$
e_k : E_k \longrightarrow E \, , \quad k \in K \, ,
$$
et des diagrammes finis connexes de ${\mathcal C}/Q$
$$
\left(X_{\bullet}^k , q_{\bullet}^k\right)
$$
indexés par des carquois finis connexes $D_{\bullet}^k$ contenant chacun deux objets
$$
d_1^k , d_2^k
$$
tels que
$$
\left(X_{d_1^k}^k , q_{d_1^k}^k \right) = (X , q_1) \, ,
$$
$$
\left(X_{d_2^k}^k , q_{d_2^k}^k \right) = (X , q_2) 
$$
et que le morphisme composé
$$
E_k \xrightarrow{ \ e_k \ } E \xrightarrow{ \ e \ } \rho (X)
$$
se factorise en
$$
E_k \longrightarrow \varprojlim_{D^k} \, \rho \left(X_{\bullet}^k\right) \, .
$$
Cela signifie que les morphismes
$$
E_k \xrightarrow{ \ e_k \ } E \longrightarrow  {\rm eg} \left({\xymatrix{ \rho \, (X) \dar[r] & \widehat\rho \, (Q) }}\right)
$$
se factorisent chacun en
$$
E_k \longrightarrow \widehat\rho \left( {\rm eg} \left( {\xymatrix{ y(X) \dar[r]^{^{^{\mbox{\footnotesize$q_1$}}}}_{q_2} & Q }} \right)\right) ,
$$
et donc que le morphisme
$$
\widehat\rho \left( {\rm eg} \left( {\xymatrix{ y(X) \dar[r]& Q }} \right) \right) \longrightarrow {\rm eg} \left({\xymatrix{ \rho \, (X) \dar[r] & \widehat\rho \, (Q) }}\right)
$$
est un épimorphisme.

On en déduit que pour toute paire de morphismes de $\widehat{\mathcal C}$
$$
{\xymatrix{ P \dar[r]& Q }}
$$
le morphisme
$$
\widehat\rho \left( {\rm eg} \left( {\xymatrix{ P \dar[r]& Q }} \right) \right) \longrightarrow {\rm eg} \left({\xymatrix{ \widehat\rho \, (P) \dar[r] & \widehat\rho \, (Q) }}\right)
$$
est un épimorphisme.

Pour terminer, il reste à montrer que $\widehat\rho$ transforme tout monomorphisme de $\widehat{\mathcal C}$
$$
P \xymatrix{\ar@{^{(}->}[r] & } Q
$$
en un monomorphisme de ${\mathcal E}$.

En effet, si $P \to Q$ est un monomorphisme de $\widehat{\mathcal C}$,
$$
P \longrightarrow {\rm eg} \left( \xymatrix{ P \times P \dar[r]& Q } \right)
$$
est un épimorphisme, donc le composé
$$
\widehat\rho \, (P) \longrightarrow \widehat\rho \left({\rm eg} \left( \xymatrix{ P \times P \dar[r]& Q } \right)\right) \longrightarrow  {\rm eg} \left( \xymatrix{ \widehat\rho \, (P) \times \widehat\rho \, (P) \dar[r]& \widehat\rho \, (Q) } \right)
$$
est aussi un épimorphisme d'après ce qu'on a montré et le fait que $\widehat\rho$ respecte les colimites, en particulier les épimorphismes.

Cela signifie comme voulu que le morphisme
$$
\widehat\rho \, (P) \longrightarrow \widehat\rho \, (Q)
$$
est un monomorphisme.

Ceci achève la démonstration de la proposition et donc aussi du théorème \ref{thm3.1.10}. 
\end{suitedemoprop}

\subsection{Produits finis d'espaces topologiques et produits finis de topos}\label{3.1.6}

Associer à tout espace topologique $X$ le topos
$$
{\mathcal E}_X = \widehat{({\mathcal C}_X)}_{J_X}
$$
des faisceaux sur $X$, c'est-à-dire sur la catégorie ${\mathcal C}_X$ de ses ouverts munie de la topologie $J_X$ définie par la notion usuelle de recouvrement, définit un foncteur de la catégorie des espaces topologiques vers celle des topos.

Ce foncteur associe à toute application continue entre deux espaces topologiques
$$
f : X \longrightarrow Y
$$
le morphisme de topos
$$
(f^* , f_*) : {\mathcal E}_X \longrightarrow {\mathcal E}_Y \, .
$$

En particulier, pour tous espaces topologiques $X_1 , \cdots , X_n$ et leur produit $X_1 \times \cdots \times X_n$, les projections continues
$$
X_1 \times \cdots \times X_n \longrightarrow X_i \, , \quad 1 \leq i \leq n \, ,
$$
définissent des morphismes de topos
$$
{\mathcal E}_{X_1 \times \cdots \times X_n} \longrightarrow {\mathcal E}_{X_i} \, , \quad 1 \leq i \leq n \, ,
$$
et donc un morphisme de topos
$$
{\mathcal E}_{X_1 \times \cdots \times X_n} \longrightarrow {\mathcal E}_{X_1} \times \cdots \times {\mathcal E}_{X_n} \, .
$$
On a:

\begin{prop}\label{prop3.1.13}

Soient des espaces topologiques $X_1 , \cdots , X_n$ qui sont localement compacts sauf éventuellement l'un d'entre eux.

Alors le morphisme de topos
$$
{\mathcal E}_{X_1 \times \cdots \times X_n} \longrightarrow {\mathcal E}_{X_1} \times \cdots \times {\mathcal E}_{X_n} 
$$
est un isomorphisme.
\end{prop}

\begin{remark}

Ici, ``localement compact'' a le sens que tout ouvert
$$
U
$$
s'écrit comme une réunion d'ouverts
$$
U_i \, , \quad i \in I \, ,
$$
dont les adhérences
$$
\overline U_i
$$
sont contenues dans $U$ et sont compactes, c'est-à-dire permettent d'extraire un sous-recouvrement fini de tout recouvrement ouvert.
\end{remark}

\begin{demosansqed}

Elle repose sur le lemme suivant:
\end{demosansqed}

\begin{lem}\label{lem3.1.14}

Soient un espace topologique $X$ et un espace topologique $Y$ qui est supposé compact.

Soit un recouvrement de l'espace topologique produit
$$
X \times Y
$$
par des ouverts $W_k$, $k \in K$.

Alors il existe un recouvrement de $X$ par des ouverts
$$
U_i \, , \quad i \in I \, ,
$$
et pour chaque indice $i \in I$, un recouvrement fini de $Y$ par des ouverts
$$
V_{i,j} \, , \quad 1 \leq j \leq n_i \, ,
$$
tels que chaque ouvert produit
$$
U_i \times V_{i,j} \, , \quad i \in I \, , \quad 1 \leq j \leq n_i \, ,
$$
soit contenu dans l'un au moins des $W_k$, $k \in K$.
\end{lem}

\begin{demolem}
Pour tout élément $(x,y)$ de $X \times Y$, choisissons deux voisinages ouverts
$$
U_{x,y} \quad \mbox{de $x$ dans $X$}
$$
et
$$
V_{x,y} \quad \mbox{de $y$ dans $Y$}
$$
tels que le produit
$$
U_{x,y} \times V_{x,y}
$$
soit contenu dans l'un au moins des $W_k$, $k \in K$.

Pour tout $x \in X$, les ouverts
$$
V_{x,y} \, , \quad y \in Y \, ,
$$
recouvrent $Y$ donc il est possible d'en extraire un sous-recouvrement fini
$$
V_{x,y_j} \, , \quad 1 \leq j \leq n_x \, .
$$

Posons alors pour tout $x \in X$
$$
U_x = \bigcap_{1 \leq j \leq n_x} U_{x,y_j} \, .
$$
C'est un voisinage ouvert de $x$ dans $X$.

On voit que l'espace $X$ est recouvert pour les ouverts
$$
U_x \, , \quad x \in X \, ,
$$
que pour tout $x \in X$, l'espace $Y$ est recouvert par la famille finie d'ouverts
$$
V_{x,y_j} \, , \quad 1 \leq j \leq n_x \, ,
$$
et que chaque produit
$$
U_x \times V_{x,y_j} \, , \quad x \in X \, , \quad 1 \leq j \leq n_x \, ,
$$
est contenu dans l'un au moins des $W_k$, $k \in K$.

C'est ce que l'on voulait. 
\end{demolem}

\begin{demoprop}
Il résulte d'abord de ce lemme que le produit de deux espaces compacts est un espace compact, donc aussi le produit d'une famille finie d'espaces compacts.

On en déduit que le produit d'une famille finie d'espaces localement compacts est localement compact.

Cela réduit la proposition au cas du produit
$$
X \times Y
$$
d'un espace topologique $X$ et d'un espace topologique $Y$ qui est supposé localement compact.

La catégorie ${\mathcal C}_{X \times Y}$ des ouverts de $X \times Y$ munie de la topologie $J_{X \times Y}$ des recouvrements au sens usuel possède comme sous-catégorie pleine et dense celle des ouverts produits de deux ouverts non vides de $X$ et $Y$
$$
U \times V \, .
$$
Donc le topos ${\mathcal E}_{X \times Y}$ des faisceaux sur $X \times Y$ s'identifie au topos
$$
\widehat{\mathcal C}_J
$$
des faisceaux sur la catégorie ${\mathcal C}$ des produits $U \times V$ d'ouverts non vides de $X$ et $Y$, munie de la topologie $J$ induite par $J_{X \times Y}$.

La topologie $J$ contient les topologies $J_X$ et $J_Y$ de ${\mathcal C}$ engendrées par les topologies de $X$ et $Y$.

Il s'agit de prouver que la topologie $J$ est engendrée par les topologies $J_X$ et $J_Y$.

Cette assertion est locale sur $X$ ou $Y$ donc il suffit de considérer le cas où $Y$ est un espace compact et de montrer que tout crible $J$-couvrant $C$ de l'espace
$$
X \times Y
$$
contient une famille couvrante pour la topologie engendrée par $J_X$ et $J_Y$.

D'après le lemme \ref{lem3.1.14}, un tel crible $C$ contient une famille de la forme
$$
U_i \times V_{i,j} \, , \quad i \in I \, , \quad 1 \leq j \leq n_i \, ,
$$
où les $U_i$ forment un recouvrement de $X$ et, pour chaque $i \in I$, les $V_{i,j}$, $1 \leq j \leq n_i$, forment un recouvrement de $X$.

On en déduit que tout crible $J_Y$-fermé qui contient $C$ comprend les ouverts
$$
U_i \times Y \, , \quad i \in I \, ,
$$
puis que tout crible $J_Y$-fermé et $J_X$-fermé qui contient $C$ comprend le produit
$$
X \times Y \, .
$$
Cela signifie comme voulu que $C$ est couvrant pour la topologie engendrée par $J_X$ et $J_Y$. 
\end{demoprop}

La proposition \ref{prop3.1.13} a la forme concrète suivante qui permet de construire le topos d'un espace topologique produit d'espaces localement compacts $X_1 , \cdots , X_n$ sans faire référence aux points de ces espaces:

\begin{cor}\label{cor3.1.15}

Soient des espaces topologiques $X_1 , \cdots , X_n$.

Soient des bases d'ouverts de $X_1 , \cdots , X_n$ et les catégories ${\mathcal C}_1 , \cdots , {\mathcal C}_n$ qu'elles constituent, munies des topologies $J_1 , \cdots , J_n$ définies par les notions usuelles de recouvrement dans $X_1 , \cdots , X_n$.

Alors:

\begin{listeimarge}

\item Le topos ${\mathcal E}_{X_1 \times \cdots \times X_n}$ des faisceaux sur l'espace produit $X_1 \times \cdots \times X_n$ s'identifie au topos des faisceaux sur la catégorie produit
$$
{\mathcal C}_1 \times \cdots \times {\mathcal C}_n
$$
munie de la topologie $J$ induite par la notion usuelle de recouvrement dans $X_1 , \cdots , X_n$.

%\smallskip

\item Cette topologie $J$ contient les topologies $J_1 , \cdots , J_n$ induites par celles de ${\mathcal C}_1 , \cdots , {\mathcal C}_n$.

%\smallskip

\item Si tous les facteurs $X_i$ sont des espaces localement compacts, à l'exception éventuelle d'un seul, la topologie $J$ est engendrée par $J_1 , \cdots , J_n$.

%\smallskip

\item Si $J$ est engendré par $J_1 , \cdots , J_n$, une famille de morphismes de ${\mathcal C}_1 \times \cdots \times {\mathcal C}_n$
$$
(U_1^i , \cdots , U_n^i) \xymatrix{\ar@{^{(}->}[r] & } (U_1 , \cdots , U_n) \, , \quad i \in I \, ,
$$
est $J$-couvrante si et seulement si le crible maximal est le seul crible qui

\smallskip

$\left\{\begin{matrix}
\bullet &\mbox{contient cette famille,} \hfill \\
\bullet &\mbox{est fermé relativement à chaque topologie $J_i$, $1 \leq i \leq n$.} \hfill
\end{matrix}\right.$
\end{listeimarge}
\end{cor}

\begin{remark}

En général, la topologie $J$ raffine la topologie ${\mathcal C}_1 \times \cdots \times {\mathcal C}_n$ engendrée par $J_1 , \cdots , J_n$ ce qui signifie que le morphisme canonique
$$
{\mathcal E}_{X_1 \times \cdots \times X_n}  \xymatrix{\ar@{^{(}->}[r] & } {\mathcal E}_{X_1} \times \cdots \times  {\mathcal E}_{X_n}
$$
est un plongement.
\end{remark}

\begin{demo}
\begin{listeisansmarge}
\item Comme l'ouvert vide de n'importe quel espace topologique est recouvert par la famille vide, on peut supposer que tous les ouverts qui composent les catégories ${\mathcal C}_1 , \cdots , {\mathcal C}_n$ sont non vides et donc que
$$
{\mathcal C}_1 \times \cdots \times {\mathcal C}_n
$$
est une sous-catégorie pleine et dense de celle des ouverts de l'espace $X_1 \times \cdots \times X_n$.

%\smallskip

\item résulte de ce que les projections
$$
X_1 \times \cdots \times X_n \longrightarrow X_i \, , \quad 1 \leq i \leq n \, ,
$$
sont des applications continues.

%\smallskip

\item est une formulation concrète de la proposition \ref{prop3.1.13}.

%\smallskip

\item résulte du corollaire \ref{cor3.1.8}.
\end{listeisansmarge} 
\end{demo}

\section{Une formule d'engendrement en termes de multi-recouvrements}\label{sec3.2}

\subsection{La notion de multi-recouvrement d'un objet}\label{3.2.1}

Si une catégorie essentiellement petite ${\mathcal C}$ est munie d'une notion de précrible couvrant qui est ``stable'' au sens du lemme \ref{lem3.1.3}, on est amené à composer des précribles couvrants pour construire des arbres de morphismes de ${\mathcal C}$ que l'on appellera alors des multi-recouvrements.

Rappelons d'abord la notion d'arbre:

\begin{defn}\label{defn3.2.1}

Un arbre $A_{\bullet}$ consiste en

\medskip

$\left\{\begin{matrix}
\bullet &\mbox{une suite d'ensemble $A_n$, $n \in {\mathbb N}$, telle que $A_0$ soit réduit à un unique élément,} \hfill \\
\bullet &\mbox{une suite d'applications entre ces ensembles} \hfill \\
&f_n : A_n \longrightarrow A_{n-1} \, , \quad n \geq 1 \, .
\end{matrix}\right.$
\end{defn}

\begin{remarksqed}

\begin{listeisansmarge}

\item On appelle ``n\oe uds'' d'un arbre $A_{\bullet}$ les paires
$$
(n,a)
$$
constituées d'un entier $n \in {\mathbb N}$ et d'un élément $a \in A_n$.

%\smallskip

\item L'ensemble des n\oe uds d'un arbre $A_{\bullet}$ est muni d'une relation d'ordre partiel définie en convenant que
$$
(n,a) \geq (m,b)
$$
si 
$$
n \geq m
$$
et
$$
b = (f_{m+1} \circ \cdots \circ f_n) (a) \, .
$$
\item Une ``branche'' d'un arbre $A_{\bullet}$ est une suite d'éléments
$$
a_0 \in A_0 , a_1 \in A_1 , \cdots , a_n \in A_n \, ,
$$
telle que
$$
a_{i-1} = f_i (a_i) \, , \quad 1 \leq i \leq n \, .
$$
\end{listeisansmarge}
\end{remarksqed}

%\medskip

On distingue parmi les arbres ceux qui sont ``à branches finies'':

\begin{defn}\label{defn3.2.2}

Un arbre
$$
A_{\bullet} = \left( A_0 \overset{f_1}{\longleftarrow} A_1 \overset{f_2}{\longleftarrow} \cdots \overset{f_n}{\longleftarrow} A_n \longleftarrow \cdots \right)
$$
est dit ``à branches finies'' si, dans l'ensemble ordonné de ses n\oe uds, toute partie non vide admet au moins un élément maximal.
\end{defn}

\begin{remarksqed}
\begin{listeisansmarge}
\item Si un arbre est ``à branches finies'', il n'admet pas de ``branche infinie'' c'est-à-dire de suite d'éléments
$$
a_n \in A_n \, , \quad n \in {\mathbb N} \, ,
$$
telle que
$$
a_{n-1} = f_n (a_n) \, , \quad \forall \, n \geq 1 \, .
$$
\item Réciproquement, si on accepte l'axiome du choix, tout arbre qui n'a pas de ``branche infinie'' est nécessairement ``à branches finies''. 
\end{listeisansmarge}
\end{remarksqed}

On appellera ``multi-recouvrements'' les composés de recouvrements, c'est-à-dire de précribles couvrants, indexés par des arbres ``à branches finies'':

\begin{defn}\label{defn3.2.3}

Soit ${\mathcal C}$ une catégorie essentiellement petite munie d'une notion $J$ de ``précrible couvrant'' qui est supposée stable.

Alors un ``multi-recouvrement'' de type $J$ d'un objet $X$ de ${\mathcal C}$ consiste en

\medskip

$\left\{\begin{matrix}
\bullet &\mbox{un arbre à branches finies $A_{\bullet} = \left( A_0 \overset{f_1}{\longleftarrow} A_1 \longleftarrow \cdots \longleftarrow A_{n-1} \overset{f_n}{\longleftarrow} A_n \longleftarrow \cdots \right)$,} \hfill \\
\bullet &\mbox{une famille d'objets de ${\mathcal C}$} \hfill \\
&X_{n,a} \\
&\mbox{indexés par les n\oe uds $(n,a)$ de l'arbre $A_{\bullet}$,} \hfill \\
\bullet &\mbox{une famille de morphismes de ${\mathcal C}$} \hfill \\
&x_{n,a} : X_{n,a} \longrightarrow X_{n-1 , f_n (a)} \, , \quad n \geq 1 \, , \quad a \in A_n \, ,
\end{matrix}\right.$

\medskip

tels que

\medskip

$\left\{\begin{matrix}
\bullet &\mbox{si $a_0$ est l'unique élément de $A_0$,} \hfill \\
&X_{0,a_0} = X \, , \\
\bullet &\mbox{pour tout n\oe ud $(n,a)$ de $A_{\bullet}$ qui n'est pas maximal, c'est-à-dire définit une fibre non vide} \hfill \\
&f_{n+1}^{-1} (a) \ne \emptyset \, , \\
&\mbox{la famille des morphismes} \hfill \\
&x_{n+1,a'} : X_{n+1,a'} \longrightarrow X_{n,a} \, , \quad a' \in f_{n+1}^{-1} (a) \, , \\
&\mbox{constitue un précrible $J$-couvrant.} \hfill
\end{matrix}\right.$

%\hfill $\Box$
\end{defn}

On observe que pour tout recouvrement d'un objet $X$ par des morphismes $x_i : X_i \to X$, $i \in I$, le choix de multi-recouvrements des objets $X_i$, $i \in I$, définit un multi-recouvrement de $X$:

\begin{lem}\label{lem3.2.4}

Soit ${\mathcal C}$ une catégorie essentiellement petite qui est munie d'une notion $J$ de ``précrible couvrant'' qui est supposée stable.

Considérons un précrible $J$-couvrant d'un objet $X$ de ${\mathcal C}$
$$
\left( x_i : X_i \longrightarrow X \right)_{i \in I}
$$
et, pour tout indice $i \in I$, un multi-recouvrement de la source $X_i$ de $x_i$ qui consiste en

\medskip

$\left\{\begin{matrix}
\bullet &\mbox{un arbre à branches finies $A_{\bullet}^i = \left( A_0^i \longleftarrow \cdots \longleftarrow A_{n-1}^i \overset{f_n^i}{\longleftarrow} A_n^i \longleftarrow \cdots \right)$,} \hfill \\
\bullet &\mbox{une famille d'objets de ${\mathcal C}$ indexés par les n\oe uds $(n,a)$ de l'arbre $A_{\bullet}^i$} \hfill \\
&X_{n,a}^i \, , \\
\bullet &\mbox{une famille de morphismes de ${\mathcal C}$} \hfill \\
&x_{n,a}^i : X_{n,a}^i \longrightarrow X^i_{n-1 , f^i_n (a)} \, , \quad n \geq 1 \, , \quad a \in A^i_n \, .
\end{matrix}\right.$

\medskip

Alors:

\begin{listeimarge}

\item On définit un arbre à branches finies $A_{\bullet}$ en posant
$$
\left\{\begin{matrix}
A_0 &= &\{\bullet \} \, , \hfill \\
A_n &= &\underset{i \in I}{\coprod} \, A_{n-1}^i \quad \mbox{si} \quad n \geq 1 \, , \hfill \\
f_1 &= &(A_1 \longrightarrow A_0 = \{\bullet \}) \, , \hfill \\
f_n &= &\underset{i \in I}{\coprod} \, f_{n-1}^i \quad \mbox{si} \quad n \geq 2 \, . \hfill
\end{matrix} \right.
$$
\item On définit un multi-recouvrement de $X$ indexé par l'arbre à branches finies $A_{\bullet}$ en posant

\medskip

$\left\{\begin{matrix}
\mbox{$X_{0,a_0} = X$ si $a_0$ est l'unique élément de $A_0$,} \hfill \\
\mbox{$X_{n,a} = X_{n-1,a}^i$ si $a \in A_{n-1}^i \subseteq A_n$, $n \geq 1$,} \hfill \\
\mbox{$x_{1,a} = (x_i : X_i \to X_0)$ si $a$ est l'unique élément de $A_0^i \subseteq A_1$,} \hfill \\
\mbox{$x_{n,a} = x_{n-1,a}^i : X_{n-1 , a}^i \to X_{n-2 , f_{n-1}^i (a)}^i$ si $a \in A_{n-1}^i \subseteq A_n$, $n \geq 2$.} \hfill
\end{matrix}\right.$
\end{listeimarge}
\end{lem}

\begin{remarks}
\begin{listeisansmarge}
\item La propriété de composition des précribles couvrants avec des multi-recouvrements pour former des multi-recouvrements ne serait pas vérifiée si on se limitait à des arbres $A_{\bullet}$ ``de hauteur bornée'' au sens que $A_n = \emptyset$ pour tout entier $n$ assez grand.

\smallskip

\item En revanche, si on part d'une notion de ``précrible couvrant'' dont tous les éléments
$$
(x_i : X_i \longrightarrow X)_{i \in I}
$$
sont des familles finies, alors tous les multi-recouvrements associés à une telle notion sont indexés par des arbres finis, donc de hauteur bornée.
\end{listeisansmarge}
\end{remarks}

\begin{demo}
\begin{listeisansmarge}
\item L'ensemble des n\oe uds de l'arbre $A_{\bullet}$ s'identifie à la somme disjointe des ensembles des n\oe uds des arbres $A_{\bullet}^i$, $i \in I$, complétée par le n\oe ud initial
$$
(0, a_0) \, .
$$
De plus, deux n\oe uds $(n,a)$ et $(m,b)$ de l'arbre $A_{\bullet}$ satisfont la relation
$$
(n,a) \geq (m,b)
$$
si et seulement si
$$
m = 0 \quad \mbox{et} \quad b = a_0
$$
ou bien si
$$
n \geq m \geq 1
$$
et qu'existe un indice $i \in I$ tel que
$$
a \in A_{n-1}^i \subseteq A_n \, , \quad b \in A_{m-1}^i \subseteq A_m
$$
et
$$
(n-1 , a) \geq (m-1 , b)
$$
pour la relation d'ordre entre les n\oe uds de l'arbre $A_{\bullet}^i$.

Toute partie non vide de l'ensemble des n\oe uds de l'arbre $A_{\bullet}$ est réduite au singleton
$$
\{(0, a_0)\}
$$
ou bien a une intersection non vide avec l'ensemble des n\oe uds de l'un au moins des arbres $A_{\bullet}^i$.

Comme une telle intersection non vide a nécessairement un élément maximal -- puisque chacun des arbres $A_{\bullet}^i$ est ``à branches finies'' --, toute partie non vide de l'ensemble des n\oe uds de $A_{\bullet}$ a un élément maximal.

Autrement dit, l'arbre $A_{\bullet}$ est ``à branches finies''.

%\smallskip

\item On sait déjà d'après (i) que $A_{\bullet}$ est un arbre ``à branches finies''. 

La fibre de $f_1 : A_1 \to A_0$ au-dessus de l'unique élément
$$
a_0 \in A_0
$$
est
$$
A_1 = \coprod_{i \in I} A_0^i = \coprod_{i \in I} \{\bullet\}
$$
et la famille correspondante de morphismes
$$
x_{1,a} : X_{1,a} \longrightarrow X_{0,a_0} \, , \quad a \in A_1 \, ,
$$
n'est autre que la famille
$$
(x_i : X_i \longrightarrow X)_{i \in I}
$$
qui est par hypothèse un précrible couvrant de l'objet $X$.

Si $n \geq 1$, tout élément
$$
a \in A_n = \coprod_{i \in I} A_{n-1}^i
$$
est élément de l'un des ensembles
$$
A_{n-1}^i \, ,
$$
et la fibre
$$
f_{n+1}^{-1} (a) \subseteq A_{n+1}
$$
se confond avec la fibre
$$
(f_n^i)^{-1} (a) \subseteq A_n^i \, .
$$
Si cette fibre est non vide, la famille de morphismes
$$
\left( x_{n+1 , a'} : X_{n+1,a'} \longrightarrow X_{n,a} \right)_{a' \in f_{n+1}^{-1} (a)}
$$
se confond avec la famille de morphismes
$$
\left( x^i_{n , a'} : X^i_{n,a'} \longrightarrow X^i_{n-1,a} \right)_{a' \in (f^i_{n})^{-1} (a)}
$$
qui est par hypothèse un précrible couvrant de l'objet
$$
X_{n-1,a}^i = X_{n,a} \, .
$$
Cela prouve (ii). 
\end{listeisansmarge}
\end{demo}

\subsection{La propriété de stabilité des multi-recouvrements}\label{ssec3.2.2}

\medskip

On rappelle qu'une notion $J$ de précrible couvrant sur une catégorie essentiellement petite ${\mathcal C}$, consistant en une collection de précribles appelés ``$J$-couvrants'', est dite ``stable'' si, pour tout morphisme de ${\mathcal C}$
$$
x : X' \longrightarrow X
$$
et tout précrible $J$-couvrant de son but $X$
$$
\left( X_i \xrightarrow{ \ x_i \ } X \right)_{i \in I} \, ,
$$
il existe un précrible $J$-couvrant de sa source $X'$
$$
\left( X'_j \xrightarrow{ \ x'_j \ } X' \right)_{j \in I'} 
$$
et une application
$$
\alpha : I' \longrightarrow I
$$
tel que chaque morphisme $X'_j \xrightarrow{ \ x'_j \ } X'$, $j \in I'$, s'inscrit dans un carré commutatif de ${\mathcal C}$:
$$
\xymatrix{
X'_j \ar[d]_{x'_j} \ar[r] &X_{\alpha(j)} \ar[d]^{x_{\alpha (j)}} \\
X' \ar[r]^x &X
}
$$

Montrons que, si une notion de précrible couvrant est stable, alors la notion associée de multi-recouvrement est également stable au sens de la proposition suivante:

\begin{prop}\label{prop3.2.5}

Considérons une catégorie essentiellement petite ${\mathcal C}$ munie d'une notion $J$ de précrible couvrant qui est stable.

Soient un morphisme de ${\mathcal C}$
$$
x : X' \longrightarrow X
$$
et un multi-recouvrement de son but $X$ de type $J$ consistant en

\smallskip

$\left\{\begin{matrix}
\bullet &\mbox{un arbre à branches finies} \hfill \\
&A_{\bullet} = \left( A_0 \longleftarrow \cdots \longleftarrow A_{n-1} \xleftarrow{ \ f_n \ } A_n \longleftarrow \cdots \right), \\
\bullet &\mbox{une famille d'objets de ${\mathcal C}$ indexés par les n\oe uds de cet arbre} \hfill \\
&X_{n,a} \, , \quad n \in {\mathbb N} \, , \ a \in A_n \, , \\
\bullet &\mbox{une famille de morphismes de ${\mathcal C}$} \hfill \\
&x_{n,a} : X_{n,a} \longrightarrow X_{n-1,f_n(a)} \, , \quad n \geq 1 \, , \ a \in A_n \, .
\end{matrix}\right.$

\smallskip

Alors il existe:

\begin{enumerate}[label=(\arabic*)]

\item un multi-recouvrement de type $J$ de la source $X'$ du morphisme $x$, consistant en

\smallskip

$\left\{\begin{matrix}
\bullet &\mbox{un arbre à branches finies} \hfill \\
&B_{\bullet} = \left( B_0 \longleftarrow \cdots \longleftarrow B_{n-1} \xleftarrow{ \ g_n \ } B_n \longleftarrow \cdots \right), \\
\bullet &\mbox{une famille d'objets de ${\mathcal C}$ indexés par les n\oe uds de cet arbre} \hfill \\
&X'_{n,b} \, , \quad n \in {\mathbb N} \, , \ b \in B_n \, , \\
\bullet &\mbox{une famille de morphismes de ${\mathcal C}$} \hfill \\
&x'_{n,b} : X'_{n,b} \longrightarrow X'_{n-1,g_n(b)} \, , \quad n \geq 1 \, , \ b \in B_n \, ,
\end{matrix}\right.$

%\medskip

\item un morphisme d'arbres
$$
B_{\bullet} \longrightarrow A_{\bullet}
$$
consistant en une famille d'applications
$$
\alpha_n : B_n \longrightarrow A_n \, , \quad n \in {\mathbb N} \, ,
$$
qui sont compatibles au sens que les carrés
$$
\xymatrix{
B_n \ar[d]_{g_n} \ar[r]^{\alpha_n} &A_n \ar[d]^{f_n} \\
B_{n-1} \ar[r]^{\alpha_{n-1}} &A_{n-1}
}
$$
sont commutatifs,

%\smallskip

\item une famille de morphismes de ${\mathcal C}$
$$
\widetilde x_{n,b} : X'_{n,b} \longrightarrow X_{n,\alpha_n (b)} \, , \quad n \in {\mathbb N} \, , \ b \in B_n \, ,
$$
telle que

\smallskip

$\left\{\begin{matrix}
\bullet &\mbox{si $n=0$, $b_0$ est l'unique élément de $B_0$ et $a_0 = \alpha_0 (b_0)$ est l'unique élément de $A_0$,} \hfill \\
&\left( X'_{0,b_0} \xrightarrow{ \ \widetilde x_{0,b_0} \ } X_{0,a_0} \right) = \left( X' \xrightarrow{ \ x \ } X \right), \\
\bullet &\mbox{si $n \geq 1$, $b$ est un élément de $B_{n-1}$ et} \hfill \\
&b' \in g_n^{-1} (b) \subseteq B_n \, , \\
&\mbox{le carré de ${\mathcal C}$} \hfill \\
&\xymatrix{
X'_{n,b'} \ar[d]_{x'_{n,b'}} \ar[rr]^{\widetilde x_{n,b'}} &&X_{n , \alpha_n (b')} \ar[d]^{x_{n , \alpha_n (b')}} \\
X'_{n-1 , b} \ar[rr]^{\widetilde x_{n-1,b}} &&X_{n-1 , \alpha_{n-1} (b)}
} \\
&\mbox{est commutatif.} \hfill
\end{matrix}\right.$
\end{enumerate}
\end{prop}

\begin{demo}

Nous allons construire par récurrence sur $n \in {\mathbb N}$ les ensembles d'indices

$$
B_n
$$
inscrits dans des carrés commutatifs

$$
\xymatrix{
B_n \ar[d]_{g_n} \ar[r]^{\alpha_n} &A_n \ar[d]^{f_n} \\
B_{n-1} \ar[r]^{\alpha_{n-1}} &A_{n-1}
}
$$
et pour chaque $a \in A_n$, des morphismes indexés par les $b \in \alpha_n^{-1} (a) \subseteq B_n$

$$
X'_{n,b} \xrightarrow{ \ \widetilde x_{n,b} \ } X_{n,a}
$$
et

$$
X'_{n,b} \xrightarrow{ \ x'_{n,b} \ } X_{n-1 , g_n (b)}
$$
tels que les carrés

$$
\xymatrix{
X'_{n,b} \ar[d]_{x'_{n,b}} \ar[rr]^{\widetilde x_{n,b}} &&X_{n , a} \ar[d]^{x_{n , a}} \\
X'_{n-1 , g_n(b)} \ar[rr]^{\widetilde x_{n-1, g_n(b)}} &&X_{n-1 , f_n (a)}
}
$$
soient commutatifs.

Pour $n=0$, on prend pour $B_0$ l'ensemble à un seul élément $b_0$ muni de l'unique application

$$
\alpha_0 : B_0 \longrightarrow \{ a_0 \} = A_0
$$
et on prend pour morphisme

$$
X'_{0,b_0} \xrightarrow{ \ \widetilde x_{0,b_0} \ } X_{0,a_0}
$$
le morphisme

$$
X' \xrightarrow{ \ x \ } X \, .
$$

Puis, pour $n \geq 1$, supposons déjà construit le diagramme commutatif
$$
\xymatrix{
B_{n-1} \ar[d]_{g_{n-1}} \ar[r]^{\alpha_{n-1}} &A_{n-1} \ar[d]^{f_{n-1}} \\
B_{n-2} \ar[d] \ar[r]^{\alpha_{n-1}} &A_{n-2} \ar[d] \\
{ \ } \ar[d] \ar @{.}[r] &{ \ } \ar[d] \\
B_0 \ar[r]^{\alpha_0} &A_0
}
$$
et les carrés commutatifs de ${\mathcal C}$
$$
\xymatrix{
X'_{m,b'} \ar[d]_{x'_{m,b'}} \ar[rr]^{\widetilde x_{m,b'}} &&X_{m , a'} \ar[d]^{x_{m , a'}} \\
X'_{m-1 , b} \ar[rr]^{\widetilde x_{m-1, b}} &&X_{m-1,a}
}
$$
indexés par les $m \in \{1 , \cdots , n-1\}$, $b' \in B_m$, $a' = \alpha_m (b')$, $b = g_m (b')$, $a = f_m (a')$.

Construire un carré commutatif
$$
\xymatrix{
B_n \ar[d]_{g_n} \ar[r]^{\alpha_n} &A_n \ar[d]^{f_n} \\
B_{n-1} \ar[r]^{\alpha_{n-1}} &A_{n-1}
}
$$
équivaut à construire, pour tout élément
$$
b \in B_{n-1} \quad \mbox{d'image} \quad \alpha_{n-1} (b) = a \in A_{n-1} \, ,
$$
un ensemble
$$
B_{n,b}
$$
muni d'une application
$$
\alpha_n : B_{n,b} \longrightarrow f_n^{-1} (a) \subseteq A_n \, .
$$
En effet, $B_n$ sera alors défini par la formule
$$
B_n = \coprod_{b \in B_{n-1}} B_{n,b} \, .
$$
Si la fibre de $f_n : A_n \to A_{n-1}$ au-dessus de $a \in \alpha_{n-1} (b)$
$$
f_n^{-1} (a)
$$
est vide, on pose
$$
B_{n,b} = \emptyset \, .
$$
Si au contraire cette fibre est non vide, la famille de morphismes
$$
\left( X_{n,a'} \xrightarrow{ \ x_{n,a'} \ } X_{n-1,a} \right)_{a' \in f_n^{-1} (a)}
$$
est un précrible $J$-couvrant.

Par hypothèse de stabilité de la notion de précrible couvrant par le morphisme
$$
\widetilde x_{n-1,b} : X'_{n-1,b} \longrightarrow X_{n-1,a} \, ,
$$
il existe un précrible $J$-couvrant indexé par un ensemble $B_{n,b}$
$$
\left( X'_{n,b'} \xrightarrow{ \ x'_{n,b'} \ } X'_{n-1,b} \right)_{b' \in B_{n,b}} \, ,
$$
une application
$$
\alpha_n : B_{n,b} \longrightarrow f_n^{-1} (a) \subseteq A_n \, ,
$$
et des morphismes
$$
\widetilde x_{n,b'} : X'_{n,b'} \longrightarrow X_{n,\alpha_n(b')} \, , \quad b' \in B_{n,b} \, ,
$$
pour lesquels les carrés
$$
\xymatrix{
X'_{n,b'} \ar[d]_{x'_{n,b'}} \ar[rr]^{\widetilde x_{n,b'}} &&X_{n , \alpha_n (b')} \ar[d]^{x_{n , \alpha_n (b')}} \\
X'_{n-1 , b} \ar[rr]^{\widetilde x_{n-1,b}} &&X_{n-1 , a}
}
$$
sont commutatifs.

Pour conclure, il reste seulement à vérifier que l'arbre $B_{\bullet}$ ainsi construit par récurrence est nécessairement à branches finies.

Cela résulte de ce qu'il est muni d'un morphisme d'arbres
$$
B_{\bullet} \longrightarrow A_{\bullet}
$$
vers l'arbre $A_{\bullet}$ qui est à branches finies.

En effet, l'application qui associe à tout n\oe ud de $B_{\bullet}$
$$
(n,b) \, , \quad b \in B_n \, ,
$$
le n\oe ud correspondant de $A_{\bullet}$
$$
(n, \alpha_n (b))
$$
respecte la relation d'ordre.

Pour toute partie non vide ${\mathcal B}$ de l'ensemble des n\oe uds de $B_{\bullet}$, son image ${\mathcal A}$ par cette application est non vide donc possède au moins un élément maximal
$$
(n,a) \, .
$$
Alors n'importe quel antécédent de cet élément dans ${\mathcal B}$
$$
(n,b)
$$
est un élément maximal de ${\mathcal B}$.

Cela termine la preuve de la proposition. 
\end{demo}

\subsection{Explicitation de l'opération de clôture des sous-préfaisceaux}\label{ssec3.2.3}

\medskip

Comme nous l'avons vu au chapitre \ref{chap2}, la relation de dualité
$$
R \xymatrix{\ar@{^{(}->}[r] & } T \times S
$$
entre la classe $T$ des cribles d'objets $X$ d'une catégorie essentiellement petite ${\mathcal C}$ et la classe $S$ des plongements $Q \hookrightarrow P$ dans un préfaisceau $P$ sur ${\mathcal C}$ d'un sous-préfaisceau $Q$ de $P$ engendre une paire adjointe d'applications respectant les relations d'ordre
$$
\xymatrix{ ({\mathcal P} (T) , \subseteq) \ar@<2pt>[r]^{F_R} & ({\mathcal P} (S) , \supseteq). \ar@<2pt>[l]^{G_R} }
$$
L'application $F_R$ associe à toute collection $J$ de cribles une notion de sous-préfaisceau $J$-fermé consistant en la sous-classe
$$
F_R (J) \subseteq {\mathcal P} (S)
$$
des sous-préfaisceaux $Q$ de préfaisceaux $P$ qui sont dits $J$-fermés. Cette notion définit à son tour une ``opération de $J$-fermeture''
$$
(Q \hookrightarrow P) \longmapsto (\overline Q \hookrightarrow P)
$$
qui associe à tout sous-préfaisceau $Q$ d'un préfaisceau $P$ le plus petit sous-préfaisceau $J$-fermé $\overline Q$ de $P$ qui contient $Q$.

Enfin, la notion de sous-préfaisceau $J$-fermé et l'opération correspondante de $J$-fermeture des sous-préfaisceaux sont les mêmes que celles définies par la topologie engendrée par $J$
$$
\overline J = G_R \circ F_R(J) \, .
$$

L'opération de $J$-fermeture des sous-préfaisceaux s'explicite de la manière suivante lorsque la notion de crible $J$-couvrant est définie par une notion ${\mathcal J}$ de précrible couvrant qui est stable au sens de la définition \ref{defn3.1.2} (ii):

\begin{thm}\label{thm3.2.6}

Soit une catégorie essentiellement petite ${\mathcal C}$.

Soit une notion ${\mathcal J}$ de précrible couvrant de ${\mathcal C}$ qui est stable.

Soit $J$ la notion associée de crible couvrant, définie en décidant qu'un crible d'un objet $X$ de ${\mathcal C}$
$$
C \xymatrix{\ar@{^{(}->}[r] & } y(X)
$$
est $J$-couvrant s'il est le crible maximal de $X$ ou s'il contient au moins un précrible élément de ${\mathcal J}$
$$
\left( X_i \xrightarrow{ \ x_i \ } X\right)_{i \in I} \, .
$$
Soit enfin $\overline J$ la topologie de ${\mathcal C}$ engendrée par $J$.

Alors, pour tout sous-préfaisceau d'un préfaisceau $P$ de ${\mathcal C}$
$$
Q \xymatrix{\ar@{^{(}->}[r] & } P
$$
et sa $J$-fermeture ou, ce qui est équivalent, sa $\overline J$-fermeture
$$
\overline Q \xymatrix{\ar@{^{(}->}[r] & } P \, ,
$$
une section du préfaisceau $P$ sur un objet $X$ de ${\mathcal C}$
$$
p \in P(X)
$$
est élément du sous-ensemble
$$
\overline Q (p) \subseteq P(X)
$$
si et seulement si il existe un multi-recouvrement de type ${\mathcal J}$ de l'objet $X$ consistant en

\medskip

$\left\{\begin{matrix}
\bullet &\mbox{un arbre à branches finies $A_{\bullet} = \left( A_0 \xleftarrow{ \, f_1 \, } A_1 \leftarrow \cdots \leftarrow A_{n-1} \xleftarrow{ \, f_n \, } A_n \leftarrow \cdots \right),$} \hfill \\
\bullet &\mbox{une famille d'objets de ${\mathcal C}$} \hfill \\
&X_{n,a} \, , \quad n \in {\mathbb N} \, , \quad a \in A_n \, , \\
\bullet &\mbox{une famille de morphismes de ${\mathcal C}$} \hfill \\
&x_{n,a} : X_{n,a} \longrightarrow X_{n-1 , f_n(a)} \, , \quad n \geq 1 \, , \quad a \in A_n \, ,
\end{matrix}\right.$

\medskip

tel que, pour tout $n \in {\mathbb N}$ et $a \in A_n$, on a l'alternative:

\medskip

$\left\{\begin{matrix}
\bullet &\mbox{Ou bien $a$ est élément de l'image de l'application $f_{n+1} : A_{n+1} \to A_n$,} \hfill \\
&\mbox{autrement dit le n\oe ud $(n,a)$ de $A_{\bullet}$ n'est pas maximal.} \hfill \\
\bullet &\mbox{Ou bien la famille vide de morphismes de but $X_{n,a}$ est élément de ${\mathcal J}$.} \hfill \\
\bullet &\mbox{Ou bien, notant} \hfill \\
&a = a_n \, , \quad f_n (a_n) = a_{n-1} , \cdots , f_1(a_1) = a_0 \, , \\
&\mbox{le morphisme composé} \hfill \\
&X_{n,a_n} \xrightarrow{ \ x_{n,a_n} \ } X_{n-1, a_{n-1}} \longrightarrow \cdots \longrightarrow X_{1,a_1} \xrightarrow{ \ x_{1,a_1} \ } X_{0,a_0} = X \\
&\mbox{envoie l'élément} \hfill \\
&p \in P(X) \\
&\mbox{dans le sous-ensemble} \hfill \\
&Q(X_{n,a_n}) \subseteq P(X_{n,a_n}) \, .
\end{matrix}\right.$
\end{thm}

\begin{remark}

Ce théorème s'applique en particulier si la catégorie essentiellement petite ${\mathcal C}$ possède des produits fibrés et si, partant d'une famille arbitraire ${\mathcal X}$ de précribles de ${\mathcal C}$, on prend pour ${\mathcal J}$ la collection stable des précribles de la forme
$$
(X_i \times_X X' \longrightarrow X')_{i \in I}
$$
déduits des précribles éléments de ${\mathcal X}$
$$
( X_i \longrightarrow X)_{i \in I}
$$
par changement de base par tous les morphismes
$$
X' \longrightarrow X \, .
$$
Dans ce cas, $\overline J$ est la topologie engendrée par la famille ${\mathcal X}$ de précribles de ${\mathcal C}$.
\end{remark}

\begin{demo}

Considérons donc un sous-préfaisceau d'un préfaisceau $P$ de ${\mathcal C}$
$$
Q \xymatrix{\ar@{^{(}->}[r] & } P \, .
$$
Pour tout objet $X$ de ${\mathcal C}$, notons
$$
\widetilde Q(X) \subseteq P(X)
$$
le sous-ensemble des éléments
$$
p \in P(X)
$$
qui satisfont la condition du théorème relativement à un multi-recouvrement $X_{\bullet}$ de type ${\mathcal J}$ de l'objet $X$ consistant en

\medskip

$\left\{\begin{matrix}
\bullet &\mbox{un arbre à branches finies $A_{\bullet}$,} \hfill \\
\bullet &\mbox{des objets $X_{n,a}$, $n \in {\mathbb N}$, $a \in A_n$, avec $X_{0,a_0} = X$,} \hfill \\
\bullet &\mbox{des morphismes $X_{n,a} \xrightarrow{ \ x_{n,a} \ } X_{n-1 , f_n(a)}$, $n \geq 1$, $a \in A_n$.} \hfill
\end{matrix}\right.$

\medskip

\noindent Montrons d'abord que $\widetilde Q$ est un sous-préfaisceau de $Q$, c'est-à-dire que tout morphisme
$$
x : X' \longrightarrow X
$$
envoie le sous-ensemble
$$
\widetilde Q (X) \subseteq P(X)
$$
dans le sous-ensemble
$$
\widetilde Q (X') \subseteq P(X') \, .
$$

Etant donné un multi-recouvrement $X_{\bullet}$ de type ${\mathcal J}$ de l'objet $X$ indexé par un arbre à branches finies $A_{\bullet}$, il existe d'après la proposition \ref{prop3.2.5} un multi-recouvrement $X'_{\bullet}$ de $X'$ de type ${\mathcal J}$ indexé par un arbre à branches finies $B_{\bullet}$, un morphisme d'arbres
$$
B_{\bullet} \longrightarrow A_{\bullet} = \left( B_n \xrightarrow{ \ \alpha_n \ } A_n \right)_{n \in {\mathbb N}}
$$
et une famille compatible de morphismes de ${\mathcal C}$
$$
\widetilde x_{n,b} : X'_{n,b} \longrightarrow X_{n,\alpha_n (b)} \, , \quad n \in {\mathbb N} \, , \quad b \in B_n \, ,
$$
dont l'unique composante de degre $n=0$
$$
\widetilde x_{0,b_0} : X'_{0,b_0} \longrightarrow X_{0,a_0}
$$
est le morphisme
$$
x : X' \longrightarrow X \, .
$$
Considérons un n\oe ud arbitraire de l'arbre $B_{\bullet}$
$$
(n,b)
$$
et son image le n\oe ud de l'arbre $A_{\bullet}$
$$
(n,a)
$$
défini par $a = \alpha_n (b)$.

Si le n\oe ud $(n,a)$ de $A_{\bullet}$ est maximal, alors a fortiori le n\oe ud $(n,b)$ de $B_{\bullet}$ est maximal.

Si la famille vide de morphismes de but $X_{n,a}$ est ${\mathcal J}$-couvrante, alors l'existence du morphisme
$$
\widetilde x_{n,b} : X'_{n,b} \longrightarrow X_{n,a}
$$
et la propriété de stabilité de ${\mathcal J}$ assurent que la famille vide de morphismes de but $X'_{n,b}$ est ${\mathcal J}$-couvrante.

Enfin, notant
$$
a = a_n , a_{n-1} , \cdots , a_1 , a_0
$$
les images successives de l'élément $a \in A_n$ dans l'arbre $A_{\bullet}$, et
$$
b = b_n , b_{n-1} , \cdots , b_1 , b_0
$$
les images successives de l'élément $b \in B_n$ dans l'arbre $B_{\bullet}$, la commutativité du diagramme
$$
\xymatrix{
X'_{n,b_n} \ar[d] \ar[rr]^{x'_{n,b_n}} &&\cdots \ar[r] &X'_{1,b_1} \ar[d] \ar[r] &X' \ar[d]^x \\
X_{n,a_n} \ar[rr]^{x_{n,a_n}} &&\cdots \ar[r] &X_{1,a_1} \ar[r] &X
}
$$
assure que si le composé
$$
X_{n,a_n} \longrightarrow \cdots \longrightarrow X_{1,a_1} \longrightarrow X
$$
envoie l'élément $p \in P(X)$ dans le sous-ensemble $Q(X_{n,a_n}) \subseteq P(X_{n,a_n})$, alors le composé
$$
X'_{n,b_n} \longrightarrow \cdots \longrightarrow X'_{1,b_1} \longrightarrow X'
$$
envoie l'élément $x^*(p) \in P(X')$ dans le sous-ensemble $Q(X'_{n,b_n}) \subseteq P(X'_{n,b_n})$.

Cela démontre comme voulu que la famille des sous-ensembles
$$
\widetilde Q (X) \subseteq P(X)
$$
définit un sous-préfaisceau de $P$
$$
\widetilde Q \xymatrix{\ar@{^{(}->}[r] & } P \, .
$$

Considérant pour tout objet $X$ de ${\mathcal C}$ son multi-recouvrement trivial réduit au seul objet $X_{0,a_0} = X$, on voit que le sous-préfaisceau
$$
\widetilde Q \xymatrix{\ar@{^{(}->}[r] & } P
$$
contient le sous-préfaisceau
$$
Q \xymatrix{\ar@{^{(}->}[r] & } P \, .
$$

Puis montrons que le sous-préfaisceau
$$
\widetilde Q \xymatrix{\ar@{^{(}->}[r] & } P
$$
est $J$-fermé.

Pour cela, il suffit de montrer que pour tout élément
$$
p \in P(X)
$$
tel qu'existe un précrible ${\mathcal J}$-couvrant
$$
(x_i : X_i \longrightarrow X)_{i \in I}
$$
dont chaque élément
$$
x_i : X_i \longrightarrow X \, , \quad i \in I \, ,
$$
envoie $p$ dans le sous-ensemble
$$
\widetilde Q (X_i) \subseteq P(X_i) \, ,
$$
alors $p$ est lui-même élément du sous-ensemble
$$
\widetilde Q (X) \subseteq P(X) \, .
$$
Par hypothèse, il existe pour tout indice $i \in I$ un multi-recouvrement de type ${\mathcal J}$ de l'objet $X_i$ consistant en

\medskip

$\left\{\begin{matrix}
\bullet &\mbox{un arbre à branches finies $A_{\bullet}^i = \left( A_0^i \longleftarrow A_1^i \longleftarrow \cdots \longleftarrow A_{n-1}^i \xleftarrow{ \ f_n^i \ } A_n^i \longleftarrow \cdots \right)$,} \hfill \\
\bullet &\mbox{des objets $X_{n,a}^i$, $n \in {\mathbb N}$, $a \in A_n^i$, avec $X_{0,a_0}^i = X_i$,} \hfill \\
\bullet &\mbox{des morphismes $X_{n,a}^i \xrightarrow{ \ x_{n,a}^i \ } X_{n-1 , f_n^i (a)}^i$, $n \geq 1$, $a \in A_n^i$,} \hfill
\end{matrix}\right.$

\medskip

\noindent qui satisfait la condition du théorème appliquée à l'élément
$$
x_i^* (p) \in P(X_i) \, .
$$
Alors la construction du lemme III.2.4 permet d'amalgamer les arbres $A_{\bullet}^i$ en un arbre à branches finies $A_{\bullet}$ et les multi-recouvrements $X_{\bullet}^i$ des $X_i$, $i \in I$, indexés par les arbres $A_{\bullet}^i$ en un multi-recouvrement $X_{\bullet}$ de type ${\mathcal J}$ de l'objet $X$ indexé par l'arbre à branches finies $A_{\bullet}$.

De plus, ce multi-recouvrement $X_{\bullet}$ de $X$ satisfait la condition du théorème relativement à l'élément
$$
p \in P(X) \, .
$$
Cela montre comme voulu que le sous-préfaisceau
$$
\widetilde Q \xymatrix{\ar@{^{(}->}[r] & } P
$$
est $J$-fermé ou, ce qui revient au même, $\overline J$-fermé.

Pour conclure, il reste à démontrer que tout élément
$$
p \in \widetilde Q (X) \subseteq P(X) \, ,
$$
c'est-à-dire qui satisfait la condition du théorème relativement à un multi-recouvrement $X_{\bullet}$ de type ${\mathcal J}$ de $X$ indexé par un arbre à branches finies $A_{\bullet}$, est aussi élément du sous-ensemble
$$
\overline Q(X) \subseteq P(X) \, .
$$
Considérons l'ensemble ${\mathcal A}$ des n\oe uds de l'arbre $A_{\bullet}$
$$
(n,a) \, , 
$$
engendrant une branche
$$
a = a_n \, , \quad f_n (a_n) = a_{n-1} , \cdots , f_2(a_2) = a_1 , f_1 (a_1) = a_0 \, ,
$$
tels que le morphisme composé
$$
X_{n,a_n} \xrightarrow{ \ x_{n,a_n} \ } X_{n-1 , a_{n-1}} \longrightarrow \cdots \longrightarrow X_{1,a_0} \xrightarrow{ \ x_{1,a_1} \ } X_{0,a_0} = X
$$
n'envoie pas l'élément
$$
p \in \widetilde Q (X) \subseteq P(X)
$$
dans le sous-ensemble
$$
\overline Q (X_{n,a_n}) \subseteq P(X_{n,a_n}) \, .
$$
Il s'agit de montrer que cet ensemble
$$
{\mathcal A}
$$
est vide.

S'il n'était pas vide, et comme l'arbre $A_{\bullet}$ et à branches finies, il contiendrait un élément maximal
$$
(n,a) \, , \quad \mbox{avec} \quad a \in A_n \, .
$$
Par définition de ${\mathcal A}$, le morphisme composé
$$
X_{n,a_n} \xrightarrow{ \ x_{n,a_n} \ } X_{n-1 , a_{n-1}} \longrightarrow \cdots \longrightarrow X_{1,a_1} \xrightarrow{ \ x_{1,a_1} \ } X_{0,a_0} = X
$$
n'envoie pas l'élément $p \in P(X)$ dans le sous-ensemble
$$
\overline Q (X_{n,a_n}) \subseteq P(X_{n,a_n}) 
$$
ni a fortiori dans le sous-ensemble
$$
Q(X_{n,a_n}) \subseteq \overline Q (X_{n,a_n}) \, .
$$
Si la fibre
$$
f_{n+1}^{-1} (a_n) = \{a_{n+1} \in A_{n+1} \mid f_{n+1} (a_{n+1}) = a_n \}
$$
est non vide, les morphismes
$$
x_{n+1 , a_{n+1}} : X_{n+1 , a_{n+1}} \longrightarrow X_{n,a_n} \, , \quad a_{n+1} \in f_{n+1}^{-1} (a_n) \, ,
$$
constituent un précrible ${\mathcal J}$-couvrant de $X_{n,a_n}$ et, par maximalité du n\oe ud $(n,a_n)$ dans la partie ${\mathcal A}$, tous les composés
$$
X_{n+1 , a_{n+1}} \xrightarrow{ \ x_{n+1 , a_{n+1}} \ } X_{n,a_n} \longrightarrow \cdots \longrightarrow X_{0,a_0} = X
$$
envoient l'élément $p \in P(X)$ dans le sous-ensemble
$$
\overline Q (X_{n+1 , a_{n+1}}) \subseteq P(X_{n+1 , a_{n+1}}) \, .
$$
Comme le sous-préfaisceau $\overline Q \hookrightarrow P$ est $J$-fermé, cela contredit l'hypothèse que le n\oe ud $(n,a_n)$ est élément de la partie ${\mathcal A}$, c'est-à-dire que le composé
$$
X_{n,a_n} \xrightarrow{ \ x_{n,a_n} \ } X_{n-1 , a_{n-1}} \longrightarrow \cdots \longrightarrow X_{0,a_0} = X
$$
n'envoie pas l'élément $p \in P(X)$ dans le sous-ensemble
$$
\overline Q (X_{n,a_n}) \subseteq P(X_{n,a_n}) \, .
$$
Si au contraire la fibre
$$
f_{n+1}^{-1} (a_n)
$$
est vide, alors le précrible vide sur l'objet $X_{n,a_n}$ est ${\mathcal J}$-couvrant. Cela impose que le composé
$$
X_{n,a_n} \xrightarrow{ \ x_{n,a_n} \ } X_{n-1 , a_{n-1}} \longrightarrow \cdots \longrightarrow X_{0,a_0} = X
$$
envoie l'élément
$$
p \in P(X)
$$
dans le sous-ensemble
$$
\overline Q (X_{n,a_n}) \subseteq P(X_{n,a_n}) \, ,
$$
ce qui contredit à nouveau l'hypothèse que le n\oe ud $(n,a_n)$ est élément de la partie ${\mathcal A}$.

En définitive, la seule possibilité est que la partie ${\mathcal A}$ est vide.

En particulier, le n\oe ud $(0,a_0)$ n'est pas élément de ${\mathcal A}$, ce qui signifie que l'élément
$$
p \in \widetilde Q(X) \subseteq P(X)
$$
est dans le sous-ensemble
$$
\overline Q (X) \subseteq P(X) \, .
$$
On a montré comme annoncé que
$$
\widetilde Q = \overline Q \, .
$$
\end{demo}

\subsection{Explicitation de l'opération d'engendrement des topologies}\label{ssec3.2.4}

\medskip

Le théorème \ref{thm3.2.6} du sous-paragraphe précédent se spécialise en particulier en une formule qui caractérise en termes de multi-recouvrements les cribles couvrants d'une topologie engendrée par une notion stable de précrible couvrant:

\begin{cor}\label{cor3.2.7}

Soit une catégorie essentiellement petite ${\mathcal C}$.

Soit une notion ${\mathcal J}$ de précrible couvrant de ${\mathcal C}$ qui est stable.

Soit $J$ la notion associée de crible couvrant, définie en décidant qu'un crible d'un objet $X$ de ${\mathcal C}$ est $J$-couvrant s'il est le crible maximal de $X$ ou s'il contient au moins un précrible élément de ${\mathcal J}$.

Soit enfin $\overline J$ la topologie de ${\mathcal C}$ engendrée par $J$. Alors un crible d'un objet $X$ de ${\mathcal C}$
$$
C \xymatrix{\ar@{^{(}->}[r] & } y(X)
$$
est couvrant pour la topologie engendrée $\overline J$ si et seulement si il existe un multi-recouvrement de type ${\mathcal J}$ de l'objet $X$ consistant en

\medskip

$\left\{\begin{matrix}
\bullet &\mbox{un arbre à branches finies $A_{\bullet} = \left( A_0 \xleftarrow{ \, f_1 \, } A_1 \leftarrow \cdots \leftarrow A_{n-1} \xleftarrow{ \, f_n \, } A_n \leftarrow \cdots \right),$} \hfill \\
\bullet &\mbox{une famille d'objets de ${\mathcal C}$} \hfill \\
&X_{n,a_n} \, , \quad n \in {\mathbb N} \, , \quad a \in A_n \, , \quad \mbox{avec $X_{0,a_0} = X$,} \\
\bullet &\mbox{une famille de morphismes de ${\mathcal C}$} \hfill \\
&x_{n,a} : X_{n,a} \longrightarrow X_{n-1 , f_n(a)} \, , \quad n \geq 1 \, , \quad a \in A_n \, ,
\end{matrix}\right.$

\bigskip

tel que, pour tout $n \in {\mathbb N}$ et tout $a \in A_n$, on a l'alternative:

\bigskip

$\left\{\begin{matrix}
\bullet &\mbox{Ou bien $a$ est élément de l'image de l'application $f_{n+1} : A_{n+1} \to A_n$,} \hfill \\
&\mbox{autrement dit le n\oe ud $(n,a)$ de $A_{\bullet}$ n'est pas maximal.} \hfill \\
\bullet &\mbox{Ou bien la famille vide de morphismes de but $X_{n,a}$ est élément de ${\mathcal J}$.} \hfill \\
\bullet &\mbox{Ou bien, notant} \hfill \\
&a = a_n \, , \quad f_n (a_n) = a_{n-1} , \cdots , f_1(a_1) = a_0 \, , \\
&\mbox{le morphisme composé} \hfill \\
&X_{n,a_n} \xrightarrow{ \ x_{n,a_n} \ } X_{n-1, a_{n-1}} \longrightarrow \cdots \longrightarrow X_{1,a_1} \xrightarrow{ \ x_{1,a_1} \ } X_{0,a_0} = X \\
&\mbox{est élément du crible $C$.} \hfill 
\end{matrix}\right.$
\end{cor}

\begin{remarks}
\begin{listeisansmarge}
\item Dans ces conditions, un précrible sur un objet $X$
$$
\left( X_i \xrightarrow{ \ x_i \ } X \right)_{i \in I}
$$
est couvrant pour la topologie engendrée $\overline J$ si et seulement si $X$ admet un multi-recouvrement $X_{\bullet}$ de type ${\mathcal J}$ indexé par un arbre à branches finies $A_{\bullet}$ dont tout n\oe ud maximal $(n,a)$ satisfait l'alternative:

\medskip

$\left\{\begin{matrix}
\bullet &\mbox{Ou bien le précrible vide sur l'objet $X_{n,a}$ est élément de ${\mathcal J}$.} \hfill \\
\bullet &\mbox{Ou bien le morphisme composé associé à la branche} \hfill \\
&\mbox{$a = a_n \, , \quad f_n (a_n) = a_{n-1} , \cdots , a_1 , f_1(a_1) = a_0  \, ,$ du n\oe ud $(n,a)$} \\
&X_{n,a_n} \xrightarrow{ \ x_{n,a_n} \ } X_{n-1, a_{n-1}} \longrightarrow \cdots \longrightarrow X_{1,a_1} \xrightarrow{ \ x_{1,a_1} \ } X_{0,a_0} = X \\
&\mbox{se factorise à travers l'un au moins des morphismes $X_i \xrightarrow{ \ x_i \ } X$, $i \in I$.} \hfill 
\end{matrix}\right.$

%\medskip

\item On rappelle que si tous les précribles éléments de ${\mathcal J}$ sont finis, alors tous les multi-recouvrements $X_{\bullet}$ de type ${\mathcal J}$ d'objets $X$ de ${\mathcal C}$ sont finis, c'est-à-dire indexés par des arbres $A_{\bullet}$ finis.

%\smallskip

\item Comme le théorème \ref{thm3.2.6}, ce corollaire s'applique en particulier si la catégorie essentiellement petite ${\mathcal C}$ possède des produits fibrés et si, partant d'une famille arbitraire ${\mathcal X}$ de précribles de ${\mathcal C}$, on prend pour ${\mathcal J}$ la collection stable des précribles qui se déduisent des précribles éléments de ${\mathcal X}$ par changement de base.

Dans ce cas, si tous les précribles éléments de ${\mathcal X}$ sont finis, il en va de même des précribles éléments de ${\mathcal J}$.
\end{listeisansmarge}
\end{remarks}

\begin{demo}

La notion stable $J$ de crible couvrant ou, ce qui revient au même, la topologie engendrée $\overline J$ définissent une notion de sous-préfaisceau fermé des préfaisceaux sur ${\mathcal C}$ et une opération de fermeture des sous-préfaisceaux
$$
(Q \hookrightarrow P) \longmapsto (\overline Q \hookrightarrow P)
$$
à laquelle s'applique le théorème \ref{thm3.2.6}.

Or, un crible d'un objet $X$ de ${\mathcal C}$ vu comme un sous-préfaisceau
$$
C \xymatrix{\ar@{^{(}->}[r] & } y(X)
$$
est couvrant pour la topologie engendrée $\overline J$ si et seulement si sa fermeture $\overline C$ est le crible maximal
$$
\overline C = y(X)
$$
engendrée par le morphisme ${\rm id}_X : X \to X$.

Le corollaire s'obtient en appliquant le théorème \ref{thm3.2.6} à l'élément
$$
{\rm id} _X \in {\rm Hom} (X,X) = y(X)(X) \, .
$$
\end{demo}

\section{Traductions topologiques des problèmes de démontrabilité}\label{sec3.3}

\subsection{Interprétation topologique des propriétés géométriques}\label{ssec3.3.1}

On a rappelé au sous-paragraphe \ref{ssec1.1.3} que toute théorie géométrique du premier ordre ${\mathbb T}$ possède un ``topos classifiant'' ${\mathcal E}_{\mathbb T}$ muni d'un ``modèle universel'' $U_{\mathbb T}$ de la théorie ${\mathbb T}$.

Toute sorte $A$ de la signature de ${\mathbb T}$ s'interprète comme un objet du topos ${\mathcal E}_{\mathbb T}$ noté
$$
U_{\mathbb T} \, A \, ,
$$
tout symbole de fonction de cette signature $f : A_1 \cdots A_n \to B$ s'interprète comme un morphisme de ${\mathcal E}_{\mathbb T}$
$$
U_{\mathbb T} \, f : U_{\mathbb T} \, A_1 \times \cdots \times U_{\mathbb T} \, A_n \longrightarrow U_{\mathbb T} \, B \, ,
$$
et tout symbole de relation $R \rightarrowtail A_1 \cdots A_n$ s'interprète comme un sous-objet
$$
U_{\mathbb T} \, R \ {\lhook\joinrel\relbar\joinrel\rightarrow} \ U_{\mathbb T} \, A_1 \times \cdots \times U_{\mathbb T} \, A_n \, .
$$

A partir des interprétations des éléments de la signature de ${\mathbb T}$ se construisent des interprétations de toutes les formules géométriques de ${\mathbb T}$
$$
\varphi (\vec x)
$$
en des suites de variables $\vec x = (x_1^{A_1} , \cdots , x_n^{A_n})$ associées à des sortes $A_1 , \cdots , A_n$, comme des sous-objets dans ${\mathcal E}_{\mathbb T}$
$$
U_{\mathbb T} \, \varphi (\vec x) \ {\lhook\joinrel\relbar\joinrel\rightarrow} \ U_{\mathbb T} \, A_1 \times \cdots \times U_{\mathbb T} \, A_n \, .
$$
On a vu dans la proposition \ref{prop1.1.30} (iii) qu'un séquent géométrique de la signature de ${\mathbb T}$
$$
\varphi (\vec x) \vdash \psi (\vec x)
$$
en des variables $\vec x = (x_1^{A_1} , \cdots , x_n^{A_n})$ est démontrable dans ${\mathbb T}$ si et seulement si les interprétations comme sous-objets des deux formules qui composent le séquent
$$
U_{\mathbb T} \, \varphi \ {\lhook\joinrel\relbar\joinrel\rightarrow} \ U_{\mathbb T} \, A_1 \times \cdots \times U_{\mathbb T} \, A_n \, ,
$$
$$
U_{\mathbb T} \, \psi \ {\lhook\joinrel\relbar\joinrel\rightarrow} \ U_{\mathbb T} \, A_1 \times \cdots \times U_{\mathbb T} \, A_n 
$$
sont reliées par la relation d'inclusion
$$
U_{\mathbb T} \, \varphi \subseteq U_{\mathbb T} \, \psi \, .
$$

Comme tout topos, le topos classifiant ${\mathcal E}_{\mathbb T}$ d'une telle théorie ${\mathbb T}$ se représente de multiples façons comme un topos de faisceaux
$$
\widehat{\mathcal C}_J \xrightarrow{ \ \sim \ } {\mathcal E}_{\mathbb T}
$$
sur une petite catégorie ${\mathcal C}$ munie d'une topologie de Grothendieck $J$.

On sait que le foncteur canonique associé à une telle représentation
$$
\ell : {\mathcal C} \longrightarrow \widehat{\mathcal C}_J \xrightarrow{ \ \sim \ } {\mathcal E}_{\mathbb T}
$$
respecte les limites finies, en particulier les produits finis d'objets, les plongements de sous-objets et leurs intersections quand ils sont bien définis dans ${\mathcal C}$.

Il est naturel et commode de poser la définition suivante:

\begin{defn}\label{defn3.3.1}
Soit ${\mathbb T}$ une théorie géométrique du premier ordre de signature $\Sigma$.

Soit une présentation de son topos classifiant ${\mathcal E}_{\mathbb T}$ comme topos des faisceaux sur une petite catégorie ${\mathcal C}$ munie d'une topologie $J$
$$
\widehat{\mathcal C}_J \xrightarrow{ \ \sim \ } {\mathcal E}_{\mathbb T} \, .
$$
\begin{listeimarge}

\item On dit qu'une sorte $A$ de $\Sigma$ est interprétable dans ${\mathcal C}$ si son interprétation comme objet de ${\mathcal E}_{\mathbb T}$
$$
U_{\mathbb T} \, A
$$
est image par le foncteur canonique
$$
\ell : {\mathcal C} \longrightarrow \widehat{\mathcal C}_J \xrightarrow{ \ \sim \ } {\mathcal E}_{\mathbb T}
$$
d'un objet de ${\mathcal C}$ noté alors
$$
U_{\mathbb T}^{\mathcal C} \, A \, .
$$
\item On dit qu'un symbole de fonction de la signature de ${\mathbb T}$
$$
f : A_1 \cdots A_n \longrightarrow B
$$
est interprétable dans ${\mathcal C}$ si

\medskip

$\left\{\begin{matrix}
\bullet &\mbox{les sortes $A_1 , \cdots , A_n,B$ sont interprétables dans ${\mathcal C}$,} \hfill \\
\bullet &\mbox{le produit $U_{\mathbb T}^{\mathcal C} \, A_1 \times \cdots \times U_{\mathbb T}^{\mathcal C} \, A_n$ des interprétations $U_{\mathbb T}^{\mathcal C} \, A_i$ des sortes $A_i$, $1 \leq i \leq n$,} \hfill \\
&\mbox{est bien défini dans ${\mathcal C}$,} \hfill \\
\bullet &\mbox{le morphisme d'interprétation de $f$ dans ${\mathcal E}_{\mathbb T}$} \hfill \\
&U_{\mathbb T} \, f : U_{\mathbb T} \, A_1 \times \cdots \times U_{\mathbb T} \, A_n \longrightarrow U_{\mathbb T} \, B \\
&\mbox{est image par le foncteur canonique $\ell$ d'un morphisme de ${\mathcal C}$ de la forme} \hfill \\
&U_{\mathbb T}^{\mathcal C} \, f : U_{\mathbb T}^{\mathcal C} \, A_1 \times \cdots \times U_{\mathbb T}^{\mathcal C} \, A_n \longrightarrow U_{\mathbb T}^{\mathcal C} \, B \, .
\end{matrix} \right.$

\medskip

\item On dit qu'un symbole de relation de la signature de ${\mathbb T}$
$$
R \ {>\joinrel\relbar\joinrel\longrightarrow} \ A_1 \cdots A_n
$$
[resp. une formule géométrique en des variables $\vec x = (x_1^{A_1} , \cdots , x_n^{A_n})$
$$
\varphi (\vec x) \ \mbox{]}
$$
est interprétable dans ${\mathcal C}$ si

\medskip

$\left\{\begin{matrix}
\bullet &\mbox{les sortes $A_1 , \cdots , A_n$ sont interprétables dans ${\mathcal C}$,} \hfill \\
\bullet &\mbox{le produit $U_{\mathbb T}^{\mathcal C} \, A_1 \times \cdots \times U_{\mathbb T}^{\mathcal C} \, A_n$ des interprétations des $A_i$, $1 \leq i \leq n$, est bien défini dans ${\mathcal C}$,} \hfill \\
\bullet &\mbox{le sous-objet de ${\mathcal E}_{\mathbb T}$} \hfill \\
&U_{\mathbb T} \, R \ {\lhook\joinrel\relbar\joinrel\rightarrow} \ U_{\mathbb T} \, A_1 \times \cdots \times U_{\mathbb T} \, A_n \\
&\mbox{[resp. $U_{\mathbb T} \, \varphi \ {\lhook\joinrel\relbar\joinrel\rightarrow} \ U_{\mathbb T} \, A_1 \times \cdots \times U_{\mathbb T} \, A_n$]} \\
&\mbox{est image par le foncteur canonique $\ell$ d'un sous-objet de ${\mathcal C}$} \hfill \\
&U_{\mathbb T}^{\mathcal C} \, R \ {\lhook\joinrel\relbar\joinrel\rightarrow} \ U_{\mathbb T}^{\mathcal C} \, A_1 \times \cdots \times U_{\mathbb T}^{\mathcal C} \, A_n \\
&\mbox{[resp. $U_{\mathbb T}^{\mathcal C} \, \varphi \ {\lhook\joinrel\relbar\joinrel\rightarrow} \ U_{\mathbb T}^{\mathcal C} \, A_1 \times \cdots \times U_{\mathbb T}^{\mathcal C} \, A_n$].} 
\end{matrix} \right.$

\medskip

\item On dit qu'un séquent géométrique de la signature de ${\mathbb T}$ en des variables $\vec x = (x_1^{A_1} , \cdots , x_n^{A_n})$
$$
\varphi (\vec x) \vdash \psi (\vec x)
$$
est interprétable dans ${\mathcal C}$ si
\medskip

$\left\{\begin{matrix}
\bullet &\mbox{les sortes $A_1 , \cdots , A_n$ sont interprétables dans ${\mathcal C}$,} \hfill \\
\bullet &\mbox{les formules $\varphi$ et $\psi$ qui composent ce séquent sont interprétables dans ${\mathcal C}$ comme des sous-objets} \hfill \\
&U_{\mathbb T}^{\mathcal C} \, \varphi \ {\lhook\joinrel\relbar\joinrel\rightarrow} \ U_{\mathbb T}^{\mathcal C} \, A_1 \times \cdots \times U_{\mathbb T}^{\mathcal C} \, A_n \, , \\
&U_{\mathbb T}^{\mathcal C} \, \psi \ {\lhook\joinrel\relbar\joinrel\rightarrow} \ U_{\mathbb T}^{\mathcal C} \, A_1 \times \cdots \times U_{\mathbb T}^{\mathcal C} \, A_n \, , \\
\bullet &\mbox{l'intersection de ces deux sous-objets de $U_{\mathbb T}^{\mathcal C} \, A_1 \times \cdots \times U_{\mathbb T}^{\mathcal C} \, A_n$} \hfill \\
&U_{\mathbb T}^{\mathcal C} \, \varphi \quad \mbox{et} \quad U_{\mathbb T}^{\mathcal C} \, \psi \\
&\mbox{est bien définie dans ${\mathcal C}$.} \hfill
\end{matrix} \right.$
\end{listeimarge}
\end{defn}

On obtient comme cas particulier du corollaire \ref{cor1.2.15}:

\begin{cor}\label{cor3.3.2}
Soit une théorie géométrique du premier ordre ${\mathbb T}$ dont le topos classifiant ${\mathcal E}_{\mathbb T}$ est présenté comme topos des faisceaux sur une petite catégorie ${\mathcal C}$ munie d'une topologie $J$
$$
\widehat{\mathcal C}_J \xrightarrow{ \ \sim \ } {\mathcal E}_{\mathbb T} \, .
$$

Considérons des séquents géométriques de la signature de ${\mathbb T}$
$$
\varphi_i \vdash \psi_i
$$
qui sont interprétables dans ${\mathcal C}$, ainsi qu'un séquent géométrique
$$
\varphi \vdash \psi \, .
$$

Soit ${\mathbb T}'$ la théorie quotient de ${\mathbb T}$ définie par les axiomes additionnels
$$
\varphi_i \vdash \psi_i \, , \quad i \in I \, .
$$

Alors le séquent
$$
\varphi \vdash \psi
$$
est démontrable dans la théorie quotient ${\mathbb T}'$ si et seulement si le monomorphisme de ${\mathcal C}$
$$
U_{\mathbb T}^{\mathcal C} \, \varphi \cap U_{\mathbb T}^{\mathcal C} \, \psi \ {\lhook\joinrel\relbar\joinrel\rightarrow} \ U_{\mathbb T}^{\mathcal C} \, \varphi 
$$
est couvrant pour la topologie $J'$ de ${\mathcal C}$ engendrée par $J$ et les monomorphismes
$$
U_{\mathbb T}^{\mathcal C} \, \varphi_i \cap U_{\mathbb T}^{\mathcal C} \, \psi_i \ {\lhook\joinrel\relbar\joinrel\rightarrow} \ U_{\mathbb T}^{\mathcal C} \, \varphi_i \, , \quad i \in I \, .
$$
\end{cor}

\subsection{Réduction aux interprétations topologiques des formules de Horn}\label{ssec3.3.2}

Dans la pratique, lorsque le topos classifiant ${\mathcal E}_{\mathbb T}$ d'une théorie géométrique du premier ordre ${\mathbb T}$ est présenté comme le topos des faisceaux sur une petite catégorie ${\mathcal C}$ munie d'une topologie $J$, il est rare que toutes les formules géométriques de la signature de ${\mathbb T}$ soient interprétables dans la catégorie ${\mathcal C}$ au sens de la définition \ref{defn3.3.1}.

Cependant, des formules complexes non interprétables dans une telle catégorie de présentation ${\mathcal C}$ peuvent souvent se ramener à des familles de formules plus simples interprétables, du fait que les symboles logiques de disjonction $\bigvee$ et de quantification existentielle $\exists$ s'interprètent en termes de recouvrements:

\begin{lem}\label{lem3.3.3}
Soient une théorie géométrique du premier ordre ${\mathbb T}$ de signature $\Sigma$, son topos classifiant ${\mathcal E}_{\mathbb T}$ et le modèle universel $U_{\mathbb T}$ de ${\mathbb T}$ dans ${\mathcal E}_{\mathbb T}$.

Considérons une formule géométrique de $\Sigma$ en une suite de variables $\vec x = (x_1^{A_1} , \cdots , x_n^{A_n})$ qui est écrite sous la forme
$$
\varphi (\vec x) = \bigvee_{k \in K} (\exists \, \vec x_k) \, \varphi_k (\vec x , \vec x_k) \, .
$$
On a:
\begin{listeimarge}
\item Dans le topos classifiant ${\mathcal E}_{\mathbb T}$, la famille des morphismes de projection
$$
U_{\mathbb T} \, \varphi_k (\vec x , \vec x_k) \longrightarrow U_{\mathbb T} \, \varphi (\vec x) \, , \quad k \in K \, ,
$$
est globalement épimorphique.

%\smallskip

\item Supposons que le topos classifiant ${\mathcal E}_{\mathbb T}$ est présenté comme topos des faisceaux
$$
\widehat{\mathcal C}_J \xrightarrow{ \ \sim \ } {\mathcal E}_{\mathbb T}
$$
sur un petit site $({\mathcal C},J)$ tel que toutes les formules
$$
\varphi_k (\vec x , \vec x_k) \, , \quad k \in K \, ,
$$
soient interprétables dans ${\mathcal C}$ par des objets
$$
U_{\mathbb T}^{\mathcal C} \, \varphi_k \, , \quad k \in K \, ,
$$
nécessairement munis de morphismes de projection
$$
U_{\mathbb T}^{\mathcal C} \, \varphi_k \longrightarrow U_{\mathbb T}^{\mathcal C} \, A_1 \times \cdots \times U_{\mathbb T}^{\mathcal C} \, A_n
$$
vers le produit dans ${\mathcal C}$ des objets $U_{\mathbb T}^{\mathcal C} \,A_i$, $1 \leq i \leq n$, qui interprètent les sortes $A_1 , \cdots , A_n$ associées aux variables $x_1^{A_1} , \cdots , x_n^{A_n}$.

Alors l'interprétation dans le topos ${\mathcal E}_{\mathbb T}$ de la formule $\varphi$
$$
U_{\mathbb T} \, \varphi
$$
s'identifie à l'image par le foncteur de faisceautisation
$$
\widehat{\mathcal C} \longrightarrow \widehat{\mathcal C}_J \xrightarrow{ \ \sim \ } {\mathcal E}_{\mathbb T}
$$
du sous-préfaisceau du préfaisceau représentable
$$
y \left( U_{\mathbb T}^{\mathcal C} \, A_1 \times \cdots \times U_{\mathbb T}^{\mathcal C} \, A_n \right) = y(U_{\mathbb T}^{\mathcal C} \, A_1) \times \cdots \times y(U_{\mathbb T}^{\mathcal C} \, A_n)
$$
défini comme la réunion des images des préfaisceaux représentables
$$
y(U_{\mathbb T}^{\mathcal C} \, \varphi_k) \, , \quad k \in K \, ,
$$
par les morphismes de projection
$$
U_{\mathbb T}^{\mathcal C} \, \varphi_k \longrightarrow U_{\mathbb T}^{\mathcal C} \, A_1 \times \cdots \times U_{\mathbb T}^{\mathcal C} \, A_n \, .
$$
\end{listeimarge}
\end{lem}

\begin{demo}
\begin{listeisansmarge}
\item En effet, pour tout modèle $M$ de la théorie ${\mathbb T}$ dans un topos ${\mathcal E}$, l'interprétation de la formule $\varphi$ en les variables $x_1^{A_1} , \cdots , x_n^{A_n}$
$$
M\varphi (\vec x)
$$
est le sous-objet de l'objet produit
$$
MA_1 \times \cdots \times MA_n
$$
qui est la réunion des images des morphismes de projection
$$
M \varphi_k (\vec x , \vec x_k) \longrightarrow MA_1 \times \cdots \times MA_n \, , \quad k \in K \, .
$$

\item résulte de (i) et de ce que le foncteur de faisceautisation
$$
\widehat{\mathcal C} \longrightarrow \widehat{\mathcal C}_J \xrightarrow{ \ \sim \ } {\mathcal E}_{\mathbb T}
$$
respecte à la fois la formation des réunions de familles de sous-objets et celle des images de morphismes. 
\end{listeisansmarge}
\end{demo}

On déduit de ce lemme:

\begin{cor}\label{cor3.3.4}
Soient une théorie géométrique du premier ordre ${\mathbb T}$ de signature $\Sigma$ dont le topos classifiant ${\mathcal E}_{\mathbb T}$ est présenté comme topos des faisceaux
$$
\widehat{\mathcal C}_J \longrightarrow {\mathcal E}_{\mathbb T}
$$
sur un petit site $({\mathcal C},J)$.

Considérons un séquent géométrique de $\Sigma$ en une suite de variables $\vec x = (x_1^{A_1} , \cdots , x_n^{A_n})$
$$
\varphi (\vec x) \vdash \psi (\vec x)
$$
dont les deux formules constitutives $\varphi$ et $\psi$ sont écrites sous la forme
$$
\varphi (\vec x) = \bigvee_{k \in K} (\exists \, \vec x_k) \, \varphi_k (\vec x , \vec x_k) \, ,
$$
$$
\psi (\vec x) = \bigvee_{k' \in K'} (\exists \, \vec y_{k'}) \, \psi_k (\vec x , \vec y_{k'}) \, .
$$
Supposons que toutes les formules
$$
\varphi_k (\vec x , \vec x_k) \, , \quad k \in K \, ,
$$
et
$$
\varphi_k (\vec x , \vec x_k) \wedge \psi_{k'} (\vec x , \vec y_{k'}) \, , \quad k \in K \, , \ k' \in K' \, ,
$$
soient interprétables dans ${\mathcal C}$ par des objets
$$
U_{\mathbb T}^{\mathcal C} \, \varphi_k \, , \quad k \in K \, ,
$$
et
$$
U_{\mathbb T}^{\mathcal C}  (\varphi_k \wedge \psi_{k'}) \, , \quad k \in K \, , \ k' \in K' \, ,
$$
munis de morphismes 
$$
U_{\mathbb T}^{\mathcal C} (\varphi_k \wedge \psi_{k'}) \longrightarrow U_{\mathbb T}^{\mathcal C} \, \varphi_k \, , \quad k \in K \, , \ k' \in K' \, ,
$$
qui relèvent les morphismes de projection $U_{\mathbb T} (\varphi_k \wedge \psi_{k'}) \to U_{\mathbb T} \, \varphi_k$.

Alors le séquent géométrique
$$
\varphi (\vec x) \vdash \psi (\vec x)
$$
est démontrable dans la théorie ${\mathbb T}$ si et seulement si, pour tout indice $k \in K$, la famille de morphismes de ${\mathcal C}$
$$
U_{\mathbb T}^{\mathcal C} (\varphi_k \wedge \psi_{k'}) \longrightarrow U_{\mathbb T}^{\mathcal C} \, \varphi_k \, , \quad k' \in K' \, ,
$$
est couvrante pour la topologie $J$.
\end{cor}

\begin{demo}

Le séquent géométrique
$$
\varphi (\vec x) \vdash \psi (\vec x)
$$
est démontrable dans la théorie ${\mathbb T}$ si et seulement si le monomorphisme de ${\mathcal E}_{\mathbb T}$
$$
U_{\mathbb T} (\varphi \wedge \psi) \ {\lhook\joinrel\relbar\joinrel\rightarrow} \ U_{\mathbb T} \, \varphi
$$
est aussi un épimorphisme.

Or, d'après le lemme \ref{lem3.3.3}, la famille des morphismes
$$
U_{\mathbb T} \, \varphi (\vec x , \vec x_k) \longrightarrow U_{\mathbb T} \, \varphi (\vec x) \, , \quad k \in K \, ,
$$
est globalement épimorphique.

Pour tout $k \in K$, le produit fibré dans ${\mathcal E}_{\mathbb T}$
$$
U_{\mathbb T} (\varphi \wedge \psi) \times_{U_{\mathbb T} \, \varphi } U_{\mathbb T} \, \varphi (\vec x , \vec x_k)
$$
s'identifie à
$$
U_{\mathbb T} (\psi (\vec x) \wedge \varphi_k (\vec x , \vec x_k)) \, ,
$$
et la formule
$$
\psi (\vec x) \wedge \varphi_k (\vec x , \vec x_k)
$$
s'écrit
$$
\bigvee_{k' \in K'} (\exists \, \vec y_{k'}) \left(\psi_{k'} (\vec x , \vec y_{k'}) \wedge \varphi_k (\vec x , \vec x_k)\right) \, .
$$

D'après le lemme \ref{lem3.3.3}, la famille de morphismes
$$
U_{\mathbb T} \left(\psi_{k'} (\vec x , \vec y_{k'}) \wedge \varphi_k (\vec x , \vec x_k)\right) \longrightarrow U_{\mathbb T} (\varphi \wedge \psi) \times_{U_{\mathbb T} \, \varphi} U_{\mathbb T} (\varphi (\vec x , \vec x_k)) \, , \quad k' \in K' \, ,
$$
est globalement épimorphique pour chaque $k \in K$.

On en déduit que le séquent géométrique
$$
\varphi (\vec x) \vdash \psi (\vec x)
$$
est démontrable dans ${\mathbb T}$ si et seulement si, pour tout indice $k \in K$, la famille des morphismes
$$
U_{\mathbb T} \left(\psi_{k'} (\vec x , \vec y_{k'}) \wedge \varphi_k (\vec x , \vec x_k)\right) \longrightarrow U_{\mathbb T} \, \varphi_k (\vec x , \vec x_k) \, , \quad k' \in K' \, ,
$$
est globalement épimorphique.

Cela revient à demander que, pour tout $k \in K$, la famille des morphismes de ${\mathcal C}$
$$
U_{\mathbb T}^{\mathcal C} \left( \varphi_k (\vec x , \vec x_k) \wedge \psi_{k'} (\vec x , \vec y_{k'})\right) \longrightarrow U_{\mathbb T}^{\mathcal C} \, \varphi_k (\vec x , \vec x_k) \, , \quad k' \in K' \, ,
$$
soit couvrante pour la topologie $J$. \hfill $\Box$
\end{demo}

On rappelle que, d'après le lemme \ref{lem1.1.27}, toute formule géométrique dans la signature d'une théorie ${\mathbb T}$
$$
\varphi (\vec x)
$$
est équivalente, selon les règles d'inférence de la logique géométrique du premier ordre, à une formule de la forme
$$
\bigvee_{k \in K} (\exists \, \vec x_k) \, \varphi_k (\vec x , \vec x_k)
$$
où les $\varphi_k (\vec x , \vec x_k)$ sont des formules de Horn, c'est-à-dire des conjonctions finies de formules atomiques.

Le corollaire \ref{cor3.3.4} s'applique donc à tous les séquents géométriques quand toutes les formules de Horn sont interprétables. Cette condition se vérifie souvent:

\begin{prop}\label{prop3.3.5}
Soient une théorie géométrique du premier ordre ${\mathbb T}$ de signature $\Sigma$, son topos classifiant ${\mathcal E}_{\mathbb T}$ muni du modèle universel $U_{\mathbb T}$ de ${\mathbb T}$, et une présentation de ce topos comme topos des faisceaux
$$
\widehat{\mathcal C}_J \xrightarrow{ \ \sim \ } {\mathcal E}_{\mathbb T}
$$
sur un petit site $({\mathcal C},J)$.

Supposons que la catégorie sous-jacente ${\mathcal C}$ soit cartésienne, c'est-à-dire admette des limites finies arbitraires, et que tous les éléments de vocabulaire de $\Sigma$ soient interprétables dans ${\mathcal C}$, plus précisément

\medskip

$\left\{\begin{matrix}
\bullet &\mbox{les sortes $A$ par les objets $U_{\mathbb T}^{\mathcal C} \, A$ de ${\mathcal C}$,} \hfill \\
\bullet &\mbox{les symboles de relation $R \rightarrowtail A_1 \cdots A_n$ par des sous-objets de ${\mathcal C}$} \hfill \\
&U_{\mathbb T}^{\mathcal C} \, R \ {\lhook\joinrel\relbar\joinrel\rightarrow} \ U_{\mathbb T}^{\mathcal C} \, A_1 \times \cdots \times U_{\mathbb T}^{\mathcal C} \, A_n \, , \\
\bullet &\mbox{les symboles de fonction $f : A_1 \cdots A_n \to B$ par des morphismes de ${\mathcal C}$} \hfill \\
&U_{\mathbb T}^{\mathcal C} \, f : U_{\mathbb T}^{\mathcal C} \, A_1 \times \cdots \times U_{\mathbb T}^{\mathcal C} \, A_n \longrightarrow U_{\mathbb T}^{\mathcal C} \, B \, .
\end{matrix}\right.$

\medskip

Alors:

\begin{listeimarge}

\item Toute formule de Horn de $\Sigma$ en une suite de variables $\vec x = (x_1^{A_1} , \cdots x_n^{A_n})$
$$
\varphi (\vec x)
$$
est interprétable dans ${\mathcal C}$ comme un sous-objet
$$
U_{\mathbb T}^{\mathcal C} \, \varphi \ {\lhook\joinrel\relbar\joinrel\rightarrow} \ U_{\mathbb T}^{\mathcal C} \, A_1 \times \cdots \times U_{\mathbb T}^{\mathcal C} \, A_n \, .
$$

\item Un séquent géométrique de la signature de ${\mathbb T}$ en une suite de variables $\vec x = (x_1^{A_1} , \cdots x_n^{A_n})$
$$
\varphi (\vec x) \vdash \psi (\vec x)
$$
est démontrable dans ${\mathbb T}$ si et seulement si, écrivant à équivalence près ses formules constituantes $\varphi$ et $\psi$ sous la forme
$$
\varphi (\vec x) = \bigvee_{k \in K} (\exists \, \vec x_k) \, \varphi_k (\vec x , \vec x_k) \, ,
$$
$$
\psi (\vec x) = \bigvee_{k' \in K'} (\exists \, \vec y_{k'}) \, \psi_{k'} (\vec x , \vec y_{k'}) 
$$
pour des formules de Horn $\varphi_k$, $k \in K$, et $\psi_{k'}$, $k' \in K'$, on a que pour tout indice $k \in K$, la famille de morphismes de ${\mathcal C}$
$$
U_{\mathbb T}^{\mathcal C} \left( \varphi_k (\vec x , \vec x_k) \wedge \psi_{k'} (\vec x , \vec y_{k'})\right) \longrightarrow U_{\mathbb T}^{\mathcal C} \, \varphi_k (\vec x , \vec x_k) \, , \quad k' \in K' \, ,
$$
est couvrante pour la topologie $J$.
\end{listeimarge}
\end{prop}

\begin{remark}

Dans la situation de (ii), les interprétations
$$
U_{\mathbb T}^{\mathcal C} \, \left(\varphi_k (\vec x , \vec x_k) \wedge \psi_{k'} (\vec x , \vec y_{k'})\right)
$$
s'identifient aux produits fibrés
$$
U_{\mathbb T}^{\mathcal C} \, \varphi_k (\vec x , \vec x_k) \times_{U_{\mathbb T}^{\mathcal C} \, A_1 \times \cdots \times U_{\mathbb T}^{\mathcal C} \, A_n} U_{\mathbb T}^{\mathcal C} \, \psi_{k'} (\vec x , \vec y_{k'}) \, .
$$
\end{remark}

\begin{demo}
\begin{listeisansmarge}
\item La définition des termes de la signature de ${\mathbb T}$ à valeurs dans une sorte $B_0$
$$
g (y_1^{B_1} , \cdots , y_m^{B_m})
$$
et la construction de leurs interprétations dans le modèle $U_{\mathbb T}$ comme des morphismes
$$
U_{\mathbb T} \, B_1 \times \cdots \times U_{\mathbb T} \, B_m \longrightarrow U_{\mathbb T} \, B_0 \, ,
$$
rappelées dans les définitions \ref{defn1.1.17} (i) et \ref{defn1.1.18} (i), montrent que ces interprétations se déduisent de celles des symboles de fonction $f : A_1 \cdots A_m \to B$
$$
U_{\mathbb T} \, f : U_{\mathbb T} \, A_1 \times \cdots \times U_{\mathbb T} A_m \longrightarrow U_{\mathbb T} \, B
$$
en formant des produits finis de morphismes et des compositions de morphismes.

Comme la catégorie ${\mathcal C}$ a des produits finis et que le foncteur canonique
$$
\ell : {\mathcal C} \longrightarrow \widehat{\mathcal C}_J \xrightarrow{ \ \sim \ } {\mathcal E}_{\mathbb T}
$$
respecte les produits finis et les compositions de morphismes, chaque terme à valeurs dans une sorte $B_0$
$$
g (y_1^{B_1} , \cdots , y_m^{B_m})
$$
s'interprète dans ${\mathcal C}$ comme un morphisme
$$
U_{\mathbb T}^{\mathcal C} \, g : U_{\mathbb T}^{\mathcal C} \, B_1 \times \cdots \times U_{\mathbb T}^{\mathcal C} \, B_m \longrightarrow U_{\mathbb T}^{\mathcal C} \, B_0 \, .
$$

Les formules atomiques de la signature de ${\mathbb T}$ ont la forme
$$
R (g_1 , \cdots , g_n)
$$
pour des symboles de relation $R \rightarrowtail A_1 \cdots A_n$ et des termes $g_1 , \cdots , g_n$ à valeurs dans les sortes $A_1 , \cdots  , A_n$, ou la forme
$$
g_1 = g_2
$$
pour des termes $g_1 , g_2$ à valeurs dans une sorte $A$.

Dans le premier cas, l'interprétation dans un modèle $U_{\mathbb T}$
$$
U_{\mathbb T} \, R (g_1 , \cdots , g_n)
$$
est l'image inverse du sous-objet
$$
U_{\mathbb T} \, R \ {\lhook\joinrel\relbar\joinrel\rightarrow} \ U_{\mathbb T} \, A_1 \times \cdots \times U_{\mathbb T} \, A_n
$$
par le morphisme produit de but $U_{\mathbb T} \, A_1 \times \cdots \times U_{\mathbb T} \, A_n$
$$
U_{\mathbb T} \, g_1 \times \cdots \times U_{\mathbb T} \, g_n \, ,
$$
et dans le second cas, l'interprétation de l'équation
$$
g_1 = g_2
$$
est l'image inverse du sous-objet diagonal
$$
U_{\mathbb T} \, A \ {\lhook\joinrel\relbar\joinrel\rightarrow} \ U_{\mathbb T} \, A \times  U_{\mathbb T} \, A
$$
par le morphisme produit de but $U_{\mathbb T} \, A \times  U_{\mathbb T} \, A$
$$
U_{\mathbb T} \, g_1 \times  U_{\mathbb T} \, g_2 \, .
$$

Comme la catégorie ${\mathcal C}$ a des produits fibrés et que le foncteur canonique
$$
\ell : {\mathcal C} \longrightarrow \widehat{\mathcal C}_J \xrightarrow{ \ \sim \ } {\mathcal E}_{\mathbb T}
$$
les respecte, les formules atomiques de la forme
$$
R(g_1 , \cdots , g_m)
$$
s'interprètent comme les produits fibrés des morphismes
$$
U_{\mathbb T}^{\mathcal C} \, R \ {\lhook\joinrel\relbar\joinrel\rightarrow} \ U_{\mathbb T}^{\mathcal C} \, A_1 \times  \cdots \times U_{\mathbb T}^{\mathcal C} \, A_n
$$
et
$$
U_{\mathbb T}^{\mathcal C} \, g_1 \times \cdots \times U_{\mathbb T}^{\mathcal C} \, g_n \, ,
$$
et les formules atomiques de la forme
$$
g_1 = g_2
$$
s'interprètent comme les produits fibrés des morphismes
$$
U_{\mathbb T}^{\mathcal C} \, A \ {\lhook\joinrel\relbar\joinrel\rightarrow} \ U_{\mathbb T}^{\mathcal C} \, A \times  U_{\mathbb T}^{\mathcal C} \, A
$$
et
$$
U_{\mathbb T}^{\mathcal C} \, g_1 \times U_{\mathbb T}^{\mathcal C} \, g_2 \, .
$$

Enfin, les formules de Horn en une suite de variables $\vec x = (x_1^{A_1} , \cdots , x_n^{A_n})$ sont les conjonctions finies
$$
\varphi (\vec x) = \varphi_1 (\vec x) \wedge \cdots \wedge \varphi_m (\vec x)
$$
de formules atomiques $\varphi_j (\vec x)$, $1 \leq j \leq m$, en ces variables.

Elles s'interprètent dans ${\mathcal C}$ comme l'intersection des sous-objets
$$
U_{\mathbb T}^{\mathcal C} \, \varphi = \bigcap_{1 \leq j \leq m} U_{\mathbb T}^{\mathcal C} \, \varphi_j
$$
c'est-à-dire le produit fibré des morphismes de plongement
$$
U_{\mathbb T}^{\mathcal C} \, \varphi_j \ {\lhook\joinrel\relbar\joinrel\rightarrow} \ U_{\mathbb T}^{\mathcal C} \, A_1 \times \cdots \times U_{\mathbb T}^{\mathcal C} \, A_n \, , \quad 1 \leq j \leq m \, .
$$
\item résulte du lemme \ref{lem1.1.27} qui affirme l'existence de telles réécritures pour les formules $\varphi$ et $\psi$ en termes de formules de Horn, de ce que les formules de Horn sont toujours interprétables dans ${\mathcal C}$ d'après (i), et du corollaire \ref{cor3.3.4}. 
\end{listeisansmarge}
\end{demo}

\subsection{Processus de traduction des problèmes géométriques vers les topologies}\label{ssec3.3.3}

La proposition \ref{prop3.3.5} amène à préciser les constructions de la définition \ref{defn1.2.14} dans le cas où le topos classifiant ${\mathcal E}_{\mathbb T}$ d'une théorie géométrique ${\mathbb T}$ est présenté comme topos des faisceaux sur une catégorie ${\mathcal C}$ supposée cartésienne et telle que les éléments de la signature de ${\mathbb T}$ soient interprétables dans ${\mathcal C}$:

\begin{defn}\label{defn3.3.6}
Soient une théorie géométrique du premier ordre ${\mathbb T}$ de signature $\Sigma$ et une présentation de son topos classifiant ${\mathcal E}_{\mathbb T}$ comme topos des faisceaux
$$
\widehat{\mathcal C}_J \xrightarrow{ \ \sim \ } {\mathcal E}_{\mathbb T}
$$
sur une petite catégorie ${\mathcal C}$ munie d'une topologie $J$.

Supposons que la catégorie ${\mathcal C}$ soit cartésienne et que les éléments de vocabulaire de $\Sigma$ soient interprétables dans ${\mathcal C}$ comme

\medskip

$\left\{\begin{matrix}
\bullet &\mbox{des objets $U_{\mathbb T}^{\mathcal C} \, A$ indexés par les sortes $A$ de $\Sigma$,} \hfill \\
\bullet &\mbox{des morphismes $U_{\mathbb T}^{\mathcal C} \, f : U_{\mathbb T}^{\mathcal C} \, A_1 \times \cdots \times U_{\mathbb T}^{\mathcal C} \, A_n \to U_{\mathbb T}^{\mathcal C} \, B$} \hfill \\
&\mbox{indexés par les symboles de fonction $f : A_1 \cdots A_n \to B$ de $\Sigma$,} \hfill \\
\bullet &\mbox{des sous-objets $U_{\mathbb T}^{\mathcal C} \, R \hookrightarrow U_{\mathbb T}^{\mathcal C} \, A_1 \times \cdots \times U_{\mathbb T}^{\mathcal C} \, A_n$} \hfill \\
&\mbox{indexés par les symboles de relation $R \rightarrowtail A_1 \cdots A_n$ de $\Sigma$,} \hfill
\end{matrix}\right.$

\medskip

si bien que toute formule de Horn de $\Sigma$ en une suite de variables $\vec x = (x_1^{A_1} , \cdots , x_n^{A_n})$
$$
\varphi (\vec x)
$$
s'interprète dans ${\mathcal C}$ comme un sous-objet
$$
U_{\mathbb T}^{\mathcal C} \, \varphi (\vec x) \ {\lhook\joinrel\relbar\joinrel\rightarrow} \ U_{\mathbb T}^{\mathcal C} \, A_1 \times \cdots \times  U_{\mathbb T}^{\mathcal C} \, A_n \, .
$$
Pour n'importe quel séquent géométrique de $\Sigma$
$$
\varphi (\vec x) \vdash \psi (\vec x)
$$
dont les formules constitutives sont écrites à équivalence près sous la forme
$$
\varphi (\vec x) = \bigvee_{k \in K} (\exists \, \vec x_k) \, \varphi_k (\vec x , \vec x_k) \, ,
$$
$$
\psi (\vec x) = \bigvee_{k' \in K'} (\exists \, \vec y_{k'}) \, \psi_{k'} (\vec x , \vec y_{k'}) \, ,
$$
on note
$$
\chi_{\varphi , \psi}
$$
la famille indexée par les éléments $k \in K$ des précribles sur les objets $U_{\mathbb T}^{\mathcal C} \, \varphi_k (\vec x , \vec x_k)$ de ${\mathcal C}$
$$
P_{\varphi , \psi , k} = \left( U_{\mathbb T}^{\mathcal C} \left( \varphi_k (\vec x , \vec x_k) \wedge \psi_{k'} (\vec x , \vec y_{k'})\right) \longrightarrow U_{\mathbb T}^{\mathcal C} \, \varphi_k (\vec x , \vec x_k) \right)_{k' \in K'} \, .
$$
\end{defn}

Ayant introduit ces notations, on déduit de la proposition \ref{prop3.3.5}:

\begin{thm}\label{thm3.3.7}
Soit comme ci-dessus une théorie géométrique ${\mathbb T}$ de signature $\Sigma$ dont le topos classifiant ${\mathcal E}_{\mathbb T}$ est présenté comme topos des faisceaux
$$
\widehat{\mathcal C}_J \xrightarrow{ \ \sim \ } {\mathcal E}_{\mathbb T}
$$
sur une petite catégorie cartésienne ${\mathcal C}$ munie d'une topologie $J$ et telle que les éléments de vocabulaire de $\Sigma$ soient interprétables dans ${\mathcal C}$.

Considérons une famille de séquents géométriques de la signature de ${\mathbb T}$
$$
\varphi_i (\vec x_i) \vdash\psi_i (\vec x_i) \, , \quad i \in I \, ,
$$
et la théorie quotient ${\mathbb T}'$ de ${\mathbb T}$ que cette famille d'axiomes supplémentaires définit.

Alors un séquent géométrique de la signature de ${\mathbb T}$
$$
\varphi (\vec x) \vdash \psi (\vec x)
$$
est démontrable dans la théorie quotient ${\mathbb T}'$ si et seulement si tous les précribles
$$
P \in \chi_{\varphi , \psi}
$$
sont couvrants pour la topologie $J'$ de ${\mathcal C}$ engendrée par $J$ et les précribles
$$
P \in \chi_{\varphi_i , \psi_i} \, , \quad i \in I \, .
$$
\end{thm}

\begin{demo}
D'après le théorème \ref{thm1.2.10}, les théories quotients ${\mathbb T}'$ de ${\mathbb T}$ correspondent aux topologies $J'$ de ${\mathcal C}$ qui raffinent la topologie $J$.

De plus, lorsque ${\mathbb T}'$ et $J'$ se correspondent, l'équivalence
$$
\widehat{\mathcal C}_J \xrightarrow{ \ \sim \ } {\mathcal E}_{\mathbb T}
$$
induit une équivalence
$$
\widehat{\mathcal C}_{J'} \xrightarrow{ \ \sim \ } {\mathcal E}_{{\mathbb T}'}
$$
telle que le modèle universel $U_{{\mathbb T}'}$ de ${\mathbb T}'$ dans ${\mathcal E}_{{\mathbb T}'}$ s'identifie au transformé du modèle universel $U_{\mathbb T}$ de ${\mathbb T}$ dans  ${\mathcal E}_{\mathbb T}$ par le foncteur de faisceautisation
$$
j^* : \widehat{\mathcal C}_J \longrightarrow \widehat{\mathcal C}_{J'} \, .
$$

Il résulte alors de la proposition \ref{prop3.3.5} que la théorie quotient ${\mathbb T}'$ définie par une famille d'axiomes additionnels
$$
\varphi_i (\vec x_i) \vdash \psi_i (\vec x_i) \, , \quad i \in I \, ,
$$
correspond à la topologie $J'$ de ${\mathcal C}$ engendrée sur $J$ par les familles de précribles
$$
P \in \chi_{\varphi_i , \psi_i} \, , \quad i \in I \, ,
$$
et qu'un séquent géométrique
$$
\varphi (\vec x) \vdash \psi (\vec x)
$$
est démontrable dans ${\mathbb T}'$ si et seulement si les précribles
$$
P \in \chi_{\varphi , \psi}
$$
sont couvrants pour cette topologie. 
\end{demo}

Il n'est pas inutile de rappeler ici comment, d'après les résultats des paragraphes \ref{sec3.1} et \ref{sec3.2}, une topologie d'une catégorie ${\mathcal C}$ supposée cartésienne est engendrée par une famille de précribles:

\begin{thm}\label{thm3.3.8}
Soit ${\mathcal C}$ une petite catégorie supposée cartésienne, c'est-à-dire qui admet des limites finies arbitraires.

Soit $\chi$ une famille de précribles ${\mathcal C}$ c'est-à-dire de familles de morphismes de ${\mathcal C}$ de même but
$$
P = \{ X' \xrightarrow{ \ x \ } X \mid x \in P \} \, .
$$
Soit ${\mathcal J}_{\chi}$ la ``notion stable de précrible couvrant'' constituée de tous les précribles
$$
x_1^* P = \{ X' \times_X X_1 \longrightarrow X_1 \mid (X' \xrightarrow{ \ x \ } X) \in P \}
$$
images réciproques des précribles éléments de $\chi$
$$
P = \{ X' \xrightarrow{ \ x \ } X \mid x \in P \} \, , \quad P \in \chi \, ,
$$
par les morphismes de ${\mathcal C}$
$$
x_1 : X_1 \longrightarrow X \, .
$$
Alors:

\begin{listeimarge}

\item Un crible $C$ d'un objet $X$ de ${\mathcal C}$ est couvrant pour la topologie $J$ engendrée par $\chi$ ou ${\mathcal J}_{\chi}$ si le crible maximal de $X$ est le seul crible qui contient $C$ et qui est ${\mathcal J}_{\chi}$-fermé au sens que

\medskip

$\left\{\begin{matrix}
\bullet &\mbox{pour tout morphisme $X_1 \xrightarrow{ \ x \ } X$ de but $X$,} \hfill \\
\bullet &\mbox{pour tout précrible $P_1$ de but $X_1$ qui est élément de ${\mathcal J}_{\chi}$} \hfill \\
&P_1 = \left\{ X'_1 \xrightarrow{ \ x_1 \ } X_1 \mid x_1 \in P_1 \right\} \, ,
\end{matrix}\right.$

\medskip

\noindent on a
$$
x \in C
$$
si
$$
x \circ x_1 \in C \, , \quad \forall \, x_1 \in P_1 \, .
$$
\item Un précrible d'un objet $X$ de ${\mathcal C}$
$$
P = \left\{ X' \xrightarrow{ \ x \ } X \mid x \in P \right\}
$$
est couvrant pour la topologie $J$ engendrée par $\chi$ ou ${\mathcal J}_{\chi}$ si et seulement si il existe un multi-recouvrement de type ${\mathcal J}_{\chi}$ de l'objet $X$ consistant en

\medskip

$\left\{\begin{matrix}
\bullet &\mbox{un arbre à branches finies $A_{\bullet} = \left( A_0 \xleftarrow{ \, f_1 \, } A_1 \leftarrow \cdots \leftarrow A_{n-1} \xleftarrow{ \, f_n \, } A_n \leftarrow \cdots \right)$ avec $A_0 = \{a_0\}$,} \hfill \\
\bullet &\mbox{une famille d'objets de ${\mathcal C}$} \hfill \\
&X_{n , a_n} \, , \quad n \in {\mathbb N} \, , \quad a \in A_n \, , \quad \mbox{avec} \ X_{0,a_0} = X \, , \\
\bullet &\mbox{une famille de morphismes de ${\mathcal C}$} \hfill \\
&x_{n,a} : X_{n,a} \longrightarrow X_{n-1,f_n(a)} \, , \quad n \geq 1 \, , \ a \in A_n \, ,
\end{matrix}\right.$

\medskip

\noindent tel que, pour tout $n \in {\mathbb N}$ et tout $a \in A_n$, on a l'alternative:

\medskip

$\left\{\begin{matrix}
\bullet &\mbox{Ou bien la fibre $f_{n+1}^{-1} (a)$ est non vide et le précrible} \hfill \\
&\left\{ X_{n+1,a'} \xrightarrow{ \ x_{n+1,a'} \ } X_{n,a} \mid a' \in f_{n+1}^{-1} (a) \right\} \\
&\mbox{est élément de ${\mathcal J}_{\chi}$.} \hfill \\
\bullet &\mbox{Ou bien le précrible vide sur l'objet $X_{n,a}$ est élément de ${\mathcal J}_{\chi}$.} \hfill \\
\bullet &\mbox{Ou bien, notant} \hfill \\
&a = a_n \, , \ f_n (a_n) = a_{n-1} , \cdots , f_1(a_1) = a_0 \, , \\
&\mbox{le morphisme composé} \hfill \\
&X_{n,a_n} \xrightarrow{ \, x_{n,a_n} \, } X_{n-1 , a_{n-1}} \longrightarrow \cdots \longrightarrow X_{1,a_1} \xrightarrow{ \, x_{1,a_1} \, } X_{0,a_0} = X \\
&\mbox{est élément du crible engendré par le précrible $P$} \hfill \\ 
&\mbox{c'est-à-dire se factorise à travers l'un au moins des morphismes} \hfill \\
&X' \xrightarrow{ \ x \ } X \, ,  \quad x \in P \, .
\end{matrix}\right.$
\end{listeimarge}
\end{thm}

\begin{demo}
Comme la catégorie ${\mathcal C}$ est cartésienne, l'image réciproque $x_1^* C$ par un morphisme $x_1 : X_1 \to X$ d'un crible $C$ de $X$ engendré par un précrible
$$
P = \left\{ X' \xrightarrow{ \ x \ } X \mid x \in P \right\}
$$
est le crible engendré par le précrible de $X_1$
$$
x_1^* P = \left\{ X' \times_X X_1 \longrightarrow X_1 \mid (X' \xrightarrow{ \ x \ } X) \in P \right\} \, .
$$

Alors la partie (i) du théorème résulte du théorème \ref{thm3.1.5} et la partie (ii) résulte du corollaire \ref{cor3.2.7}. 
\end{demo}

\subsection{Le problème de présentation des topos classifiants}\label{ssec3.3.4}

Etant donnée une théorie géométrique du premier ordre ${\mathbb T}$, on se demande comment présenter son topos classifiant ${\mathcal E}_{\mathbb T}$ comme topos des faisceaux pour une certaine topologie sur une petite catégorie ${\mathcal C}$ qui soit cartésienne et telle que les éléments de vocabulaire de ${\mathbb T}$ s'interprètent dans ${\mathcal C}$.

Remarquons d'abord que toute telle présentation définit une théorie cartésienne de même signature que ${\mathbb T}$:

\begin{defn}\label{defn3.3.9}
	Soit une théorie géométrique du premier ordre ${\mathbb T}$ de signature $\Sigma$ dont le topos classifiant ${\mathcal E}_{\mathbb T}$ est présenté comme topos des faisceaux
	$$
	\widehat{\mathcal C}_J \xrightarrow{\ \sim \ } {\mathcal E}_{\mathbb T}
	$$
	sur un petit site $({\mathcal C},J)$.
	
	On suppose que ${\mathcal C}$ est cartésienne et que
	
	\medskip
	
	$\left\{\begin{matrix}
		\bullet &\mbox{les sortes $A$ de $\Sigma$ s'interprètent comme des objets $U_{\mathbb T}^{\mathcal C} A$ de ${\mathcal C}$,} \hfill \\
		\bullet &\mbox{les symboles de relation $R \rightarrowtail A_1 \cdots A_n$ de $\Sigma$ s'interprètent comme des sous-objets} \hfill \\ 
		&U_{\mathbb T}^{\mathcal C} R \hookrightarrow U_{\mathbb T}^{\mathcal C} A_1 \times \cdots \times U_{\mathbb T}^{\mathcal C}  A_n \, , \\
		\bullet &\mbox{les symboles de fonction $f : A_1 \cdots A_n \to B$ de $\Sigma$ s'interprètent comme des morphismes de ${\mathcal C}$} \hfill \\
		&U_{\mathbb T}^{\mathcal C}  f : U_{\mathbb T}^{\mathcal C} A_1 \times \cdots \times U_{\mathbb T}^{\mathcal C}  A_n \longrightarrow U_{\mathbb T}^{\mathcal C} B \, ,
	\end{matrix}\right.$
	
	\medskip
	
	\noindent si bien que toute formule de Horn de $\Sigma$ en une suite de variables $\vec x = (x_1^{A_1} , \cdots , x_n^{A_n})$
	$$
	\varphi (\vec x) = \varphi (x_1^{A_1} , \cdots , x_n^{A_n})
	$$
	s'interprète comme un sous-objet dans ${\mathcal C}$
	$$
	U_{\mathbb T}^{\mathcal C} \varphi (\vec x) \ {\lhook\joinrel\relbar\joinrel\rightarrow} \ U_{\mathbb T}^{\mathcal C} A_1 \times \cdots \times U_{\mathbb T}^{\mathcal C}  A_n \, .
	$$
	
	Alors on appelle théorie cartésienne de l'interprétation de $\Sigma$ dans ${\mathcal C}$ et on note ${\mathbb T}_{\mathcal C}$ la théorie de signature $\Sigma$ définie en les deux étapes suivantes:
	
	\begin{listeimarge}
	
	\item[(1)] Pour toute famille de Horn $\varphi (\vec x , \vec y)$ de $\Sigma$ en deux suites de variables $\vec x = (x_1^{A_1} , \cdots , x_n^{A_n})$ et $\vec y$, telle que le morphisme de ${\mathcal C}$
	$$
	U_{\mathbb T}^{\mathcal C} \varphi (\vec x , \vec y) \longrightarrow U_{\mathbb T}^{\mathcal C} A_1 \times \cdots \times U_{\mathbb T}^{\mathcal C}  A_n 
	$$
	est un monomorphisme, on décide que le séquent de Horn
	$$
	\varphi (\vec x , \vec y) \wedge \varphi (\vec x , \vec y\,') \vdash \vec y = \vec y\, '
	$$
	est un axiome de ${\mathbb T}_{\mathcal C}$ et donc que la formule
	$$
	(\exists \, \vec y) \, \varphi (\vec x , \vec y)
	$$
	est une formule ${\mathbb T}_{\mathcal C}$-cartésienne qui s'interprète dans ${\mathcal C}$ comme le sous-objet de
	$$
	U_{\mathbb T}^{\mathcal C} A_1 \times \cdots \times U_{\mathbb T}^{\mathcal C}  A_n 
	$$
	défini par le monomorphisme de ${\mathcal C}$
	$$
	U_{\mathbb T}^{\mathcal C} \varphi (\vec x , \vec y) \ {\lhook\joinrel\relbar\joinrel\rightarrow} \ U_{\mathbb T}^{\mathcal C} A_1 \times \cdots \times U_{\mathbb T}^{\mathcal C}  A_n \, .
	$$
	
	\item[(2)] Pour toutes formules ${\mathbb T}_{\mathcal C}$-cartésiennes $\varphi (\vec x)$ et $\psi (\vec x)$ en une suite de variables $\vec x = (x_1^{A_1} , \cdots , x_n^{A_n})$ de la forme de (1), on décide que le séquent
	$$
	\varphi (\vec x) \vdash \psi (\vec x)
	$$
	est démontrable dans ${\mathbb T}_{\mathcal C}$ si et seulement si les interprétations dans ${\mathcal C}$ de $\varphi$ et $\psi$ satisfont la relation
	$$
	U_{\mathbb T}^{\mathcal C} \varphi (\vec x) \subseteq U_{\mathbb T}^{\mathcal C} \psi (\vec x) \, .
	$$
\end{listeimarge}
\end{defn}

Cette définition permet de se ramener au cas des catégories syntactiques cartésiennes comme catégories de présentation des topos classifiants:

\begin{thm}\label{thm3.3.10}
	Soit une théorie géométrique du premier ordre ${\mathbb T}$ de signature $\Sigma$ dont le topos classifiant ${\mathcal E}_{\mathbb T}$ est présenté comme topos des faisceaux pour une topologie $J$
	$$
	\widehat{\mathcal C}_J \xrightarrow{\ \sim \ } {\mathcal E}_{\mathbb T}
	$$
	sur une petite catégorie ${\mathcal C}$ supposée cartésienne et dans laquelle les éléments du vocabulaire $\Sigma$ de ${\mathbb T}$ s'interprètent.
	
	Soit ${\mathbb T}_{\mathcal C}$ la théorie cartésienne de l'interprétation de $\Sigma$ dand ${\mathcal C}$.
	
	Soit ${\mathcal C}_c$ la catégorie syntactique cartésienne de ${\mathbb T}_{\mathcal C}$.
	
	Alors:
	
	\begin{listeimarge}
	\item Les catégories ${\mathcal C}$ et ${\mathcal C}_c$ sont reliées par un foncteur canonique
	$$
	{\mathcal C}_c \longrightarrow {\mathcal C}
	$$
	qui est fidèle et respecte les limites finies.
	
	\item Il existe une unique topologie $J_c$ sur ${\mathcal C}_c$ telle que le foncteur composé
	$$
	{\mathcal C}_c \ {\lhook\joinrel\relbar\joinrel\rightarrow} \ {\mathcal C} \longrightarrow \widehat{\mathcal C}_J \xrightarrow{ \ \sim \ } {\mathcal E}_{\mathbb T}
	$$
	induise une équivalence
	$$
	\widehat{({\mathcal C}_c)}_{J_c} \xrightarrow{ \ \sim \ } {\mathcal E}_{\mathbb T} \, .
	$$
	
	\item Un précrible de ${\mathcal C}_c$
	$$
	\left( X_i \xrightarrow{ \ x_i \ } X \right)_{i \in I}
	$$
	est couvrant pour la topologie $J_c$ si et seulement si son image par le foncteur de plongement
	$$
	{\mathcal C}_c \ {\lhook\joinrel\relbar\joinrel\rightarrow} \ {\mathcal C}
	$$
	est couvrant pour la topologie $J$.
	
	\item Tout objet $X$ de ${\mathcal C}$ admet un précrible $J$-couvrant
	$$
	\left( X_i \xrightarrow{ \ x_i \ } X \right)_{i \in I}
	$$
	dont les sources $X_i$, $i \in I$, sont des objets de ${\mathcal C}_c$.
	
	\item Pour tout morphisme de ${\mathcal C}$ reliant deux objets de ${\mathcal C}_c$
	$$
	X \xrightarrow{ \ y \ } Y \, ,
	$$
	il existe un précrible $J_c$-couvrant de $X$ dans ${\mathcal C}_c$
	$$
	\left( X_i \xrightarrow{ \ x_i \ } X \right)_{i \in I}
	$$
	tel que tous les morphismes composés
	$$
	X_i \xrightarrow{ \ y \, \circ \, x_i \ } Y \, , \quad i \in I \, ,
	$$
	soient des morphismes de ${\mathcal C}_c$.
	\end{listeimarge}
\end{thm}

\begin{demo}
	\begin{listeisansmarge}
	\item Les objets de ${\mathcal C}_c$ sont par définition les formules ${\mathbb T}_{\mathcal C}$-cartésiennes en des suites de variables $\vec x = (x_1^{A_1} , \cdots , x_n^{A_n})$.
	
	Ces formules sont de la forme
	$$
	(\exists \, \vec y) \, \varphi (\vec x , \vec y)
	$$
	pour une formule de Horn $\varphi (\vec x , \vec y)$ telle que le séquent
	$$
	\varphi (\vec x , \vec y) \wedge \varphi (\vec x, \vec y\, ') \vdash \vec y = \vec y\, '
	$$
	soit ${\mathbb T}_{\mathcal C}$-démontrable, c'est-à-dire telle que le morphisme de ${\mathcal C}$
	$$
	U_{\mathbb T}^{\mathcal C} \varphi (\vec x , \vec y) \longrightarrow U_{\mathbb T}^{\mathcal C} A_1 \times \cdots \times U_{\mathbb T}^{\mathcal C}  A_n 
	$$
	soit un monomorphisme et donc définisse un sous-objet de $U_{\mathbb T}^{\mathcal C} A_1 \times \cdots \times U_{\mathbb T}^{\mathcal C}  A_n$.
	
	A chaque objet de ${\mathcal C}_c$ est ainsi associé un objet de ${\mathcal C}$.
	
	Les morphismes de ${\mathcal C}_c$ sont les formules ${\mathbb T}_{\mathcal C}$-cartésiennes
	$$
	\varphi (\vec x) \xrightarrow{ \ \theta (\vec x , \vec y) \ } \psi (\vec y)
	$$
	qui sont démontrablement fonctionnelles, c'est-à-dire telles que les séquents
	$$
	\left\{ \begin{matrix}
		\theta (\vec x , \vec y) \vdash \varphi (\vec x) \wedge \psi (\vec y) \hfill \\
		\theta (\vec x , \vec y) \wedge \theta (\vec x , \vec y\, ') \vdash \vec y = \vec y\, ' \, , \hfill \\
		\varphi (\vec x) \vdash (\exists \, \vec y) \, \theta (\vec x , \vec y) \hfill
	\end{matrix}\right.
	$$
	soient ${\mathbb T}_{\mathcal C}$-démontrables.
	
	Cela revient à demander que la première projection définit un isomorphisme de ${\mathcal C}$
	$$
	U_{\mathbb T}^{\mathcal C} \theta (\vec x , \vec y) \xrightarrow{ \ \sim \ } U_{\mathbb T}^{\mathcal C} \varphi (\vec x)
	$$
	et donc qu'au morphisme de ${\mathcal C}_c$
	$$
	\varphi (\vec x) \xrightarrow{ \ \theta (\vec x , \vec y) \ } \psi (\vec y)
	$$
	est associé un morphisme de ${\mathcal C}$
	$$
	U_{\mathbb T}^{\mathcal C} \varphi (\vec x) \longrightarrow U_{\mathbb T}^{\mathcal C} \psi (\vec y) \, .
	$$
	Il y a compatibilité avec la composition des morphismes de ${\mathcal C}_c$
	$$
	\varphi (\vec x) \xrightarrow{ \ \theta (\vec x , \vec y) \ } \psi (\vec y) \xrightarrow{ \ \theta' (\vec y , \vec z) \ } \chi (\vec z)
	$$
	définie par les formules
	$$
	(\exists \, \vec y) \, (\theta (\vec x , \vec y) \wedge \theta' (\vec y , \vec z)) \, .
	$$
	On a ainsi défini un foncteur
	$$
	{\mathcal C}_c \longrightarrow {\mathcal C} \, .
	$$
	Ce foncteur est fidèle car deux morphismes de ${\mathcal C}_c$
	$$
	\xymatrix{ \varphi (\vec x) \dar[rr]^{^{^{\mbox{\scriptsize $\theta_1 (\vec x , \vec y)$}}}}_{\theta_2 (\vec x , \vec y)} && \psi (\vec y) }
	$$
	ont même image dans ${\mathcal C}$ si et seulement si les deux formules 
	$$
	\theta_1 (\vec x , \vec y) \quad \mbox{et} \quad \theta_2 (\vec x , \vec y)
	$$
	sont ${\mathbb T}_{\mathcal C}$-démontrablement équivalentes.
	
	Enfin, le foncteur 
	$$
	{\mathcal C}_c \longrightarrow {\mathcal C}
	$$
	respecte les limites finies car
	
	\medskip
	
	$\left\{\begin{matrix}
		\bullet &\mbox{le produit de deux objets $\varphi (\vec x)$ et $\psi (\vec y)$ de ${\mathcal C}_c$} \hfill \\
		&\varphi (\vec x) \wedge \psi (\vec y) \\
		&\mbox{est transformé en le produit} \hfill \\
		&U_{\mathbb T}^{\mathcal C} (\varphi (\vec x) \wedge \psi (\vec y)) = U_{\mathbb T}^{\mathcal C} \varphi (\vec x) \times U_{\mathbb T}^{\mathcal C} \psi (\vec y) \, , \\
		\bullet &\mbox{l'objet terminal de ${\mathcal C}_c$, qui est la formule sans variable} \hfill \\
		&\top \, , \\
		&\mbox{est transformé en un objet terminal de ${\mathcal C}$,} \hfill \\
		\bullet &\mbox{l'égalisateur de deux morphismes de ${\mathcal C}_c$} \hfill \\
		&\xymatrix{ \varphi (\vec x) \dar[rr]^{^{^{\mbox{\scriptsize $\theta_1 (\vec x , \vec y)$}}}}_{\theta_2 (\vec x , \vec y)} && \psi (\vec y) } , \\
		&\mbox{qui s'identifie à la formule ${\mathbb T}_{\mathcal C}$-cartésienne} \hfill \\
		&(\exists \, \vec y) (\theta_1 (\vec x ,\vec y) \wedge \theta_2 (\vec x ,\vec y)) \, , \\
		&\mbox{est transformé en l'égalisateur des morphismes de ${\mathcal C}$ associés.} \hfill
	\end{matrix} \right.$
	
	\medskip
	
	\item Comme le foncteur
	$$
	{\mathcal C} \longrightarrow \widehat{\mathcal C}_J \xrightarrow{ \ \sim \ } {\mathcal E}_{\mathbb T}
	$$
	respecte les limites finies, tous les axiomes de la théorie cartésienne ${\mathbb T}_{\mathcal C}$ sont vérifiés par le modèle universel $U_{\mathbb T}$ de la théorie ${\mathbb T}$ dans le topos classifiant ${\mathcal E}_{\mathbb T}$, et donc ils sont démontrables dans ${\mathbb T}$.
	
	Ainsi, la théorie ${\mathbb T}$ appara{\^\i}t comme un quotient de la théorie ${\mathbb T}_{\mathcal C}$.
	
	Comme la théorie ${\mathbb T}_{\mathcal C}$ est cartésienne, elle admet pour topos classifiant le topos de préfaisceaux
	$$
	\widehat{\mathcal C}_c \, .
	$$
	Sa théorie quotient ${\mathbb T}$ est classifiée par un sous-topos de $\widehat{\mathcal C}_c$, nécessairement de la forme
	$$
	\widehat{({\mathcal C}_c)}_{J_c}
	$$
	pour une unique topologie $J_c$ sur la catégorie ${\mathcal C}_c$.
	
	Le foncteur canonique
	$$
	{\mathcal C}_c \longrightarrow \widehat{({\mathcal C}_c)}_{J_c}\xrightarrow{ \ \sim \ } {\mathcal E}_{\mathbb T}
	$$
	se factorise à travers le foncteur
	$$
	{\mathcal C}_c \longrightarrow {\mathcal C}
	$$
	puisqu'il consiste à associer à toute formule ${\mathbb T}_{\mathcal C}$-cartésienne
	$$
	\varphi(\vec x)
	$$
	son interprétation dans le topos classifiant ${\mathcal E}_{\mathbb T}$
	$$
	U_{\mathbb T} \, \varphi (\vec x)
	$$
	laquelle est image de l'interprétation dans ${\mathcal C}$
	$$
	U_{\mathbb T}^{\mathcal C} \, \varphi (\vec x)
	$$
	par le foncteur canonique
	$$
	{\mathcal C} \longrightarrow \widehat{\mathcal C}_J \xrightarrow{ \ \sim \ } {\mathcal E}_{\mathbb T} \, .
	$$
	\item résulte de ce que les deux conditions équivalent à ce que le transformé du précrible
	$$
	\left( X_i \xrightarrow{ \ x_i \ } X \right)_{i \in I}
	$$
	par le foncteur
	$$
	{\mathcal C}_c \ {\lhook\joinrel\relbar\joinrel\rightarrow} \ {\mathcal C} \longrightarrow \widehat{\mathcal C}_J \xrightarrow{ \ \sim \ } {\mathcal E}_{\mathbb T}
	$$
	soit une famille globalement épimorphique de ${\mathcal E}_{\mathbb T}$.
	
	\item Pour tout objet $X$ de ${\mathcal C}$, son image par le foncteur
	$$
	{\mathcal C} \xrightarrow{ \ \ell \ } \widehat{\mathcal C}_J \xrightarrow{ \ \sim \ } \widehat{({\mathcal C}_c)}_{J_c}
	$$
	s'identifie avec son image par le foncteur
	$$
	{\mathcal C} \xrightarrow{ \ y \ } \widehat{\mathcal C} \longrightarrow \widehat{\mathcal C}_c \longrightarrow \widehat{({\mathcal C}_c)}_{J_c}
	$$
	où $\widehat{\mathcal C} \to \widehat{\mathcal C}_c$ désigne le foncteur de restriction des préfaisceaux sur ${\mathcal C}$ à la sous-catégorie ${\mathcal C}_c \hookrightarrow {\mathcal C}$.
	
	Il en résulte que l'objet
	$$
	\ell (X) \quad \mbox{de} \quad \widehat{\mathcal C}_J
	$$
	admet une famille globalement épimorphique de morphismes
	$$
	\ell (X_i) \longrightarrow \ell (X) \, , \quad i \in I \, ,
	$$
	dont les sources sont les images par
	$$
	{\mathcal C}_c \longrightarrow \widehat{({\mathcal C}_c)}_{J_c} \xrightarrow{ \ \sim \ } \widehat{\mathcal C}_J
	$$
	d'objets $X_i$ de ${\mathcal C}_c$, puis que pour chaque morphisme
	$$
	\ell (X_i) \longrightarrow \ell (X) \, , \quad i \in I \, ,
	$$
	il existe un précrible $J_c$-couvrant de $X_i$ dans ${\mathcal C}_c$
	$$
	(X_{i,j} \longrightarrow X_i)_{j \in I_i} 
	$$
	tel que tous les composés
	$$
	\ell (X_{i,j}) \longrightarrow \ell (X_i) \longrightarrow \ell (X) \, , \quad i \in I \, , \ j \in I_i \, ,
	$$
	soient les images par $\ell$ de morphismes de ${\mathcal C}$
	$$
	X_{i,j}  \longrightarrow X \, .
	$$
	\item Le morphisme $X \xrightarrow{ \ y \ } Y$ de ${\mathcal C}$ induit un morphisme
	$$
	\ell (X) \xrightarrow{ \ \ell (y) \ } \ell (Y) \quad \mbox{de} \quad \widehat{\mathcal C}_J \xrightarrow{ \ \sim \ } \widehat{({\mathcal C}_c)}_{J_c} \, .
	$$
	Il existe un précrible $J_c$-couvrant de $X$ dans ${\mathcal C}_c$
	$$
	\left( X_i \xrightarrow{ \ x_i \ } X \right)_{i \in I}
	$$
	tel que les morphismes composés
	$$
	\ell (y) \circ \ell (x_i) : \ell (X_i) \longrightarrow \ell (X) \longrightarrow \ell (Y) \, , \quad i \in I \, ,
	$$
	se relèvent en des morphismes de ${\mathcal C}_c$
	$$
	y_i : X_i \longrightarrow Y \, , \quad i \in I \, .
	$$
	Pour chaque indice $i \in I$, les deux morphismes de ${\mathcal C}$
	$$
	\xymatrix{ X_i \dar[rr]^{^{^{\mbox{\scriptsize $y \circ x_i$}}}}_{y_i} && Y}
	$$
	ont même image par le foncteur
	$$
	{\mathcal C} \longrightarrow \widehat{\mathcal C} \longrightarrow \widehat{\mathcal C}_J \xrightarrow{ \ \sim \ } \widehat{({\mathcal C}_c)}_{J_c}
	$$
	ou, ce qui revient au même, par le foncteur
	$$
	{\mathcal C} \longrightarrow \widehat{\mathcal C} \longrightarrow \widehat{\mathcal C}_c \longrightarrow \widehat{({\mathcal C}_c)}_{J_c} \, .
	$$
	Par conséquent, chaque objet $X_i$ de ${\mathcal C}_c$ admet un précrible $J_c$-couvrant dans ${\mathcal C}_c$
	$$
	\left( X_{i,j} \xrightarrow{ \ x_{i,j} \ } X_i \right)_{j \in I_i}
	$$
	tel que les deux morphismes
	$$
	\xymatrix{ X_i \dar[rr]^{^{^{\mbox{\scriptsize $y \circ x_i$}}}}_{y_i} && Y}
	$$
	aient le même composé avec chaque $X_{i,j} \xrightarrow{ \, x_{i,j} \, } X_i$, $j \in I_i$.
	
	Cela termine la démonstration du théorème.
\end{listeisansmarge}
\end{demo}

Réciproquement, on a:

\begin{prop}\label{prop3.3.11}
	Soit ${\mathbb T}$ une théorie géométrique du premier ordre de signature $\Sigma$.
	
	Soit ${\mathbb T}_c$ une théorie cartésienne de même signature $\Sigma$ et dont tous les axiomes sont démontrables dans ${\mathbb T}$.
	
	Soit ${\mathcal C}$ la théorie syntactique cartésienne de la théorie ${\mathbb T}_c$.
	
	Alors:
	
	\begin{listeimarge}
	\item La catégorie ${\mathcal C}$ est petite et cartésienne.
	
	\item Les éléments du vocabulaire $\Sigma$ s'interprètent dans ${\mathcal C}$, ainsi que les formules de Horn (ou plus généralement ${\mathbb T}_c$-cartésiennes)
	$$
	\varphi (\vec x) \quad \mbox{en des suites de variables $\vec x = (x_1^{A_1} , \cdots , x_n^{A_n})$,}
	$$
	sous la forme de sous-objets
	$$
	U_{{\mathbb T}_c}^{\mathcal C} \, \varphi (\vec x) \ {\lhook\joinrel\relbar\joinrel\rightarrow} \ U_{{\mathbb T}_c}^{\mathcal C} \, A_1 \times \cdots \times U_{{\mathbb T}_c}^{\mathcal C} \, A_n \, .
	$$
	\item Le topos classifiant ${\mathcal E}_{\mathbb T}$ de ${\mathbb T}$ se présente comme topos des faisceaux
	$$
	\widehat{\mathcal C}_J \xrightarrow{ \ \sim \ } {\mathcal E}_{\mathbb T}
	$$
	sur la catégorie ${\mathcal C}$ munie de la topologie $J$ engendrée par les familles de précribles
	$$
	P_{\varphi , \psi , k} = \left( U_{{\mathbb T}_c}^{\mathcal C} (\varphi_k (\vec x , \vec x_k) \wedge \psi_{k'} (\vec x , \vec y_{k'})) \longrightarrow U_{{\mathbb T}_c}^{\mathcal C} \, \varphi_k (\vec x , \vec x_k) \right)_{k' \in K'} \, , \quad k \in K \, ,
	$$
	associés aux axiomes de définition de la théorie ${\mathbb T}$ comme quotient de ${\mathbb T}_c$
	$$
	\varphi (\vec x) \vdash \psi (\vec x)
	$$
	et à leurs décompositions en termes de formules de Horn $\varphi_k$ et $\psi_{k'}$
	$$
	\varphi (\vec x) = \bigvee_{k \in K} (\exists \, \vec x_k) \, \varphi_k (\vec x , \vec x_k) \, ,
	$$
	$$
	\psi (\vec x) = \bigvee_{k' \in K'} (\exists \, \vec y_{k'}) \, \psi_{k'} (\vec x , \vec y_{k'}) \, .
	$$
\end{listeimarge}
\end{prop}

\begin{demo}
	\begin{listeisansmarge}
	\item reprend le contenu des remarques qui suivent la définition \ref{defn1.1.33}.
	
	Pour la théorie cartésienne ${\mathbb T}_c$, on a des plongements entre catégories syntactiques
	$$
	{\mathcal C} = {\mathcal C}_{{\mathbb T}_c}^{\rm car} \ {\lhook\joinrel\relbar\joinrel\rightarrow} \ {\mathcal C}_{{\mathbb T}_c}^{\rm reg} \ {\lhook\joinrel\relbar\joinrel\rightarrow} \ {\mathcal C}_{\mathbb T}^{\rm coh}
	$$
	qui sont toutes petites puisque les formules cohérentes de $\Sigma$, considérées à substitution près des variables, forment un ensemble.
	
	Toutes les catégories syntactiques ont des limites finies arbitraires toujours construites de la manière rappelée dans la preuve du lemme \ref{lem1.1.26} (ii).
	
	\item Les éléments de vocabulaire de $\Sigma$ s'interprètent dans toutes les catégories syntactiques de toutes les théories de signature $\Sigma$.
	
	\item D'après le théorème \ref{thm1.1.35}, la théorie cartésienne ${\mathbb T}_c$ admet pour topos classifiant le topos des préfaisceaux sur la catégorie syntactique cartésienne ${\mathcal C} = {\mathcal C}_{{\mathbb T}_c}^{\rm car}$
	$$
	{\mathcal E}_{{\mathbb T}_c} = \widehat{\mathcal C} \, .
	$$
	
	D'après le théorème \ref{thm1.2.10}, la théorie ${\mathbb T}$ vue comme une théorie quotient de ${\mathbb T}_c$ a pour topos classifiant un sous-topos
	$$
	{\mathcal E}_{\mathbb T} \ {\lhook\joinrel\relbar\joinrel\rightarrow} \ {\mathcal E}_{{\mathbb T}_c}  = \widehat{\mathcal C}
	$$
	nécessairement défini d'après le théorème \ref{thm1.2.8} par une topologie $J$ de ${\mathcal C}$.
	
	D'après le théorème \ref{thm3.3.7}, la topologie $J$ de ${\mathcal C}$ est engendrée par les familles de précribles
	$$
	P_{\varphi , \psi , k}
	$$
	associées aux axiomes
	$$
	\varphi (\vec x) \vdash \psi (\vec x)
	$$
	qui définissent ${\mathbb T}$ comme quotient de ${\mathbb T}_c$ et à leurs décompositions en termes de formules de Horn.
\end{listeisansmarge}
\end{demo}

\subsection{Interprétation sémantique des catégories syntactiques cartésiennes}\label{ssec3.3.5}

On rappelle le théorème suivant:

\begin{thm}\label{thm3.3.12}
	Soit ${\mathbb T}$ une théorie géométrique du premier ordre dont le topos classifiant se présente comme topos des préfaisceaux sur une catégorie ${\mathcal C}$ essentiellement petite
	$$
	\widehat{\mathcal C} \xrightarrow{ \ \sim \ } {\mathcal E}_{\mathbb T} \, .
	$$
	Alors:
	
	\begin{listeimarge}
	\item La catégorie des modèles ensemblistes de ${\mathbb T}$
	$$
	{\mathbb T}\mbox{-{\rm mod}} \, ({\rm Ens})
	$$
	est équivalente à la catégorie
	$$
	{\rm Ind} ({\mathcal C}^{\rm op})
	$$
	constituée des préfaisceaux sur ${\mathcal C}^{\rm op}$
	$$
	\rho : {\mathcal C} \longrightarrow {\rm Ens}
	$$
	qui s'écrivent comme des colimites filtrantes de préfaisceaux représentables
	$$
	y(X) = {\rm Hom}_{\mathcal C} (X,\bullet)
	$$
	associés à des objets $X$ de ${\mathcal C}$.
	
	\item Un modèle ensembliste $M$ de ${\mathbb T}$ est ``finiment présentable'', au sens que le foncteur
	$$
	{\rm Hom}(M,\bullet) : {\mathbb T}\mbox{-{\rm mod}} \, ({\rm Ens}) \longrightarrow {\rm Ens}
	$$
	respecte les colimites filtrantes, si et seulement si le préfaisceau
	$$
	\rho_M : {\mathcal C} \longrightarrow {\rm Ens}
	$$
	est défini par la formule
	$$
	\rho_M = {\rm eg} \left( \xymatrix{ y(X) \dar[r]^{^{^{\mbox{\scriptsize $y(e)$}}}}_{\rm id} & y(X)}\right)
	$$
	pour un objet $X$ de ${\mathcal C}$ muni d'un ``idempotent'' c'est-à-dire un endomorphisme
	$$
	e : X \longrightarrow X
	$$
	tel que
	$$
	e = e \circ e \, .
	$$
	\item La sous-catégorie pleine de
	$$
	{\mathbb T}\mbox{-{\rm mod}} \, ({\rm Ens})
	$$
	constituée des modèles qui sont ``finiment présentables'' est équivalente à la ``complétion Karoubienne'' de ${\mathcal C}^{\rm op}$
	$$
	{\rm Kar} ({\mathcal C}^{\rm op}) = {\rm Kar} ({\mathcal C})^{\rm op}
	$$
	dont
	
	\medskip
	
	$\left\{\begin{matrix}
		\bullet &\mbox{les objets sont les paires $(X,e)$ constituées d'un objet $X$ de ${\mathcal C}^{\rm op}$} \hfill \\
		&\mbox{muni d'un endomorphisme idempotent $e : X \to X$,} \hfill \\
		\bullet &\mbox{les morphismes $(X',e') \to (X,e)$ sont les morphismes $u : X' \to X$} \hfill \\
		&\mbox{de ${\mathcal C}^{\rm op}$ tels que $e \circ u = u = u \circ e'$.} \hfill
	\end{matrix}\right.$
	
	\medskip
	
	\item On a une équivalence de topos
	$$
	\widehat{{\rm Kar} ({\mathcal C})} \xrightarrow{ \ \sim \ } \widehat{\mathcal C}
	$$
	si bien que le topos classifiant ${\mathcal E}_{\mathbb T}$ de ${\mathbb T}$ se présente comme le topos des préfaisceaux
	$$
	\widehat{{\mathcal M}^{\rm op}}
	$$
	sur l'opposée de la catégorie des modèles finiment présentables de ${\mathbb T}$
	$$
	{\mathcal M} = {\mathbb T}\mbox{-{\rm mod}} \, ({\rm Ens})_{\rm fp} \, .
	$$
\end{listeimarge}
\end{thm}

\begin{remark}
	Si la catégorie ${\mathcal C}$ est ``Karoubienne'' c'est-à-dire si tout idempotent $e : X \to X$ d'un objet $X$ de ${\mathcal C}$ définit un sous-objet
	$$
	X_e = {\rm eg} \left( \xymatrix{ X \dar[r]^{^{^{\mbox{\scriptsize $e$}}}}_{\rm Id} & X}\right),
	$$
	alors le plongement
	$$
	{\mathcal C} \ {\lhook\joinrel\relbar\joinrel\rightarrow} \ {\rm Kar} ({\mathcal C})
	$$
	est une équivalence et la catégorie ${\mathcal C}^{\rm op}$ est équivalente à la catégorie ${\mathbb T}\mbox{-mod} ({\rm Ens})_{\rm fp}$ des modèles finiment présentables de ${\mathcal C}$.
\end{remark}

\begin{esquissedemo}
	\begin{listeisansmarge}
	\item D'après l'équivalence
	$$
	\widehat{\mathcal C} \xrightarrow{ \ \sim \ } {\mathcal E}_{\mathbb T} \, ,
	$$
	la catégorie des modèles ensemblistes de ${\mathbb T}$ est equivalente à celle des points du topos $\widehat{\mathcal C}$ c'est-à-dire des morphismes de topos
	$$
	{\rm Ens} \longrightarrow \widehat{\mathcal C} \, .
	$$
	D'après le théorème \ref{thm1.1.23}, elle est encore équivalente à la catégorie des foncteurs
	$$
	\rho : {\mathcal C} \longrightarrow {\rm Ens}
	$$
	qui sont plats, c'est-à-dire dont l'unique prolongement respectant les colimites
	$$
	\widehat\rho : \widehat{\mathcal C} \longrightarrow {\rm Ens}
	$$
	respecte aussi les limites finies.
	
	Pour tout objet $X$ de ${\mathcal C}$, le foncteur d'évaluation des préfaisceaux en $X$
	$$
	\begin{matrix}
		\widehat{\mathcal C} &\longrightarrow &{\rm Ens}, \\
		P &\longmapsto &P(X)
	\end{matrix}
	$$
	respecte toutes les colimites et toutes les limites, donc définit un foncteur plat
	$$
	{\rm Hom} \, (X,\bullet) : {\mathcal C} \longrightarrow {\rm Ens}.
	$$
	Les colimites filtrantes dans ${\rm Ens}$ respectent les limites finies, donc toutes les colimites filtrantes de foncteurs de la forme
	$$
	{\rm Hom} \, (X,\bullet) : {\mathcal C} \longrightarrow {\rm Ens}
	$$
	sont encore des foncteurs plats.
	
	Inversement, tout préfaisceau sur ${\mathcal C}^{\rm op}$
	$$
	\rho : {\mathcal C} \longrightarrow {\rm Ens}
	$$
	s'écrit d'après le lemme \ref{lem1.1.14} comme une colimite
	$$
	\rho = \varinjlim_{(X,x) \in {\mathcal C}^{\rm op} / \rho} {\rm Hom} (X,\bullet)
	$$
	indexée par la catégorie ${\mathcal C}^{\rm op} / \rho$ des éléments de $\rho$ c'est-à-dire des paires $(X,x)$ constituées d'un objet $X$ de ${\mathcal C}$ et d'un élément $x \in \rho (X)$.
	
	On montre que si un tel préfaisceau $\rho : {\mathcal C} \to {\rm Ens}$ est plat, alors la catégorie de ses éléments
	$$
	{\mathcal C}^{\rm op} / \rho
	$$
	est filtrante et donc $\rho$ s'écrit comme une colimite filtrante de préfaisceaux représentables
	$$
	{\rm Hom} (X,\bullet) \, .
	$$
	\item Si $e : X \to X$ est un idempotent d'un objet $X$ de ${\mathcal C}$, la catégorie réduite à $X$ et ses deux endomorphismes
	$$
	\xymatrix{ X \dar[r]^{^{^{\mbox{\scriptsize $e$}}}}_{\rm id} & X}
	$$
	est filtrante, et sa colimite est le foncteur plat
	$$
	\begin{matrix}
		\rho : {\mathcal C} &\longrightarrow &{\rm Ens}, \hfill \\
		\hfill X' &\longmapsto &\{x \in {\rm Hom} (X,X') \mid x \circ e = x \} \, .
	\end{matrix}
	$$
	Pour tout préfaisceau
	$$
	P : {\mathcal C} \longrightarrow {\rm Ens},
	$$
	on a
	$$
	{\rm Hom} (\rho , P) = \{p \in P(X) \mid P (e)(p) = p \}
	$$
	et donc le foncteur 
	$$
	{\rm Hom} (\rho , \bullet)
	$$
	respecte les colimites filtrantes dans ${\rm Ind} ({\mathcal C}^{\rm op})$.
	
	Réciproquement, si $\rho : {\mathcal C} \to {\rm Ens}$ est une colimite filtrante de préfaisceaux représentables
	$$
	{\rm Hom} (X,\bullet)
	$$
	et si le foncteur
	$$
	{\rm Hom} (\rho , \bullet)
	$$
	respecte les colimites filtrantes, alors le morphisme
	$$
	{\rm id} : \rho \longrightarrow \rho
	$$
	se factorise nécessairement sous la forme
	$$
	\rho \xrightarrow{ \ s \ } {\rm Hom} (X,\bullet) \xrightarrow{ \ p \ } \rho
	$$
	pour un objet $X$ de ${\mathcal C}$.
	
	Alors l'équation
	$$
	p \circ s = {\rm id}_{\rho}
	$$
	implique
	$$
	(s \circ p) \circ (s \circ p) = s \circ p
	$$
	et l'idempotent
	$$
	s \circ p : {\rm Hom} (X, \bullet) \longrightarrow {\rm Hom} (X, \bullet)
	$$
	provient d'un endomorphisme idempotent
	$$
	e : X \longrightarrow X
	$$
	tel que
	$$
	\rho = {\rm eg} \left( \xymatrix{ {\rm Hom} (X, \bullet) \dar[r]^{^{^{\mbox{\scriptsize $y(e)$}}}}_{\rm id} & {\rm Hom} (X, \bullet)} \right).
	$$
	\item résulte de (ii).
	
	\item résulte de (iii) puisque tout préfaisceau
	$$
	P : {\mathcal C}^{\rm op} \longrightarrow {\rm Ens}
	$$
	se prolonge canoniquement en un préfaisceau
	$$
	\begin{matrix}
		P : {\rm Kar} ({\mathcal C})^{\rm op} &\longrightarrow &{\rm Ens}, \hfill \\
		\hfill (X,e) &\longmapsto &\{p \in P (X) \mid P(e)(p) = p \} \, .
	\end{matrix}
	$$
\end{listeisansmarge}
\end{esquissedemo}

On déduit de ce théorème:

\begin{cor}\label{cor3.3.13}
	Soit ${\mathbb T}$ une théorie du premier ordre qui est cartésienne.
	
	Soit ${\mathcal C} = {\mathcal C}_{\mathbb T}^{\rm car}$ la catégorie syntactique cartésienne ${\mathbb T}$.
	
	Soit ${\mathcal M} = {\mathbb T}\mbox{-{\rm mod}} \, ({\rm Ens})_{\rm fp}$ la catégorie de ses modèles ensemblistes finiment présentables.
	
	Alors:
	
	\begin{listeimarge}
	\item La catégorie ${\mathcal M}$ est équivalente à la catégorie opposée ${\mathcal C}^{\rm op}$.
	
	\item L'équivalence
	$$
	{\mathcal C}^{\rm op} \longrightarrow {\mathcal M}
	$$
	associe à tout objet de ${\mathcal C}$, c'est-à-dire à toute formule cartésienne $\varphi (\vec x)$ en une suite de variables $\vec x = (x_1^{A_1} , \cdots , x_n^{A_n})$, un modèle
	$$
	M_{\varphi}
	$$
	qui est ``présenté'' par la formule $\varphi$ au sens que, pour tout modèle
	$$
	M \quad \mbox{dans} \quad {\mathbb T}\mbox{-{\rm mod}} \, ({\rm Ens}) \, ,
	$$
	l'ensemble
	$$
	{\rm Hom} (M_{\varphi} , M)
	$$
	s'identifie au  sous-ensemble défini par la formule $\varphi (\vec x)$ dans $M$
	$$
	M\varphi (\vec x) \subseteq MA_1 \times \cdots \times MA_n \, .
	$$
	\end{listeimarge}
\end{cor}

%\newpage

\begin{demo}
	\begin{listeisansmarge}
	\item La catégorie syntactique cartésienne ${\mathcal C} = {\mathcal C}_{\mathbb T}^{\rm car}$ est Karoubienne puisqu'elle est cartésienne.
	
	D'après la remarque qui suit le théorème \ref{thm3.3.12}, on a donc une équivalence naturelle
	$$
	{\mathcal C} \xrightarrow{ \ \sim \ } {\mathcal M}^{\rm op} \, .
	$$
	\item Les modèles ensemblistes de ${\mathbb T}$ correspondent aux morphismes de topos
	$$
	{\rm Ens} \longrightarrow \widehat{\mathcal C}
	$$
	et donc aux foncteurs plats
	$$
	{\mathcal C} \longrightarrow {\rm Ens} \, .
	$$
	Dans cette correspondance, les modèles $M$ sont vus comme les foncteurs
	$$
	\begin{matrix}
		{\mathcal C} = {\mathcal C}_{\mathbb T}^{\rm car} &\longrightarrow &{\rm Ens} \, , \hfill \\
		\hfill \varphi (\vec x) &\longmapsto &M\varphi (\vec x) \, .
	\end{matrix}
	$$
\end{listeisansmarge}
\end{demo}

Ce corollaire implique que les catégories syntactiques cartésiennes des théories cartésiennes ont des ensembles de morphismes entre paires d'objets qui ne peuvent pas être calculés en général, même dans le cas de théories définies simplement.

Par exemple, si ${\mathbb A}$ est la théorie algébrique (et donc cartésienne) des anneaux commutatifs, alors pour tout modèle $A$ présenté par des équations polynomiales en des variables $X_1 , \cdots , X_n$
$$
P_i (X_1 , \cdots , X_n) = 0 \, , \quad 1 \leq i \leq m \, ,
$$
et tout modèle $A'$, l'ensemble
$$
{\rm Hom} (A,A')
$$
s'identifie à l'ensemble des solutions
$$
\left\{ (a_1 , \cdots , a_n) \in A'^n \mid P_i (a_1 , \cdots , a_n) = 0 \, , \ 1 \leq i \leq m \right\}
$$
qui est rarement calculable.

\subsection{Les catégories syntactiques cartésiennes des théories sans axiomes}\label{ssec3.3.6}

D'après la définition \ref{defn1.1.17} (i), un terme à valeurs dans une sorte $B$ en une suite de variables $\vec x = (x_1^{A_1} , \cdots , x_n^{A_n})$ associées à des sortes $A_1 , \cdots , A_n$ d'une signature $\Sigma$ peut être

\medskip

$\left\{\begin{matrix}
	\bullet &\mbox{ou bien l'une des variables} \hfill \\
	&x_i^{A_i} \\
	&\mbox{pour un indice $i \in \{1, \cdots , n \}$ tel que $A_i = B$,} \hfill \\
	\bullet &\mbox{ou bien de la forme} \hfill \\
	&f \left(g_1^{B_1} (\vec x) , \cdots , g_m^{B_m} (\vec x) \right) \\
	&\mbox{pour un symbole de fonction de $\Sigma$ à valeurs dans la sorte $B$} \hfill \\
	&f : B_1 \cdots B_m \longrightarrow B \\
	&\mbox{et des termes en $\vec x = (x_1^{A_1} , \cdots , x_n^{A_n})$ déjà définis} \hfill \\
	&g_i^{B_i} (\vec x) \, , \quad 1 \leq i \leq m \, , \\
	&\mbox{à valeurs dans les sortes $B_i$ de $\Sigma$.} \hfill
\end{matrix}\right.$

\medskip

Toute famille d'équations entre des termes en une suite de variables $\vec x$ définit une relation d'équivalence entre les termes en $\vec x$ ainsi qu'entre les formules atomiques en $\vec x$:

\begin{defn}\label{defn3.3.14}
	Soit $\Sigma$ une signature constituée de sortes, de symboles de fonctions et de symboles de relations.
	
	Soit $\vec x = (x_1^{A_1} , \cdots , x_n^{A_n})$ une suite de variables associées à des sortes $A_1 , \cdots  , A_n$ de $\Sigma$.
	
	Soit une famille d'équations
	$$
	f_i^{B_i} (\vec x) = f'^{B_i}_i (\vec x) \, , \quad i \in I \, ,
	$$
	entre des paires de termes en $\vec x$ à valeurs dans les mêmes sortes $B_i$ de $\Sigma$.
	
	Alors:
	
	\begin{listeimarge}
	\item On appelle relation d'équivalence entre termes en $\vec x$ à valeurs dans les mêmes sortes
	$$
	f^B (\vec x) \sim f'^{B} (\vec x)
	$$
	la plus petite relation d'équivalence telle que
	
	\medskip
	
	$\left\{\begin{matrix}
		\bullet &\mbox{pour tout $i \in I$, on a} \hfill \\
		&f_i^{B_i} (\vec x) \sim f'^{B_i}_i (\vec x) \, , \\
		\bullet &\mbox{pour tout symbole de fonction de $\Sigma$} \hfill \\
		&f : B_1 \cdots B_m \longrightarrow B \\
		&\mbox{et toute paire de familles de termes en $\vec x$ à valeurs dans les sortes $B_i$, $1 \leq i \leq m$,} \hfill \\
		&g_1^{B_i} (\vec x) \, , \quad 1 \leq i \leq m \, , \\
		&\mbox{et} \hfill \\
		&g'^B_i (\vec x) \, , \quad 1 \leq i \leq m \, , \\
		&\mbox{on a} \hfill \\
		&f \left( g_1^{B_1} (\vec x) , \cdots , g_m^{B_m} (\vec x) \right) \sim \left( g'^{B_1}_1 (\vec x) , \cdots , g'^{B_m}_m (\vec x) \right) \\
		&\mbox{si} \hfill \\
		&g_i^{B_i} (\vec x) \sim g'^{B_i}_i (\vec x) \, , \quad 1 \leq i \leq m \, .
	\end{matrix}\right.$
	
	\medskip
	
	\item Pour n'importe quel symbole de relation de $\Sigma$
	$$
	R \rightarrowtail B_1 \cdots B_m \, ,
	$$
	deux formules atomiques en $\vec x$ de la forme
	$$
	R \left( g_1^{B_1} (\vec x) , \cdots , g_m^{B_m} (\vec x) \right) \quad \mbox{et} \quad R \left( g'^{B_1}_1 (\vec x) , \cdots , g'^{B_m}_m (\vec x) \right)
	$$
	sont dites équivalentes si et seulement si on a
	$$
	g_i^{B_i} (\vec x) \sim g'^{B_i}_i (\vec x) \, , \quad 1 \leq i \leq m \, .
	$$
\end{listeimarge}
\end{defn}

Cette définition permet de décrire les catégories syntactiques cartésiennes des théories sans axiomes:

\begin{prop}\label{prop3.3.15}
	Soit $\Sigma$ une théorie géométrique du premier ordre réduite à sa signature, sans axiome.
	
	Alors la catégorie syntactique cartésienne de $\Sigma$
	$$
	{\mathcal C}_{\Sigma}^{\rm car}
	$$
	se décrit de la manière suivante:
	
	\begin{listeimarge}
	\item Ses objets sont les formules de Horn en des suites de variables $\vec x = (x_1^{A_1} , \cdots , x_n^{A_n})$, considérées à substitution près des variables $x_i^{A_i}$, $1 \leq i \leq n$, par des variables de mêmes sortes associées $A_i$.
	
	Autrement dit, ce sont les conjonctions finies
	$$
	\varphi (\vec x) = \bigwedge_{1 \leq j \leq m} \left( f_j^{B_j} (\vec x) = f'^{B_j}_j (\vec x)\right) \wedge \bigwedge_{1 \leq k \leq \ell} R_k \left( f_{k,1}^{B_{k,1}} (\vec x) , \cdots ,  f^{B_{k,m_k}}_{k,m_k} (\vec x)\right)
	$$
	d'équations entre des paires de termes en $\vec x$ à valeurs dans la même sorte et de symboles de relations de $\Sigma$ appliqués à des suites de termes en $\vec x$.
	
	\item Les morphismes de ${\mathcal C}_{\Sigma}^{\rm car}$
	$$
	\varphi (\vec x) \longrightarrow \psi (\vec y)
	$$
	entre deux formules de Horn $\varphi$ et $\psi$ en
	$$
	\vec x = (x_1^{A_1} , \cdots , x_n^{A_n})
	$$
	et
	$$
	\vec y = (y_1^{A'_1} , \cdots , y_{n'}^{A'_{n'}} )
	$$
	qui ont les formes
	$$
	\varphi (\vec x) = \bigwedge_{1 \leq j \leq m} \left( f_j^{B_j} (\vec x) = f'^{B_j}_j (\vec x)\right) \wedge \bigwedge_{1 \leq k \leq \ell} R_k \left( f_{k,1}^{B_{k,1}} (\vec x) , \cdots ,  f^{B_{k,m_k}}_{k,m_k} (\vec x)\right) ,
	$$
	$$
	\psi (\vec y) = \bigwedge_{1 \leq j \leq m'} \left( g_j^{C_j} (\vec y) = g'^{C_j}_j (\vec y)\right) \wedge \bigwedge_{1 \leq k \leq \ell'} S_k \left( g_{k,1}^{C_{k,1}} (\vec y) , \cdots , g^{C_{k,m'_k}}_{k,m'_k} (\vec y)\right)
	$$
	sont les classes d'équivalence de suites de termes en $\vec x$ à valeurs dans les sortes $A'_1 , \cdots , A'_{n'}$
	$$
	\left( h_1^{A'_1} (\vec x) , \cdots , h_{n'}^{A'_n} (\vec x) \right)
	$$
	telles que
	
	\medskip
	
	$\left\{\begin{matrix}
		\bullet &\mbox{pour tout indice $j$, $1\leq j \leq m'$, il y a équivalence des termes en $\vec x$} \hfill \\
		&g_j^{C_j} \left( h_1^{A'_1} (\vec x) , \cdots , h_{n'}^{A'_{n'}} (\vec x) \right) \sim g'^{C_j}_j \left( h_1^{A'_1} (\vec x) , \cdots , h_{n}^{A'_{n'}} (\vec x) \right), \\
		\bullet &\mbox{pour tout indice $k$, $1 \leq k \leq \ell'$, la formule atomique en $\vec x$ obtenue par substitution de} \hfill \\
		&\vec y = \left( h_1^{A'_1} (\vec x) , \cdots , h_{n'}^{A'_{n'}} (\vec x) \right) \\
		&\mbox{dans} \hfill \\
		&S_k \left( g_{k,1}^{C_{k,1}} (\vec y) , \cdots ,  g^{C_{k,m'_k}}_{k,m'_k} (\vec y)\right) \\
		&\mbox{est équivalente à l'une au moins des formules atomiques} \hfill \\
		&R_{k'} \left( f_{k',1}^{B_{k',1}} (\vec x) , \cdots ,  f^{B_{k',m_{k'}}}_{k',m_{k'}} (\vec x)\right) , \quad 1 \leq k' \leq \ell \, ,
	\end{matrix}\right.$
	
	\medskip
	
	\noindent pour la relation d'équivalence entre termes en $\vec x$ ou entre formules atomiques en $\vec x$ qui est définie par les équations
	$$
	f_j^{B_j} (\vec x) = f'^{B_j}_j (\vec x) \, , \quad 1 \leq j \leq m \, .
	$$
	(iii) La composition des morphismes de ${\mathcal C}_{\Sigma}^{\rm car}$
	$$
	\varphi (\vec x) \xrightarrow{ \ \left( h_1^{A'_1} (\vec x) , \cdots , h_{n'}^{A'_{n'}} (\vec x) \right) \ } \psi (\vec y) \xrightarrow{ \ \left( g_1^{A''_1} (\vec y) , \cdots , g_{n''}^{A''_{n''}} (\vec y) \right) \ } \chi (\vec z)
	$$
	est définie par la substitution
	$$
	\vec y = \left( h_1^{A'_1} (\vec x) , \cdots , h_{n'}^{A'_{n'}} (\vec x) \right)
	$$
	dans les termes 
	$$
	g_i^{A''_i} (\vec y) \, , \quad 1 \leq i \leq n'' \, .
	$$
\end{listeimarge}
\end{prop}

\begin{demo}
	Comme la théorie $\Sigma$ est sans axiome, ses seules formules cartésiennes sont les formules de Horn et ses seules formules démontrablement fonctionnelles
	$$
	\varphi (\vec x) \xrightarrow{ \ \theta (\vec x , \vec y) \ } \psi (\vec y)
	$$
	reliant deux formules de Horn
	$$
	\varphi (\vec x) = \varphi (x_1^{A_1} , \cdots , x_n^{A_n})
	$$
	et
	$$
	\psi (\vec y) = \psi (y_1^{A'_1} , \cdots , y_{n'}^{A'_{n'}})
	$$
	sont les formules de la forme
	$$
	\bigwedge_{1 \leq i \leq n'} \left( y_i^{A'_i} = h_i^{A'_i} (\vec x) \right)
	$$
	pour des termes $h_i^{A'_i} (\vec x)$ en $\vec x$ à valeurs dans les sortes $A'_i$ tels que le séquent
	$$
	\varphi (\vec x) \vdash \psi \left( h_1^{A'_1} (\vec x) , \cdots , h_{n'}^{A'_n} (\vec x) \right)
	$$
	soit démontrable.
	
	Cela revient à demander que chaque équation entre termes constitutives de $\psi \left( h_1^{A'_1} (\vec x) , \cdots , h_{n'}^{A'_n} (\vec x) \right)$
	$$
	g_j^{C_j} \left( h_1^{A'_1} (\vec x) , \cdots , h_{n'}^{A'_{n'}} (\vec x) \right) = g'^{C_j}_j \left( h_1^{A'_1} (\vec x) , \cdots , h_{n'}^{A'_{n'}} (\vec x) \right)
	$$
	relie deux termes en $\vec x$ équivalents au sens de la relation engendrée par les équations constitutives de $\varphi (\vec x)$
	$$
	f_j^{B_j} (\vec x) = f'^{B_j}_j (\vec x) \, , \quad 1 \leq j \leq m
	$$
	et que chaque formule atomique en $\vec x$ obtenue par substitution de 
	$$
	\vec y  = \left( h_1^{A'_1} (\vec x) , \cdots , h_{n'}^{A'_n} (\vec x) \right)
	$$
	dans une relation constitutive de $\psi (\vec y)$
	$$
	S_k \left( g_{k,1}^{C_{k,1}} (\vec y) , \cdots ,  g^{C_{k,m'_k}}_{k,m'_k} (\vec y)\right)
	$$
	figure à équivalence près parmi les relations constitutives de $\varphi (\vec x)$
	$$
	R_{k'} \left( f_{k',1}^{B_{k',1}} (\vec x) , \cdots , f_{k',m_{k'}}^{B_{k' , m_{k'}}} (\vec x) \right) .
	$$
	
	Enfin, deux morphismes
	$$
	\xymatrix{
		\varphi (\vec x) \dar[r] &\psi (\vec y)
	}
	$$
	définis par deux familles de termes
	$$
	\left(h_1^{A'_1} (\vec x) , \cdots , h_{n'}^{A'_{n'}} (\vec x) \right) \quad \mbox{et} \quad \left( h'^{A'_1}_1 (\vec x) , \cdots , h'^{A'_{n'}}_{n'} (\vec x) \right)
	$$
	sont démontrablement égaux si et seulement si ces familles sont équivalentes au sens induit par les équations
	$$
	f_j^{B_j} (\vec x) = f'^{B_j}_j (\vec x) \, , \quad 1 \leq j \leq m \, .
	$$
	(iii) est un cas particulier de la règle générale de composition des morphismes des catégories syntactiques définis par des formules démontrablement fonctionnelles.
	
\end{demo}

Dans la pratique, calculer des classes d'équivalence est difficile.

Cependant, les relations d'équivalence disparaissent de l'énoncé de la proposition précédente dans le cas de signatures $\Sigma$ qui n'ont pas de symboles de fonctions:

\begin{cor}\label{cor3.3.16}
	Soit $\Sigma$ une signature réduite à des sortes et des symboles de relations, sans symboles de fonctions.
	
	Alors la catégorie syntactique cartésienne de la théorie sans axiome de signature $\Sigma$
	$$
	{\mathcal C}_{\Sigma}^{\rm car}
	$$
	se décrit de la manière suivante:
	
	\begin{listeimarge}
	\item Ses objets sont indexés par les données combinatoires constituées de
	
	$\left\{\begin{matrix}
		\bullet &\mbox{un sous-ensemble fini ${\mathcal A}$ de l'ensemble des sortes de $\Sigma$,} \hfill \\
		\bullet &\mbox{une application $A \mapsto m_A$ qui associe à chaque sorte $A \in {\mathcal A}$ une multiplicité $m_A \geq 1$,} \hfill \\
		\bullet &\mbox{un sous-ensemble fini ${\mathcal R}$ de l'ensemble des paires $(R,\alpha)$} \hfill \\
		&\mbox{consistant en un symbole de relation de $\Sigma$} \hfill \\
		&R \rightarrowtail A_1 \cdots A_{n_R} \quad \mbox{avec} \quad A_i \in {\mathcal A} \, , \ 1 \leq i \leq n_R \, , \\
		&\mbox{et une application} \hfill \\
		&\alpha : \{1 , \ldots , n_R \} \longrightarrow {\mathbb N} \\
		&\mbox{telle que} \hfill \\
		&1 \leq \alpha (i) \leq m_{A_i} \, , \quad 1 \leq i \leq n_R \, .
	\end{matrix}\right.$
	
	Choisissant des noms de variables en nombre $m_A$ pour chaque $A \in {\mathcal A}$
	$$
	x_i^A \, , \quad 1 \leq i \leq m_A \, , \quad A \in {\mathcal A} \, ,
	$$
	l'objet correspondant est la formule de Horn en la famille finie de variables $\vec x_{{\mathcal A} , m_{\bullet}} = (x_i^A)_{A \in {\mathcal A} , 1 \leq i \leq m_A}$ définie comme la conjonction
	$$
	\varphi_{\mathcal R} (\vec x_{{\mathcal A},m_{\bullet}}) = \bigwedge_{(R,\alpha) \in {\mathcal R}} R \left(x_{\alpha (1)}^{A_1} , \cdots , x_{\alpha (n_R)}^{A_{n_R}} \right).
	$$
	\item Les morphismes d'un objet de la forme
	$$
	\varphi_{\mathcal R} (\vec x_{{\mathcal A},m_{\bullet}}) = \bigwedge_{(R,\alpha) \in {\mathcal R}} R \left(x_{\alpha (1)}^{A_1} , \cdots , x_{\alpha (n_R)}^{A_{n_R}} \right)
	$$
	dans un objet de la forme d'une formule de Horn
	$$
	\varphi_{{\mathcal R}'} (\vec y_{{\mathcal A}',m'_{\bullet}}) = \bigwedge_{(R,\beta) \in {\mathcal R}'} R \left(y_{\beta (1)}^{A_1} , \cdots , y_{\beta (n_R)}^{A_{n_R}} \right)
	$$
	en une famille finie de variables $\vec y_{{\mathcal A}',m'_{\bullet}} = (y_i^A)_{A \in {\mathcal A}' , 1 \leq i \leq m'_A}$, ne peuvent exister que si ${\mathcal A}' \subseteq {\mathcal A}$.
	
	Si cette condition est satisfaite, ils sont indexés par les familles finies d'applications
	$$
	f_A : \{ 1,2,\cdots , m'_A \} \longrightarrow \{1,2,\cdots , m_A \} \, , \quad A \in {\mathcal A}' \subseteq {\mathcal A} \, ,
	$$
	telles que, pour tout élément
	$$
	(R \rightarrowtail A_1 \cdots A_{n_R} , \beta  : \{ 1, \cdots , n_R \} \longrightarrow {\mathbb N}) \in {\mathcal R}'
	$$
	avec donc $A_1 , \cdots , A_{n_R} \in {\mathcal A}' \subseteq {\mathcal A}$, l'application
	$$
	\begin{matrix}
		\alpha : \{ 1 , \cdots , n_R \} &\longrightarrow &{\mathbb N} \, , \hfill \\
		\hfill i &\longmapsto &f_{A_i} (\beta (i))
	\end{matrix}
	$$
	forme avec le symbole de relation $R \rightarrowtail A_1 \cdots A_{n_R}$ un élément
	$$
	(R,\alpha) \in {\mathcal R} \, .
	$$
	
	Le morphisme indexé par une telle famille d'applications
	$$
	\left( f_A : \{1,2,\cdots , m'_A\} \longrightarrow \{1,2,\cdots ,m_A\}\right)_{A \in {\mathcal A}}
	$$
	est défini par les formules
	$$
	y_i^A = x^A_{f_A(i)} \, , \quad A \in {\mathcal A}' \subseteq {\mathcal A} \, , \quad 1 \leq i \leq m'_A \, .
	$$
	\item Le composé de deux morphismes
	$$
	\varphi_{\mathcal R} (\vec x_{{\mathcal A},m_{\bullet}}) \longrightarrow \varphi_{{\mathcal R}'} (\vec y_{{\mathcal A}',m'_{\bullet}}) \longrightarrow \varphi_{{\mathcal R}''} (\vec z_{{\mathcal A}'',m''_{\bullet}})
	$$
	indexés par deux familles d'applications
	$$
	\left( f_A : \{1,2,\cdots , m'_A\} \longrightarrow \{1,2,\cdots ,m_A\}\right)_{A \in {\mathcal A}' \subseteq {\mathcal A}}
	$$
	et
	$$
	\left( g_A : \{1,2,\cdots , m''_A\} \longrightarrow \{1,2,\cdots ,m'_A\}\right)_{A \in {\mathcal A}'' \subseteq {\mathcal A}'}
	$$
	est le morphisme
	$$
	\varphi_{\mathcal R} (\vec x_{{\mathcal A},m_{\bullet}}) \longrightarrow \varphi_{{\mathcal R}''} (\vec z_{{\mathcal A}'',m''_{\bullet}})
	$$
	indexé par la famille des applications composées
	$$
	\left( f_A \circ g_A : \{1,2,\cdots , m''_A\} \longrightarrow \{1,2,\cdots ,m_A \} \right)_{A \in {\mathcal A}'' \subseteq {\mathcal A}} \, .
	$$
	\end{listeimarge}
\end{cor}

\begin{remark}
	Ainsi, les objets d'une telle catégorie syntactique cartésienne ${\mathcal C}_{\Sigma}^{\rm car}$ sont indexés par des données combinatoires finies $({\mathcal A} , m_{\bullet} , {\mathcal R})$ d'un type entièrement explicite.
	
	De plus, les morphismes entre deux objets indexés par des données $({\mathcal A} , m_{\bullet} , {\mathcal R})$ et $({\mathcal A}' , m'_{\bullet} , {\mathcal R}')$ forment un ensemble fini facile à calculer à partir de ces données.
\end{remark}

\begin{demo}
	Ce corollaire est un cas particulier de la proposition précédente puisque, si $\Sigma$ n'a pas de symbole de fonction, les seuls termes en une famille de variables
	$$
	\vec x = (x_1^{A_1} , \cdots , x_n^{A_n})
	$$
	sont les variables $x_i^{A_i}$ qui composent cette famille, et donc que les équations entre termes ont toujours la forme
	$$
	x_i^{A_i} = x_j^{A_j}
	$$
	pour des indices $i,j$ tels que $A_i = A_j$.
	
	Toute formule de Horn en $\vec x = (x_1^{A_1} , \cdots , x_n^{A_n})$ dont l'une des composantes atomiques a la forme
	$$
	x_i^{A_i} = x_j^{A_j} \quad \mbox{avec} \quad i < j \quad \mbox{et} \quad A_i = A_j 
	$$
	c'est-à-dire s'écrit comme une conjonction
	$$
	\varphi (\vec x) = \varphi_1 (\vec x) \wedge (x_i^{A_i} = x_j^{A_j})
	$$
	est canoniquement isomorphe dans ${\mathcal C}_{\Sigma}^{\rm car}$ a la formule en les $n-1$ variables  $(x_1^{A_1} , \cdots , x_{j-1}^{A_{j-1}} ,  x_{j+1}^{A_{j+1}} , \cdots , x_n^{A_n})$
	$$
	\varphi_1 (x_1^{A_1} , \cdots , x_i^{A_i} , \cdots , x_{j-1}^{A_{j-1}} , x_i^{A_i} ,  x_{j+1}^{A_{j+1}} , \cdots , x_n^{A_n}) \, .
	$$
	
	Après un nombre fini de ces réductions, tout tel objet $\varphi (\vec x)$ de ${\mathcal C}_{\Sigma}^{\rm car}$ se réécrit sous la forme de (i).
	
	La partie (ii) du corollaire se déduit alors de la partie (ii) de la proposition puisque
	
	\medskip
	
	$\left\{\begin{matrix}
		\bullet &\mbox{les termes définis sur un objet de ${\mathcal C}_{\Sigma}^{\rm car}$ de la forme} \hfill \\
		&\varphi_{\mathcal R} (\vec x_{{\mathcal A},m_{\bullet}}) = \underset{(R,\alpha) \in {\mathcal R}}{\bigwedge} R (x_{\alpha(1)}^{A_1} , \cdots , x_{\alpha(n_R)}^{A_{n_R}}) \\
		&\mbox{sont les variables} \hfill \\
		&x_i^A \, , \quad A \in {\mathcal A} \, , \quad 1 \leq i \leq m_A \, , \\
		&\mbox{qui composent la famille $\vec x_{{\mathcal A} , m_{\bullet}}$,} \hfill \\
		\bullet &\mbox{deux tels termes sont équivalents si et seulement si ils sont égaux,} \hfill \\
		\bullet &\mbox{deux formules atomiques en $\vec x_{{\mathcal A} , m_{\bullet}}$,} \hfill \\
		&R (x_{\alpha(1)}^{A_1} , \cdots , x_{\alpha(n_R)}^{A_{n_R}}) \quad \mbox{et} \quad R' (x_{\beta(1)}^{B_1} , \cdots , x_{\beta(n_{R'})}^{B_{n_{R'}}}) \\
		&\mbox{associées à deux symboles de relation} \hfill \\
		&R \rightarrowtail A_1 \cdots A_{n_R} \quad \mbox{et} \quad R' \rightarrowtail B_1 \cdots B_{n_{R'}} \\
		&\mbox{complétés par deux applications} \hfill \\
		&\alpha : \{ 1, \cdots , n_R\} \longrightarrow {\mathbb N} \quad \mbox{et} \quad \beta : \{1, \cdots , n_{R'}\} \longrightarrow {\mathbb N} \\
		&\mbox{sont  équivalentes si et seulement si $R=R'$ dans $\Sigma$ et $\alpha = \beta$.} \hfill
	\end{matrix}\right.$
	
	\medskip
	
	Enfin, (iii) est immédiat.
\end{demo}

%SubsecIII.3.g
\subsection{Réduction aux théories sans symboles de fonction}\label{ssec3.3.7}

Pour une théorie géométrique du premier ordre ${\mathbb T}$ de signature $\Sigma$, l'interprétation dans un modèle $M$ de ${\mathbb T}$ à valeurs dans un topos ${\mathcal E}$ de n'importe quel symbole de fonction de $\Sigma$
$$
f : A_1 \cdots A_n \longrightarrow B
$$
est un morphisme
$$
Mf : MA_1 \times \cdots \times MA_n \longrightarrow MB \, .
$$

Or, se donner un tel morphisme équivaut à se donner son graphe c'est-à-dire un sous-objet
$$
R_f M \ {\lhook\joinrel\relbar\joinrel\rightarrow} \ MA_1 \times \cdots \times MA_n \times MB
$$
dont la projection
$$
R_f M \longrightarrow MA_1 \times \cdots \times MA_n
$$
est un isomorphisme ou, ce qui revient au même, est à la fois un monomorphisme et un épimorphisme.

Cette remarque conduit à introduire le procédé général suivant de remplacement des symboles de fonctions par des symboles de relations:

\begin{defn}\label{defn3.3.17}
	Soit $\Sigma$ une signature qui comprend des symboles de fonctions.
	
	Alors:
	
	\begin{listeimarge}
	\item On appelle signature relationnelle associée à $\Sigma$ et on note
	$$
	\Sigma_r
	$$
	la signature déduite de $\Sigma$ en remplaçant chacun de ses symboles de fonctions
	$$
	f : A_1 \cdots A_n \longrightarrow B
	$$
	par un symbole de relation
	$$
	R_f \rightarrowtail  A_1 \cdots A_n \, B \, .
	$$
	\item On appelle théorie relationnelle associée à $\Sigma$ et on note
	$$
	{\mathbb T}_r^{\Sigma}
	$$
	la théorie cartésienne de signature $\Sigma_r$ définie par les axiomes
	$$
	\left\{ \begin{matrix}
		R_f (x_1^{A_1} , \cdots , x_n^{A_n} , x_B) \wedge R_f (x_1^{A_1} , \cdots , x_n^{A_n} , x'_B) \vdash x_B = x'_B \, , \hfill \\
		\top (x_1^{A_1} , \cdots , x_n^{A_n}) \vdash (\exists \, x_B) \, R_f (x_1^{A_1} , \cdots , x_n^{A_n} , x_B) \hfill
	\end{matrix} \right.
	$$
	associés aux symboles de fonction de $\Sigma$
	$$
	f : A_1 \cdots A_n \longrightarrow B \, .
	$$
	\item On appelle traduction de $\Sigma$ à ${\mathbb T}_r^{\Sigma}$ le processus suivant de remplacement des formules géométriques de $\Sigma$ par des formules géométriques de ${\mathbb T}_r^{\Sigma}$:
	\begin{enumerate}
	\item[(1)] Pour tout terme de $\Sigma$ en une suite de variables $\vec x = (x_1^{A_1} , \cdots , x_n^{A_n})$ qui est écrit sous la forme
	$$
	g(\vec x) = f(g_1 (\vec x) , \cdots , g_m (\vec x))
	$$
	pour un symbole de fonction de $\Sigma$
	$$
	f : B_1\cdots B_m \longrightarrow B
	$$
	et des termes plus simples $g_1 (\vec x) , \cdots , g_m (\vec x)$ à valeurs dans les sortes $B_1 , \cdots , B_m$, on introduit des variables
	$$
	z_g^B \, , \ z_{g_1}^{B_1} , \cdots , z_{g_m}^B
	$$
	et la formue atomique
	$$
	R_f (z_{g_1}^{B_1} , \cdots , z_{g_m}^{B_m} , z_g^B)
	$$
	éventuellement complétée par des formules
	$$
	z_{g_j}^{B_j} = x_i^{A_i}
	$$
	pour chaque terme $g_j (\vec x)$ qui se confond avec l'une des variables $x_i^{A_i}$ composant la suite $\vec x = (x_1^{A_1} , \cdots , x_n^{A_n})$.
	
	\item[(2)] Pour toute formule atomique de $\Sigma$ nécessairement de la forme
	$$
	g_1^{B_1} (\vec x) = g_2^{B_1} (\vec x) \quad \mbox{ou} \quad R(g_1^{B_1} (\vec x) , \cdots , g_n^{B_n} (\vec x))
	$$
	pour des termes $g_1^{B_1} , g_2^{B_1}$ ou $g_1^{B_1} , \cdots , g_n^{B_n}$ en une suite de variables $\vec x = (x_1^{A_1} , \cdots , x_n^{A_n})$ et éventuellement un symbole de relation $R \rightarrowtail B_1 \cdots B_m$ de $\Sigma$, on introduit une variable
	$$
	z_g^B
	$$
	pour chaque terme en $\vec x$ à valeur dans une sorte $B$ qui apparaît dans le processus fini de construction des termes $g_1^{B_1} , g_2^{B_2}$ ou $g_1^{B_1} , \cdots , g_n^{B_n}$ à partir de symboles de fonction de $\Sigma$, et on lit ces variables par les relations de (1)
	$$
	R_f (z_{g_1} , \cdots , z_{g_m} , z_g)
	$$
	associées aux formules de construction de ces termes
	$$
	g(\vec x) = f (g_1(\vec x) , \cdots , g_m (\vec x)) \, .
	$$
	
	On forme alors la conjonction de ces formules en nombre fini et de la formule
	$$
	z_{g_1}^{B_1} = z_{g_2}^{B_1} \quad \mbox{ou} \quad R(z_{g_1}^{B_1} , \cdots ,  z_{g_n}^{B_n}) \, ,
	$$
	et on applique à la formule atomique obtenue en les variables
	$$
	x_1^{A_1}  \cdots  x_n^{A_n} \quad \mbox{et} \quad z_g^B
	$$
	le symbole 
	$$
	\exists
	$$
	appliqué à toutes les variables $z_g^B$, ce qui définit une formule ${\mathbb T}_r^{\Sigma}$-cartésienne en la suite de variables $\vec x = (x_1^{A_1} , \cdots , x_n^{A_n})$.
	
	\item[(3)] Pour toute formule géométrique $\varphi (\vec x)$ de $\Sigma$ construite à partir de formules atomiques en appliquant des symboles $\exists$, $\vee$ ou $\bigvee$, on remplace chaque formule atomique qui entre dans la construction de $\varphi$ par la formule ${\mathbb T}_r^{\Sigma}$-cartésienne qui lui est associée à la manière de (2), de façon à obtenir une formule géométrique de ${\mathbb T}_r^{\Sigma}$ en la même suite de variables $\vec x$.
\end{enumerate}
\end{listeimarge}
\end{defn}

\begin{remark}
	Le procédé de (3) remplace toute formule géométrique de $\Sigma$ qui est cohérente [resp. régulière] par une formule géométrique de ${\mathbb T}_r^{\Sigma}$ qui est cohérente [resp. régulière].
\end{remark}

Le procédé de traduction de $\Sigma$ à ${\mathbb T}_r^{\Sigma}$ de cette définition permet d'obtenir:

\begin{prop}\label{prop3.3.18}
	Soit ${\mathbb T}$ une théorie géométrique du premier ordre de signature $\Sigma$.
	
	Soient $\Sigma_r$ la signature relationnelle associée à $\Sigma$ et ${\mathbb T}_r^{\Sigma}$ la théorie cartésienne de signature $\Sigma_r$ obtenue en remplaçant chaque symbole de fonction de $\Sigma$
	$$
	f : A_1 \cdots A_n \longrightarrow B
	$$
	par un symbole de relation de $\Sigma_r$
	$$
	R_f \rightarrowtail A_1 \cdots A_n \, B
	$$
	complété par les deux axiomes cartésiens de ${\mathbb T}_r^{\Sigma}$
	$$
	\left\{ \begin{matrix}
		R_f (x_1^{A_1} , \cdots , x_n^{A_n} , x_B) \wedge R_f (x_1^{A_1} , \cdots , x_n^{A_n} , x'_B) \vdash x_B = x'_B \, , \hfill \\
		\top (x_1^{A_1} , \cdots , x_n^{A_n}) \vdash (\exists \, x_B) \, R_f (x_1^{A_1} , \cdots , x_n^{A_n} , x_B) \, . \hfill
	\end{matrix} \right.
	$$
	
	Soit ${\mathbb T}_r$ la théorie géométrique de signature $\Sigma_r$ qui est déduite de ${\mathbb T}_r^{\Sigma}$ en lui ajoutant les axiomes
	$$
	\varphi_i^r (\vec x_i) \vdash \psi_i^r (\vec x_i)
	$$
	obtenus à partir des axiomes de ${\mathbb T}$
	$$
	\varphi_i (\vec x_i) \vdash \psi_i (\vec x_i) \, , \quad i \in I \, ,
	$$
	en remplaçant leurs formules géométriques constitutives
	$$
	\varphi_i (\vec x_i) \, , \ \psi_i (\vec x_i)
	$$
	par les formules géométriques de $\Sigma_r$ qui leur sont associées par le procédé de traduction de $\Sigma$ à ${\mathbb T}_r^{\Sigma}$.
	
	Alors:
	
	\begin{listeimarge}
	\item Le procédé de traduction de $\Sigma$ à ${\mathbb T}_r^{\Sigma}$ définit un foncteur entre les catégories syntactiques
	$$
	{\mathcal C}_{\mathbb T} \longrightarrow {\mathcal C}_{{\mathbb T}_r}
	$$
	qui est une équivalence de catégories.
	
	\item Si ${\mathbb T}$ est une théorie cohérente [resp. régulière, resp. cartésienne] sa traduction relationnelle ${\mathbb T}_r$ est également cohérente [resp. régulière, resp. cartésienne] et l'équivalence de (i) se restreint en une équivalence
	$$
	{\mathcal C}_{\mathbb T}^{\rm coh} \longrightarrow {\mathcal C}_{{\mathbb T_r}}^{\rm coh}
	$$
	[resp.
	$$
	{\mathcal C}_{\mathbb T}^{\rm reg} \longrightarrow {\mathcal C}_{{\mathbb T_r}}^{\rm reg} \, ,
	$$
	resp.
	$$
	{\mathcal C}_{\mathbb T}^{\rm car} \longrightarrow {\mathcal C}_{{\mathbb T_r}}^{\rm car} \, \mbox{].}
	$$
\end{listeimarge}
\end{prop}

\begin{demo}
	\begin{listeisansmarge}
	\item Le procédé de traduction des formules géométriques de $\Sigma$ en des formules géométriques de ${\mathbb T}_r^{\Sigma}$ transforme tout séquent
	$$
	\varphi (\vec x) \vdash \psi (\vec x)
	$$
	démontrable à partir d'axiomes
	$$
	\varphi_i (\vec x_i) \vdash \psi_i (\vec x_i) \, , \quad i \in I \, ,
	$$
	en un séquent
	$$
	\varphi^r (\vec x) \vdash \psi^r (\vec x)
	$$
	démontrable à partir des axiomes
	$$
	\varphi_i^r (\vec x_i) \vdash \psi_i^r (\vec x_i) 
	$$
	complétés par les axiomes de ${\mathbb T}_r^{\Sigma}$.
	
	En particulier, il transforme les formules démontrablement fonctionnelles
	$$
	\varphi (\vec x) \xrightarrow{ \ \theta (\vec x , \vec y) \ } \psi (\vec y)
	$$
	en des formules démontrablement fonctionnelles
	$$
	\varphi^r (\vec x) \xrightarrow{ \ \theta^r (\vec x , \vec y) \ } \psi^r (\vec y)
	$$
	et respecte les équivalences démontrables.
	
	De plus, il transforme toute formule de la forme
	$$
	(\exists \, \vec y) \, (\theta_1 (\vec x , \vec y) \wedge \theta_2 (\vec y , \vec z))
	$$
	en la formule
	$$
	(\exists \, \vec y) \, (\theta_1^r (\vec x , \vec y) \wedge \theta_2^r (\vec x , \vec y)) \, .
	$$
	Ainsi, il définit un foncteur
	$$
	{\mathcal C}_{\mathbb T} \longrightarrow {\mathcal C}_{{\mathbb T}_r} \, .
	$$
	
	Il résulte de l'interprétation sémantique des morphismes entre des objets de topos
	$$
	Mf : MA_1 \times \cdots \times MA_n \longrightarrow MB
	$$
	comme des graphes c'est-à-dire des sous-objets
	$$
	R_f \, M \ {\lhook\joinrel\relbar\joinrel\rightarrow} \ MA_1 \times \cdots \times MA_n \times MB
	$$
	dont la projection
	$$
	R_f \, M \longrightarrow MA_1 \times \cdots \times MA_n
	$$
	est un isomorphisme, que le foncteur ${\mathcal C}_{\mathbb T} \to {\mathcal C}_{{\mathbb T}_r}$ induit une équivalence des topos classifiants
	$$
	{\mathcal E}_{\mathbb T} \xrightarrow{ \ \sim \ } {\mathcal E}_{{\mathbb T}_r} \, .
	$$
	
	Comme ${\mathcal C}_{\mathbb T}$ et ${\mathcal C}_{{\mathbb T}_r}$ sont des sous-catégories pleines de ${\mathcal E}_{\mathbb T}$ et ${\mathcal E}_{{\mathbb T}_r}$, le foncteur
	$$
	{\mathcal C}_{\mathbb T} \longrightarrow {\mathcal C}_{{\mathbb T}_r}
	$$
	est pleinement fidèle.
	
	Il est essentiellement surjectif, et donc une équivalence, car toute formule gémétrique de ${\mathbb T}_r$ provient à équivalence démontrable près de la formule géométrique de ${\mathbb T}$ qui s'en déduit en remplaçant chaque formule de la forme
	$$
	R_f (x_1^{A_1} , \cdots , x_n^{A_n} , x^B)
	$$
	indexée par un symbole de fonction $f : A_1 \cdots A_n \to B$ de $\Sigma$, par la formule
	$$
	f (x_1^{A_1} , \cdots , x_n^{A_n}) = x^B \, .
	$$
	\item résulte de (i) car le procédé de traduction de $\Sigma$ à ${\mathbb T}_r^{\Sigma}$ transforme les formules atomiques [resp. régulières, resp. cohérentes, resp. ${\mathbb T}$-cartésiennes] en formules ${\mathbb T}_r^{\Sigma}$-cartésiennes [resp. régulières, resp. cohérentes, resp. ${\mathbb T}_r$-cartésiennes], et que réciproquement toute formule de ${\mathbb T}_r^{\Sigma}$ qui est atomique [resp. régulière, resp. cohérente, resp. ${\mathbb T}_r$-cartésienne] provient à équivalence démontrable près d'une formule de $\Sigma$ qui est atomique [resp. régulière, resp. cohérente, resp. ${\mathbb T}$-cartésienne].
	\end{listeisansmarge}
\end{demo}

	%Subsec IV.1.a et b
\chapter[Morphismes localement connexes et images extraordinaires]{Morphismes localement connexes et images extraordinaires des sous-topos}\label{chap4}

\section{Fibrations, topologies de Giraud et morphismes localement connexes}\label{sec4.1}

\subsection{La notion de morphisme horizontal}\label{ssec4.1.a}

On rappelle que tout foncteur $F : {\mathcal C} \to {\mathcal B}$ d'une catégorie ${\mathcal C}$ dans une catégorie ${\mathcal B}$ définit dans ${\mathcal C}$ une notion de morphisme  ``horizontal'' ou ``cartésien'' relativement au foncteur $F$:

\begin{defn}\label{defn4.1.1}
Soit un foncteur 
$$
F : {\mathcal C} \longrightarrow {\mathcal B}
$$
d'une catégorie ${\mathcal C}$ sur une catégorie de base ${\mathcal B}$.

On dit qu'un morphisme de ${\mathcal C}$
$$
X' \xrightarrow{ \ x \ } X
$$
est ``$F$-cartésien'' ou ``horizontal'' si, pour tout morphisme de ${\mathcal C}$ de même but $X$
$$
X'' \xrightarrow{ \ x_1 \ } X \, ,
$$
toute factorisation dans ${\mathcal B}$
$$
\xymatrix{
F(X'') \ar[rr]^{b'} \ar[rd]_{F(x_1)} &&F(X') \ar[ld]^{F(x)} \\
&F(X)
}
$$
se relève dans ${\mathcal C}$ en une unique factorisation
$$
\xymatrix{
X'' \ar[rr]^{x'} \ar[rd]_{x_1} &&X' \ar[ld]^{x} \\
&X
}
$$
avec donc $F(x') = b'$.
\end{defn}

\newpage

\begin{exqed}
Soit ${\mathcal C} = {\mathcal B} / P$ la ``catégorie des éléments'' d'un préfaisceau
$$
P : {\mathcal B}^{\rm op} \longrightarrow {\rm Ens}
$$
sur une catégorie ${\mathcal B}$, dont les objets sont les paires $(Y,p)$ constituées d'un objet $Y$ de ${\mathcal B}$ et d'un élément $p \in P(Y)$.

Alors les morphismes de ${\mathcal C} = {\mathcal B}/P$
$$
(Y',p') \longrightarrow (Y,p) \, ,
$$
qui sont constitués des morphismes $Y' \xrightarrow{ \, y \, } Y$ tels que
$$
P(y)(p) = p' \, ,
$$
sont tous horizontaux relativement au foncteur d'oubli ${\mathcal C} = {\mathcal B}/P \to {\mathcal B}$.

En effet, pour tout triangle commutatif de ${\mathcal B}$
$$
\xymatrix{
Y'' \ar[rr]^{y'} \ar[rd]_{y_1} &&Y' \ar[ld]^{y} \\
&Y
}
$$
et tous éléments $p \in P(Y)$, $p' \in P(Y')$, $p'' \in P(Y'')$, les conditions
$$
p' = P(y)(p) \quad \mbox{et} \quad p'' = P(y_1)(p)
$$
impliquent
$$
P(y')(p') = p'' \, .
$$
\end{exqed}

On remarque:

\begin{lem}\label{lem4.1.2}
Considérons deux catégories localement petites ${\mathcal C}$ et ${\mathcal B}$ reliées par une paire de foncteurs adjoints
$$
\left( {\mathcal C} \xrightarrow{ \ F \ } {\mathcal B} , {\mathcal B} \xrightarrow{ \ G \ } {\mathcal C} \right).
$$

Alors un morphisme de ${\mathcal C}$
$$
X' \xrightarrow{ \ x \ } X
$$
est $F$-cartésien si et seulement si le carré commutatif de ${\mathcal C}$
$$
\xymatrix{
X' \ar[d] \ar[r]^x &X \ar[d] \\
G \circ F (X') \ar[r] &G \circ F (X)
}
$$
est un carré cartésien.
\end{lem}

\newpage

\begin{ex}
Si ${\mathcal C}$ est une catégorie localement petite, considérons la catégorie
$$
{\mathcal C} / {\mathcal C}
$$
dont les objets sont les morphismes de ${\mathcal C}$
$$
\xymatrix{
X \ar[d] \\
Y
}
$$
et dont les morphismes
$$
\left(
     \raisebox{0.5\depth}{
     \xymatrix{X' \ar[d] \\ Y'}} \right) \longrightarrow \left(
     \raisebox{0.5\depth}{
     \xymatrix{X \ar[d] \\ Y}} \right) 
$$
sont les carrés commutatifs de ${\mathcal C}$:
$$
\xymatrix{
X' \ar[d] \ar[r] &X \ar[d] \\
Y' \ar[r] &Y
}
$$
On observe que le foncteur d'oubli
$$
\begin{matrix}
\hfill {\mathcal C} / {\mathcal C} &\longrightarrow &{\mathcal C} \, , \hfill \\
\hfill \left(
     \raisebox{0.5\depth}{
     \xymatrix{X \ar[d] \\ Y}} \right) &\longmapsto &Y \, , \hfill \\
\left(
     \raisebox{0.5\depth}{ 
     \xymatrix{
X' \ar[d] \ar[r]^x &X \ar[d] \\
Y' \ar[r]^y &Y
}} \right) &\longmapsto &\left(Y' \xrightarrow{ \, y' \, }Y \right)
\end{matrix}
$$
admet pour adjoint à droite le foncteur
$$
\begin{matrix}
{\mathcal C} &\longrightarrow &{\mathcal C} / {\mathcal C} \, , \hfill \\
X &\longmapsto &\left(
     \raisebox{0.5\depth}{
     \xymatrix{
X \ar[d]^{{\rm id}_X} \\ X}}\right). \\
\end{matrix}
$$
Donc un morphisme de ${\mathcal C} / {\mathcal C}$
$$
\xymatrix{
X' \ar[d] \ar[r] &X \ar[d] \\
Y' \ar[r] &Y
}
$$
est ``horizontal'' relativement au foncteur d'oubli ${\mathcal C} / {\mathcal C} \to {\mathcal C}$ si et seulement si le carré commutatif de ${\mathcal C}$ qu'il forme est cartésien.
\end{ex}

\newpage

\begin{demolem}
Considérons deux morphismes de ${\mathcal C}$
$$
\xymatrix{
X'' \ar[rd]_{x_1} &&X' \ar[ld]^x \\
&X
}
$$
et une factorisation dans ${\mathcal B}$
$$
\xymatrix{
F(X'') \ar[rr]^{b'} \ar[rd]_{F(x_1)} &&F(X') \ar[ld]^{F(x)} \\
&F(X)
}
$$
laquelle correspond par adjonction à un triangle commutatif de ${\mathcal C}$
$$
\xymatrix{
X'' \ar[r] \ar[rd] &G \circ F (X') \ar[d]^{G \circ F(x)} \\
&G \circ F(X)
}
$$
dont la composante $X'' \to G \circ F(X)$ est le composé
$$
X'' \xrightarrow{ \ x_1 \ } X \longrightarrow G \circ F(X) \, .
$$

Ainsi, la factorisation en $F(X'') \xrightarrow{ \, b' \, } F(X')$ correspond à un carré commutatif de ${\mathcal C}$
$$
\xymatrix{
X'' \ar[d] \ar[r] &G \circ F(X') \ar[d] \\
X \ar[r] &G \circ F(X)
}
$$
et il y a existence et unicité d'une factorisation en un $X'' \xrightarrow{ \, x' \, } X'$ relevant $F(X'') \xrightarrow{ \, b' \, } F(X')$ si et seulement si le carré commutatif de ${\mathcal C}$
$$
\xymatrix{
X' \ar[d] \ar[r] &G \circ F(X') \ar[d] \\
X \ar[r] &G \circ F(X)
}
$$
est cartésien.
\end{demolem}

\subsection{La notion de fibration}\label{4.1.2}

On rappelle encore qu'un foncteur $F : {\mathcal C} \to {\mathcal B}$ est dit une ``fibration'' quand sa catégorie de définition ${\mathcal C}$ a suffisamment de morphismes horizontaux relativement à $F$, au sens suivant:

\newpage

\begin{defn}\label{def4.1.3}
Un foncteur entre deux catégories
$$
F : {\mathcal C} \longrightarrow {\mathcal B}
$$ est appelé une ``fibration'' si, pour tout objet $X$ de ${\mathcal C}$ et tout morphisme de ${\mathcal B}$ vers $F(X)$
$$
Y' \xrightarrow{ \ b \ } F(X) \, ,
$$
il existe dans ${\mathcal C}$ un morphisme horizontal ou $F$-cartésien
$$
X' \xrightarrow{ \ x \ } X
$$
et dans ${\mathcal B}$ un isomorphisme
$$
F(X') \xrightarrow{ \ \sim \ } Y'
$$
qui rende commutatif le triangle:
$$
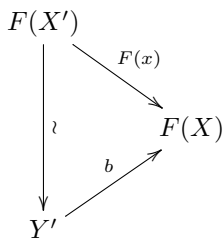
\xymatrix{
F(X') \ar[dd]^{\wr} \ar[rd]^{F(x)} \\
&F(X) \\
Y' \ar[ru]^b
}
$$
\end{defn}

\begin{remark}
Si $X'_1 \xrightarrow{ \, x_1 \, } X$ et $X'_2 \xrightarrow{ \, x_2 \, } X$ sont deux morphismes horizontaux de ${\mathcal C}$ complétés par deux isomorphismes qui rendent commutatif le diagramme
$$
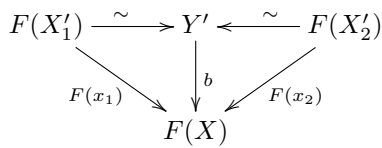
\xymatrix{
F(X'_1) \ar[r]^{\sim} \ar[rd]_{F(x_1)} &Y' \ar[d]^b &F(X'_2) \ar[l]_{\sim} \ar[ld]^{F(x_2)} \\
&F(X)
}
$$
alors il résulte de la définition des morphismes horizontaux que l'isomorphisme composé
$$
F(X'_1) \xrightarrow{ \ \sim \ } F(X'_2)
$$
se relève en un unique isomorphisme de ${\mathcal C}$
$$
X'_1 \xrightarrow{ \ \sim \ } X'_2
$$
qui rende commutatif le triangle:
$$
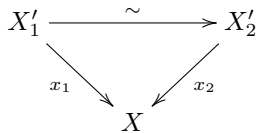
\xymatrix{
X'_1 \ar[rr]^{\sim} \ar[rd]_{x_1} &&X'_2 \ar[ld]^{x_2} \\
&X
}
$$
Autrement dit, la paire $\left( X' \xrightarrow{ \, x \, } X , F(X') \xrightarrow{ \, \sim \, } Y \right)$ qui relève $\left( Y \xrightarrow{ \, b \, } F(X) \right)$ est unique à unique isomorphisme près.
\end{remark}

\newpage

\begin{exsqed}
\noindent (i) Pour tout préfaisceau $P$ sur une catégorie ${\mathcal B}$
$$
P : {\mathcal B}^{\rm op} \longrightarrow {\rm Ens} \, ,
$$
le foncteur d'oubli de la catégorie de ses éléments
$$
{\mathcal C} = {\mathcal B}/P \longrightarrow {\mathcal B}
$$
est une fibration.

En effet, pour tout objet $(Y,p)$ de ${\mathcal C}$ composé d'un objet $Y$ de ${\mathcal B}$ et d'un élément $p \in P(Y)$, tout morphisme de ${\mathcal B}$ vers $Y$
$$
Y' \xrightarrow{ \ y \ } Y
$$
se relève en un morphisme de ${\mathcal C}$ nécessairement horizontal
$$
(Y',p') \longrightarrow (Y,p) \quad \mbox{en posant} \quad p' = P(y)(p) \, .
$$ 
(ii) Pour toute catégorie ${\mathcal C}$ qui a des produits fibrés, le foncteur
$$
\begin{matrix}
{\mathcal C} / {\mathcal C} &\longrightarrow &{\mathcal C} \, , \\
\left(
     \raisebox{0.5\depth}{
      \xymatrix{X \ar[d] \\Y}} \right) &\longmapsto &Y
\end{matrix}
$$
est une fibration.

En effet, tout morphisme de ${\mathcal C}$ vers l'image $Y$ d'un objet {\small $\left(
     \raisebox{0.5\depth}{
      \xymatrix{X \ar[d] \\Y}} \right)$} de ${\mathcal C} / {\mathcal C}$
$$
Y' \longrightarrow Y
$$
se relève en le morphisme défini par le carré cartésien de ${\mathcal C}$
$$
\xymatrix{
Y' \times_Y X \ar[d] \ar[r] &X \ar[d] \\
Y' \ar[r] &Y
}
$$
qui est un morphisme horizontal. 
\end{exsqed}

Dans la situation du lemme \ref{lem4.1.2} d'un foncteur $F : {\mathcal C} \to {\mathcal B}$ qui admet un adjoint à droite $G : {\mathcal B} \to {\mathcal C}$, on a:

\begin{prop}\label{prop4.1.4}
Considérons deux catégories localement petites ${\mathcal C}$ et ${\mathcal B}$ reliées par une paire de foncteurs adjoints
$$
\left( {\mathcal C} \xrightarrow{ \ F \ } {\mathcal B} , {\mathcal B} \xrightarrow{ \ G \ } {\mathcal C} \right).
$$

Supposons que les catégories ${\mathcal C}$ et ${\mathcal B}$ ont des produits fibrés.

Alors le foncteur
$$
F : {\mathcal C} \longrightarrow {\mathcal B}
$$
est une fibration si et seulement si pour tout objet $X$ de ${\mathcal C}$ et toute paire de morphisme de ${\mathcal B}$
$$
\xymatrix{
&F(X) \ar[d] \\
Y' \ar[r] &Y
}
$$
le morphisme naturel
$$
F(X \times_{G(Y)} G(Y')) \longrightarrow F(X) \times_Y Y'
$$
induit par les morphismes
$$
F(X \times_{G(Y)} G(Y')) \longrightarrow F(X) \qquad \mbox{(image par $F$ de $X \times_{G(Y)} G(Y')\longrightarrow X$)}
$$
et
$$
F(X \times_{G(Y)} G(Y')) \longrightarrow Y' \qquad \mbox{(adjoint de $X \times_{G(Y)} G(Y')\longrightarrow G(Y')$)}
$$
au-dessus du morphisme
$$
F(X \times_{G(Y)} G(Y')) \longrightarrow Y \qquad \mbox{(adjoint de $X \times_{G(Y)} G(Y')\longrightarrow G(Y)$),}
$$
est un isomorphisme
$$
F(X \times_{G(Y)} G(Y')) \xrightarrow{ \ \sim \ } F(X) \times_Y Y' \, .
$$
\end{prop}

\begin{demo}
Supposons d'abord que $F$ soit une fibration et considérons deux morphismes de ${\mathcal B}$ de la forme
$$
\xymatrix{ &F(X) \ar[d] \\
Y' \ar[r] &Y
}
$$
pour un objet $X$ de ${\mathcal C}$.

D'après la définition des fibrations, il existe un morphisme horizontal de ${\mathcal C}$
$$
X' \xrightarrow{ \ x \ } X
$$
et un isomorphisme
$$
F(X') \xrightarrow{ \ \sim \ } F(X) \times_Y Y'
$$
dont le composé avec la projection $F(X) \times_Y Y' \to F(X)$ soit le morphisme image
$$
F(x) : F(X') \longrightarrow F(X) \, .
$$

D'après le lemme \ref{lem4.1.2}, le morphisme horizontal $x$ définit un carré cartésien de ${\mathcal C}$:
$$
\xymatrix{
X' \ar[d] \ar[r]^x &X \ar[d] \\
G \circ F(X') \ar[r] &G \circ F(X)
}
$$
\newpage
\noindent Or, le carré cartésien de ${\mathcal B}$
$$
\xymatrix{
F(X') \ar[d] \ar[r] &F(X) \ar[d] \\
Y' \ar[r] &Y
}
$$
induit un carré cartésien de ${\mathcal C}$
$$
\xymatrix{
G \circ F(X') \ar[d] \ar[r] &G \circ F(X) \ar[d] \\
G(Y') \ar[r] &G (Y)
}
$$
puisque le foncteur $G$, qui est un adjoint à droite, respecte les produits fibrés.

On conclut comme voulu que le carré commutatif de ${\mathcal C}$
$$
\xymatrix{
X' \ar[d] \ar[r] &X \ar[d] \\
G (Y') \ar[r] &G (Y)
}
$$
est cartésien.

Réciproquement, supposons que la condition de la proposition soit vérifiée, et considérons un morphisme de ${\mathcal B}$ de la forme
$$
Y' \xrightarrow{ \ b \ } F(X)
$$
pour un objet $X$ de ${\mathcal C}$.

La paire  de morphismes de ${\mathcal C}$
$$
\xymatrix{
&X \ar[d] \\
G(Y') \ar[r]^-{G(b)} &G \circ F(X)
}
$$
se complète en un carré cartésien
$$
\xymatrix{
X' \ar[d] \ar[r]^x&X \ar[d] \\
G(Y') \ar[r]^-{G(b)} &G \circ F(X)
}
$$
et, d'après l'hypothèse, le carré commutatif correspondant de ${\mathcal B}$
$$
\xymatrix{
F(X') \ar[d] \ar[r]^-{F(x)} &F(X) \ar[d]^{\rm id} \\
Y' \ar[r]^{b} &F(X)
}
$$
est également cartésien, ce qui signifie que le morphisme 
$$
F(X') \longrightarrow Y'
$$
est un isomorphisme de ${\mathcal B}$.

On en déduit que le carré commutatif de ${\mathcal C}$
$$
\xymatrix{
X' \ar[d] \ar[r]^x &X \ar[d] \\
G \circ F(X') \ar[r] &G \circ F(X)
}
$$
est cartésien, ce qui signifie d'après le lemme \ref{lem4.1.2} que le morphisme
$$
X' \xrightarrow{ \ x \ } X
$$
est horizontal.

Il relève le morphisme composé
$$
F(X') \xrightarrow{ \ \sim \ } Y' \xrightarrow{ \ b \ } F(X) \, .
$$

Cela montre comme voulu que le foncteur
$$
F : {\mathcal C} \longrightarrow {\mathcal B}
$$
est une fibration.
\end{demo}

%Subsec IV.1.c et d
\subsection{Morphismes localement connexes et fibrations}\label{4.1.3}

Pour les morphismes de topos, la propriété de connexité locale se caractérise en termes de la notion de fibration:

\begin{prop}\label{prop4.1.5}
Considérons un morphisme de topos
$$
f : {\mathcal E}' \longrightarrow {\mathcal E}
$$
qui est ``essentiel'' c'est-à-dire dont la composante d'image réciproque
$$
f^* : {\mathcal E} \longrightarrow {\mathcal E}'
$$
a non seulement un adjoint à droite $f_* : {\mathcal E}' \to {\mathcal E}$ mais aussi un adjoint à gauche
$$
f_! : {\mathcal E}' \longrightarrow {\mathcal E} \, .
$$

Alors le morphisme $f$ est ``localement connexe'' si et seulement si ce foncteur
$$
f_! : {\mathcal E}' \longrightarrow {\mathcal E}
$$
est une fibration.
\end{prop}

\begin{remark}
C'est la proposition 5.2.8 de l'article \cite{Localfibrations} de L. Bertoli et O.C.
\end{remark}

\begin{demo}
Le foncteur $f_! : {\mathcal E}' \to {\mathcal E}$ admet pour adjoint à droite le foncteur
$$
f^* : {\mathcal E} \longrightarrow {\mathcal E}' \, .
$$

Les catégories ${\mathcal E}$ et ${\mathcal E}'$, qui sont des topos, sont localement  petites et admettent tous les produits fibrés.

D'après la proposition \ref{prop4.1.4}, le foncteur $f_!$ est une fibration si et seulement si, pour tout morphisme de ${\mathcal E}$
$$
E_2 \xrightarrow{ \ e \ } E_1 \, ,
$$
tout objet $E'_1$ de ${\mathcal E}'$ et tout morphisme de ${\mathcal E}$
$$
f_! \, E'_1 \longrightarrow E_1
$$
correspondant à un morphisme de ${\mathcal E}'$
$$
E'_1 \longrightarrow f^* \, E_1 \, ,
$$
le morphisme naturel
$$
f_! \, (E' \times_{f^* E_1} f^* \, E_2) \longrightarrow f_! \, E' \times_{E_1} E_2
$$
est un isomorphisme.

C'est la seconde des trois propriétés équivalentes qui définissent la notion de ``morphisme localement connexe'' de topos, telles que la définition \ref{defn1.3.18} les rappelle.
\end{demo}

On rappelle d'autre part le résultat suivant (qui est un cas particulier du Corollaire 4.58 de \cite{Denseness}):

\begin{prop}\label{prop4.1.6}
Considérons un foncteur entre deux catégories essentiellement petites
$$
p : {\mathcal C} \longrightarrow {\mathcal B}
$$
qui est une fibration.

Alors le morphisme de topos de préfaisceaux induit par ce foncteur
$$
(p^* , p_*) : \widehat{\mathcal C} \longrightarrow \widehat{\mathcal B}
$$
est localement connexe.
\end{prop}

\begin{remark}
On se souvient que, pour tout foncteur entre des catégories essentiellement petites
$$
p : {\mathcal C} \longrightarrow {\mathcal B} \, ,
$$
le foncteur de composition des préfaisceaux avec $p$
$$
p^* = (\bullet) \circ p : \widehat{\mathcal B} \longrightarrow \widehat{\mathcal C}
$$
possède à la fois un adjoint à droite $p_*$ et un adjoint à gauche $p_!$, donc définit un morphisme de topos qui est essentiel
$$
(p^* , p_*) : \widehat{\mathcal C} \longrightarrow \widehat{\mathcal B} \, .
$$

Les adjoints $p_!$ et $p_*$ associent à tout préfaisceau $P$ sur ${\mathcal C}$ les préfaisceaux sur ${\mathcal B}$ définis par les formules
$$
Y \longmapsto (p_! \, P)(Y) = \varprojlim_{(X,Y \to p(X)) \in Y \backslash {\mathcal C}} P(X) \, ,
$$
$$
Y \longmapsto (p_* \, P)(Y) = \varinjlim_{(X, p(X) \to Y) \in {\mathcal C}/Y} P(X) \, ,
$$
où $Y \backslash {\mathcal C}$ [resp. ${\mathcal C}/Y$] désigne la catégorie essentiellement petite dont

\medskip

$\left\{\begin{matrix}
\bullet &\mbox{les objets sont les paires constituées d'un objet $X$ de ${\mathcal C}$ et d'un morphisme de ${\mathcal B}$} \hfill \\
&Y \longrightarrow p(X) \quad \mbox{[resp.} \ p(X) \longrightarrow Y \, \mbox{]}, \\
\bullet &\mbox{les morphismes sont les morphismes de ${\mathcal C}$} \hfill \\
&X_2 \xrightarrow{ \ x \ } X_1 \\
&\mbox{qui rendent commutatif le triangle} \hfill \\
&\xymatrix{
&p(X_2) \ar[dd]^{p(x)} \\
Y \ar[ru] \ar[rd] \\
&p(X_1)
} \qquad \qquad \xymatrix{
p(X_2) \ar[dd]^{p(x)}_{\mbox{[resp.} \qquad \qquad} \ar[rd] \\
&Y \quad \mbox{].} \\
p(X_1) \ar[ru]
} 
\end{matrix}\right.$
\end{remark}

\begin{demo}
Nous devons montrer que si le foncteur
$$
p : {\mathcal C} \longrightarrow {\mathcal B}
$$
est une fibration, alors pour tout préfaisceau $P$ sur ${\mathcal C}$ et toute paire de morphismes de $\widehat{\mathcal B}$
$$
\xymatrix{
&p_! \, P \ar[d] \\
Q_2 \ar[r] &Q_1
}
$$
le morphisme canonique
$$
p_! (P \times_{p^* Q_1} p^* Q_2) \longrightarrow p_! P \times_{Q_1} Q_2
$$
est un isomorphisme, autrement dit que en tout objet $Y$ de ${\mathcal B}$ l'application
$$
p_! (P \times_{p^* Q_1} p^* Q_2) (Y) \longrightarrow p_! P(Y) \times_{Q_1 (Y)} Q_2 (Y)
$$
est une bijection.

Or, comme $p$ est une fibration, tout objet de $Y/{\mathcal C}$
$$
(X , Y \longrightarrow p(X))
$$
définit un morphisme horizontal de ${\mathcal C}$
$$
X' \xrightarrow{ \ x \ } X
$$
et un isomorphisme $Y \xrightarrow{ \, \sim \, } p(X')$ qui fait commuter le triangle:
$$
\xymatrix{
Y \ar[rr]^{\sim} \ar[rd] &&p(X') \ar[ld]^{p(x)} \\
&p(X)
}
$$
De plus, tout morphisme de $Y/{\mathcal C}$
$$
(X_2 , Y \to p(X_2)) \longrightarrow (X_1 , Y \to p(X_1))
$$
\newpage
\noindent définit un carré commutatif de ${\mathcal C}$
$$
\xymatrix{
X'_2 \ar[d]^{x'} \ar[r] &X_2 \ar[d]^x \\
X'_1 \ar[r]&X_1
}
$$
dont l'image par $p$ s'inscrit dans un diagramme commutatif:
$$
\xymatrix{
&p(X'_2) \ar[dd]^{p(x')} \ar[r] &p(X_2) \ar[dd]^{p(x)} \\
Y \ar[ru]^{\sim} \ar[rd]^{\sim} \\
&p(X'_1) \ar[r]&p(X_1)
}
$$
Par conséquent, le calcul de la limite
$$
p_! (P \times_{p^* Q_1} p^* Q_2)(Y) = \varinjlim_{(X,Y \to p(X)) \in Y \backslash {\mathcal C}} P(X) \times_{Q_1 (p(X))} Q_2 (p(X))
$$
peut être restreint à la sous-catégorie pleine de
$$
Y \backslash {\mathcal C}
$$
constituée des objets
$$
(X,Y \longrightarrow p(X))
$$
pour lesquels le morphisme de ${\mathcal B}$
$$
Y \longrightarrow p(X)
$$
est un isomorphisme, donc induit des bijections
$$
Q_1 (p(X)) \xrightarrow{ \ \sim \ } Q_1 (Y) \, ,
$$
$$
Q_2 (p(X)) \xrightarrow{ \ \sim \ } Q_2 (Y) \, .
$$

Cela montre que l'application
$$
p_! (P \times_{p^* Q_1} p^* Q_2) \longrightarrow p_! P(Y) \times_{Q_1(Y)} Q_2 (Y)
$$
est une bijection, comme annoncé.
\end{demo}

\subsection{Topologies de Giraud}\label{ssec4.1.4}

Comme on a vu dans la proposition \ref{prop1.3.12}, l'expression des images réciproques de sous-topos par un morphisme de topos
$$
f : {\mathcal E}' \longrightarrow {\mathcal E}
$$
se fait en général en termes de topologies engendrées.

Rappelons que les images réciproques de topologies sont données par des formules explicites dans le cas de morphismes de topos de préfaisceaux
$$
\widehat{\mathcal C} \longrightarrow \widehat{\mathcal B}
$$
induits par des fibrations $p : {\mathcal C} \to {\mathcal B}$.
\newpage
\begin{thm}[\textbf{(J. Giraud)}]\label{thm4.1.7} 
Soit un foncteur entre deux catégories essentiellement petites
$$
p : {\mathcal C} \longrightarrow {\mathcal B}
$$
qui est une fibration et induit donc un morphisme localement connexe de topos de préfaisceaux
$$
(p^* , p_*) : \widehat{\mathcal C} \longrightarrow \widehat{\mathcal B} \, .
$$

Soient un sous-topos du topos de base défini par une topologie $J$ de ${\mathcal B}$
$$
\widehat{\mathcal B}_J \ {\lhook\joinrel\relbar\joinrel\rightarrow} \ \widehat{\mathcal B}
$$
et son image réciproque par le morphisme $\widehat{\mathcal C} \to \widehat{\mathcal B}$, qui est un sous-topos
$$
\widehat{\mathcal C}_{J'} \ {\lhook\joinrel\relbar\joinrel\rightarrow} \ \widehat{\mathcal C} \, .
$$

Pour tout crible d'un objet $X$ de ${\mathcal C}$
$$
C \ {\lhook\joinrel\relbar\joinrel\rightarrow} \ y(X) \, ,
$$
notons
$$
C_p \ {\lhook\joinrel\relbar\joinrel\rightarrow} \ y(p(X))
$$
le crible de l'objet image $p(X)$ dans ${\mathcal B}$ constitué des morphismes
$$
Y \longrightarrow p(X)
$$
dont un relèvement horizontal dans ${\mathcal C}$
$$
X' \longrightarrow X
$$
est élément du crible $C$.

Alors un tel crible dans ${\mathcal C}$
$$
C \ {\lhook\joinrel\relbar\joinrel\rightarrow} \ y(X)
$$
est $J'$-couvrant si et seulement si son crible associé dans ${\mathcal B}$
$$
C_p \ {\lhook\joinrel\relbar\joinrel\rightarrow} \ y(p(X))
$$
est $J$-couvrant.
\end{thm}

\begin{remarks}
\noindent (i) Dans le cas général d'un foncteur
$$
p : {\mathcal C} \longrightarrow {\mathcal B}
$$
qui n'est pas nécessairement une fibration, la topologie $J'$ de ${\mathcal C}$ qui définit le sous-topos
$$
\widehat{\mathcal C}_{J'} \ {\lhook\joinrel\relbar\joinrel\rightarrow} \ \widehat{\mathcal C} 
$$
image réciproque par $\widehat{\mathcal C} \to \widehat{\mathcal B}$ d'un sous-topos défini par une topologie $J$ de ${\mathcal B}$
$$
\widehat{\mathcal B}_J \ {\lhook\joinrel\relbar\joinrel\rightarrow} \ \widehat{\mathcal B}
$$
est la plus petite topologie qui contient les cribles
$$
y(X) \times_{p^* y(X)} p^* C \ {\lhook\joinrel\relbar\joinrel\rightarrow} \ y(X)
$$
associés aux cribles $J$-couvrants des objets $Y$ de ${\mathcal B}$
$$
C \ {\lhook\joinrel\relbar\joinrel\rightarrow} \ y(Y)
$$
et aux morphismes de ${\mathcal B}$ définis sur les images par $p$ des objets $X$ de ${\mathcal C}$
$$
p(X) \xrightarrow{ \ b \ } Y
$$
vus comme des morphismes $y(X) \to p^* y(Y)$ de $\widehat{\mathcal C}$.

Ces cribles s'écrivent concrètement
$$
y(X) \times_{p^* y(Y)} p^* C = \left\{ X' \xrightarrow{ \ x \ } X \mid p(x) \in b^* C \right\}.
$$
(ii) Dans le cas de fibrations $p : {\mathcal C} \to {\mathcal B}$, les topologies $J'$ de ${\mathcal C}$ qui sont les images réciproques de topologies $J$ de ${\mathcal B}$ sont appelées les topologies de Giraud.
\end{remarks}

\begin{demo}
Pour tout objet $X$ de ${\mathcal C}$ et tout morphisme de ${\mathcal B}$
$$
Y \longrightarrow p(X) \, ,
$$
deux relèvements horizontaux de ce morphisme
$$
X'_1 \longrightarrow X \quad \mbox{et} \quad X'_2 \longrightarrow X
$$
sont reliés par un unique isomorphisme $X'_1 \xrightarrow{ \, \sim \, } X'_2$ qui fait commuter le triangle:
$$
\xymatrix{
X'_1 \ar[rr]^{\sim} \ar[rd] &&X'_2 \ar[ld] \\
&X
}
$$
Donc $X'_1 \to X$ est élément d'un crible $C$ de $X$ si et seulement si $X'_2 \to X$ est élément de $C$.

De plus, pour tous morphisme de ${\mathcal B}$
$$
Y_2 \longrightarrow Y_1 \longrightarrow p(X) \, ,
$$
tout relèvement horizontal de $Y_1 \to p(X)$ en un morphisme
$$
X_1 \longrightarrow X
$$ 
dont l'image par $p$ est munie d'une factorisation
$$
p(X_1) \xrightarrow{ \ \sim \ } Y_1 \longrightarrow X
$$
et tout relèvement horizontal du morphisme composé
$$
Y_2 \longrightarrow Y_1 \xrightarrow{ \ \sim \ } p(X_1)
$$
en un morphisme de ${\mathcal C}$
$$
X_2 \longrightarrow X_1 \, ,
$$
le morphisme de ${\mathcal C}$ composé
$$
X_2 \longrightarrow X_1 \longrightarrow X
$$
est un relèvement horizontal du morphisme de ${\mathcal B}$ composé
$$
Y_2 \longrightarrow Y_1 \longrightarrow p(X) \, .
$$

Il en résulte que les morphismes de ${\mathcal B}$ de but $p(X)$
$$
Y \longrightarrow p(X)
$$
dont un relèvement horizontal est élément du crible $C$ de $X$ constituent bien un crible $C_p$ de l'objet $p(X)$ de ${\mathcal B}$.

D'après la remarque (i), la topologie $J'$ de ${\mathcal C}$ qui définit le sous-topos de $\widehat{\mathcal C}$ image réciproque de $\widehat{\mathcal B}_J \hookrightarrow \widehat{\mathcal B}$ est  la topologie engendrée par les cribles d'objets $X$ de ${\mathcal C}$ de la forme
$$
C' = \left\{ X' \xrightarrow{ \ x \ } X \mid p(x) \in C \right\} = p^{-1} C
$$
pour des cribles $J$-couvrants $C$ des objets $p(X)$.

Comme tout morphisme $Y \to p(X)$ se relève en un morphisme horizontal dans ${\mathcal C}$, on voit que pour un tel crible
$$
C' = p^{-1} C \quad \mbox{d'un objet $X$},
$$
on a
$$
C'_p = C \, .
$$
Donc tous les générateurs $C' = p^{-1} C$ de la topologie $J'$ satisfont la condition que le crible associé
$$
C'_p
$$
est $J$-couvrant.

D'autre part, si $C$ est un crible d'un objet $X$ de ${\mathcal C}$, tout morphisme de ${\mathcal C}$
$$
X'' \xrightarrow{ \ x \ } X
$$
se factorise en
$$
X'' \longrightarrow X' \longrightarrow X
$$
pour un relèvement horizontal
$$
X' \longrightarrow X
$$
du morphisme $p(x) : p(X'') \to p(X')$.

Cela signifie que pour tout crible $C$ d'un objet $X$ de ${\mathcal C}$, on a
$$
C \supseteq p^{-1} C_p \, .
$$

En particulier, le crible $C$ est $J'$-couvrant si le crible associé $C_p$ est $J$-couvrant.

Pour conlure, il reste à montrer que les cribles d'objets $X$ de ${\mathcal C}$
$$
C \ {\lhook\joinrel\relbar\joinrel\rightarrow} \ y(X)
$$
qui satisfont la condition que le crible associé
$$
C_p \ {\lhook\joinrel\relbar\joinrel\rightarrow} \ y(p(X))
$$
soit $J$-couvrant forment une topologie de ${\mathcal C}$, c'est-à-dire satisfont les axiomes (1), (2), (3) du théorème \ref{thm2.2.2} (i).

Pour (1), il est évident que le crible maximal d'un objet $X$ de ${\mathcal C}$ satisfait la condition.

Pour (2), considérons un morphisme de ${\mathcal C}$
$$
X' \xrightarrow{ \ x \ } X
$$
un crible $C \hookrightarrow y(X)$ de $X$ dont le crible associé
$$
C_p \ {\lhook\joinrel\relbar\joinrel\rightarrow} \ y(p(X))
$$
est $J$-couvrant, et le crible de $X'$
$$
x^* C = C' \ {\lhook\joinrel\relbar\joinrel\rightarrow} \ y(X') \, .
$$
Pour tout morphisme de ${\mathcal B}$ de but $p(X')$
$$
Y \longrightarrow p(X')
$$ 
qui est élément de $p(x)^* C_p$, c'est-à-dire tel que le composé
$$
Y \longrightarrow p(X') \xrightarrow{ \ p(x) \ } p(X)
$$
soit élément de $C_p$, son relèvement horizontal
$$
X'' \xrightarrow{ \ x' \ } X'
$$
s'inscrit dans un diagramme commutatif
$$
\xymatrix{
X'' \ar[d] \ar[r]^{x'} &X' \ar[d] \ar[rd]^x \\
\overline X'' \ar[r] &\overline X' \ar[r] &X
}
$$
où $\overline X' \to X$ désigne un relèvement horizontal de
$$
p(X') \xrightarrow{ \ p(x) \ } p(X)
$$
et $\overline X'' \to \overline X'$ désigne un relèvement horizontal de
$$
Y \longrightarrow p(X') \xrightarrow{ \ \sim \ } p(\overline X') \, .
$$
Cela signifie que
$$
x \circ x' \in C
$$
soit
$$
x' \in C' \, .
$$
Ainsi, on a une inclusion
$$
p(x)^* C_p \subseteq C'_p
$$
et donc le crible $C'_p$ de $p(X')$ est $J$-couvrant.

Pour (3) enfin, il suffit de considérer un crible $C$ d'un objet $X$ de ${\mathcal C}$ tel qu'existe une famille de morphismes horizontaux 
$$
X_i \xrightarrow{ \ x_i \ } X \, , \quad i \in I \, ,
$$
tels que les cribles des objets $p(X_i)$
$$
(x_i^* C)_p \ {\lhook\joinrel\relbar\joinrel\rightarrow} \ y(p(X_i))
$$
soient tous $J$-couvrants.

Comme les composés de morphismes horizontaux
$$
X'_i \longrightarrow X_i
$$
avec les morphismes
$$
X_i \xrightarrow{ \ x_i \ } X 
$$
sont des morphismes horizontaux, on obtient que le crible
$$
C_p \ {\lhook\joinrel\relbar\joinrel\rightarrow} \ y(p(X))
$$
est $J$-couvrant.

Donc la condition (3) aussi est vérifiée.

Cela achève la démonstration du théorème.
\end{demo}

%Subsec IV.1.e
\subsection{Images réciproques des sous-topos par des morphismes localement connexes}\label{4.1.5}

La caractérisation des topologies de Giraud a une variante plus abstraite et plus générale dans le contexte des morphismes localement connexes de topos
$$
f : {\mathcal E}' \longrightarrow {\mathcal E} \, .
$$

On a vu dans le théorème \ref{thm2.3.2} et son corrolaire \ref{cor2.3.3} que les sous-topos d'un topos ${\mathcal E}$
$$
{\mathcal E}_J \ {\lhook\joinrel\relbar\joinrel\rightarrow} \ {\mathcal E}
$$
correspondent bijectivement aux ``topologies'' $J$ de ${\mathcal E}$ c'est-à-dire aux classes des monomorphismes de ${\mathcal E}$
$$
C \ {\lhook\joinrel\relbar\joinrel\rightarrow} \ E
$$
qui satisfont les axiomes (1), (2) et (3) du théorème \ref{thm2.3.2} (i).

La topologie $J$ qui correspond à un sous-topos
$$
\xymatrix{
{\mathcal E}_J \,  \ar@{^{(}->}[rr]^{(j^*,j_*)} &&{\mathcal E}
}
$$
est constituée des monomorphismes de ${\mathcal E}$
$$
C  \ {\lhook\joinrel\relbar\joinrel\rightarrow} \ E
$$
dont le transformé par le foncteur $j^* : {\mathcal E} \to {\mathcal E}_J$ est un isomorphisme.

Démontrons la variante suivante du théorème \ref{thm4.1.7} de Giraud:

\begin{thm}\label{thm4.1.8}
Soit un morphisme de topos
$$
f : {\mathcal E}' \longrightarrow {\mathcal E} 
$$
qui est localement connexe, c'est-à-dire dont la composante d'image réciproque
$$
f^* : {\mathcal E} \longrightarrow {\mathcal E}'
$$
admet non seulement un adjoint à droite $f_*$ mais aussi un adjoint à gauche
$$
f_! : {\mathcal E}' \longrightarrow {\mathcal E} 
$$
qui est une fibration.

On a:

\noindent (i) Pour tout monomorphisme de ${\mathcal E}'$
$$
C \ {\lhook\joinrel\relbar\joinrel\rightarrow} \ X
$$
il existe un plus grand sous-objet de $f_! X$, noté
$$
C_f \ {\lhook\joinrel\relbar\joinrel\rightarrow} \ f_! X \, ,
$$
dont le relèvement horizontal
$$
f^* C_f \times_{f^* f_! X} X \ {\lhook\joinrel\relbar\joinrel\rightarrow} \ X
$$
se factorise à travers $C \hookrightarrow X$.

\noindent (ii) Considérons un sous-topos du topos de base
$$
{\mathcal E}_J \ {\lhook\joinrel\relbar\joinrel\rightarrow} \ {\mathcal E}
$$
qui est défini par une topologie $J$ de ${\mathcal E}$, puis la topologie $J'$ de ${\mathcal E}'$ qui définit le sous-topos
$$
f^{-1} {\mathcal E}_J = {\mathcal E}'_{J'} \ {\lhook\joinrel\relbar\joinrel\rightarrow} \ {\mathcal E}'
$$
image réciproque de ${\mathcal E}_J \hookrightarrow {\mathcal E}$ par le morphisme $f : {\mathcal E}' \to {\mathcal E}$.

Alors un monomorphisme de ${\mathcal E}'$
$$
C \ {\lhook\joinrel\relbar\joinrel\rightarrow} \ X
$$
est élément de la topologie $J'$ si et seulement si le monomorphisme de ${\mathcal E}$ qui lui est associé
$$
C_f \ {\lhook\joinrel\relbar\joinrel\rightarrow} \ f_! X
$$
est élément de la topologie $J$.
\end{thm}

\begin{remarks}
\noindent (i) Dans le cas général d'un morphisme de topos
$$
f : {\mathcal E}' \longrightarrow {\mathcal E} 
$$
qui n'est pas nécessairement localement connexe, la topologie $J'$ de ${\mathcal E}'$ qui définit le sous-topos
$$
{\mathcal E}'_{J'} \ {\lhook\joinrel\relbar\joinrel\rightarrow} \ {\mathcal E}'
$$
image réciproque par $f$ d'un sous-topos
$$
{\mathcal E}_{J} \ {\lhook\joinrel\relbar\joinrel\rightarrow} \ {\mathcal E}
$$
défini par une topologie $J$, est la plus petite topologie qui contient les monomorphismes
$$
f^* S \ {\lhook\joinrel\relbar\joinrel\rightarrow} \ f^* Y
$$
qui sont les transformés par $f^* : {\mathcal E} \to {\mathcal E}'$ des monomorphismes de ${\mathcal E}$
$$
S \ {\lhook\joinrel\relbar\joinrel\rightarrow} \ Y
$$
qui sont éléments de la topologie $J$.

\noindent (ii) Si $f = (f_! , f^* , f_*) : {\mathcal E}' \to {\mathcal E}$ est un morphisme de topos localement connexe, le relèvement horizontal dans ${\mathcal E}'$ de tout morphisme de ${\mathcal E}$ de la forme
$$
Y \longrightarrow f_! X
$$
s'identifie à
$$
f^* Y \times_{f^* f_! X} X \longrightarrow X \, .
$$
En effet, si $X' \xrightarrow{ \, x \, } X$ désigne un tel relèvement, il est muni d'un isomorphisme
$$
f_! X' \xrightarrow{ \ \sim \ } Y
$$
dont le composé avec $Y \to f_! X$ est $f_! (x) : f_! X' \to f_! X$ et, d'après le lemme \ref{lem4.1.2}, le carré commutatif
$$
\xymatrix{
X' \ar[d] \ar[r] &X \ar[d] \\
f^* f_! X' \ar[r] &f^* f_! X
}
$$
est cartésien.
\end{remarks}

\begin{demo}
\noindent (i) Le foncteur
$$
{\mathcal E} / f_! X \longrightarrow {\mathcal E}'/X \, ,
$$
$$
(Y \to f_! X) \longmapsto (f^* Y \times_{f^* f_! X} X \to X)
$$
respecte les limites et les colimites arbitraires.

En particulier, il respecte les réunions arbitraires de sous-objets.

Pour un sous-objet $C \hookrightarrow X$, la réunion des sous-objets
$$
S \ {\lhook\joinrel\relbar\joinrel\rightarrow} \ f_! X
$$
tels que
$$
f^* S \times_{f^* f_! X} X \ {\lhook\joinrel\relbar\joinrel\rightarrow} \ X
$$
se factorise à travers $C \hookrightarrow X$, répond donc à la question posée.

\noindent (ii) D'après la remarque (i), la topologie $J'$ de ${\mathcal E}'$ qui définit le sous-topos image réciproque par $f : {\mathcal E}' \to {\mathcal E}$ du sous-topos ${\mathcal E}_J \hookrightarrow {\mathcal E}$ défini par la topologie $J$ est engendrée par les éléments
$$
C = f^* S \ {\lhook\joinrel\relbar\joinrel\rightarrow} \ f^* Y = X
$$
qui sont les transformés par $f^*$ des éléments de $J$
$$
S \ {\lhook\joinrel\relbar\joinrel\rightarrow} \ Y \, .
$$

Or, pour tout tel élément de $J$
$$
S \ {\lhook\joinrel\relbar\joinrel\rightarrow} \ Y \, ,
$$
son produit fibré avec le morphisme d'adjonction
$$
f_! f^* Y \longrightarrow Y
$$
est un élément de $J$
$$
S \times_Y f_! f^* Y \ {\lhook\joinrel\relbar\joinrel\rightarrow} \ f_! f^* Y = f_! X \, ,
$$
et le transformé de celui-ci par le foncteur
$$
f^* (\bullet) \times_{f^* f_! X} X =  f^* (\bullet) \times_{f^* f_! f^* Y} f^* Y
$$
est
$$
(f^* S \times_{f^* Y} f^* f_! f^* Y) \times_{f^* f_! f^* Y} f^* Y = f^* S = C \ {\lhook\joinrel\relbar\joinrel\rightarrow} \ X \, .
$$
\newpage
Donc le sous-objet
$$
C_f \ {\lhook\joinrel\relbar\joinrel\rightarrow} \ f_! X \, ,
$$
contient le sous-objet
$$
S \times_Y f_! f^* Y \ {\lhook\joinrel\relbar\joinrel\rightarrow} \ f_! f^* Y = f_! X 
$$
et il est élément de $J$.

Ainsi, tous les générateurs de la topologie $J'$ satisfont la condition du théorème.

Réciproquement, considérons un monomorphisme de ${\mathcal E}'$
$$
C \ {\lhook\joinrel\relbar\joinrel\rightarrow} \ X 
$$
tel que 
$$
C_f \ {\lhook\joinrel\relbar\joinrel\rightarrow} \ f_! X 
$$
soit élément de $J$.

Alors
$$
f^* C_f \ {\lhook\joinrel\relbar\joinrel\rightarrow} \ f^* f_! X 
$$
est élément de $J'$ donc aussi le monomorphisme
$$
f^* C_f \times_{f^* f_! X} X \ {\lhook\joinrel\relbar\joinrel\rightarrow} \ X \, .
$$
Comme il se factorise à travers le monomorphisme
$$
C \ {\lhook\joinrel\relbar\joinrel\rightarrow} \ X \, ,
$$
celui-ci est élément de la topologie $J'$.

Pour conclure, il reste à montrer que les monomorphismes de ${\mathcal E}'$
$$
C \ {\lhook\joinrel\relbar\joinrel\rightarrow} \ X 
$$
qui satisfont la condition que le monomorphisme de ${\mathcal E}$ associé
$$
C_f \ {\lhook\joinrel\relbar\joinrel\rightarrow} \ f_! X 
$$
soit élément de $J$, forment une topologie de ${\mathcal E}'$, c'est-à-dire satisfont les axiomes (1), (2), (3) du théorème \ref{thm2.3.2} (i)

La propriété (1) résulte de ce que $C_f \hookrightarrow f_! X$ est un isomorphisme si $C \hookrightarrow X$ est un isomorphisme.

Pour les axiomes (2) et (3), on a besoin du lemme suivant:

\begin{lem}\label{lem4.1.9}
Considérons un morphisme localement connexe de topos
$$
(f_! , f^* , f_*) : {\mathcal E}' \longrightarrow {\mathcal E}
$$
et, pour tout monomorphisme de ${\mathcal E}'$
$$
C \ {\lhook\joinrel\relbar\joinrel\rightarrow} \ X \, ,
$$
le plus grand sous-objet dans ${\mathcal E}$
$$
C_f \ {\lhook\joinrel\relbar\joinrel\rightarrow} \ f_! X 
$$
dont le relèvement horizontal est contenu dans $C \hookrightarrow X$.

Alors, pour tout morphisme de ${\mathcal E}'$
$$
X' \xrightarrow{ \ x \ } X \, ,
$$
tout monomorphisme
$$
C \ {\lhook\joinrel\relbar\joinrel\rightarrow} \ X 
$$
et son transformé par changement de base
$$
C' = C \times_X X' \ {\lhook\joinrel\relbar\joinrel\rightarrow} \ X' \, ,
$$
on a:

\noindent (i) Le sous-objet
$$
C'_f \ {\lhook\joinrel\relbar\joinrel\rightarrow} \ f_! X' 
$$
contient le sous-objet
$$
C_f \times_{f_! X} f_! X' \ {\lhook\joinrel\relbar\joinrel\rightarrow} \ f_! X' \, . 
$$
(ii) Il y a même égalité de sous-objets
$$
C'_f = C_f \times_{f_! X} f_! X'
$$
si le morphisme
$$
X' \xrightarrow{ \ x \ } X 
$$
est un épimorphisme.
\end{lem}

\begin{demolem}
\noindent (i) Pour tout sous-objet
$$
S \ {\lhook\joinrel\relbar\joinrel\rightarrow} \ f_! X 
$$
et son image réciproque
$$
S' = S \times_{f_! X} f_! X' \ {\lhook\joinrel\relbar\joinrel\rightarrow} \ f_! X' \, , 
$$
le relèvement horizontal
$$
f^* S' \times_{f^* f_! X'} X' \ {\lhook\joinrel\relbar\joinrel\rightarrow} \ X'
$$
s'identifie à
$$
(f^* S \times_{f^* f_! X} f^* f_! X') \times_{f^* f_! X'} X' \ {\lhook\joinrel\relbar\joinrel\rightarrow} \ X '
$$
et donc à
$$
(f^* S \times_{f^* f_! X} X) \times_X X' \ {\lhook\joinrel\relbar\joinrel\rightarrow} \ X ' \, .
$$
Il en résulte que le relèvement horizontal de
$$
S' \ {\lhook\joinrel\relbar\joinrel\rightarrow} \ f_! X'
$$
est contenu dans
$$
C' = C \times_X X' \ {\lhook\joinrel\relbar\joinrel\rightarrow} \ X' 
$$
si le relèvement horizontal de
$$
S \ {\lhook\joinrel\relbar\joinrel\rightarrow} \ f_! X 
$$
est contenu dans $C \hookrightarrow X$.

\noindent (ii) Si $X' \xrightarrow{ \, x \, } X$ est un épimorphisme, l'objet $X$ s'écrit dans ${\mathcal E}'$ comme une colimite
$$ 
X = \varinjlim \left(\xymatrix{X' \times_X X' \dar[r]^-{^{^{\mbox{\scriptsize $p_1$}}}}_-{p_2} &X'} \right) .
$$
Comme le foncteur $f_! : {\mathcal E}' \to {\mathcal E}$ est un adjoint à gauche, il respecte les colimites, donc l'objet $f_! X$ de ${\mathcal E}$ s'écrit comme la colimite
$$
f_! X = \varinjlim \left(\xymatrix{f_! (X' \times_X X') \dar[rr]^-{^{^{\mbox{\scriptsize $f_!(p_1)$}}}}_-{f_!(p_2)} &&f_! X'} \right) .
$$
Pour tout sous-objet
$$
S \ {\lhook\joinrel\relbar\joinrel\rightarrow} \ f_! X'
$$
qui satisfait la condition
$$
f^* S \times_{f^* f_! X'} X' \subseteq C \times_X X' \, ,
$$
les deux sous-objets
$$
S_1 = f_! (p_1)^{-1} S \quad \mbox{et} \quad S_2 = f_! (p_2)^{-1} S \quad \mbox{de} \quad f_! (X' \times_X X')
$$
satisfont la condition
$$
f^* S_i \times_{f^* f_! (X' \times_X X')} X' \times_X X' \subseteq C \times_X X' \times_X X'
$$
donc leurs images par les projections
$$
S_{1,2} = f_! (p_2)_* (S_1)
$$
$$
S_{2,1} = f_! (p_1)_* (S_2)
$$
sont encore deux sous-objets de $f_! X'$ dont les relèvements horizontaux sont contenus dans $C' = C \times_X X'$.

Par maximalité du sous-objet
$$
C'_f \ {\lhook\joinrel\relbar\joinrel\rightarrow} \ f_! X' \, ,
$$
il a donc la même image réciproque par les deux morphismes
$$
\xymatrix{f_! (X' \times_X X') \dar[rr]^-{^{^{\mbox{\scriptsize $f_!(p_1)$}}}}_-{f_!(p_2)} &&f_! X'}  .
$$
Il provient donc d'un sous-objet de la colimite
$$
S \ {\lhook\joinrel\relbar\joinrel\rightarrow} \ f_! X
$$
au sens que
$$
S \times_{f_! X} f_! X' = C'_f \, .
$$
Cela implique
$$
f^* C'_f \times_{f^* f_! X'} X' = (f^* S \times_{f^* f_! X} X) \times_X X' \, .
$$
Comme $X' \xrightarrow{ \, x \, } X$ est un épimorphisme, la condition
$$
f^* C'_f \times_{f^* f_! X} X' \subseteq C' = C \times_X X'
$$
implique
$$
f^* S \times_{f^* f_! X} X \subseteq C
$$
et donc
$$
S \subseteq C_f
$$
puis
$$
C'_f \subseteq C_f \times_{f_! X} f_! X' \, .
$$
Cette inclusion est en fait une égalité puisque l'inclusion en sens inverse est connue d'après (i).

Cela termine la démonstration du lemme.
\end{demolem}

\noindent {\bf Fin de la démonstration du théorème \ref{thm4.1.8} (ii).}

Il reste à vérifier les axiomes (2) et (3) des topologies sur ${\mathcal E}'$.

Pour (2), considérons un morphisme de ${\mathcal E}'$
$$
X' \longrightarrow X
$$
et un monomorphisme
$$
C \ {\lhook\joinrel\relbar\joinrel\rightarrow} \ X
$$
tel que le monomorphisme $C_f \hookrightarrow f_! X$ de ${\mathcal E}$ soit élément de $J$.

D'après le lemme \ref{lem4.1.9} (i), le sous-objet
$$
C'_f \ {\lhook\joinrel\relbar\joinrel\rightarrow} \ f_! X'
$$
associé au monomorphisme
$$
C' = C \times_X X' \ {\lhook\joinrel\relbar\joinrel\rightarrow} \ X'
$$
contient le sous-objet
$$
C_f \times_{f_! X} f_! X' \ {\lhook\joinrel\relbar\joinrel\rightarrow} \ f_! X'
$$
et donc il est élément de $J$.

Cela prouve l'axiome (2).

Pour l'axiome (3) enfin, considérons un monomorphisme de ${\mathcal E}'$
$$
C \ {\lhook\joinrel\relbar\joinrel\rightarrow} \ X
$$
tel qu'existent un monomorphisme
$$
C' \ {\lhook\joinrel\relbar\joinrel\rightarrow} \ X
$$
et une famille globalement épimorphique de morphismes
$$
X_k \longrightarrow C' \, , \quad k \in K \, ,
$$
pour lequels les morphismes de ${\mathcal E}$ associés
$$
C'_f \ {\lhook\joinrel\relbar\joinrel\rightarrow} \ f_! X
$$
et
$$
(C \times_X X_k)_f \ {\lhook\joinrel\relbar\joinrel\rightarrow} \ f_! X_k \, , \quad k \in K \, ,
$$
sont éléments de $J$.

Quitte à remplacer $C' \hookrightarrow X$ par un sous-objet
$$
C'' = f^* C'_f \times_{f^* f_! X} X \ {\lhook\joinrel\relbar\joinrel\rightarrow} \ X
$$
qui est le relèvement horizontal de
$$
C'_f \ {\lhook\joinrel\relbar\joinrel\rightarrow} \ f_! X \, ,
$$
et la famille globalement épimorphique
$$
(X_k \longrightarrow C')_{k \in K}
$$
par sa restriction au-dessus de $C''$
$$
(X_k \times_{C'} C'' \longrightarrow C'')_{k \in K} \, ,
$$
on peut supposer que
$$
C' \ {\lhook\joinrel\relbar\joinrel\rightarrow} \ X
$$
est le relèvement horizontal de
$$
C'_f \ {\lhook\joinrel\relbar\joinrel\rightarrow} \ f_! X \, .
$$

Alors, pour tout sous-objet
$$
S \ {\lhook\joinrel\relbar\joinrel\rightarrow} \ C'_f \ {\lhook\joinrel\relbar\joinrel\rightarrow} \ f_! X 
$$
le relèvement horizontal de
$$
S \ {\lhook\joinrel\relbar\joinrel\rightarrow} \ f_! X 
$$
est le composé de $C' \hookrightarrow X$ et du relèvement horizontal de
$$
S \ {\lhook\joinrel\relbar\joinrel\rightarrow} \ C'_f = f_! C' \, .
$$
De plus, $S \hookrightarrow C'_f$ est dans $J$ si et seulement si $S \hookrightarrow f_! X$ est dans $J$.

Cela nous ramène au cas où
$$
C'_f = f_! X \, ,
$$
$$
C' = X
$$
et la famille
$$
(X_k \longrightarrow X)_{k \in K}
$$
est globlament épimorphique.

Posant
$$
X' = \coprod_{k \in K} X_k \, ,
$$
on a
$$
C \times _X X' = \coprod_{k \in K} C \times_X X_k \, ,
$$
$$
f_! X' = \coprod_{k \in K}  f_! X_k
$$
et, notant $C' = C \times_X X'$,
$$
C'_f = \coprod_{k \in K}  (C \times_X X_k)_f \, .
$$
Comme  le morphisme
$$
X' \longrightarrow X
$$
est un épimorphisme, on obtient d'après le lemme \ref{lem4.1.9} (ii) la formule
$$
C'_f = C_f \times_{f_! X} f_! X' \, .
$$
Or, le monomorphisme
$$
C'_f \ {\lhook\joinrel\relbar\joinrel\rightarrow} \ f_! X' 
$$
est élément de $J$ puisqu'il est somme des monomorphismes
$$
(C \times_X X_k)_f \ {\lhook\joinrel\relbar\joinrel\rightarrow} \ f_! X_k
$$
qui sont éléments de $J$.

Comme le morphisme de ${\mathcal E}$
$$
f_! X' \longrightarrow f_! X
$$
est un épimorphisme, on conclut comme voulu que le monomorphisme
$$
C_f \ {\lhook\joinrel\relbar\joinrel\rightarrow} \ f_! X 
$$
est élément de $J$.

Ainsi, l'axiome (3) est bien vérifié. 
\end{demo}

Le théorème \ref{thm4.1.8} implique en particulier le résultat suivant (dont le point (i) coincide avec le Corollaire 4.58 de \cite{Denseness}):

\begin{cor}\label{cor4.1.10}
Soient un morphisme localement connexe de topos
$$
f = (f_! , f^* , f_*) : {\mathcal E}' \longrightarrow {\mathcal E} \, ,
$$
un sous-topos du topos de base ${\mathcal E}$ défini par une topologie $J$
$$
(j^* , j_*) : {\mathcal E}_J \ {\lhook\joinrel\relbar\joinrel\rightarrow} \ {\mathcal E} 
$$
et son image réciproque le sous-topos défini par une topologie $J'$
$$
(j'^* , j'_*) : {\mathcal E}'_{J'} \ {\lhook\joinrel\relbar\joinrel\rightarrow} \ {\mathcal E}' \, .
$$
Alors:

\noindent (i) Le morphisme $f : {\mathcal E}' \to {\mathcal E}$ induit un morphisme localement connexe
$$
f_J = (f_! , f^* , f_*) : {\mathcal E}'_{J'} \longrightarrow {\mathcal E}_J \, .
$$
(ii) Les deux carrés
$$
\xymatrix{
{\mathcal E}'_{J'} \ar@{^{(}->}[r]^{j'_*} &{\mathcal E}' \\
{\mathcal E}_J \ar[u]_{f^*} \ar@{^{(}->}[r]^{j_*} &{\mathcal E} \ar[u]_{f^*}
} \qquad  \xymatrix{
{\mathcal E}'_{J'} \ar[d]_{\mbox{et} \qquad \qquad f_!} &{\mathcal E}' \ar[l]_{j'^*} \ar[d]^{f_!} \\
{\mathcal E}_J &{\mathcal E} \ar[l]_{j^*}
}
$$
sont commutatifs.
\end{cor}

\begin{demo}
\noindent (ii) Montrons d'abord que le foncteur
$$
f^* : {\mathcal E} \longrightarrow {\mathcal E}'
$$
envoie la sous-catégorie pleine
$$
\xymatrix{
{\mathcal E}_J \, \ar@{^{(}->}[r]^{j_*} &{\mathcal E} 
}
$$
dans la sous-catégorie pleine
$$
\xymatrix{
{\mathcal E}'_{J'} \ar@{^{(}->}[r]^{j_*} &{\mathcal E}' \, .
}
$$
Considérons pour cela la classe $J''$ des monomorphismes de ${\mathcal E}'$
$$
C \ {\lhook\joinrel\relbar\joinrel\rightarrow} \ X 
$$
tels que, pour tout objet $Y$ de ${\mathcal E}_J \hookrightarrow {\mathcal E}$ et tout morphisme de ${\mathcal E}'$
$$
X' \longrightarrow X \, ,
$$
l'application
$$
{\rm Hom} (X' , f^* Y) \longrightarrow {\rm Hom} (C \times_X X' , f^* Y)
$$
est une bijection.

On sait d'après le théorème \ref{thm2.3.2} que $J''$ est une topologie.

Il s'agit de prouver que $J''$ contient $J'$.

Il suffit pour cela de montrer que $J''$ contient une famille de générateurs de $J'$.

Or, d'après le théorème \ref{thm4.1.8} (ii), la topologie $J'$ est engendrée par les monomorphismes
$$
C = f^* S \times_{f^* f_! X} X \\ {\lhook\joinrel\relbar\joinrel\rightarrow} \ X 
$$
qui sont les relèvements horizontaux de monomorphismes
$$
S \ {\lhook\joinrel\relbar\joinrel\rightarrow} \ f_! X 
$$
éléments de $J$.

Pour tout morphisme $X' \to X$,
$$
C' = C \times_X X' \ {\lhook\joinrel\relbar\joinrel\rightarrow} \ X'
$$
est le relèvement horizontal de
$$
S' = S \times_{f_! X} f_! X' \ {\lhook\joinrel\relbar\joinrel\rightarrow} \ f_! X'
$$
si $C \hookrightarrow X$ est le relèvement horizontal de $S \hookrightarrow f_! X$, et on a pour tout objet $Y$ de ${\mathcal E}$
$$
{\rm Hom} (X' , f^* Y) = {\rm Hom} (f_! X' , Y) \, ,
$$
$$
{\rm Hom} (C' , f^* Y) = {\rm Hom} (f_! C' , Y) 
$$
avec
$$
f_! C' = S' \, .
$$
D'où il résulte comme voulu que l'application
$$
{\rm Hom} (X' , f^* Y) \longrightarrow {\rm Hom} (C' , f^* Y)
$$
est une bijection si $Y$ est un objet de
$$
\xymatrix{
{\mathcal E}_J \, \ar@{^{(}->}[r]^{j_*} &{\mathcal E} \, .
}
$$

Ainsi, le foncteur d'image réciproque 
$$
f^* : {\mathcal E} \longrightarrow {\mathcal E}'
$$
se restreint en un foncteur
$$
f^* : {\mathcal E}_J \longrightarrow {\mathcal E}'_{J'}
$$
qui est nécessairement la composante d'image réciproque du morphisme de topos induit
$$
f_J : {\mathcal E}'_{J'} \longrightarrow {\mathcal E}_J \, .
$$
La commutativité du carré
$$
\xymatrix{
{\mathcal E}'_{J'} \ar@{^{(}->}[r]^{j'_*} &{\mathcal E}' \\
{\mathcal E}_J \ar[u]_{f^*} \ar@{^{(}->}[r]^{j_*} &{\mathcal E} \ar[u]_{f^*}
}
$$
implique que le foncteur
$$
f_! : {\mathcal E}'_{J'} \longrightarrow {\mathcal E}_J 
$$
défini par restriction à
$$
\xymatrix{
{\mathcal E}'_{J'} \ar@{^{(}->}[r]^{j'_*} &{\mathcal E}'
}
$$
du foncteur composé
$$
{\mathcal E}' \xrightarrow{ \ f_! \ } {\mathcal E} \xrightarrow{ \ j^* \ } {\mathcal E}_J
$$
est un adjoint à gauche du foncteur
$$
f^* : {\mathcal E}_J \longrightarrow {\mathcal E}'_{J'} \, .
$$

Par passage aux adjoints à gauche, on obtient finalement un carré commutatif:
$$
\xymatrix{
{\mathcal E}'_{J'} \ar[d]_{f_!} &{\mathcal E}' \ar[l]_-{j'^*} \ar[d]^{f_!} \\
{\mathcal E}_J &{\mathcal E} \ar[l]_-{j^*}
}
$$
(ii) La composante d'image réciproque
$$
f^* : {\mathcal E}_J \longrightarrow {\mathcal E}'_{J'} 
$$
du morphisme de topos
$$
f_J : {\mathcal E}'_{J'} \longrightarrow {\mathcal E}_J
$$
admet un adjoint à gauche
$$
f_! : {\mathcal E}'_{J'} \longrightarrow {\mathcal E}_J \, .
$$

\newpage

\noindent Il reste seulement à montrer que, pour tout objet $X$ de ${\mathcal E}'_{J'}$ et tous morphismes de ${\mathcal E}_J$
$$
\xymatrix{ &f_! X \ar[d] \\
Y' \ar[r] &Y
}
$$
le morphisme canonique
$$
f_! (X \times_{f^* Y} f^* Y') \longrightarrow f_! X \times_Y Y'
$$
est un isomorphisme.

Cela résulte de ce que

\medskip

$\left\{\begin{matrix}
\bullet &\mbox{cette propriété est vérifiée par la paire de morphismes adjoints} \hfill \\
&\left( {\mathcal E}' \xrightarrow{ \ f_! \ } {\mathcal E} , {\mathcal E} \xrightarrow{ \ f^* \ } {\mathcal E}' \right), \\
\bullet &\mbox{le foncteur} \hfill \\
&f^* : {\mathcal E}_J \longrightarrow {\mathcal E}'_{J'} \\
&\mbox{est la restriction du foncteur} \hfill \\
&f^* : {\mathcal E} \longrightarrow {\mathcal E}' \, , \\
\bullet &\mbox{le foncteur} \hfill \\
&f_! : {\mathcal E}'_{J'} \longrightarrow {\mathcal E}_J \\
&\mbox{est la restriction à ${\mathcal E}'_{J'} \hookrightarrow {\mathcal E}'$ du foncteur composé} \hfill \\
&{\mathcal E}' \xrightarrow{ \ f_! \ } {\mathcal E} \xrightarrow{ \ f^* \ } {\mathcal E}_J \, , \\
\bullet &\mbox{les produits fibrés sont respectés par} \hfill \\
&j'_* : {\mathcal E}'_{J'} \ {\lhook\joinrel\relbar\joinrel\rightarrow} \ {\mathcal E}' \\
&j_* : {\mathcal E}_J \ {\lhook\joinrel\relbar\joinrel\rightarrow} \ {\mathcal E} \\
&\mbox{et} \hfill \\
&j^* : {\mathcal E} \longrightarrow {\mathcal E}_J \, . 
\end{matrix} \right.$

\medskip

Ainsi, le morphisme de topos
$$
f_J : {\mathcal E}'_{J'} \ {\lhook\joinrel\relbar\joinrel\rightarrow} \ {\mathcal E}_J
$$
est bien localement connexe.
\end{demo}

%Subsec IV.2.a.b et c.
\section{Images directes extraordinaires et images réciproques des réunions}\label{sec4.2}

\subsection{Images directes extraordinaires par un morphisme localement connexe}\label{4.2.1}

On déduit du théorème \ref{thm4.1.8} le corollaire suivant:

\begin{cor}\label{cor4.2.1}
Soit un morphisme localement connexe de topos
$$
f = (f_! , f^* , f_*) : {\mathcal E} ' \longrightarrow {\mathcal E} \, .
$$

Considérons le morphisme d'image réciproque entre les ensembles ordonnés des sous-topos de ${\mathcal E}$ et ${\mathcal E}'$
$$
f^{-1} : {\rm ST} ({\mathcal E}) \longrightarrow {\rm ST} ({\mathcal E}') \, .
$$

Alors:

\noindent (i) Ce morphisme $f^{-1}$ respecte non seulement les intersections arbitraires de sous-topos mais aussi les réunions arbitraires.

\noindent (ii) Il possède non seulement un adjoint à droite qui est le morphisme d'image directe des sous-topos
$$
f_* : {\rm ST} ({\mathcal E}') \longrightarrow {\rm ST} ({\mathcal E})
$$
mais aussi un adjoint à gauche, appelé le morphisme d'image directe extraordinaire des sous-topos
$$
f_! : {\rm ST} ({\mathcal E}') \longrightarrow {\rm ST} ({\mathcal E}) \, ,
$$
qui est caractérisé par la propriété que, pour tous sous-topos
$$
{\mathcal E}_1 \ {\lhook\joinrel\relbar\joinrel\rightarrow} \ {\mathcal E} \quad \mbox{et} \quad {\mathcal E}'_1 \ {\lhook\joinrel\relbar\joinrel\rightarrow} \ {\mathcal E}' \, ,
$$
on a
$$
f_! ({\mathcal E}'_1) \supseteq {\mathcal E}_1
$$
si et seulement si
$$
{\mathcal E}'_1 \supseteq f^{-1} {\mathcal E}_1 \, .
$$
(iii) Pour tout sous-topos de ${\mathcal E}'$
$$
{\mathcal E}'_{J'} \ {\lhook\joinrel\relbar\joinrel\rightarrow} \ {\mathcal E}'
$$
qui est défini par une topologie de ${\mathcal E}'$ engendrée par une famille de monomorphismes
$$
(C_i \ {\lhook\joinrel\relbar\joinrel\rightarrow} \ X_i)_{i \in I} \, ,
$$
son image directe extraordinaire par $f$
$$
f_! ({\mathcal E}'_{J'}) = {\mathcal E}_J \ {\lhook\joinrel\relbar\joinrel\rightarrow} \ {\mathcal E}
$$
est définie par la topologie $J$ de ${\mathcal E}$ engendrée par la famille des monomorphismes
$$
(C_i)_f \ {\lhook\joinrel\relbar\joinrel\rightarrow} \ f_! X_i \, , \quad i \in I \, ,
$$
où, pour tout monomorphisme de ${\mathcal E}'$
$$
C \ {\lhook\joinrel\relbar\joinrel\rightarrow} \ X \, ,
$$
on note
$$
C_f \ {\lhook\joinrel\relbar\joinrel\rightarrow} \ f_! X
$$
le plus grand sous-objet
$$
S \ {\lhook\joinrel\relbar\joinrel\rightarrow} \ f_! X
$$
dont le relèvement horizontal dans ${\mathcal E}'$
$$
f^* S \times_{f^* f_! X} X \ {\lhook\joinrel\relbar\joinrel\rightarrow} \ X
$$
est contenu dans $C \hookrightarrow X$.
\end{cor}

\begin{demo}
Dans la situation de (iii) d'un sous-topos de ${\mathcal E}'$
$$
{\mathcal E}'_{J'} \ {\lhook\joinrel\relbar\joinrel\rightarrow} \ {\mathcal E}'
$$
défini par une topologie $J'$ de ${\mathcal E}'$ engendrée par une famille de monomorphismes,
$$
C_i \ {\lhook\joinrel\relbar\joinrel\rightarrow} \ X_i \, , \quad i \in I \, ,
$$
considérons donc le sous-topos de ${\mathcal E}$
$$
{\mathcal E}_J \ {\lhook\joinrel\relbar\joinrel\rightarrow} \ {\mathcal E}
$$
défini par la topologie $J$ engendrée par les monomorphismes
$$
(C_i)_f \ {\lhook\joinrel\relbar\joinrel\rightarrow} \ f_! X_i \, , \quad i \in I \, .
$$
\newpage
Pour tout sous-topos de ${\mathcal E}$ défini par une topologie $K$
$$
{\mathcal E}_K \ {\lhook\joinrel\relbar\joinrel\rightarrow} \ {\mathcal E} \, ,
$$
on sait d'après le théorème \ref{thm4.1.8} (ii) que son image réciproque par $f : {\mathcal E}' \to {\mathcal E}$
$$
f^{-1} {\mathcal E}_K = {\mathcal E}'_{K'} \ {\lhook\joinrel\relbar\joinrel\rightarrow} \ {\mathcal E}'
$$
est définie par la topologie $K'$ de ${\mathcal E}'$ constituée des monomorphismes de ${\mathcal E}'$
$$
C \ {\lhook\joinrel\relbar\joinrel\rightarrow} \ X
$$
dont le monomorphisme de ${\mathcal E}$ associé
$$
C_f \ {\lhook\joinrel\relbar\joinrel\rightarrow} \ f_! X
$$
est élément de la topologie $K$.

Il en résulte que la condition
$$
J' \subseteq K'
$$
est équivalente à la condition
$$
J \subseteq K \, .
$$
Autrement dit, la condition
$$
{\mathcal E}'_{J'} \supseteq f^{-1} {\mathcal E}_K
$$
est équivalente à la condition
$$
f_! \, {\mathcal E}'_{J'} \supseteq {\mathcal E}_K
$$
si l'on note ${\mathcal E}_J = f_! \, {\mathcal E}'_{J'}$.

Comme toute topologie $J'$ de ${\mathcal E}'$ est engendrée par des familles de monomorphismes, la formule de (iii) définit un foncteur
$$
f_! : ({\rm ST} ({\mathcal E}') , \supseteq) \longrightarrow ({\rm ST} ({\mathcal E}) , \supseteq)
$$
qui est adjoint à gauche de $f^{-1}$.

C'est le contenu de (ii).

Enfin, (i) est conséquence formelle de (ii).
\end{demo}

\subsection{Le cas particulier des images directes extraordinaires le long de fibrations}\label{4.2.2}

Dans le cas particulier de morphismes de topos induits par des fibrations de catégories de présentation et par des topologies de Giraud, on a:

\begin{cor}\label{cor4.2.2}
Soit un foncteur entre deux catégories essentiellement petites
$$
p : {\mathcal C} \longrightarrow {\mathcal B}
$$
qui est une fibration et induit donc un morphisme localement connexe de topos de préfaisceaux
$$
(p^* , p_*) : \widehat{\mathcal C} \longrightarrow \widehat{\mathcal B} \, .
$$

Soient un sous-topos du topos de base défini par une topologie $J$ de ${\mathcal B}$
$$
\widehat{\mathcal B}_J \ {\lhook\joinrel\relbar\joinrel\rightarrow} \ \widehat{\mathcal B}
$$
et son image réciproque par le morphisme $\widehat{\mathcal C} \to \widehat{\mathcal B}$ qui est un sous-topos
$$
\widehat{\mathcal C}_{J'} \ {\lhook\joinrel\relbar\joinrel\rightarrow} \ \widehat{\mathcal C}
$$
défini par la topologie de Giraud $J'$ de ${\mathcal C}$ associée à $J$.

Alors:

\noindent (i) Le morphisme de topos induit 
$$
(p^* , p_*) : \widehat{\mathcal C}_{J'} \longrightarrow \widehat{\mathcal C}_J
$$
est localement connexe.

\noindent (ii) Le morphisme induit d'image réciproque des sous-topos
$$
p^{-1} : ({\rm ST} (\widehat{\mathcal B}_J) , \supseteq) \longrightarrow ({\rm ST} (\widehat{\mathcal C}_{J'}), \supseteq)
$$
respecte les réunions arbitraires commes les intersections arbitraires, et il possède à la fois un adjoint à droite
$$
p_* : ({\rm ST} (\widehat{\mathcal C}_{J'}) , \supseteq) \longrightarrow ({\rm ST} (\widehat{\mathcal B}_J), \supseteq)
$$
et un adjoint à gauche, le foncteur d'image directe extraordinaire,
$$
p_! : ({\rm ST} (\widehat{\mathcal C}_{J'}) , \supseteq) \longrightarrow ({\rm ST} (\widehat{\mathcal B}_J), \supseteq) \, .
$$

\noindent (iii) Pour tout sous-topos de $\widehat{\mathcal C}_{J'}$
$$
\widehat{\mathcal C}_{K'} \ {\lhook\joinrel\relbar\joinrel\rightarrow} \ \widehat{\mathcal C}_{J'}
$$
défini par une topologie $K' \supseteq J'$ de ${\mathcal C}$ qui est engendrée sur $J'$ par une famille de cribles d'objets $X_i$ de ${\mathcal C}$
$$
C_i \ {\lhook\joinrel\relbar\joinrel\rightarrow} \ y(X_i) \, , \quad i \in I \, ,
$$
son image directe extraordinaire est le sous-topos
$$
\widehat{\mathcal B}_K \ {\lhook\joinrel\relbar\joinrel\rightarrow} \ \widehat{\mathcal B}_J
$$
défini par la topologie $K$ de ${\mathcal B}$ qui est engendrée sur $J$ par la famille des cribles
$$
(C_i)_p \ {\lhook\joinrel\relbar\joinrel\rightarrow} \ y(p(X_i)) \, , \quad i \in I \, ,
$$
où, pour tout crible d'un objet $X$ de ${\mathcal C}$
$$
C \ {\lhook\joinrel\relbar\joinrel\rightarrow} \ y(X) \, ,
$$
on a noté
$$
C_p \ {\lhook\joinrel\relbar\joinrel\rightarrow} \ y(p(X))
$$
le crible de l'objet $p(X)$ de ${\mathcal B}$ constitué des morphismes
$$
Y \longrightarrow p(X)
$$
dont un relèvement horizontal
$$
X' \longrightarrow X
$$
est élément du crible $C$.
\end{cor}

\begin{demo}
\noindent (i) est un cas particulier du corollaire \ref{cor4.1.10} (i) puisque, d'après la proposition \ref{prop4.1.6}, le morphisme de topos de préfaisceaux
$$
\widehat{\mathcal C} \longrightarrow \widehat{\mathcal B}
$$
induit par une fibration $p : {\mathcal C} \to {\mathcal B}$ est localement connexe.

\noindent (ii) résulte de (i) et du corollaire \ref{cor4.2.1}.

\noindent (iii) Considérons un sous-topos
$$
\widehat{\mathcal B}_{J_1} \ {\lhook\joinrel\relbar\joinrel\rightarrow} \ \widehat{\mathcal B}_J
$$
défini par une topologie $J_1$ de ${\mathcal B}$ qui contient $J$, et son image réciproque
$$
p^{-1} \widehat{\mathcal B}_{J_1} = \widehat{\mathcal C}_{J'_1} \ {\lhook\joinrel\relbar\joinrel\rightarrow} \ \widehat{\mathcal C}_{J'}
$$
qui est définie par la topologie de Giraud $J'_1$ de ${\mathcal C}$ associée à $J_1$.

Cela signifie qu'un crible d'un objet $X$ de ${\mathcal C}$
$$
C \ {\lhook\joinrel\relbar\joinrel\rightarrow} \ y(X)
$$
est élément de $J'_1$ si et seulement si le crible associé de l'objet $p(X)$ de ${\mathcal B}$
$$
C_p \ {\lhook\joinrel\relbar\joinrel\rightarrow} \ y(p(X))
$$
est élément de $J_1$.

Ainsi, on a
$$
K' \subseteq J'_1
$$
si et seulement si
$$
K \subseteq J_1 \, .
$$
Autrement dit, on a la relation entre sous-topos
$$
\widehat{\mathcal C}_{K'} \supseteq \widehat{\mathcal C}_{J'_1} = p^{-1} \widehat{\mathcal B}_{J_1}
$$
si et seulement si
$$
\widehat{\mathcal B}_K \supseteq \widehat{\mathcal B}_{J_1} \, .
$$
On conclut comme voulu que
$$
\widehat{\mathcal B}_K = p_! \, \widehat{\mathcal C}_{K'} \, .
$$
\end{demo}

\subsection{Factorisation des morphismes de topos et images réciproques des réunions de sous-topos}\label{4.2.3}

Nous allons déduire du corollaire \ref{cor4.2.2} le résultat suivant déjà annoncé comme corollaire \ref{cor1.3.13}:

\begin{cor}\label{cor4.2.3}
Pour tout morphisme de topos
$$
f = (f^* , f_*) : {\mathcal E}' \longrightarrow {\mathcal E} \, ,
$$
le morphisme induit des images réciproques des sous-topos
$$
f^{-1} : ({\rm ST} ({\mathcal E}), \supseteq) \longrightarrow ({\rm ST} ({\mathcal E}'), \supseteq)
$$
respecte les réunions finies de sous-topos.
\end{cor}

\begin{demo}
On sait déjà que $f^{-1}$ respecte les réunions arbitraires de sous-topos si le morphisme $f$ est de la forme
$$
{\mathcal E}' = \widehat{\mathcal C}_{J'} \xrightarrow{ \ (p^*\!, \, p_*) \ } \widehat{\mathcal B}_J = {\mathcal E}
$$
pour une fibration $p : {\mathcal C} \to {\mathcal B}$ et une topologie de Giraud $J'$ sur ${\mathcal C}$ associée à une topologie $J$ sur ${\mathcal B}$, ou plus généralement si
$$
f : {\mathcal E}' \longrightarrow {\mathcal E}
$$
est un morphisme de topos localement connexe.

On sait d'autre part d'après la proposition \ref{prop1.3.2} (ii) que si
$$
f = j : {\mathcal E}' \ {\lhook\joinrel\relbar\joinrel\rightarrow} \ {\mathcal E}
$$
est un morphisme de plongement d'un sous-topos, l'application d'image réciproque par $j$, qui n'est autre que l'application d'intersection avec le sous-topos ${\mathcal E}'$ de ${\mathcal E}$
$$
j^{-1} = {\mathcal E}' \wedge (\bullet) : ({\rm ST} ({\mathcal E}) , \supseteq) \longrightarrow ({\rm ST} ({\mathcal E}') , \supseteq) \, ,
$$
respecte les réunions finies de sous-topos.
\end{demo}

La conclusion résulte du théorème de factorisation suivant:

\begin{thm}\label{thm4.2.3}
Soit un morphisme de topos
$$
f : {\mathcal E}' \longrightarrow {\mathcal E} \, .
$$
Alors il admet une factorisation de la forme
$$
{\mathcal E}' \xrightarrow{ \ \sim \ } {\mathcal C}'_{K'} \ {\lhook\joinrel\relbar\joinrel\rightarrow} \ {\mathcal C}'_{J'} \longrightarrow {\mathcal B}_J \xrightarrow{ \ \sim \ } {\mathcal E}
$$
où

\medskip

$\left\{\begin{matrix}
\bullet &\mbox{${\mathcal C}'$ et ${\mathcal B}$ sont des petites catégories munies de topologies $K'$ contenant $J'$, et $J$,} \hfill \\
\bullet &\mbox{les morphismes de topos} \hfill \\
&{\mathcal E}' \xrightarrow{ \ \sim \ } \widehat{\mathcal C}'_{K'} \quad \mbox{et} \quad \widehat{\mathcal B}_J \xrightarrow{ \ \sim \ } {\mathcal E} \\
&\mbox{sont des équivalences,} \hfill \\
\bullet &\mbox{le morphisme de topos} \hfill \\
&\widehat{\mathcal C}'_{K'} \ {\lhook\joinrel\relbar\joinrel\rightarrow} \ \widehat{\mathcal C}'_{J'} \\
&\mbox{est le plongement de sous-topos induit par la relation $K' \supseteq J'$,} \hfill \\
\bullet &\mbox{le morphisme de topos} \hfill \\
&\widehat{\mathcal C}'_{J'} \longrightarrow \widehat{\mathcal B}_J \\
&\mbox{est induit par un foncteur} \hfill \\
&{\mathcal C}' \longrightarrow {\mathcal B} \\
&\mbox{qui est une fibration et par la topologie de Giraud $J'$ de ${\mathcal C}'$ associée à la topologie $J$ de ${\mathcal B}$, } \hfill \\
&\mbox{donc est localement connexe.} \hfill
\end{matrix}\right.$
\end{thm}

\begin{remark}
Ainsi, tout morphisme de topos peut s'écrire comme le composé d'un morphisme de plongement de sous-topos et d'un morphisme de topos localement connexe.

Le corollaire \ref{cor4.2.3} en résulte.
\end{remark}

\begin{demo}
On commence par rappeler:

\begin{prop}\label{prop4.2.4}
\noindent (i) Soit ${\mathcal E}$ un topos.

Soit ${\mathcal B}$ une petite sous-catégorie pleine de ${\mathcal E}$ qui contient une famille séparante d'objets de ${\mathcal E}$.

Soit $J$ la topologie de ${\mathcal B}$ induite par la topologie canonique de ${\mathcal E}'$, pour laquelle une famille de morphismes
$$
(E_i \longrightarrow E)_{i \in I}
$$
est couvrante si et seulement si elle est globalement épimorphique dans ${\mathcal E}$. Alors le foncteur de plongement
$$
{\mathcal C} \ {\lhook\joinrel\relbar\joinrel\rightarrow} \ {\mathcal E}
$$
induit deux équivalences de topos réciproques l'une de l'autre
$$
{\mathcal E} \xrightarrow{ \sim \ } \widehat{\mathcal B}_J \quad \mbox{et} \quad {\mathcal B}_J \xrightarrow{ \sim \ } {\mathcal E} \, .
$$
(ii) Soit un morphisme de topos
$$
{\mathcal E}' \xrightarrow{ \ f \ } {\mathcal E} \, .
$$

Soeint ${\mathcal B}$ et ${\mathcal C}$ deux petites sous-catégories pleines de ${\mathcal E}$ et ${\mathcal E}'$ telles que le foncteur
$$
f^* : {\mathcal E} \longrightarrow {\mathcal E}'
$$
se restreigne en un foncteur
$$
\rho : {\mathcal B} \longrightarrow {\mathcal C} \, .
$$

Soient $J$ et $K$ les topologies de ${\mathcal B}$ et ${\mathcal C}$ induites par les topologies canoniques de ${\mathcal E}$ et ${\mathcal E}'$.

Alors le foncteur
$$
{\mathcal B} \xrightarrow{ \ \rho \ } {\mathcal C} \xrightarrow{ \ \ell \ } \widehat{\mathcal C}_K
$$
est plat et $J$-continu, donc définit un morphisme de topos
$$
\widehat{\mathcal C}_K \longrightarrow \widehat{\mathcal B}_J
$$
et le morphisme
$$
f : {\mathcal E}' \longrightarrow {\mathcal E}
$$
se factorise en
$$
{\mathcal E}' \xrightarrow{ \ \sim \ } \widehat{\mathcal C}_K \longrightarrow \widehat{\mathcal B}_J \xrightarrow{ \ \sim \ } {\mathcal E} \, .
$$
(iii) Dans la représentation de (ii), on peut supposer que les petites sous-catégories pleines ${\mathcal B}$ et ${\mathcal C}$ de ${\mathcal E}$ et ${\mathcal E}'$ sont cartésiennes, c'est-à-dire ont des limites finies arbitraires.

Dans ce cas, le foncteur induit
$$
\rho : {\mathcal B} \longrightarrow {\mathcal C}
$$
est nécessairement cartésien.
\end{prop}

\noindent {\bf Esquisse de démonstration:}

\noindent (i) est le contenu de la remarque (ii) qui suit le théorème \ref{thm1.1.7}. Sa démonstration s'applique à toute catégorie ${\mathcal E}$ qui possède les propriétés (0), (1), (2), (3), (5), (7) et (10) du théorème \ref{thm1.1.7}.

C'est la démonstration du théorème de Giraud rappelé dans la remarque (i) qui suit ce théorème \ref{thm1.1.7}.

\noindent (ii) est conséquence de (i).

Il suffit de prendre pour ${\mathcal B}$ et ${\mathcal C}$ n'importe quelles petites cous-catégories pleines de ${\mathcal E}$ et ${\mathcal E}'$ qui contiennent des familles séparantes d'objets pour obtenir des équivalences de topos
$$
{\mathcal E} \xrightarrow{ \ \sim \ } \widehat{\mathcal B}_J \, , \quad \widehat{\mathcal B}_J \xrightarrow{ \ \sim \ } {\mathcal E} \, ,
$$
et
$$
{\mathcal E}' \xrightarrow{ \ \sim \ } \widehat{\mathcal C}_K \, , \quad \widehat{\mathcal C}_K \xrightarrow{ \ \sim \ } {\mathcal E}' \, .
$$
Si l'on impose de plus que le foncteur
$$
{\mathcal B} \ {\lhook\joinrel\relbar\joinrel\rightarrow} \ {\mathcal E} \xrightarrow{ \ f^* \ } {\mathcal E}'
$$
se factorise à travers la petite sous-catégorie pleine
$$
{\mathcal C} \ {\lhook\joinrel\relbar\joinrel\rightarrow} \ {\mathcal E}' \, ,
$$
on obtient par restriction un foncteur
$$
{\mathcal B} \xrightarrow{ \ \rho \ } {\mathcal C}
$$
qui a les propriétés requises puisque le foncteur
$$
f^* : {\mathcal E} \longrightarrow {\mathcal E}'
$$
respecte les limites finies ainsi que les familles globalement épimorphiques.

\noindent (iii) Partant de deux familles séparantes ${\mathcal J}$ et ${\mathcal J}'$ d'objets de ${\mathcal E}$ et ${\mathcal E}'$ on peut prendre

\medskip

$\left\{\begin{matrix}
\bullet &\mbox{pour ${\mathcal B}$ la petite sous-catégorie pleine de ${\mathcal E}$ fermée sur les limites finies d'objets de ${\mathcal J}$,} \hfill \\
\bullet &\mbox{pour ${\mathcal C}$ la petite sous-catégorie pleine de ${\mathcal E}'$ fermée sur les limites finies d'objets qui sont dans ${\mathcal J}$} \hfill \\
&\mbox{ou sont des images par $f^* : {\mathcal E} \to {\mathcal E}'$ d'objets de ${\mathcal B} \hookrightarrow {\mathcal E}$.} \hfill
\end{matrix}\right.$

\medskip

\noindent Alors les petites sous-catégories pleines
$$
{\mathcal B} \ {\lhook\joinrel\relbar\joinrel\rightarrow} \ {\mathcal E} \quad \mbox{et} \quad {\mathcal C} \ {\lhook\joinrel\relbar\joinrel\rightarrow} \ {\mathcal E}'
$$
sont cartésiennes, et les limites finies se calculent dans ${\mathcal B}$ ou ${\mathcal C}$ comme dans ${\mathcal E}$ ou ${\mathcal E}'$.

Le foncteur
$$
f^* : {\mathcal E} \longrightarrow {\mathcal E}'
$$
respecte les limites finies, donc aussi sa restriction
$$
\rho : {\mathcal B} \longrightarrow {\mathcal C} \, .
$$
\end{demo}

Puis on a le lemme suivant:

\begin{lem}\label{lem4.2.5}
Soient ${\mathcal B}$ et ${\mathcal C}$ deux catégories localement petites reliées par un foncteur
$$
\rho : {\mathcal B} \longrightarrow {\mathcal C} \, .
$$

Considérons la catégorie
$$
{\mathcal C}' = {\mathcal C} / {\mathcal B}
$$
dont

\medskip

$\left\{\begin{matrix}
\bullet &\mbox{les objets sont les triplets $(X,Y,X \to \rho (Y))$ constitués de} \hfill \\
&\mbox{ \quad -- un objet $X$ de ${\mathcal C}$,} \hfill \\
&\mbox{ \quad -- un objet $Y$ de ${\mathcal B}$,} \hfill \\
&\mbox{ \quad -- un morphisme $X \to \rho (Y)$ de ${\mathcal C}$,} \hfill \\
\bullet &\mbox{les morphismes} \hfill \\
&(X',Y',X' \to \rho (Y')) \longrightarrow (X,Y,X \to \rho (Y)) \\
&\mbox{sont les paires constituées de} \hfill \\
&\mbox{ \quad -- un morphisme $X' \xrightarrow{ \, x \, } X$ de ${\mathcal C}$,} \hfill \\
&\mbox{ \quad -- un morphisme $Y' \xrightarrow{ \, y \, } Y$ de ${\mathcal B}$,} \hfill \\
&\mbox{rendant commutatif le carré de ${\mathcal C}$:} \hfill \\
&\xymatrix{
X' \ar[d] \ar[r]^x &X \ar[d] \\
\rho (Y') \ar[r]^{\rho (y)} &\rho (Y)
}
\end{matrix}\right.$

\medskip

\noindent Alors:

\noindent (i) Le foncteur d'oubli vers ${\mathcal B}$
$$
\begin{matrix}
p &: &{\mathcal C}' = {\mathcal C}/{\mathcal B} \hfill &\longrightarrow &{\mathcal B} \, , \\
&&(X,Y,X \to \rho(Y)) &\longmapsto &Y
\end{matrix}
$$
admet pour adjoint à droite le foncteur
$$
\begin{matrix}
s : {\mathcal B} &\longrightarrow &{\mathcal C}' = {\mathcal C}/{\mathcal B} \, , \hfill \\
\hfill Y &\longmapsto &\left(\rho (Y) , Y , \rho(Y) \xrightarrow{ \, {\rm id} \, } \rho(Y)\right) .
\end{matrix}
$$
(ii) Si la catégorie ${\mathcal C}$ est cartésienne, le foncteur
$$
p : {\mathcal C}' = {\mathcal C}/{\mathcal B} \longrightarrow {\mathcal B}
$$
est une fibration.

Plus précisément, pour tout objet $(X,Y,X \to \rho(Y))$ de ${\mathcal C}' = {\mathcal C}/{\mathcal B}$, tout morphisme de ${\mathcal B}$
$$
Y' \longrightarrow Y
$$
admet pour relèvement horizontal le morphisme de ${\mathcal C}' = {\mathcal C} / {\mathcal B}$
$$
(X \times_{\rho (Y)} \rho (Y') , Y' , X \times_{\rho (Y)} \rho(Y') \to \rho (Y')) \longrightarrow (X,Y,X \to \rho(Y)) \, .
$$
(iii) Si la catégorie ${\mathcal B}$ admet un objet terminal $1_{\mathcal B}$, et
$$
\rho (1_{\mathcal B}) = 1_{\mathcal C}
$$
est un objet terminal de ${\mathcal C}$, le foncteur d'oubli vers ${\mathcal C}$
$$
\begin{matrix}
q &: &{\mathcal C}' = {\mathcal C}/{\mathcal B} \hfill &\longrightarrow &{\mathcal C} \, , \\
&&(X,Y,X \to \rho(Y)) &\longmapsto &X
\end{matrix}
$$
admet pour adjoint à droite le foncteur
$$
\begin{matrix}
r : {\mathcal C} &\longrightarrow &{\mathcal C}' = {\mathcal C}/{\mathcal B} \, , \hfill  \\
\hfill X &\longmapsto &(X , 1_{\mathcal B} , X \to \rho (1_{\mathcal B})) \, .
\end{matrix}
$$
\end{lem}

\begin{demo}
\noindent (i) En effet, pour tous objets
$$
(X,Y,X \to \rho (Y)) \quad \mbox{de} \quad {\mathcal C}' = {\mathcal C}/{\mathcal B}
$$
et
$$
Y' \quad \mbox{de} \quad {\mathcal B} \, ,
$$
se donner un morphisme de ${\mathcal B}$
$$
Y \xrightarrow{ \ y \ } Y'
$$
équivaut à se donner un morphisme de ${\mathcal C}' = {\mathcal C}/{\mathcal B}$
$$
(X,Y,X \to \rho (Y)) \longrightarrow (\rho (Y') , Y' , \rho (Y') \xrightarrow{ \, {\rm id} \, } \rho(Y'))
$$
puisque le carré de ${\mathcal C}$
$$
\xymatrix{
X \ar[d] \ar[r] &\rho(Y) \ar[d]^{\rho (Y)} \\
\rho (Y') \ar[r]^{{\rm id}} &\rho (Y')
}
$$
est commutatif si et seulement si le morphisme
$$
X \longrightarrow \rho (Y')
$$
est le composé
$$
X \longrightarrow \rho (Y) \xrightarrow{ \ \rho (y) \ } \rho(Y') \, .
$$
(ii) D'après le lemme \ref{lem4.1.2}, un morphisme de ${\mathcal C}' = {\mathcal C}/{\mathcal B}$
$$
(X',Y',X' \to \rho (Y')) \longrightarrow (X,Y,X \to \rho (Y))
$$
est horizontal relativement au foncteur $p : {\mathcal C}' \to {\mathcal B}$ si et seulement si le carré commutatif de ${\mathcal C}'$
$$
\xymatrix{
(X',Y',X' \to \rho (Y)) \ar[d] \ar[r] &(X,Y,X \to \rho (Y)) \ar[d] \\
(\rho (Y') , Y' , \rho (Y') \xrightarrow{ \, {\rm id} \, } \rho (Y')) \ar[r] &(\rho (Y) , Y , \rho (Y) \xrightarrow{ \, {\rm id} \, } \rho (Y)) 
}
$$
est cartésien, c'est-à-dire si et seulement si le carré commutatif de ${\mathcal C}$
$$
\xymatrix{
X' \ar[d] \ar[r] &X \ar[d] \\
\rho (Y') \ar[r] &\rho (Y)
}
$$
est cartésien.

Si les produits fibrés sont bien définis dans ${\mathcal C}$, alors pour tout objet $(X,Y,X \to \rho(Y))$ de ${\mathcal C}' = {\mathcal C}/{\mathcal B}$, tout morphisme de ${\mathcal B}$
$$
Y' \longrightarrow Y
$$ se relève en le morphisme horizontal de ${\mathcal C}' = {\mathcal C}/{\mathcal B}$
$$
(X \times_{\rho (Y)} \rho (Y') , Y' , X \times_{\rho (Y)} \rho(Y') \to \rho (Y')) \longrightarrow (X,Y,X \to \rho (Y)) \, .
$$
Ainsi, le foncteur
$$
\begin{matrix}
\rho &: &{\mathcal C}' = {\mathcal C} / {\mathcal B} \hfill &\longrightarrow {\mathcal B} \, , \\
&&(X,Y,X \to \rho (Y)) &\longrightarrow Y
\end{matrix}
$$
est bien une fibration.

\noindent (iii) En effet, si $(X,Y,X \to \rho (Y))$ est un objet de ${\mathcal C}' = {\mathcal C}/{\mathcal B}$ et $X'$ est un objet de ${\mathcal C}$, tout morphisme de ${\mathcal C}$
$$
X \longrightarrow X'
$$
complété par l'unique morphisme
$$
Y \longrightarrow 1_{\mathcal B}
$$
définit un morphisme de ${\mathcal C}' = {\mathcal C}/{\mathcal B}$
$$
(X,Y,X \to \rho(Y)) \longrightarrow (X' , 1_{\mathcal B} , X' \to \rho (1_{\mathcal B}))
$$
vers l'objet défini par l'unique morphisme $X' \to \rho (1_{\mathcal B})$, puisque le carré
$$
\xymatrix{
X \ar[d] \ar[r] &\rho(Y) \ar[d] \\
X' \ar[r] &\rho (1_{\mathcal B})
}
$$
est nécessairement commutatif.
\end{demo}

Puis on a le théorème de factorisation suivant tiré de l'article \cite{Denseness}:

\begin{thm}\label{thm4.2.6}
Soient $({\mathcal B} , J)$ et $({\mathcal C},K)$ deux sites reliés par un foncteur
$$
\rho : {\mathcal B} \longrightarrow {\mathcal C}
$$
qui définit un morphisme de topos
$$
\widehat{\mathcal C}_K \longrightarrow \widehat{\mathcal B}_J
$$
c'est-à-dire est tel que le foncteur composé
$$
{\mathcal B} \xrightarrow{ \ \rho \ } {\mathcal C} \xrightarrow{ \ \ell \ } \widehat{\mathcal C}_K
$$
est plat et $J$-continu.

Considérons la catégorie
$$
{\mathcal C}' = {\mathcal C}/{\mathcal B}
$$
munie de ses deux foncteurs d'oubli vers ${\mathcal B}$ et ${\mathcal C}$
$$
\begin{matrix}
p &: &{\mathcal C}' = {\mathcal C} / {\mathcal B} \hfill &\longrightarrow &{\mathcal B} \, , \\
&&(X,Y,X \to \rho (Y)) &\longmapsto &Y \, , 
\end{matrix}
$$
et
$$
\begin{matrix}
q &: &{\mathcal C}' = {\mathcal C} / {\mathcal B} \hfill &\longrightarrow &{\mathcal C}  \\
&&(X,Y,X \to \rho (Y)) &\longmapsto &X \, .
\end{matrix}
$$
Alors:

\noindent (i) Il existe une topologie $K'$ sur ${\mathcal C}'$ pour laquelle une famille de morphismes
$$
(X_i,Y_i ,X_i \to \rho(Y_i)) \longrightarrow (X , Y , X \to \rho (Y)) \, , \quad i \in I \, ,
$$
est couvrante si et seulement si son image dans ${\mathcal C}$
$$
X_i \longrightarrow X \, , \quad i \in I \, ,
$$
est couvrante pour la topologie $K$.

\noindent (ii) Le foncteur d'oubli de ${\mathcal C}'$ vers ${\mathcal C}$
$$
\begin{matrix}
q &: &{\mathcal C}' = {\mathcal C} / {\mathcal B} \hfill &\longrightarrow &{\mathcal C} \, , \\
&&(X,Y,X \to \rho (Y)) &\longmapsto &X 
\end{matrix}
$$
définit deux morphismes de topos
$$
\widehat{\mathcal C}_K \longrightarrow \widehat{\mathcal C}'_{K'} \quad \mbox{et} \quad \widehat{\mathcal C}'_{K'} \longrightarrow \widehat{\mathcal C}_K
$$
qui sont des équivalences quasi-inverses l'une de l'autre.

\noindent (iii) Le foncteur d'oubli de ${\mathcal C}'$ vers ${\mathcal B}$
$$
\begin{matrix}
p &: &{\mathcal C}' = {\mathcal C} / {\mathcal B} \hfill &\longrightarrow &{\mathcal B} \, , \\
&&(X,Y,X \to \rho (Y)) &\longmapsto &Y 
\end{matrix}
$$
définit un carré commutatif de morphismes de topos:
$$
\xymatrix{
\, \widehat{\mathcal C}'_{K'} \ar@{_{(}->}[d] \ar[r] &\widehat{\mathcal B}_J \ar@{_{(}->}[d] \\
\widehat{\mathcal C}' \ar[r] &\widehat{\mathcal B}
}
$$
(iv) Le morphisme de topos induit par $\rho$
$$
\widehat{\mathcal C}_K \longrightarrow \widehat{\mathcal B}_J
$$
se factorise comme le composé
$$
\widehat{\mathcal C}_K \xrightarrow{ \ \sim \ } \widehat{\mathcal C}'_{K'} \longrightarrow \widehat{\mathcal B}_J \, .
$$
\end{thm}

\noindent {\bf Démonstration dans un cas suffisant:}

Donnons la démonstration dans le cas plus simple, suffisant pour nous, où les catégories ${\mathcal B}$ et ${\mathcal C}$ sont cartésiennes et où le foncteur
$$
\rho : {\mathcal B} \longrightarrow {\mathcal C}
$$
est lui-même cartésien.

\noindent (i) Il s'agit de vérifier que cette définition de $K'$ satisfait les conditions (0), (1), (2), (3) de la définition \ref{defn1.1.1}, qui sont déjà satisfaites par $K$.

C'est immédiat pour les conditions (0), (1) et (3).

Pour la condition (2), considérons un morphisme de ${\mathcal C}'$
$$
(X',Y',X' \to \rho(Y')) \longrightarrow (X,Y,X \to \rho(Y))
$$
et une famille de morphismes
$$
(X_i,Y_i,X_i \to \rho(Y_i)) \longrightarrow (X,Y,X \to \rho(Y)) \, , \quad i \in I \, ,
$$
qui est $K'$-couvrante, c'est-à-dire telle que la famille
$$
X_i \longrightarrow X \, , \quad i \in I \, ,
$$
est $K$-couvrante.

Si les catégories ${\mathcal B}$ et ${\mathcal C}$ sont cartésiennes et le foncteur
$$
\rho : {\mathcal B} \longrightarrow {\mathcal C}
$$
respecte les produit fibrés, on peut former la famille de morphismes de ${\mathcal C}'$
$$
(X_i \times_X X',  Y_i \times_Y Y', X_i \times_X X' \to \rho(Y_i \times_Y Y')) \longrightarrow (X',Y',X' \to \rho(Y')) \, , \quad i \in I \, .
$$
Elle est $K'$-couvrante puisque la famille de morphismes de ${\mathcal C}$
$$
X_i \times_X X' \longrightarrow X' \, , \quad i \in I \, ,
$$
est $K$-couvrante.

Ainsi, $K'$ satisfait la condition (2).

C'est une topologie.

\noindent (ii) Si les catégories ${\mathcal B}$ et ${\mathcal C}$ sont cartésiennes et que le foncteur $\rho$ respecte les limites finies, la catégorie
$$
{\mathcal C}' = {\mathcal C}/{\mathcal B}
$$
est cartésienne et le foncteur d'oubli
$$
\begin{matrix}
q &: &{\mathcal C}' = {\mathcal C} / {\mathcal B} \hfill &\longrightarrow &{\mathcal C} \, , \\
&&(X,Y,X \to \rho (Y)) &\longmapsto &X 
\end{matrix}
$$
respecte les limites finies.

A fortiori, le foncteur composé
$$
{\mathcal C}' \xrightarrow{ \ q \ } {\mathcal C} \xrightarrow{ \ \ell \ } \widehat{\mathcal C}_K
$$
est plat. Il est $K'$-continu car le foncteur
$$
q : {\mathcal C}' \longrightarrow {\mathcal C}
$$
transforme toute famille $K'$-couvrante en une famille $K$-couvrante.

Ainsi, $q$ définit un morphisme de topos
$$
g=(g^*,g_*) : \widehat{\mathcal C}_K \longrightarrow \widehat{\mathcal C}'_{K'}
$$
dont la composante d'image réciproque $g^*$ s'inscrit dans un carré commutatif:
$$
\xymatrix{
{\mathcal C}' \ar[d]_-{\ell} \ar[r]^q &{\mathcal C} \ar[d]^-{\ell} \\
\widehat{\mathcal C}'_{K'} \ar[r]^{g^*} &\widehat{\mathcal C}_K
}
$$

D'autre part, le foncteur d'oubli
$$
q : {\mathcal C}' = {\mathcal C}/{\mathcal B} \longrightarrow {\mathcal C}
$$
définit le foncteur de composition des préfaisceaux avec $q$
$$
q^* : \widehat{\mathcal C} \longrightarrow \widehat{\mathcal C}' \, .
$$

\newpage

\noindent Comme, d'après le lemme \ref{lem4.2.5} (iii), le foncteur admet pour adjoint à droite le foncteur
$$
\begin{matrix}
r : {\mathcal C} &\longrightarrow &{\mathcal C}' = {\mathcal C}/{\mathcal B} \, , \hfill \\
\hfill X &\longmapsto &(X , 1_{\mathcal B} , X \to \rho (1_{\mathcal B})) \, ,
\end{matrix}
$$
on a un carré commutatif de foncteurs:
$$
\xymatrix{
{\mathcal C}_{ \ } \ar@{_{(}->}[d]_y \ar[r]^r &{\mathcal C}'_{ \ }  \ar@{_{(}->}[d]^y \\
\widehat{\mathcal C} \ar[r]^{q^*} &\widehat{\mathcal C}'
}
$$
Le foncteur composé
$$
{\mathcal C} \ {\lhook\joinrel\relbar\joinrel\rightarrow} \ \widehat{\mathcal C} \longrightarrow \widehat{\mathcal C}' \longrightarrow \widehat{\mathcal C}'_{K'}
$$
est plat. Il est aussi $K$-continu puisque le foncteur $r : {\mathcal C} \to {\mathcal C}'$ transforme toute famille $K$-couvrante en une famille $K'$-couvrante.

Donc il définit un morphisme de topos
$$
h = (h^* , h_*) : \widehat{\mathcal C}'_{K'} \longrightarrow \widehat{\mathcal C}_K
$$
dont la composante d'image réciproque $h^*$ s'inscrit dans un carré commutatif:
$$
\xymatrix{
{\mathcal C} \ar[d]_-{\ell} \ar[r]^r &{\mathcal C}' \ar[d]^-{\ell} \\
\widehat{\mathcal C}_K \ar[r]^{h^*} &\widehat{\mathcal C}'_{K'}
}
$$
Nous devons montrer que les deux foncteurs
$$
g^* : \widehat{\mathcal C}'_{K'} \longrightarrow \widehat{\mathcal C}_K \quad \mbox{et} \quad h^* : \widehat{\mathcal C}_K \longrightarrow \widehat{\mathcal C}'_{K'}
$$
sont quasi-inverses l'un de l'autre.

On a
$$
q \circ r = {\rm id}_{\mathcal C}
$$
d'où on tire un isomorphisme canonique
$$
g^* \circ h^* \xrightarrow{ \ \sim \ } {\rm id}_{\widehat{\mathcal C}_K} \, .
$$
Pour obtenir également un isomorphisme canonique
$$
h^* \circ g^* \xrightarrow{ \ \sim \ } {\rm id}_{\widehat{\mathcal C}'_{K'}} \, ,
$$
il suffit de prouver que, pour tout objet de ${\mathcal C}' = {\mathcal C}/{\mathcal B}$
$$
{\mathcal X} = (X,Y,X \to \rho (Y)) \, ,
$$
le morphisme canonique vers son image par $r \circ q$
$$
r \circ q({\mathcal X}) = (X , 1_{\mathcal B} , X \to \rho(1_{\mathcal B}))
$$
est transformé par le foncteur
$$
\ell : {\mathcal C}' \longrightarrow \widehat{\mathcal C}'_{K'}
$$
en un isomorphisme.

Le morphisme
$$
{\mathcal X} \longrightarrow r \circ q ({\mathcal X})
$$
est $K'$-couvrant, donc il est transformé par $\ell$ en un épimorphisme.

D'autre part, on a
$$
{\mathcal X} \times_{r \circ q ({\mathcal X})} {\mathcal X} = (X,Y \times Y,X \to \rho (Y \times Y))
$$
et le morphisme diagonal
$$
(X,Y,X \to \rho (Y)) = {\mathcal X} \longrightarrow {\mathcal X} \times_{r \circ q ({\mathcal X})} {\mathcal X}
$$
est $K'$-couvrant, donc le morphisme
$$ 
{\mathcal X} \longrightarrow r \circ q ({\mathcal X})
$$
est transformé par $\ell$ en un monomorphisme.

Ainsi, chaque morphisme du topos $\widehat{\mathcal C}'_{K'}$
$$
\ell ({\mathcal X}) \longrightarrow \ell (r \circ q ({\mathcal X}))
$$
est à  la fois un épimorphisme et un monomorphisme. C'est donc un isomorphisme.

On a bien un isomorphisme canonique
$$
h^* \circ g^* \xrightarrow{ \ \sim \ } {\rm id}_{\widehat{\mathcal C}'_{K'}} \, ,
$$
ce qui termine la démonstration de (ii).

\noindent (iii) Le foncteur d'oubli de ${\mathcal C}' = {\mathcal C}/{\mathcal B}$ vers ${\mathcal B}$
$$
\begin{matrix}
p &: &{\mathcal C}' = {\mathcal C} / {\mathcal B} \hfill &\longrightarrow &{\mathcal B} \, , \\
&&(X,Y,X \to \rho (Y)) &\longmapsto &Y \, ,
\end{matrix}
$$
définit un morphisme de topos de préfaisceaux
$$
(p^* , p_*) : \widehat{\mathcal C}' \longrightarrow \widehat{\mathcal B}
$$
dont la composante d'image réciproque, qui est le foncteur de composition des préfaisceaux avec $p$
$$
p^* : \widehat{\mathcal B} \longrightarrow \widehat{\mathcal C}' \, ,
$$
s'inscrit dans un carré commutatif de foncteurs:
$$
\xymatrix{
{\mathcal B}_{ \ } \ar@{_{(}->}[d]_y \ar[r]^s &{\mathcal C}'_{ \ } \ar@{_{(}->}[d]^y \\
\widehat{\mathcal B} \ar[r]^{p^*} &\widehat{\mathcal C}'
}
$$
Cela résulte en effet de ce que, d'après le lemme \ref{lem4.2.5} (i), le foncteur d'oubli
$$
p : {\mathcal C}' = {\mathcal C} / {\mathcal B} \longrightarrow {\mathcal B}
$$
admet pour adjoint à droite le foncteur 
$$
\begin{matrix}
s : {\mathcal B}  &\longrightarrow &{\mathcal C}' \, , \hfill \\
\hfill Y &\longmapsto &(\rho (Y) , Y , \rho (Y) \xrightarrow{ \, {\rm id} \, } \rho (Y)) \, .
\end{matrix}
$$

Le foncteur $\rho$ transforme les familles $J$-couvrantes de ${\mathcal B}$ en des familles $K$-couvrantes de ${\mathcal C}$, donc le foncteur
$$
s : {\mathcal B} \longrightarrow {\mathcal C}'
$$
transforme les familles $J$-couvrantes en familles $K'$-couvrantes.

Ainsi, le foncteur composé
$$
\xymatrix{
{\mathcal B} \, \ar@{^{(}->}[r]^y &\widehat{\mathcal B} \ar[r]^{p^*} &\widehat{\mathcal C}' \ar[r] &\widehat{\mathcal C}'_{K'}
}
$$
est plat et $J$-continu.

Il définit un morphisme de topos
$$
\widehat{\mathcal C}'_{K'} \longrightarrow \widehat{\mathcal B}_J
$$
qui s'inscrit dans un carré commutatif:
$$
\xymatrix{
\widehat{\mathcal C}'_{K'} \ar@{_{(}->}[d] \ar[r] &\widehat{\mathcal B}_J \ar@{_{(}->}[d] \\
\widehat{\mathcal C}' \ar[r]^{(p^*\!, \, p_*)} &\widehat{\mathcal B}
}
$$
(iv) La composante d'image réciproque du morphisme de topos
$$
\widehat{\mathcal C}'_{K'} \longrightarrow \widehat{\mathcal B}_J
$$
s'inscrit dans un carré commutatif:
$$
\xymatrix{
{\mathcal B} \ar[d]_-{\ell} \ar[r]^s &{\mathcal C}' \ar[d]^-{\ell} \\
\widehat{\mathcal B}_J \ar[r] &\widehat{\mathcal C}'_{K'}
}
$$

D'autre part, la composante d'image réciproque du morphisme de topos
$$
(g^* , g_*) : \widehat{\mathcal C}_K \xrightarrow{ \ \sim \ } \widehat{\mathcal C}'_{K'}
$$
s'inscrit dans un caré commutatif:
$$
\xymatrix{
{\mathcal C}' \ar[d]_-{\ell} \ar[r]^q &{\mathcal C} \ar[d]^-{\ell} \\
\widehat{\mathcal C}'_{K'} \ar[r]^{g^*} &\widehat{\mathcal C}_K
}
$$

Or, le composé du foncteur
$$
\begin{matrix}
s : {\mathcal B}  &\longrightarrow &{\mathcal C}' \, , \hfill \\
\hfill Y &\longmapsto &(\rho (Y) , Y , \rho (Y) \xrightarrow{ \, {\rm id} \, } \rho (Y)) 
\end{matrix}
$$
et du foncteur
$$
\begin{matrix}
q &: &{\mathcal C}'  \hfill &\longrightarrow &{\mathcal C} \, , \hfill \\
&&(X,Y,X \to \rho(Y)) &\longmapsto &X \hfill
\end{matrix}
$$
est le foncteur
$$
\rho : {\mathcal B} \longrightarrow {\mathcal C} \, .
$$

On en déduit comme annoncé que le foncteur 
$$
\widehat{\mathcal C}_K \longrightarrow \widehat{\mathcal B}_J
$$
induit par $\rho$ se factorise comme le composé de l'équivalence
$$
(g^* , g_*) : \widehat{\mathcal C}_K \xrightarrow{ \ \sim \ } \widehat{\mathcal C}'_{K'}
$$
et du morphisme de topos
$$
\widehat{\mathcal C}'_{K'} \longrightarrow \widehat{\mathcal B}_J \, .
$$

Cela termine la démonstration du théorème \ref{thm4.2.6} dans le cas où les catégories ${\mathcal B}$ et ${\mathcal C}$ sont cartésiennes et le foncteur
$$
\rho : {\mathcal B} \longrightarrow {\mathcal C}
$$
est cartésien. \hfill $\Box$

\bigskip

\noindent {\bf Fin de la démonstration du théorème \ref{thm4.2.3}:}

Il résulte de la proposition \ref{prop4.2.4} que tout morphisme de topos
$$
f : {\mathcal E}'  \longrightarrow {\mathcal E}
$$
se représente sous la forme d'un composé
$$
{\mathcal E}' \xrightarrow{ \ \sim \ } \widehat{\mathcal C}_K \longrightarrow \widehat{\mathcal B}_J \xrightarrow{ \ \sim \ } {\mathcal E}
$$
où

$\left\{\begin{matrix}
\bullet &\mbox{les équivalences de topos} \hfill \\
&{\mathcal E}' \xrightarrow{ \ \sim \ } \widehat{\mathcal C}_K \quad \mbox{et} \quad \widehat{\mathcal B}_J \xrightarrow{ \ \sim \ } {\mathcal E} \\
&\mbox{sont définies par deux petites sous-catégories pleines} \hfill \\
&{\mathcal C} \hookrightarrow {\mathcal E}' \quad \mbox{et} \quad {\mathcal B} \hookrightarrow {\mathcal E} \\
&\mbox{que l'on peut supposer cartésiennes,} \hfill \\
\bullet &\mbox{le morphisme de topos} \hfill \\
&\widehat{\mathcal C}_K \longrightarrow \widehat{\mathcal B}_J \\
&\mbox{est induit par un foncteur} \hfill \\
&\rho : {\mathcal B} \longrightarrow {\mathcal C} \\
&\mbox{que l'on peut supposer cartésien, et qui transforme les familles $J$-couvrantes} \hfill \\
&\mbox{en familles $K$-couvrantes.} \hfill
\end{matrix} \right.$

\medskip

D'après le théorème \ref{thm4.2.6}, le morphisme de topos induit par $\rho$
$$
\widehat{\mathcal C}_K \longrightarrow \widehat{\mathcal B}_J
$$
se factorise comme le composé de l'équivalence
$$
\widehat{\mathcal C}_K \xrightarrow{ \ \sim \ } \widehat{\mathcal C}'_{K'}
$$
et d'un morphisme de topos
$$
\widehat{\mathcal C}'_{K'} \longrightarrow \widehat{\mathcal B}_J
$$
qui s'inscrit dans un carré commutatif:
$$
\xymatrix{
\widehat{\mathcal C}'_{K'} \ar@{_{(}->}[d] \ar[r] &\widehat{\mathcal B}_J \ar@{_{(}->}[d] \\
\widehat{\mathcal C}' \ar[r]^{(p^*\!, \, p_*)} &\widehat{\mathcal B}
}
$$
D'après le lemme \ref{lem4.2.5} (ii), le foncteur d'oubli
$$
\begin{matrix}
{\mathcal C}' = {\mathcal C}/{\mathcal B} &\longrightarrow &{\mathcal B} \, , \\
(X,Y,X \to \rho(Y)) &\longmapsto &Y
\end{matrix}
$$
est une fibration.

L'image réciproque du sous-topos
$$
\widehat{\mathcal B}_J \ {\lhook\joinrel\relbar\joinrel\rightarrow} \ \widehat{\mathcal B}
$$
par le morphisme de topos
$$
\widehat{\mathcal C}' \xrightarrow{ \ (p^*\!, \, p_*) \ } \widehat{\mathcal B}
$$
est le sous-topos
$$
\widehat{\mathcal C}'_{J'} \ {\lhook\joinrel\relbar\joinrel\rightarrow} \ \widehat{\mathcal C}'
$$
défini par la topologie de Giraud $J'$ de ${\mathcal C}'$ associée à la topologie $J$ de ${\mathcal B}$ via la fibration
$$
p : {\mathcal C}' \longrightarrow {\mathcal B} \, .
$$
La commutativité du carré
$$
\xymatrix{
\widehat{\mathcal C}'_{K'} \ar@{_{(}->}[d] \ar[r] &\widehat{\mathcal B}_J \ar@{_{(}->}[d] \\
\widehat{\mathcal C}' \ar[r]^{(p^*\!, \, p_*)} &\widehat{\mathcal B}
}
$$
signifie que $K'$ contient la topologie de Girand $J'$ et que le morphisme de topos
$$
\widehat{\mathcal C}'_{K'}  \longrightarrow \widehat{\mathcal B}_J
$$
se factorise comme le composé du plongement
$$
\widehat{\mathcal C}'_{K'} \ {\lhook\joinrel\relbar\joinrel\rightarrow} \ \widehat{\mathcal C}'_{J'}
$$
et du morphisme de topos
$$
\widehat{\mathcal C}'_{J'} \longrightarrow \widehat {\mathcal B}_J 
$$
qui, d'après le corollaire \ref{cor4.2.2} (i), est localement connexe. \hfill $\Box$

\end{document}